\documentclass[10pt,a4paper]{report}
\usepackage{t1enc}

\usepackage{makeidx}
\makeindex
\usepackage[utf8]{inputenc}
\usepackage{lmodern}
\usepackage[frenchb]{babel}
\usepackage{graphicx}
\usepackage{asymptote}
\usepackage{amsmath}
\usepackage{amssymb}
\usepackage{amsthm}
\usepackage{mathrsfs}
\usepackage[all]{xy}
\usepackage{enumerate}

\theoremstyle{plain}
\newtheorem{thm}{Théorème}[section]
\newtheorem{prop}[thm]{Proposition}

\newtheorem{cor}[thm]{Corollaire}
\newtheorem{lem}[thm]{Lemme}

\theoremstyle{definition}
\newtheorem{dfn}{Définition}[section]
\newtheorem*{ex}{Exemple}
\newtheorem*{exs}{Exemples}

\theoremstyle{remark}
\newtheorem*{rem}{Remarque}

\DeclareMathOperator{\An}{An}
\DeclareMathOperator{\Out}{Out}
\DeclareMathOperator{\Aut}{Aut}

\DeclareMathOperator{\End}{End}
\DeclareMathOperator{\Gal}{Gal}

\DeclareMathOperator{\Spec}{Spec}
\DeclareMathOperator{\Spf}{Spf}
\DeclareMathOperator{\gf}{\pi_1}
\DeclareMathOperator{\gpd}{\pi_1^{pd}}
\DeclareMathOperator{\get}{\pi_1^{\acute{e}t}}
\DeclareMathOperator{\ga}{\pi_1^{alg}}
\DeclareMathOperator{\gtop}{\pi_1^{top}}
\DeclareMathOperator{\gctop}{\pi_1^{\underline{top}}}
\DeclareMathOperator{\gt}{\pi_1^{t}}
\DeclareMathOperator{\gctemp}{\pi_1^{\underline{temp}}}
\DeclareMathOperator{\gorb}{\pi_1^{orb}}
\DeclareMathOperator{\gtemp}{\pi_1^{temp}}
\DeclareMathOperator{\glog}{\pi_1^{log}}
\DeclareMathOperator{\ggeom}{\pi_1^{{log-g\acute{e}om}}}
\DeclareMathOperator{\gtempgeom}{\pi_1^{temp-geom}}
\DeclareMathOperator{\g0}{\pi_0}
\newcommand{\temp}{\mathrm{temp}}
\DeclareMathOperator{\Norm}{Norm}
\DeclareMathOperator{\Cov}{Cov}
\DeclareMathOperator{\KCov}{KCov}
\DeclareMathOperator{\KCovgeom}{KCov_{g\acute{e}om}}
\DeclareMathOperator{\GalKCov}{GalKCov}
\newcommand{\LGalKCov}{\mbb L\text{-}\GalKCov}
\DeclareMathOperator{\Covalg}{Cov^{alg}}
\DeclareMathOperator{\Covtop}{Cov^{top}}
\DeclareMathOperator{\Covtemp}{Cov^{temp}}
\DeclareMathOperator{\ket}{{k\acute{e}t}}

\DeclareMathOperator{\Ba}{\mcal B^{alg}}
\DeclareMathOperator{\Bpd}{\mcal B^{pd}}
\DeclareMathOperator{\Btop}{\mcal B^{top}}
\DeclareMathOperator{\Bctop}{\mcal B^{\underline{top}}}
\DeclareMathOperator{\Bcov}{\mathcal B^{cov}}
\DeclareMathOperator{\Btemp}{\mathcal B^{temp}}
\DeclareMathOperator{\BtempL}{\mathcal B^{temp,\mbb L}}
\DeclareMathOperator{\Bctemp}{\mathcal B^{\underline{temp}}}
\DeclareMathOperator{\lcs}{lcs}
\DeclareMathOperator{\rk}{rk}
\newcommand{\gp}{{\mathrm{gp}}}
\DeclareMathOperator{\Ker}{Ker}
\DeclareMathOperator{\Coker}{Coker}
\DeclareMathOperator{\red}{red}
\DeclareMathOperator{\Str}{Str}

\DeclareMathOperator{\DD}{DD}
\DeclareMathOperator{\DDtemp}{DD_{temp}}
\DeclareMathOperator{\Dtop}{\mcal D_{top}}
\DeclareMathOperator{\Dtopan}{\mcal D_{top}^{an}}
\DeclareMathOperator{\Dtops}{\mcal D_{top}^{sp}}
\DeclareMathOperator{\Dtopgeom}{\mcal D_{top-g\acute{e}om}}
\DeclareMathOperator{\Dtemp}{\mcal D_{temp}}
\DeclareMathOperator{\Cart}{Cart}
\DeclareMathOperator{\Ob}{Ob}
\DeclareMathOperator{\Hom}{Hom}
\newcommand{\ab}{{\mathrm{ab}}}
\newcommand{\et}{{\mathrm{\acute{e}t}}}
\newcommand{\etale}{\et}
\newcommand{\pd}{\mathrm{pd}}
\DeclareMathOperator{\Et}{\acute{E}t}
\newcommand{\qet}{{\mathrm{q\acute{e}t}}}
\DeclareMathOperator{\Qet}{Q\acute{e}t}
\newcommand{\alg}{\mathrm{alg}}
\DeclareMathOperator{\Zar}{Zar}
\DeclareMathOperator{\Ens}{Ens}
\DeclareMathOperator{\Set}{Ens}
\DeclareMathOperator{\tSet}{\text{-}Ens}
\DeclareMathOperator{\tEns}{\text{-}Ens}

\DeclareMathOperator{\fSet}{fEns}
\DeclareMathOperator{\Poset}{Poset }
\DeclareMathOperator{\Stab}{Stab}
\DeclareMathOperator{\Astt}{\mbox{\huge{\textasteriskcentered}}}
\DeclareMathOperator{\Ind}{Ind}
\DeclareMathOperator{\Indob}{Ind-}
\newcommand{\an}{{\mathrm{an}}}
\DeclareMathOperator{\Image}{Im}
\DeclareMathOperator{\codim}{codim}
\DeclareMathOperator{\Irr}{Irr}
\newcommand{\sq}{\square}
\newcommand{\sm}{{\mathrm{sm}}}
\newcommand{\tps}{{\mathrm{tps}}}
\newcommand{\op}{{\mathrm{op}}}
\DeclareMathOperator{\Lim}{Lim}
\DeclareMathOperator{\projLim}{\underset{\longleftarrow}{\Lim}}
\DeclareMathOperator{\injLim}{\underset{\longrightarrow}{\Lim}}
\DeclareMathOperator{\id}{id}
\DeclareMathOperator{\Tors}{Tors}
\DeclareMathOperator{\diametre}{diam}

\DeclareMathOperator{\rklog}{rk^{log}}
\DeclareMathOperator{\Ke}{\mcal{K}e}
\DeclareMathOperator{\Top}{\mcal{T}op}
\newcommand{\mring}{\mathring}
\newcommand{\triv}{\mathrm{tr}}
\newcommand{\tr}{\triv}

\newcommand{\Gm}{\mbf G_m}
\newcommand{\fGm}{\fk G_m}
\DeclareMathOperator{\OutGptop}{OutGp_{top}}
\DeclareMathOperator{\Pt}{Pt}
\newcommand{\nd}{\mathrm{nd}}
\DeclareMathOperator{\pt}{pt}
\DeclareMathOperator{\Graph}{Graph}
\DeclareMathOperator{\GenGraph}{GenGraph}
\DeclareMathOperator{\Frac}{Frac}

\DeclareMathOperator{\C}{C}
\DeclareMathOperator{\Cgeom}{C_{g\acute{e}om}}
\DeclareMathOperator{\Ctop}{C_{top}}
\DeclareMathOperator{\cusp}{Cusp}

\DeclareMathOperator{\PGL}{PGL}

\newcommand{\ie}{c'est-à-dire }
\newcommand{\dem}{\begin{proof}}
\newcommand{\findem}{\end{proof}}

\newcommand{\da}{\begin{displaystyle}}
\newcommand{\db}{\end{displaystyle}}
\newcommand{\dar}{\downarrow}
\newcommand{\uar}{\uparrow}

\newcommand{\mcal}{\mathcal}
\newcommand{\mbf}{\mathbf}
\newcommand{\bb}{\mathbb}
\newcommand{\mbb}{\mathbb}
\newcommand{\fk}{\mathfrak}

\newcommand\sommet{*-[o]{\bullet}}
\graphicspath{{figures/}}
\title{Géométrie anabélienne tempérée}
\author{Emmanuel Lepage}

\begin{document}
{\thispagestyle{empty}
\setlength{\topmargin}{-1.2cm} 
\setlength{\headheight}{0pt}
\setlength{\headsep}{0pt}
\setlength{\textheight}{12in}
\pagestyle{empty}
\setlength{\parindent}{0pt}
\setlength{\parskip}{\baselineskip}
\setlength{\oddsidemargin}{0in}
\setlength{\evensidemargin}{0in}
\setlength{\textwidth}{6.5in}

\begin{minipage}{25mm}
\vspace{1mm}
\end{minipage}
\begin{center}
UNIVERSITE PARIS DIDEROT (PARIS 7)\\
Ecole Doctorale : Paris Centre
    \end{center}

\vspace{1.5cm}
\begin{center}
{\bf DOCTORAT
}\\
de Mathématiques
\end{center}
\vspace{1.5cm}

\begin{center}
 {\large \bf {Emmanuel LEPAGE}}
\vspace{0.1cm}

{\Large Géométrie anabélienne tempérée}

\end{center}
\vspace{1cm}

\begin{center}
{\bf Thèse dirigée par Yves {\sc André}}

Soutenue le 04/12/2009
\end{center}
\vspace{2cm}
\begin{center}
JURY~:
\end{center}
\begin{tabbing}
\hspace{2cm}\=\hspace{8cm}\=\kill

\> M. Yves {\sc André}     \>Directeur de thèse \\
\> M. Francesco {\sc Baldassarri}    \\

\> M. Jean-Benoît {\sc Bost}          \\
\> M. Antoine {\sc Chambert-Loir}     \>Rapporteur (absent) \\
\> M. Antoine {\sc Ducros}            \\
\> M. Shinichi {\sc Mochizuki}    \>Rapporteur (absent)\\
\> Mme Leila {\sc Schneps}             \\
\> M. Jakob {\sc Stix}         \\

\end{tabbing}

}

\setcounter{page}{0}
\tableofcontents
\section*{Remerciements}
Je tiens à remercier tous ceux qui m'ont aidé dans la réalisation de cette
thèse~: Yves {\sc André} pour m'avoir guidé dans ce
travail de recherche et pour les multiples relectures qu'il a fait de mes
différents travaux~; Antoine {\sc Chambert-Loir} et
Shinichi {\sc Mochizuki} pour avoir accepté la tâche ingrate de rapporteur
et pour leurs commentaires et remarques~;
Framcesco {\sc Baldassarri}, Jean-Benoît {\sc Bost}, Antoine {\sc Ducros},
Leila {\sc Schneps} et
Jakob {\sc Stix} qui ont accepté de faire partie du jury~; Luc {\sc Illusie} pour son
intérêt pour mes questions sur le groupe fondamental logarithmique~;
Fumiharu {\sc Kato} pour m'avoir invité pendant deux mois à l'Université de Kyoto
(et la JSPS pour avoir financé le séjour) et Akio {\sc Tamagawa} pour les
discussions que j'ai eues avec lui lors de ce séjour.

\chapter*{Introduction}

\section*{Probl\'ematique}
\subsection*{Groupes fondamentaux en g\'eom\'etrie analytique $p$-adique~:
  groupe fondamental temp\'er\'e}
Pour un espace topologique connexe localement contractile, le groupe fondamental de Poincar\'e
classifie les revêtements de cet espace.\\
A. Grothendieck a développé une théorie similaire en géométrie algébrique~:
il associe à un schéma connexe un groupe profini qui classifie les revêtements
étales finis~(\cite[V.7]{sga}).\\
Dans le cas d'une variété algébrique complexe, le groupe fondamental
algébrique de Grothendieck s'identifie au complété profini du groupe fondamental
topologique de Poincar\'e de la variété analytique associée~(\cite[cor. XII.5.2]{sga}).\\

On voudrait une théorie analogue en $p$-adique. Plus précisément, on
voudrait associer à une variété analytique $p$-adique un "groupe fondamental" dont le complété
profini soit le groupe fondamental algébrique si la vari\'et\'e est
alg\'ebrique, mais qui mette aussi en \'evidence les diff\'erentes
\emph{uniformisations} dans la g\'eom\'etrie analytique rigide. Supposons momentan\'ement
pour simplifier que le corps de base est $\mbf C_p$. J. Tate a
montr\'e que, en termes de sa g\'eom\'etrie rigide, une courbe
elliptique $E$ ayant mauvaise r\'eduction (une telle courbe est appel\'ee
une courbe de Tate\index{Courbe de Tate})
peut \^etre d\'ecrite sous la forme
$\Gm/q^{\mbf Z}$ avec $|q|<1$. Nous voudrions que le groupe fondamental
fasse appara\^itre $\Gm$ comme un rev\^etement de $E$. Cette uniformisation peut se g\'en\'eraliser
aux courbes de Mumford. Par d\'efinition, une courbe de Mumford est une courbe dont toutes les composantes
irr\'eductibles de
la r\'eduction stable sont rationnelles\index{Courbes de Mumford}. Ceci \'equivaut \`a dire que sa
jacobienne a r\'eduction multiplicative. Si $X$ est une courbe de Mumford,
alors il existe un sous-groupe $\Gamma$
discret de type fini  de $\PGL_2(\mbf C_p)$ dont tous les \'el\'ements sont
hyperboliques et tel que $X=\Omega/\Gamma$, o\`u $\Omega$ est le
compl\'ementaire dans $\mbf P^1$ des points d'accumulation d'orbites par
l'action de $\Gamma$. On voudrait que $\Omega$ soit
un rev\^etement de $X$.\\

Le cadre topologico-géométrique que nous utiliserons est celui de
V. Berkovich~(\cite{berk}). Une définition naïve de la géométrie analytique
conduit à des espaces totalement discontinus. L'approche de Berkovich
consiste à rajouter des points, correspondant localement aux semi-normes
sur des algèbres de séries convergentes. Ainsi, pour une variété sur un
corps complet non archimédien $K$, Berkovich ne considère plus seulement,
comme Tate, des points à valeurs dans une extension finie de $K$, mais dans
une extension isométrique complète quelconque de $K$. Dans le cas de la droite, on
ajoute ainsi entre autres des points correspondant à toutes les boules de la droite (la
seminorme correspondante est la norme infinie sur cette boule dans
une extension assez grosse de $K$). Ce faisant, les espaces de Berkovich sont 
localement connexes.\\
La géométrie analytique non archimédienne développée par Berkovich
dans~\cite{berk} permet d'avoir une bonne notion de revêtement
topologique~; en effet, l'espace de Berkovich associé à une variété
algébrique lisse est localement contractile~(\cite{berk2}). De plus, si $X$
est une courbe de Mumford, l'uniformisation $\Omega\to X$ est un
rev\^etement topologique. En particulier l'uniformisation $\mbf G_m\to\mbf
G_m/q^{\mbf Z}$ d'une courbe de Tate est un revêtement topologique infini. Cependant, contrairement \`a la situation complexe, les revêtements
étales finis n'induisent pas nécessairement des revêtements
topologiques. Par exemple, le revêtement kummérien $\mbf G_m\to\mbf G_m$ qui
envoie $x$ sur $x^n$ est un revêtement étale fini d'ordre $n$ mais n'est pas
un revêtement topologique car
la préimage d'un point de Berkovich correspondant à une boule de centre $0$
a un unique élément.\\
Pour combler cette lacune, on est donc obligé de considérer une plus grosse
catégorie de revêtements.\\
Dans cette optique, A.J. de Jong a développé dans~\cite{dJ1} une notion de
revêtement étale (non n\'ecessairement fini) d'un espace de Berkovich lisse. Ces revêtements étales sont
classifiés par un groupe topologique, dont le complété profini est le
groupe fondamental algébrique dans le cas d'une variété algébrique lisse.\\
Cependant, ce groupe fondamental \'etale est en général ``trop gros'' pour pouvoir être décrit
raisonnablement. Ainsi, le groupe fondamental étale de la droite projective
n'est m\^eme pas prodiscret (en tant que groupe topologique), alors que, dans le cas complexe, le groupe
fondamental topologique est trivial.\\
On pourrait aussi consid\'erer la cat\'egorie des faisceaux localement constants sur le site \'etale de Berkovich. Cette cat\'egorie est une sous-cat\'egorie pleine de la cat\'egorie des rev\^etements \'etales. Le groupe fondamental associ\'e est le compl\'et\'e prodiscret du groupe fondamental \'etale. Mais il est encore trop difficile \`a d\'ecrire.\\

Pour obtenir un groupe plus maniable, Y. André a défini une notion de
\emph{revêtement tempéré}. Un revêtement tempéré est un revêtement étale
(au sens de de Jong) qui devient un revêtement topologique après pullback
par un revêtement étale fini. Ces revêtements tempérés sont classifiés par
un groupe topologique, le \emph{groupe fondamental tempéré} dont le
complété profini est encore le groupe fondamental algébrique. Mais tout
revêtement tempéré est dominé par le revêtement topologique universel d'un
revêtement étale fini galoisien (ce qui assure que le groupe fondamental
tempéré est un groupe prodiscret). On est donc assez proche de la situation complexe. Dans le cas d'une courbe, le groupe
fondamental tempéré est résiduellement fini, ce qui permet de borner
sa taille par celle du groupe fondamental algébrique.\\

\subsection*{Applications du groupe fondamental temp\'er\'e}
Les groupes fondamentaux $p$-adiques, temp\'er\'es ou \'etales au sens de
de Jong, ont aussi pour but de donner un analogue $p$-adique (au moins
partiel) de la correspondance de Riemann-Hilbert complexe entre
repr\'esentations du groupe fondamental et fibr\'es vectoriels munis d'une
connexion int\'egrable.  Pour l'analogue $p$-adique, il faut cependant
rajouter une hypoth\`ese sur la constance locale du faisceau des sections
horizontales de la connexion~(voir~\cite[th. 4.2]{dJ1} et~\cite[th.~III.3.4.6, \S{}
III.3.5]{andre1} pour des \'enonc\'es pr\'ecis).\\

Les groupes fondamentaux temp\'er\'es apparaissent aussi dans la recherche
d'un analogue $p$-adique \`a la th\'eorie des groupes triangulaires de
Schwarz de la g\'eom\'etrie complexe.  
Les groupes triangulaires $p$-adiques sont d\'efinis dans~\cite[def.
III.5.2.]{andre1} comme les images (quand elles sont discr\`etes) de
l'action de monodromie du groupe fondamental \'etale sur $\PGL_2(\mbf C_p)$
d\'efini par certaines \'equations diff\'erentielles hyperg\'eom\'etriques quand
le faisceau de leurs solutions est localement
constant pour la topologie \'etale. Contrairement au cas complexe, le
faisceau des solutions n'est a priori pas n\'ecessairement localement
constant et le groupe de monodromie n'est a priori pas n\'ecessairement
discret. La question se pose donc de savoir quand les groupes
triangulaires existent.\\
On peut, de façon similaire au groupe fondamental temp\'er\'e, d\'efinir des groupes fondamentaux temp\'er\'e d'orbifolds~(\cite[§ III.4.4]{andre1}).
Les groupes fondamentaux temp\'er\'es de l'orbifold $\mbf P_1$ avec
ramification en trois points donnent une interpr\'etation des groupes
triangulaires  $p$-adiques (\cite[prop. III.5.2.4]{andre1}). Les groupes
triangulaires de type Mumford correspondent \`a de tels orbifolds qui
admettent une uniformisation par une courbe de
Mumford. Dans~\cite{katotriangular}, F. Kato classifie ces groupes
triangulaires de type Mumford. Andr\'e montre qu'essentiellement tous les groupes
triangulaires proviennent d'orbifolds de type Mumford (\cite[th. III.5.3.7]{andre1}).\\
 P.E. Bradley se sert de ces groupes triangulaires pour d\'emontrer que tout groupe fini appara\^it comme groupe de Galois d'un rev\^etement ramifi\'e de $\mbf P^1$ par une courbe de Mumford de genre $\geq 2$~(\cite[th. 6]{bradley}).\\

Le groupe fondamental tempéré (tout comme le groupe fondamental
algébrique) est défini pour des corps de base non n\'ecessairement
alg\'ebriquement clos. On a alors une suite fondamentale de type usuel
reliant groupes fondamentaux temp\'er\'es ``arithm\'etique'' et ``g\'eom\'etrique''.  Ainsi, si $X$ est une variété sur $\mbf Q_p$, le
groupe fondamental tempéré de $X_{\mbf C_p}$ est muni d'une action
extérieure du groupe de Galois $\Gal(\mbf C_p/\mbf Q_p)$. Ces actions
galoisiennes donnent au groupe fondamental tempéré un intérêt
arithmétique.\\
Ainsi, le groupe fondamental tempéré permet de construire un analogue $p$-adique de la théorie de Grothendieck-Teichm\"uller. On a un morphisme du groupe de Galois $G_{\mbf Q}$ de $\mbf Q$ dans le groupe $\Out\ga(\mbf P^1\backslash\{0,1,\infty\})$ des automorphismes ext\'erieurs du groupe fondamental profini g\'eom\'etrique de $\mbf P^1\backslash\{0,1,\infty\}$ induit par l'action de $G_{\mbf Q}$ sur $\mbf P^1_{\overline{\mbf Q}}\backslash\{0,1,\infty\}$. Ce morphisme est injectif~(\cite{belyi}). On peut alors caract\'eriser l'image du groupe de Galois $G_{\mbf Q_p}$ de $\mbf Q_p$ en termes du groupe fondamental temp\'er\'e de $\mbf P^1_{\mbf C_p}\backslash\{0,1,\infty\}$~(\cite[th. 7.2.1]{andre2}). Plus pr\'ecis\'ement, $G_{\mbf Q_p}$ est l'intersection de $G_{\mbf Q}$ et de $\Out\ga(\mbf P_{\mbf C_p}^1\backslash\{0,1,\infty\})$. Ceci permet, en consid\'erant des groupes d'automorphismes des groupes fondamentaux temp\'er\'es des espaces de modules de courbes de type $(0,n)$, de d\'efinir une variante $p$-adique du groupe de Grothendieck-Teichm\"uller~(\cite[\S 8]{andre2}).\\

La g\'eom\'etrie anab\'elienne s'int\'eresse \`a ce que l'on peut retrouver
de vari\'et\'es \`a partir de groupes fondamentaux. Cette \'etude a \'et\'e
lanc\'ee par Grothendieck, qui a conjectur\'e que l'on pouvait reconstuire
une courbe hyperbolique sur un corps de nombres $K$ \`a partir de son groupe
fondamental profini g\'eom\'etrique muni de l'action ext\'erieure du groupe
de Galois $G_K$. S. Mochizuki a r\'esolu ce problème et m\^eme le problème
analogue sur une extension finie $K$ de $\mbf Q_p$
(\cite[th. A]{mochi4}). Dans cette optique anab\'elienne, on peut se demander ce qu'on
peut reconstruire d'une courbe hyperbolique sur $\mbf C_p$ \`a partir de
son groupe fondamental temp\'er\'e, malgr\'e l'absence de l'action galoisienne.
 Mochizuki s'est int\'eress\'e \`a cette question et a prouvé qu'on pouvait
reconstruire le graphe de sa réduction stable à partir de son groupe
fondamental tempéré~(\cite{mochi}). Les sommets du graphe correspondent aux sous-groupes compacts maximaux d'un quotient caract\'eristique du groupe fondamental temp\'er\'e (une variante $(p')$ du groupe), et les ar\^etes aux intersections non triviales de tels sous-groupes compacts maximaux.

\section*{Contenu de la th\`ese}
\subsection*{Vers une description du groupe fondamental temp\'er\'e}

Le groupe fondamental temp\'er\'e est en g\'en\'eral tr\`es difficile \`a
calculer explicitement. De par sa d\'efinition m\^eme, le groupe
fondamental temp\'er\'e d\'epend non seulement de la structure topologique de la vari\'et\'e (et donc de la structure combinatoire d'une r\'eduction), mais aussi de celle de tous ses rev\^etements \'etales finis. Sauf dans de tr\`es rares cas (comme les vari\'et\'es ab\'eliennes~; cf. \S{} \ref{varab}), on ne sait pas en donner une description explicite. Il n'est en g\'en\'eral m\^eme pas localement compact (par exemple pour $\mbf P^1\backslash\{0,1,\infty\}$~; cf. \cite[prop. III.2.3.12]{andre1}).\\

\subsubsection*{Propri\'et\'es g\'en\'erales analogues \`a celles des groupes fondamentaux usuels}
Cependant, je d\'emontre pour le groupe fondamental tempéré un certain
nombre de propriétés g\'en\'erales analogues à celles bien connues pour le groupe fondamental
topologique ou pour le groupe fondamental algébrique. Je d\'emontre par exemple que le
groupe fondamental tempéré (sur le corps de base $\mbf C_p$) est~:
\begin{itemize}
\item \emph{invariant par un morphisme
birationnel de variétés algébriques propres et lisses}
(proposition~\ref{birat})~;
\item \emph{compatible
aux produits} (``formule de Künneth''~; proposition~\ref{prod})~;
\item \emph{invariant par
extension isométrique algébriquement close du corps de base}
(proposition~\ref{invariance}).
\end{itemize}
La preuve de chacun de ces r\'esultats repose sur l'\'etude du comportement topologique (pour la topologie de Berkovich) des rev\^etements \'etales finis. L'invariance birationnelle utilise simplement le fait que le groupe fondamental topologique d'une vari\'et\'e lisse est le m\^eme pour tous les ouverts de Zariski denses. Les deux autres r\'esultats reposent sur la description du type d'homotopie d'un espace de Berkovich muni d'une r\'eduction semistable en termes de la combinatoire de cette r\'eduction. On se ram\`ene au cas d'espaces ayant r\'eduction semistable en utilisant les th\'eor\`emes d'alt\'eration de J. de Jong.\\
Je d\'emontre \'egalement un \emph{th\'eor\`eme de type Van Kampen} (les rev\^etements sont locaux pour la topologie \'etale), mais uniquement \emph{pour les courbes} (corollaire~\ref{courbestempet}). Il serait int\'eressant de savoir si ce r\'esultat est encore vrai en dimension sup\'erieure.\\
J'\'etablis aussi un \emph{isomorphisme entre l'abélianisé du groupe
fondamental d'une courbe et le groupe fondamental tempéré de sa jacobienne}~(théorème~\ref{abel}). Dans le cas d'une \emph{vari\'et\'e ab\'elienne}, je donne une \emph{description explicite du groupe fondamental temp\'er\'e} (\S{}~\ref{varab}).

\subsubsection*{Variantes $(p')$}
La difficult\'e de d\'ecrire le comportement combinatoire de la
r\'eduction des rev\^etements \'etales finis diminue sensiblement si
l'on se restreint aux  rev\^etements pro-$(p')$ (\ie d'ordre premier \`a la
caract\'eristique r\'esiduelle $p$). Pour cette raison, S. Mochizuki a
introduit, dans le cas des courbes, une variante $(p')$ du groupe fondamental
temp\'er\'e qui classifie les rev\^etements qui deviennent topologique
apr\`es changement de base par un rev\^etement \'etale fini pro-$(p')$
(cette d\'efinition se g\'en\'eralise sans difficult\'e au cas
g\'en\'eral). Mochizuki d\'ecrit alors cette variante $(p')$ en termes de
la combinatoire de la r\'eduction stable de la courbe. Cette description se
fait \`a travers un graphe de groupes.\\
Nous \'etudierons plus g\'en\'eralement ce groupe fondamental temp\'er\'e
pro-$(p')$ en dimension sup\'erieure.

\paragraph*{Pont entre groupe fondamental temp\'er\'e pro-$(p')$ et g\'eom\'etrie logarithmique.}
Un outil qui nous sera indispensable pour d\'ecrire le groupe fondamental
$(p')$-temp\'er\'e d'un espace de Berkovich sera la log g\'eom\'etrie, et
plus sp\'ecifiquement le groupe fondamental logarithmique.\\
 L'int\'er\^et du groupe fondamental logarithmique dans notre contexte est qu'il permet de d\'ecrire le groupe fondamental pro-$(p')$ de la fibre g\'en\'erique en termes de la fibre sp\'eciale dans un cadre plus g\'en\'eral que celui offert par le th\'eor\`eme de sp\'ecialisation de Grothendieck. Ainsi, pour une vari\'et\'e alg\'ebrique lisse et propre \`a r\'eduction semistable, le groupe fondamental pro-$(p')$ de la vari\'et\'e est isomorphe au groupe fondamental logarithmique pro-$(p')$ de la r\'eduction semistable (munie d'une structure logarithmique). La r\'eduction semistable pr\'esente une structure combinatoire int\'eressante et Berkovich a fait le lien entre la structure combinatoire d'une telle r\'eduction semistable et la structure topologique de l'espace de Berkovich associ\'e \`a la vari\'et\'e~(\cite{berk2}). Nous prolongerons cette \'etude en faisant le lien entre la structure combinatoire d'un rev\^etement logarithmique de la r\'eduction et la structure topologique de l'espace de Berkovich du rev\^etement correspondant de la vari\'et\'e  (proposition~\ref{skelretract}).\\
Ceci nous conduira \`a \emph{d\'efinir un groupe fondamental temp\'er\'e d'une
vari\'et\'e alg\'ebrique semistable munie d'une structure logarithmique
compatible}. Nous le comparerons avec le groupe fondamental temp\'er\'e de
notre espace de Berkovich~:\\

\emph{Soit $\underline X:X\to\cdots\to \Spec(O_K)$ une log fibration polystable sur
$\Spec(O_K)$.\\
Soit $\bar x$ un point géométrique de $X^{\an}_\eta$, et soit $\tilde x$
sa log réduction.
On a un morphisme $\gtemp(X^{\an}_\eta,\bar x)^{\mbb L}\to\gtemp(X_s,\tilde
x)^{\mbb
  L}$ qui est un isomorphisme si $p\notin\mbb L$} (th\'eor\`eme~\ref{isomfondtemp}).

\paragraph*{Complexes de groupes.}
Pour donner une description combinatoire du groupe fondamental
$(p')$-temp\'er\'e g\'en\'eralisant celle de Mochizuki pour les courbes,
nous serons amen\'e \`a introduire une notion purement cat\'egorique de
\emph{complexes classifiants} (\S~\ref{cplxclassifiants}) et de groupes
fondamentaux de tels complexes classifiants. Cette notion vise principalement \`a g\'en\'eraliser la notion de complexes de groupes d\'efinis par A. H\ae fliger dans~\cite{haef}. L'id\'ee des complexes de groupes consiste en la donn\'ee suivante~:
\begin{itemize}
\item
un espace de base sous forme d'une donn\'ee combinatoire,
\item
pour chaque composante de l'espace de base, un groupe (qui, intuitivement, classifie les rev\^etements de cette composantes).\end{itemize}
L'int\'er\^et est alors de construire \`a partir de ces donn\'ees un groupe qui, intuitivement, classifie les rev\^etements de l'espace.\\
Dans la th\'eorie d'H\ae fliger, l'espace de base est un ensemble simplicial, et les groupes des composantes sont discrets.\\
Nous \'elargirons le cadre~:
\begin{itemize}\item l'espace de base ne sera plus un ensemble simplicial mais une cat\'egorie. Ceci nous permettra en particulier de passer sans difficult\'e \`a un espace de base polysimplicial, qui constitue le cadre de la th\'eorie introduite par Berkovich dans~\cite{berk2}.
\item Les groupes en les diff\'erentes composantes de l'espace ne seront plus suppos\'es discrets~(voir~\cite{stix} pour une th\'eorie des complexes de groupes profinis). Ce seront des groupes topologiques quelconques.
\end{itemize}
La diff\'erence principale entre l'\'etude ici faite et la th\'eorie de H\ae
fliger consiste en la d\'efinition des groupes fondamentaux. Nous
d\'efinirons ici le groupe fondamental d'un complexe de groupes en terme
d'une cat\'egorie de rev\^etements, comme Mochizuki dans~\cite{mochi}. A
d\'efaut de donner une description explicite du groupe fondamental, une
telle d\'efinition met plus en \'evidence le lien avec la situation
g\'eom\'etrique qui nous int\'eresse o\`u tous les groupes fondamentaux
sont d\'efinis en termes de rev\^etements. Ceci nous permettra en
particulier, en consid\'erant des sous-cat\'egories de rev\^etements, de
d\'efinir  plusieurs types de groupes fondamentaux, notamment un groupe
fondamental temp\'er\'e (et une variante $p'$).\\

Nous associerons alors \`a un log sch\'ema semistable sur un log point un
complexe classifiant, dont les groupes des composantes sont des groupes
fondamentaux logarithmiques. On construira alors un \emph{isomorphisme entre le
groupe fondamental $(p')$-temp\'er\'e du log sch\'ema semistable et le
groupe fondamental $(p')$-temp\'er\'e de son complexe classifiant}
(\S~\ref{cclogschtemp}).

\subsection*{Propri\'et\'es anab\'eliennes du groupe fondamental temp\'er\'e}
La complication du groupe fondamental temp\'er\'e traduit aussi la richesse d'information qu'il contient. Le groupe fondamental tempéré géométrique d'une variété dépend ainsi
beaucoup plus de la variété elle-même que ne le font le groupe fondamental
topologique complexe ou m\^eme le groupe fondamental algébrique. Par exemple, le
groupe fondamental tempéré géométrique d'une droite projective privée de
quatre points dépend de la position des quatre points. Plus généralement,
Mochizuki a prouvé que, pour une courbe hyperbolique, on pouvait
reconstruire le graphe de groupe pro-$(p')$ de sa réduction stable à partir de son groupe
fondamental tempéré (et m\^eme de son groupe fondamental $(p')$-temp\'er\'e~: les sommets du graphe
correspondent aux classes de conjugaison de sous-groupes compacts maximaux
du groupe fondamental $(p')$-tempéré et les arêtes aux classes de
conjugaison d'intersections non triviales de deux sous-groupes compacts
maximaux distincts~; \cite{mochi})).\\ 
On peut alors se demander à quel point ce groupe fondamental tempéré dépend
de la courbe elle-même.\\
On prouvera , \`a ce propos, le r\'esultat suivant~:\\

\emph{L'isomorphisme entre les graphes de la réduction stable de deux courbes de
Mumford induit par un isomorphisme des groupes fondamentaux tempérés de ces
deux courbes conserve les métriques naturelles sur ces graphes}
(th\'eor\`eme~\ref{mumford}).\\

La preuve de ce r\'esultat repose principalement sur l'\'etude du
comportement des rev\^etements ab\'eliens, et de leurs descriptions en
termes de fonctions th\'eta et de courants sur le graphe de la r\'eduction
semi-stable. Contrairement au r\'esultat de Mochizuki, ce th\'eor\`eme
repose essentiellement sur les rev\^etements sauvagement ramifi\'es. Il
devient faux si l'on remplace le groupe fondamental temp\'er\'e par sa
variante $(p')$, puisqu'on peut reconstruire le groupe fondamental
$(p')$-temp\'er\'e \`a partir du graphe de groupes pro-$(p')$ de sa
r\'eduction stable (dans le cas de courbes de Mumford, ce graphe de groupes pro-$(p')$ peut m\^eme
\^etre reconstruit, \`a isomorphisme pr\`es, \`a partir du seul graphe de la r\'eduction stable).\\
Il serait int\'eressant de g\'en\'eraliser ce r\'esultat \`a une courbe quelconque.

\subsection*{Variation du groupe fondamental temp\'er\'e en famille et cosp\'ecialisation}
Les propri\'et\'es anab\'eliennes pr\'ec\'edentes ont la cons\'equence suivante~: si l'on a un morphisme propre et lisse,
on ne peut pas espérer que le groupe fondamental tempéré des fibres soit
localement constant.

On peut y rem\'edier partiellement en consid\'erant uniquement le groupe
fondamental $(p')$-tempér\'e, mais il n'est pas non plus constant sur les
fibres d'un morphisme propre et lisse (en fait, pour une courbe
hyperbolique, on peut encore reconstruire le graphe de la réduction stable
à partir du groupe fondamental $(p')$-tempéré).\\
Une importante partie de cette th\`ese est d\'edi\'ee \`a comprendre comment varie le groupe fondamental temp\'er\'e  pour une famille de vari\'et\'es.\\
Notre principal r\'esultat dans cette direction sera la construction de \emph{morphismes
de cospécialisation} pour le groupe fondamental $(p')$-tempéré, dans la situation d'un
morphisme propre et lisse ayant réduction semistable
(th\'eor\`emes~\ref{cospcourbes} et~\ref{cospthm} ; l'appellation
``morphismes de cosp\'ecialisation'' fait r\'ef\'erence aux morphismes de sp\'ecialisation
de Grothendieck pour les groupes fondamentaux profinis, mais ici les morphismes vont dans l'autre sens).\\
Le lien \'etabli entre le groupe fondamental temp\'er\'e d'une vari\'et\'e
et le groupe fondamental temp\'er\'e d'une log r\'eduction semistable
nous conduira \`a une reformulation purement log g\'eom\'etrique de la
question. Plus pr\'ecis\'ement, nous serons amen\'es \`a construire des
morphismes de cosp\'ecialisation de groupes fondamentaux topologiques. Nous
aurons aussi à \'etudier l'invariance du groupe fondamental logarithmique
g\'eom\'etrique par changement de log point de base.

\section*{Questions  en suspens}
\begin{itemize}
\item
L'\'enonc\'e~\ref{courbestempet} de type Van Kampen sur les revêtements temp\'er\'es pour la topologie
\'etale au sens de Berkovich pour les courbes propres est-il encore valable en
dimension sup\'erieure~?  A d\'efaut d'un tel r\'esultat, on peut
s'int\'eresser \`a cette question pour des topologies moins
fines comme la topologie usuel de Berkovich, la topologie \'etale
alg\'ebrique, la topologie de Zariski...
\item L'\'enonc\'e de type Van Kampen~\ref{courbestempet} permet aussi de se poser la question
  suivante sur la variation du groupe fondamental temp\'er\'e des fibres
  d'une famille de courbes.\\
Si $Y\to X$ est une famille alg\'ebrique propre et lisse de courbes, on peut
  consid\'erer le champ $\mcal C$ sur la topologie \'etale (au sens de
  Berkovich) de $X^{\an}$
  d\'efini en appelant $\mcal C_U$ la cat\'egorie des faisceaux localement
  constants pour la topologie \'etale de $Y^{\an}_U$. Si $x$ est un point
  g\'eom\'etrique de $X^{\an}$ (qui d\'efinit donc un point du topos
  \'etale de $X^{\an}$), le foncteur naturel de la fibre $\mcal C_x$ vers
  la cat\'egorie des rev\^etements temp\'er\'es de $Y_x$ est-il une
  \'equivalence de cat\'egories~?\\
Notre th\'eor\`eme de Van Kampen permet essentiellement de se ramener \`a
l'\'etude du groupe fondamental topologique de voisinages d'une fibre d'une famille propre et lisse.
\item Les morphismes de cosp\'ecialisation du groupe fondamental
  $(p')$-temp\'er\'e pour une fibration $Y\to X$ \`a r\'eduction polystable
  donnent la constance locale du groupe fondamental $(p')$-temp\'er\'e sur les
  strates, tant que l'on se restreint aux points de $X^{\an}$ \`a valuation
  discr\`ete. Il est naturel de se demander si l'on peut \'etendre
  cette description aux points de
  valuation non discr\`ete.
\item Dans une optique anab\'elienne, il serait int\'eressant de g\'en\'eraliser le
  th\'eor\`eme~\ref{mumford} qui dit qu'on peut reconstruire la m\'etrique
  du graphe de la r\'eduction stable d'une courbe de Mumford \`a partir de
  son groupe fondamental temp\'er\'e. Ainsi le r\'esultat reste-t-il vrai
  si l'on ne se restreint plus aux courbes de Mumford~?
  Il serait \'egalement int\'eressant,
  m\^eme dans le cas des courbes de Mumford, de donner une description plus
  explicite des longueurs des arêtes du graphe de la r\'eduction stable en
  terme du groupe fondamental temp\'er\'e.\\
Il serait int\'eressant de d\'eterminer quelles donn\'ees
  suppl\'ementaires sur les courbes on peut reconstruire \`a partir du
  groupe fondamental temp\'er\'e. De fa\c con plus drastique, une courbe hyperbolique
  est-elle uniquement d\'etermin\'ee, \`a isomorphisme pr\`es,
  par son groupe fondamental temp\'er\'e g\'eom\'etrique~?
\item Il serait \'egalement int\'eressant de classifier les sous-groupes
  compacts maximaux du groupe fondamental temp\'er\'e d'une courbe
  hyperbolique, en \'etendant les r\'esultats de Mochizuki pour le groupe
  fondamental $(p')$-temp\'er\'e. Plus pr\'ecis\'ement, il serait
  int\'er\'essant de les relier aux groupes de d\'ecomposition des points
  de l'espace de Berkovich de la courbe. Ces groupes de d\'ecomposition
  sont bien des sous-groupes compacts du groupe fondamental temp\'er\'e et
  on peut, par des arguments similaires \`a ceux de Mochizuki dans le cas
  pro-$(p')$ montrer que tout sous-groupe compact du groupe fondamental
  temp\'er\'e est, \`a conjugaison pr\`es dans le compl\'et\'e profini,
  inclus dans un tel sous-groupe de d\'ecomposition. La question se pose de
  savoir quand ces groupes de d\'ecompositions sont \'egaux.
\end{itemize}  

\section*{Plan}
Le premier chapitre sera consacré aux préliminaires. Nous commencerons par
rappeler la définition du groupe fondamental tempéré et les résultats de
bases le concernant. Nous rappellerons également le
lien entre le graphe de la réduction stable d'une courbe hyperbolique et le
groupe fondamental tempéré.\\
Nous nous orienterons ensuite vers la log géométrie et plus spécifiquement
les groupes fondamentaux logarithmiques dont nous aurons besoin dans les
deux derniers chapitres.\\
Finalement, nous revisiterons l'article~\cite{berk2}, où Berkovich construit
une rétraction de la fibre générique (en tant qu'espace de Berkovich) d'un
schéma formel pluristable sur un fermé qui se d\'ecrit en termes de la
structure combinatoire de la réduction.\\

Le deuxième chapitre sera consacré à \'etablir, pour le groupe fondamental
tempéré, des propriétés analogues à celles connues pour le groupe
fondamental algébrique cités plus haut~: invariance birationnelle,
invariance par extension algébriquement close du corps de base, formule de
Künneth et abélianisation du groupe fondamental tempéré des courbes.\\

Le troisième chapitre sera consacré à prouver que l'isomorphisme entre les graphes
des réductions semi-stables de deux  courbes de Mumford induit par un
isomorphisme de leurs groupes fondamentaux tempérés conserve la métrique de
ces graphes. Nous commencerons par l'étude des cas particuliers d'une
droite épointée et d'une courbe elliptique épointée.\\

La quatri\`eme partie, purement cat\'egorique, sera consacr\'ee \`a l'\'etude des complexes classifiants et de leurs groupes fondamentaux.\\  

La cinqui\`eme partie consistera à définir un groupe fondamental tempéré pour
un schéma pluristable muni d'une structure logarithmique compatible. Nous
comparerons ensuite pour un schéma pluristable sur l'anneau des entiers
$O_K$ le groupe fondamental tempéré de sa fibre géométrique et celui de sa
fibre spéciale. Nous donnerons également une description en termes de complexes classifiants
de ce groupe fondamental tempéré d'un log schéma pluristable.\\

La sixi\`eme partie sera consacrée à la définition de morphismes de
cospécialisation pour une famille propre et lisse avec réduction
pluristable.\\

\chapter{Préliminaires}
\section{Espaces de Berkovich}
Rappelons ici la d\'efinition des espaces de Berkovich, tels qu'ils sont d\'efinis dans~\cite{berketale}.

\subsection{Spectre d'un anneau de Banach}
Un \emph{anneau de Banach}\index{Anneau de Banach} $A$ est un anneau
(commutatif unitaire) muni d'une norme sous-multiplicative $\|\ \|$ pour
laquelle $A$ est complet.\\ 
Soit $(A,\|\ \|)$ un anneau de Banach.
Le \emph{spectre}\index{Spectre} $\mcal M(A)$ de $A$ est l'ensemble des
seminormes $|\ |:A\to\mbf R_+$ multiplicatives et born\'ees (\ie pour tout
$f,g\in A$, $|fg|=|f||g|$ et $|f|\leqslant\|f\|$) muni de la topologie la
plus faible rendant continues les applications $|\ |\mapsto|f|$. $\mcal
M(A)$ est alors un espace topologique s\'epar\'e compact et non vide
(\cite[th. 1.2.1]{berk}).\\ 
Si $x$ est un point de $\mcal M(A)$ correspondant \`a une seminorme $|\ |$,
le noyau $\fk p_x$ de $|\ |$ est un id\'eal premier ferm\'e de $A$. $|\ |$
d\'efinit une valuation sur le corps des fractions de $A/\fk p_x$. Notons
$\mcal H(x)$ le compl\'et\'e de ce corps pour cette valuation~: $\mcal
H(x)$ est un corps valu\'e complet, appel\'e le \emph{corps
  r\'esiduel}\index{Corps r\'esiduel} de $x$.\\ 
Tout morphisme born\'e $A\to B$ d'anneaux de Banach induit une application continue $\mcal M(B)\to\mcal M(A)$.

\subsection{Espaces affino\"ides}
Soit $K$ un corps complet non archim\'edien. Soient $O_K$ son anneau des
entiers, $\fk m_K$ l'id\'eal maximal de $O_K$ et $k$ son corps
r\'esiduel.\\ 
Pour $r=(r_1,\dots,r_n)\in (\mbf R^*_+)^n$, notons 
\[K\{r^{-1}T\}:=\{f=\sum_\nu a_\nu T^\nu|a_\nu \in K,|a_\nu|r^\nu\stackrel{|\nu|\to\infty}{\to}0\}\]
o\`u $\nu=(\nu_1,\dots,\nu_n)\in \mbf N^n$, $|\nu|=\nu_1+\cdots +\nu_n$,
$T^\nu=T_1^{\nu_1}\cdots T_n^{\nu_n}$, et $r^\nu=r_1^{\nu_1}\dots
r_n^{\nu_n}$.\\ 
C'est une $K$-alg\`ebre de Banach pour la norme multiplicative
$\|f\|=\max_\nu |a_\nu|r^\nu$. Pour $r=(1,\dots,1)$, cette alg\`ebre est
l'alg\`ebre de Tate\index{Alg\`ebre de Tate} $K\{T_1,\dots,T_n\}$. 
\begin{dfn}[{\cite[def. 2.1.1]{berk}}]
Une $K$-alg\`ebre de Banach $A$ est \emph{$K$-affino\"ide}\index{Alg\`ebre
  affino\"ide} si il existe un \'epimorphisme $K\{r^{-1}T\}\to A$ pour un
certain $r$ qui est \emph{admissible}\index{Morphisme!admissible} (\ie la
norme de $A$ est \'equivalente \`a la norme quotient). Si l'on peut prendre
$r=(1,\dots,1)$, $A$ est dite strictement $K$-affino\"ide. 
\end{dfn}
Une alg\`ebre affino\"ide est n\oe th\'erienne et tous ses id\'eaux sont ferm\'es.
\begin{dfn}[{\cite[def. 2.2.1]{berk}}]
Un ferm\'e V de $\mcal M(A)$ est un \emph{domaine
  affino\"ide}\index{Domaine!affino"ide} si il existe un morphisme born\'e
de $K$-alg\`ebres affino\"ides $A\to A_V$ satisfaisant la propri\'et\'e
universelle suivante~: pour tout morphisme born\'e d'alg\`ebres
$K$-affino\"ides  $A\to B$ tel que l'image de $\mcal B$ dans $\mcal A$ soit
contenue dans $V$, il existe un unique morphisme born\'e $A_V\to B$ rendant
commutatif le diagramme suivant~: 
\[\xymatrix{A\ar[r]\ar[dr] & A_V\ar[d] \\ & B}\] 
\end{dfn}
L'intersection de deux domaines affino\"ides est encore un domaine
affino\"ide ($A_{U\cap V}=A_U\otimes_AA_V$).\\ 
Si $V$ est un domaine affino\"ide $\mcal M(A_V)=V$ et donc $A\to A_V$ d\'etermine $V$.\\ 
Par exemple, si $g, f_1,\dots, f_n\in A$ n'ont pas de z\'eros en commun
dans $X:=\mcal M(A)$ et si $p=(p_1,\dots,p_n)\in (\mbf R_+^*)^n$, alors le
ferm\'e 
\[X(p^{-1}\frac{f}{g}):=\{x\in X||f_i(x)|\leqslant p_i|g(x)|\}\]
est un domaine affino\"ide repr\'esent\'e par
\[A\to A\{p^{-1}\frac{f}{g}\}:=A\{p^{-1}T\}/(gT_i-f_i).\]
Un domaine affino\"ide de cette forme est appel\'e \emph{domaine rationnel}\index{Domaine!rationnel}.\\
\begin{thm}[Th\'eor\`eme d'acyclicit\'e de Tate] 
Si $\mcal V=\{V_i\}_{i\in I}$ est un recouvrement fini de $X:=\mcal M(A)$
par des domaines affino\"ides et si $M$ est un $A$-module, le complexe de
$\breve{C}ech$  
\[0\to M\to \prod_i M\otimes_AA_{V_i}\to\prod_{i,j}M\otimes_AA_{V_i\cap V_j}\to\cdots\]
est exact.
\end{thm}
Un sous-ensemble $V$ de $X:=\mcal M(A)$ est dit sp\'ecial si c'est l'union
$V=\bigcup V_i$ d'un nombre fini de domaines affino\"ides $V_i$. On pose
alors $A_V=\Ker(\prod_i A_{V_i}\to\prod_{i,j}A_{V_i\cap V_j})$. $A_V$ ne
d\'epend pas du choix du recouvrement $(V_i)_{i\in I}$.\\ 
Si $U$ est un ouvert de $X$, on note \[O_X(U)=\varprojlim A_V,\] o\`u $V$
d\'ecrit l'ensemble des sous-ensembles sp\'eciaux contenus dans $U$. $O_X$
est alors un faisceau d'anneau. $O_{X,x}$ est un anneau local n\oe
th\'erien. L'espace localement annel\'e ainsi obtenu est appel\'e
\emph{espace $K$-affino\"ide}\index{Espace affino\"ide}. 
La cat\'egorie des espaces affino\"ides est la cat\'egorie oppos\'ee de la
cat\'egorie des alg\`ebres $K$-affino\"ides.\\ 
Un morphisme $\mcal M(B)\to\mcal M(A)$ est
\emph{fini}\index{Morphisme!fini} (resp. une \emph{immersion
  ferm\'ee}\index{Immersion ferm\'ee}) si $B$ est un $A$-module de Banach
fini (resp. $A\to B$ est un \'epimorphisme admissible).\\  

D\'efinissons la \emph{norme spectrale}\index{Norme spectrale} $\rho$ sur
$A$ par \[\rho(f):=\lim \sqrt[n]{\|f^n\|}=\max_{x\in\mcal M(A)} |f|_x.\] 
$A^\circ:=\{f\in A|\rho(f)\leqslant 1\}$ est un sous-anneau de $A$ et
$A^{\circ\circ}:=\{f\in A|\rho(f)<1\}$ en est un id\'eal. Notons
$\widetilde{A}:=A^\circ/A^{\circ\circ}$.\\ 
Il existe une application $\pi:\mcal M(A)\to\Spec(\widetilde A)$.\\
SI $A$ est strictement affino\"ide, $\pi$ est surjective et tout point
g\'en\'erique de $\Spec(\widetilde A)$ a une unique pr\'eimage par
$\pi$~(\cite[prop. 2.4.4]{berk}). 

\subsection{Espaces analytiques}
\begin{dfn}
Soit $X$ un espace topologiques. Une famille $\tau$ de sous-ensembles de $X$ est
\begin{itemize}
\item un quasifilet de $X$ si pour tout $x\in X$, il existe
  $V_1,\dots,V_n\in\tau$ tels que $x\in V_1\cap\cdots\cap V_n$ et
  $V_1\cup\cdots\cup V_n$ est un voisinage de $x$, 
\item un \emph{filet}\index{Filet} de $X$ si c'est un quasifilet de $X$ et
  si pour tout $U,V\in \tau$, $\tau_{|U\cap V}$ est un quasifilet de $U\cap
  V$. 
\end{itemize}
\end{dfn}
\begin{dfn}[{\cite[def. 1.2.3]{berketale}}]
Un \emph{espace $K$-analytique}\index{Espace analytique} est donn\'e par~:
\begin{itemize}
\item un espace topologique $|X|$ localement s\'epar\'e,
\item un filet $\tau$ sur $|X|$ de sous-ensembles compacts de $|X|$,
\item pour tout $V\in \tau$, une alg\`ebre $K$-affino\"ide $A_V$ et un
  hom\'eomorphisme $V\simeq \mcal M(A_V)$, 
\item pour tous $U,V\in\tau$ tels que $U\subset V$, un morphisme born\'e
  $A_V\to A_U$ qui identifie $(U,A_U)$ avec un domaine affino\"ide de
  $(V,A_V)$. 
\end{itemize}
\end{dfn}
Un morphisme fort d'espaces $K$-analytiques $\phi:X\to X'$ est une
application continue $\phi:|X|\to |X'|$, telle que pour tout $V\in \tau$ il
existe $V'\in \tau'$ avec $\phi(V)\subset V'$, et un syst\`eme compatible
de morphismes d'espaces $K$-affino\"ides $\phi_{V/V'}:(V, A_V)\to
(V',A_{V'})$ pour tout $V\in\tau$, $V'\in \tau'$ tels que $\phi(V)\subset
V'$. On notera $K-\An$ la cat\'egorie des espaces $K$-analytiques pour
laquelle les morphismes sont les morphismes forts.\\ 
Un morphisme fort $\phi: X\to X'$ est un quasi-isomorphisme si c'est un
hom\'eomorphisme et si pour tous $V\in\tau$, $V'\in\tau'$ tels que
$\phi(V)\subset V'$, $\phi_{V/V'}$ identifie $V$ \`a un domaine analytique
de $V'$.\\ 
Le syst\`eme des quasi-isomorphismes v\'erifie le calcul des fractions \`a
droite.\\ 
La cat\'egorie des espaces $K$-analytiques est par d\'efinition la
cat\'egorie $K\text{-}\An$ des fractions de $K\text{-}\widetilde{\An}$ pour
le syst\`eme des quasi-isomorphismes.\\ 
Un espace $K$-analytique $X$ admet un unique filet maximal $\hat \tau$
contenant $\tau$ et des alg\`ebres $\widehat A_U$ pour $U\in \hat \tau$ qui
d\'efinisse un espace analytique et tel que $(|X|,\tau,(A_U)_U)\to
(|X|,\hat \tau, (\widehat A_U)_U)$ soit un isomorphisme. Les \'el\'ements
de $\tau$ sont appel\'es \emph{domaines affino\"ides}\index{Domaine
  affino\"ide} de  $X$. Un sous-ensemble de $X$ est un \emph{domaine
  sp\'ecial}\index{Domaine sp\'ecial} s'il en existe un recouvrement fini
$(U_i)_{i\in I}$ par des domaines affino\"ides tels que $A_{U_i}\otimes
A_{U_j}\to A_{U_I\cap U_j} $ soit admissible. Un sous-ensemble $Y$ de $X$
est un domaine analytique si la restriction \`a $Y$ du filet des domaines
affino\"ides est encore un filet. Si $Y$ est un domaine analytique de $X$
la restriction du filet des domaines affino\"ides de $X$ d\'efinit une
structure d'espace $K$-analytique sur $Y$.\\ 
On peut construire des espaces $K$-analytiques par recollement d'espaces
$K$-analytiques le long d'ouverts, mais aussi par recollement d'un nombre
fini (ou localement fini) d'espaces $K$-analytiques le long de domaines
analytiques ferm\'es~(\cite[prop. 1.3.3]{berketale}).\\ 
La $G$-topologie $X_G$ sur $X$ est une topologie de Grothendieck
engendr\'ee par la pr\'etopologie sur la cat\'egorie des domaines
analytiques de $X$, o\`u les recouvrements d'un domaine analytique $Y$ sont
les familles $(Y_i)_I$ de domaines analytiques de $Y$ qui forment un
quasifilet sur $Y$. Si $X'$ est un espace $K$-analytique, le pr\'efaisceau
qui \`a $Y$ associe $\Hom(Y,X')$ est un faisceau pour la $G$-topologie
$X_G$ de $X$ (et donc aussi pour la topologie de $|X|$). Pour $X'=\mbf
A^1$, on obtient le faisceau structural $O_{X_G}$ de $X_G$. Le
\emph{faisceau structural} 
\index{Faisceau structural} de $X$ est le pouss\'e en avant
$O_X:=\pi_*O_{X_G}$ de $O_{X_G}$ le long de $\pi:X_G\to |X|$.\\  

Si $x$ est un point d'un espace analytique $X$, il est inclus dans un domaine affino\"ide. Le corps non archim\'edien $\mcal H(x)$ ne d\'epend pas du choix d'un tel affino\"ide. On l'appelle encore corps r\'esiduel de $x$.\\
$O_{X,x}$ est un anneau local, son corps r\'esiduel est naturellement
valu\'e et $\mcal H(x)$ est le compl\'et\'e de ce corps r\'esiduel.\\ 

Dans~\cite{berk}, Berkovich d\'efinissait une notion plus restrictive
d'espaces $K$-analytiques, qu'il appelle dans~\cite{berketale} bon espace
$K$-analytique. Ce sont les espaces $K$-analytiques dont tout point admet
un voisinage affino\"ide.\\ 

La cat\'egorie des espaces K-analytiques admet des produits
fibr\'es~(\cite[prop. 1.4.1]{berketale}). Si $L$ est une extension
isom\'etrique de $K$, il existe aussi un foncteur naturel de changement de
base $\text{-}\otimes_KL:K\text{-}\An\to L\text{-}\An$ et une application
continue fonctorielle $|X\otimes_KL|\to |X|$ (si $X=\mcal M(A)$ est
affino\"ide, $X\otimes_KL=\mcal M(A\hat{\otimes}_KL)$).\\ 
Soient $\phi:Y\to X$ un morphisme d'espaces $K$-analytiques et $x$ un point
de $X$. $x$ induit un morphisme naturel $\mcal M(\mcal H(x))\to
X\otimes_K\mcal H(x)$ d'espaces $\mcal H(x)$-analytiques. L'espace $\mcal
H(x)$-analytique $Y_x:=(Y\otimes_K\mcal H(x))\times_{X\otimes_K\mcal
  H(x)}\mcal M(\mcal H(x))$ est la \emph{fibre}\index{Fibre} de $\phi$ en
$x$. L'application $Y_x\to Y$ induit un hom\'eomorphisme
$Y_x\to\phi^{-1}(x)$. 

\subsection{Espace analytique associ\'e \`a une vari\'et\'e alg\'ebrique}
Soit $X$ un sch\'ema de type fini sur $K$. 
Soit $F$ le foncteur de la cat\'egorie oppos\'ee \`a celle des bons espaces
$K$-analytiques dans la cat\'egorie des ensembles qui \`a $Z$ associe
l'ensemble $\Hom(Z,X)$ des morphismes d'espaces localement annel\'es
au-dessus de $\Spec K$. Ce foncteur est
repr\'esentable~(\cite[th. 3.4.1]{berk}). On note $X^{\an}$ le bon espace
analytique correspondant.\\ 

D\'ecrivons plus pr\'ecis\'ement $X^{\an}$.\\
Si $X=\mbf A^n$, $X^{\an}$ est l'union croissante des boules $D(0,r)=\mcal
M(K\{r^{-1}T_1,\dots,r^{-1}T_n \})$ de multirayon $(r,\dots,r)$.\\ 
Si $X\to Y$ est une immersion ouverte, notons $\phi$ le morphisme canonique
d'espaces localement annel\'es $Y^{\an}\to Y$. Alors $X^{\an}$ est juste
l'ouvert $\phi^{-1}(X)$.\\ 
Si $X\to Y$ est une immersion ferm\'ee d\'efinie par un id\'eal $\mcal J$
de $O_Y$, alors $\mcal JO_{Y^{\an}}$ d\'efinit un sous-espace analytique
ferm\'e de $Y^{\an}$ qui n'est autre que $X^{\an}$.\\ 
Si $(X_i)_{i\in I}$ est un recouvrement ouvert de $X$, $X^{\an}$ est obtenu
par recollement des $X^{\an}_{i}$ le long des $(X_i\cap X_j)^{\an}$.\\

\subsection{Espace analytique associ\'e \`a un sch\'ema formel}
Un sch\'ema formel sur $\Spf O_K$ est localement topologiquement de
pr\'esentation finie si il est recouvert par des sch\'emas formels de la
forme $\Spf O_K\{T_1,\dots,T_n\}/I$ o\`u $O_K\{T_1,\dots,T_n\}$ est
l'anneau des s\'eries enti\`eres $\sum a_\nu T^\nu$ avec $a_\nu\in O_K$ et
$a_\nu\to 0$ quand $|\nu|\to \infty$ et $I$ est un id\'eal de type fini de
$O_K\{T_1,\dots,T_n\}$.\\ 
A un sch\'ema formel localement topologiquement de pr\'esentation finie
$\fk X$ sur $\Spf O_K$, on associe un espace $K$-analytique $\fk X_\eta$,
appel\'e sa fibre g\'en\'erique. Il existe alors une application
fonctorielle $\pi:\fk X_\eta\to\fk X_s$, appel\'ee \emph{application de
  r\'eduction}\index{R\'eduction} (cette application n'est pas
continue). Rappelons-en bri\`evement la
construction~(\cite[sect. 1]{berkvc1}).\\ 

Si $\fk X=Spf(A)$, $\mcal A=A\otimes_{O_K}K$ est une alg\`ebre affino\"ide
et $\fk X_\eta:=\mcal M(\mcal A)$. On a un morphisme $A\to\mcal A^\circ$
qui induit un morphisme $A\otimes k\to \widetilde{\mcal A}$. L'application
de r\'eduction est alors donn\'ee par la compos\'ee $\fk X_\eta=\mcal
M(A)\to\Spec(\widetilde{\mcal A})\to\Spec(A\times O_K)=\fk X_s$.\\ 
Si $\fk Y$ est un ouvert de $\fk X$, $\fk Y_\eta:=\pi^{-1}(\fk Y)$ est un
domaine analytique ferm\'e de $\fk X_{\eta}$.\\ 
Si $(\fk X_i)$ est un recouvrement localement fini de $\fk X$, $\fk X_\eta$
est obtenu par recollement des $\fk X_{i,\eta}$ le long des $(\fk
X_{i}\cap\fk X_j)_\eta$ (il est \`a remarquer que le recollement se fait le
long de domaine analytique ferm\'e et non d'ouvert comme dans le cas de
l'espace analytique associ\'e \`a une vari\'et\'e alg\'ebrique, et donc
$\fk X_\eta$ n'a pas de raison a priori d'\^etre un bon espace
analytique).\\ 

Si $X$ est un sch\'ema localement de pr\'esentation fini sur $O_K$, son compl\'et\'e formel $\fk X$ le long de la fibre sp\'eciale est un sch\'ema formel localement topologiquement de type fini.
Alors, $\fk X_\eta$ s'identifie naturellement \`a un domaine analytique
ferm\'e de $X^{\an}$. L'immersion $\fk X_\eta\to X^{\an}$ est un
isomorphisme si $X\to O_K$ est propre~(\cite[sect. 5]{berkvc1}).  

\subsection{Morphismes \'etales et morphismes lisses}
Un morphisme $\phi:Y\to X$ est \emph{quasifini} si pour tout point $y$ de
$Y$, il existe des voisinages ouverts $V$ de $y$ et $U$ de $\phi(y)$
contenant $\phi(V)$ tel que $\phi:V\to U$ soit fini (\ie pour tout domaine
affino\"ide $U_0$ de $U$, $V\cap\phi^{-1}(U_0)\to U_0$ est un morphisme
fini d'espaces $K$-affino\"ides).\\ 
Un morphisme $Y\to X$ est plat si pour tout $y\in Y$, $O_{Y,y}$ est une
$O_{X,\phi(y)}$-alg\`ebre plate.\\ 
Un morphisme $\phi:Y\to X$ est $G$-localement une immersion ferm\'ee si il
existe un quasifilet $\tau$ sur $Y$ de domaines analytiques et pour tout
$V\in\tau$ un domaine analytique $U$ de $X$ contenant l'image de $V$ tel
que $\phi_{V/U}:V\to U$ soit une immersion ferm\'ee (\ie pour tout domaine
affino\"ide $U_0$ de $U$, $V\cap\phi^{-1}(U_0)\to U_0$ est une immersion
ferm\'ee d'espaces $K$-affino\"ides). Si $X\to Y$ est $G$-localement une
immersion ferm\'ee, et $V\to U$ sont comme pr\'ec\'edemment, soit $\mcal J$
l'id\'eal de $O_{U_G}$ correspondant \`a $V$. $\mcal J/\mcal J^2$ est un
$O_{V_G}$-module, qui se recollent en un $O_{Y_G}$-module, appel\'e
\emph{faisceau conormal} de $\phi$.\\ 
Si $\phi:Y\to X$ est maintenant un morphisme d'espaces $K$-analytiques, le
morphisme diagonal $\Delta : Y\to Y\times_XY$ est $G$-localement une
immersion ferm\'ee, et son faisceau conormal est appel\'e \emph{faisceau
  des diff\'erentielles}\index{Faisceau des diff\'erentielles} de $\phi$,
et not\'e $\Omega_{Y_G/X_G}$.\\  
\begin{rem}
Dans les d\'efinitions pr\'ec\'edentes, quand $X$ et $Y$ sont de bons
espaces $K$-analytiques (ce qui sera en fait quasiment toujours le cas dans
cette th\`ese), on peut faire les constructions pr\'ec\'edentes pour la
topologie usuelle de $Y$, obtenant ainsi un faisceau $\Omega_{Y/X}$ sur
$Y$. Par pullback le long de $Y_G\to Y$, on retrouve bien
$\omega_{Y_G/X_G}$. 
\end{rem}
Un morphisme $Y\to X$ est \emph{\'etale}\index{Morphisme!d'espaces de
  Berkovich!\'etale} si il est quasifini, plat et $\Omega_{Y_G/X_G}$ est
nul.\\ 
La cat\'egorie des germes de rev\^etements \'etales finis en un point $x$
est \'equivalente \`a la cat\'egorie des rev\^etements \'etales finis de
$\mcal H(x)$.\\ 
Ceci permet \`a Berkovich de d\'efinir en~\cite[sect. 4.1]{berketale} une
\emph{topologie \'etale}\index{Topologie!\'etale} sur un espace
$K$-analytique $X$~: c'est la topologie de Grothendieck sur la cat\'egorie
$\Et(X)$ des espaces $K$-analytiques \'etales sur $X$ o\`u les
recouvrements d'un objets $U$ sont les familles $(f_i:U_i\to U)$ telles que
$U=\bigcup_i f_i(U_i)$. On note $X_{\et}$ le topos associ\'e.\\ 

Berkovich d\'efinit \'egalement dans~\cite[sect. 3]{berkvc1} une topologie
\emph{quasi-\'etale}\index{Topologie!quasi-\'etale} qui pour la topologie
\'etale joue le role que la $G$-topologie joue pour la topologie
usuelle. On dit qu'un morphisme $\phi:Y\to X$ est \emph{quasi-\'etale} si
l'ensemble des domaines affino\"ides $V$ de $Y$ tels que le morphisme induit
$V\to X$ puisse se factoriser en $V\to U\to X$ où $U\to X$ est étale et $V\to
U$ est l'immersion d'un domaine affinoïde forme un quasifilet. La topologie
quasi-\'etale sur $X$ est alors la topologie de Grothendieck sur la
cat\'egorie $\Qet(X)$ des morphismes quasi-\'etales $U\to X$ engendr\'ee
par la pr\'etopologie o\`u une famille $(f_i:U_i\to U)_{i\in I}$ est un
recouvrement si $(f_i(U_i))_{i\in I}$ est un quasifilet.\\ 
Si $\fk Y\to \fk X$ est un morphisme \'etale de sch\'emas formels
localement topologiquement de pr\'esentation finie sur $\Spf O_K$, $\fk
Y_\eta\to\fk X_\eta$ est quasi-\'etale~(\cite[prop. 2.3]{berkvc1}). On en
d\'eduit un morphisme de topos ${\fk X_{\eta}}_\qet\to {\fk X_s}_\et$.\\ 

Un morphisme $Y\to X$ est \emph{lisse}\index{Morphisme!d'espaces de
  Berkovich!lisse} si localement sur $Y$ (pour la topologie usuelle), il
peut se factoriser en $Y\to \mbf A^d\times X\to X$ avec $Y\to\mbf A^d\times
X$ \'etale.

\section{Groupe fondamental tempéré}

\subsection{Définition}

Soit $K$ un corps complet non archimédien.\\
Soit $\mbb L$ un ensemble de nombres premiers (dans le cas de l'ensemble des nombres premiers autres que la caractéristique résiduelle $p$ de
$K$, nous \'ecrirons souvent $(p')$ au lieu de $\mbb L$). Un $\mbb L$-entier est un entier qui s'écrit comme produit d'éléments
de $\mbb L$.\\
Par commodit\'e, nous traduirons provisoirement par "$K$-vari\'et\'e" ce qui est appel\'e "$K$-manifold" dans~\cite[§4]{andre2}, \ie une $K$-\emph{variété} est un
$K$-espace strictement analytique lisse et paracompact au sens de
Berkovich. Par exemple, si $X$ est un $K$-schéma lisse, $X^{\an}$ est une
$K$-variété (en pratique, nous nous intéresserons principalement à ces
espaces). D'après~\cite{berk2}, toute $K$-variété est localement
contractile (nous expliquerons plus en détail les résultats de~\cite{berk2}
dans la section~\ref{berkspaces}). En particulier, elle a un revêtement universel.\\
On dit qu'un morphisme $f:S'\to S$ de $K$-variétés est un \emph{revêtement
  étale}\index{Rev\^etement!\'etale} si $S$ est recouvert par des ouverts $U$ tel que
$f^{-1}(U)=\coprod V_j$ et $V_j\to U$ est étale fini~(\cite{dJ1}).\\
Par exemple, les revêtements $\mbb L$-finis (\ie les revêtements étales finis
dont l'ordre d'une clôture galoisienne est un $\mbb L$-entier), que nous
appellerons aussi \emph{revêtements $\mbb
  L$-algébriques} (dans le cas o\`u les vari\'et\'es proviennent de vari\'et\'es alg\'ebriques), et les revêtements au sens topologique usuel pour la
topologie de Berkovich, que nous appellerons aussi \emph{revêtements
  topologiques}\index{Rev\^etement!topologique}, sont des 
revêtements étales.\\
André définit alors la notion de revêtement tempéré dans
\cite[def. 2.1.1]{andre1}. Nous généralisons ici la définition aux revêtements
$\mbb L$-tempérés (quand $\mbb L$ est l'ensemble de tous les nombres
premiers, on retrouve la définition de revêtements tempérés)~:
\begin{dfn} \label{def:rvt:temp}
Un revêtement étale $S' \to S$ est \emph{$\mbb L$-tempéré}\index{Rev\^etement!temp\'er\'e} si c'est le quotient
du composé d'un revêtement topologique $T' \to T$ et d'un revêtement
étale $\mbb L$-fini
 $T \to S$.
\end{dfn}
Ici, un quotient $T'\to S'$ est simplement un morphisme surjectif de revêtements étales.\\
Ceci \'equivaut à dire que $S'\to S$ devient un revêtement topologique après
pullback par un revêtement étale $\mbb L$-fini.\\
Nous noterons $\Covtemp(X)^{\mbb L}$ (resp. $\Covalg(X)^{\mbb L}$,
$\Covtop(X)$) la catégorie des revêtements $\mbb L$-tempérés (resp. $\mbb
L$-algebriques, topologiques) de $X$ (muni des morphismes évidents ).\\

Un \emph{point géométrique}\index{Point!g\'eom\'etrique} d'une $K$-variété $X$ est un morphisme d'espaces de
Berkovich $\mcal M(\Omega)\to X$ o\`u $\Omega$ est une extension isométrique
s\'eparablement close de $K$.\\
Soit $\bar x$ un point géométrique de $X$. On a un foncteur \[F^{\mbb L}_{\bar
x}:\Covtemp(X)^{\mbb L}\to\Set\] qui envoie un revêtement $S\to X$ sur l'ensemble
$S_{\bar x}$. Si $\bar x$ et $\bar x'$ sont deux points géométriques d'une
même composante connexe de $X$, alors
$F^{\mbb L}_{\bar x}$ et $F^{\mbb L}_{\bar x'}$ sont (non canoniquement) isomorphes (\cite[prop. 2.9]{dJ1}).\\
Le groupe fondamental tempéré\index{Groupe fondamental!tempéré} de $X$ pointé en $\bar x$ est par définition
\[\gtemp(X,\bar x)^{\mbb L}=\Aut F^{\mbb L}_{\bar x}.\]
Si $X$ est un $K$-schéma lisse, on notera simplement 
$\Covtemp(X)^{\mbb L}$ et $\gtemp(X,\bar x)^{\mbb L}$ à la place de $\Covtemp(X^{\an})^{\mbb L}$
et $\gtemp(X^{\an},\bar x)^{\mbb L}$.\\
En considérant les stabilisateurs $(\Stab_{F(S)}(s))_{S\in \Covtemp(X)^{\mbb
L},s\in F_{\bar x}(S)}$ comme une base de voisinages ouverts de $\gtemp(X,\bar
x)^{\mbb L}$, $\gtemp(X,\bar x)^{\mbb L}$ devient un groupe
topologique. C'est un groupe topologique prodiscret.\\
\begin{rem}
Si $X$ est une variété algébrique, le foncteur
$\Covalg(X^{\an})^{\mbb L}\to\Covalg(X)^{\mbb L}$ est une équivalence de
catégories dans chacun des
cas suivants~:
\begin{itemize}
\item K est de caractéristique nulle~:
\item $X\to K$ est propre~;
\item la caractéristique de $K$ n'est pas dans $\mbb L$.
\end{itemize}
Si $X$ est algébrique, $K$ est de caractéristique 0 et a seulement un
nombre fini d'extensions d'ordre fixé dans une clôture algébrique $\overline K$,
$\gtemp(X,\bar x)^{\mbb L}$ a un système dénombrable de voisinages de $1$
et tout quotient discret est finiment engendré~(\cite[prop. 2.1.7]{andre1}).\\
\end{rem}

Si $\bar x$ et $\bar x'$ sont deux points géométriques, alors $F^{\mbb L}_{\bar
x}$ et $F^{\mbb L}_{\bar x'}$ sont (non canoniquement) isomorphes
(\cite[prop. 2.9]{dJ1}). Comme d'habitude, le groupe fondamental tempéré
dépend du point-base uniquement à automorphisme intérieur près (ce groupe,
considéré à conjugaison près, est alors simplement noté
$\gtemp(X)^{\mbb L}$).\\
La sous-catégorie pleine des revêtements $\mbb L$-tempérés $S$ pour lesquels  $F^{\mbb L}_{\bar
x}(S)$ est $\mbb L$-fini est équivalente à $\Covalg(X)^{\mbb L}$, d'où l'isomorphisme
\[\big(\gtemp(X,\bar x)^{\mbb L}\big)^{\mbb L}=\ga(X,\bar x)^{\mbb L}\] (où
 $(\ )^{\mbb L}$ note le complété pro-$\mbb L$).\\
Pour tout morphisme $X\to Y$, le changement de base définit un foncteur
$\Covtemp(Y)^{\mbb L}\to\Covtemp(X)^{\mbb L}$. Si $\bar x$ est un point
géométrique de $X$ d'image $\bar y$ dans $Y$, on obtient un morphisme continu \[\gtemp(X,\bar x)^{\mbb L}\to\gtemp(Y,\bar y)^{\mbb L}\]
(d'où un morphisme "extérieur" $\gtemp(X)^{\mbb L}\to\gtemp(Y)^{\mbb L}$).\\
On a un analogue de la correspondance de Galois usuelle~:
\begin{thm}[{\cite[th. 1.4.5]{andre1}}]  \label{galcorr} $F^{\mbb
L}_{\bar x}$ induit une équivalence de cat\'egories entre la cat\'egorie des
sommes directes de revêtements $\mbb L$-tempérés de $X$ et la catégorie
$\gtemp(X,\bar x)^{\mbb L}\tSet$ des ensembles discrets munis d'une action
continue à gauche de $\gtemp(X,\bar x)^{\mbb L}$.\end{thm}
Une somme directe de revêtements tempérés n'est pas nécessairement un
revêtement tempéré. Ainsi $f:S=\coprod_{n\in\mbf N^*}\mbf G_m\to\mbf G_m$, où la
$n$\ieme{} copie de $\mbf G_m$ s'envoie sur $\mbf G_m$ par $x\mapsto x^n$,
n'est pas un revêtement tempéré (les degrés $[\mcal H(s):\mcal H(f(s))]$ ne
sont pas bornés quand $s$ décrit $S$).  \\

Si $S$ est un revêtement galoisien $\mbb L$-fini de $X$, son revêtement
topologique universel $S^{\infty}$ est encore galoisien et tout revêtement $\mbb
L$-tempéré est dominé par un tel revêtement $\mbb L$-tempéré galoisien.\\
Si $((S_i,\bar s_i))_{i\in I}$ est un système projectif cofinal de
revêtements $\mbb L$-finis galoisiens géométriquement pointés de $(X,\bar
x)$, soit $((S^{\infty}_i,\bar s^{\infty}_i))_{i\in I}$ le système
projectif, muni des morphismes $f_{ij}^\infty$ pour $i\geq j$, des
revêtements topologiques pointés universels. Alors $F^{\mbb L}_{\bar
x}(S^{\infty}_i)=\gtemp(X,\bar x)^{\mbb L}/\Stab_{F(S^\infty_i)}(\bar
s^\infty_i)$ est naturellement un groupe quotient $G$ de $\gtemp(X,\bar x)^{\mbb
L}$  pour lequel $s^{\infty}_i$ est l'élément neutre. De plus $G$ agit par
$G$-automorphismes sur $F^{\mbb L}_{\bar x}(S^{\infty}_i)$ par translation
à droite (et donc aussi sur $S^{\infty}_i$ grâce à la correspondance
galoisienne~(théorème \ref{galcorr})). On obtient un morphisme
$\gtemp(X,\bar x)^{\mbb L}\to\Gal(S^{\infty}_i/X)$. Comme
$f_{ij}^\infty(s_i^\infty)=s_j^\infty$, ces morphismes sont compatibles avec $\Gal(S^{\infty}_i/X)\to \Gal(S^{\infty}_j/X)$.\\
Alors, grâce à~\cite[lem. III.2.1.5]{andre1},
\begin{prop}\label{limproj} \[\gtemp(X,\bar x)^{\mbb L}\to\varprojlim
  \Gal(S^{\infty}_i/X)\] est un isomorphisme.\end{prop}
Donnons maintenant une description de la catégorie des revêtements tempérés qui fait
uniquement intervenir la catégorie des revêtements topologiques de tous les
revêtements étale finis.\\

On dispose d'une catégorie fibrée
$\Dtop(X)\to\Covalg(X)$, où la fibre $\Dtop(X)_S$ en un revêtement
algébrique $S$ de $X$ est $\Covtop(X)$.\\
 $\Dtop(X)$ est un préchamp sur $\Covalg(X)$ muni de sa topologie canonique
 au sens de~\cite[déf. II.1.2.1]{giraud}.\\
Soit $\Dtemp(X)\to\Covalg(X)$ la catégorie fibrée dont la fibre en $S$ est
la catégorie $\Covtemp(U)$ des revêtements tempérés de $U$. $\Dtemp(X)$
est un champ. Le foncteur cartésien pleinement fidèle de préchamps
$\Dtop(X)\to\Dtemp(X)$ induit un foncteur cartésien pleinement fidèle de
champs $\Dtop(X)^a\to\Dtemp(X)$  où $\Dtop(X)^a$ est le champ associé à $\Dtop(X)$. Puisqu'un revêtement tempéré est topologique
localement sur
$\Covalg(X)$, ce foncteur est en fait une équivalence~(\cite[th. II.2.1.3]{giraud}).\\
De même le champ $(\Dtop(X)|_{\Covalg(X)^{\mbb
    L}})^a$ est naturellement équivalent au champ $\Dtemp(X)^{\mbb L}$ des revêtements $\mbb L$-tempérés
sur $\Covalg(X)^{\mbb L}$.\\
 
Comme les revêtements galoisiens sont cofinaux parmi les recouvrements de
l'objet final de $\Covalg(X)^{\mbb L}$, on a la description suivante.\\
 
Soit $S$ un revêtement  $\mbb L$-fini de $X$. Posons 
$\DDtemp_S$ la catégorie des données de descente dans la catégorie fibrée
$\Dtop(X)$ par rapport à $S\to
X$. $\DDtemp_S$ est naturellement équivalente à la sous-catégorie pleine de
$\Covtemp(X)$ des revêtements tempérés qui deviennent topologiques après le
changement de base $S\to X$. \\
Si $``\varprojlim"\ S_i$ est un pro-revêtement $\mbb L$-fini universel
$(X,x)$, on obtient une équivalence naturelle \[\Covtemp(X)^{\mbb L}=\projLim_i \DDtemp_{S_i}.\]

Nous utiliserons également les deux résultats suivants~:
\begin{prop}\label{andre114}\emph{(\cite[prop. III.1.1.4]{andre1})}
Soient $\overline S$ une $K$-variété et $S$ un ouvert de Zariski
dense. Alors tout revêtement topologique de $S$ s'étend
de façon unique en un revêtement topologique de $\overline S$. Ainsi
$\gtop(S,s)\to\gtop(\overline S,s)$ est un isomorphisme.\end{prop}

\begin{prop}[{\cite[th. III.2.1.11,
    prop. III.2.1.13]{andre1}}] \label{andre2111}Supposons $K$ alg\'ebriquement clos de
  caractéristique nulle. Soit
  $\overline S$ une $K$-variété et soit $S$ un ouvert de Zariski 
  dense. Alors le foncteur de la catégorie des revêtements tempérés
  $\overline S$ vers la catégorie des revêtements tempérés de $S$ est
  pleinement fidèle. Si $Z:=\overline S\backslash S$ est de codimension $\geq 2$,
  ce foncteur est une équivalence de catégories.
\end{prop}
On en déduit immédiatement un résultat analogue pour les revêtements $\mbb
L$-tempérés.\\

\subsection{Résultats de Mochizuki sur le groupe fondamental
  $(p')$-tempéré d'une courbe}\label{mochipart}

S. Mochizuki décrit dans~\cite{mochi} le groupe fondamental $(p')$-tempéré d'une
courbe algébrique $X$ en terme du graphe de la réduction stable de $X$. Il
prouve en particulier qu'on peut reconstruire le graphe de la réduction
stable de $X$ à partir de $\gtemp(X,\bar x)$.\\
Expliquons ici certains r\'esultats majeurs de son article.\\

Un \emph{semigraphe}\index{Semigraphe} $\mbb G$ est donné
par un ensemble $\mcal V$ de "sommets", un ensemble $\mcal E$ d'"arêtes" et,
pour tout $e\in\mcal E$, un ensemble $\mcal B_e$ de
"branches"\index{Branche (d'un semigraphe)} ayant au plus $2$ éléments et muni d'une fonction $\zeta_e:\mcal B_e\to\mcal V$~(\cite[Appendix]{mochi3}).\\
Nous dirons qu'une branche $b$ de $e$ \emph{aboutit} en $v$ si
$\zeta_e(b)=v$. Nous dirons qu'un semigraphe est un \emph{graphe} si pour
tout $e\in\mcal E$, $e$ a exactement deux branches.\\
Rappelons que si $\mbb G$ est un semigraphe, une structure de \emph{semigraphe
d'anabélioïdes}\index{Semigraphe!d'anabélioïdes} sur ce graphe correspond \`a la donnée suivante~:
\begin{itemize}
\item pour tout sommet ou arête $x$, une catégorie galoisienne (aussi nommée
anabélioïde connexe dans~\cite{mochi}) $\mcal G_x$;
\item pour toute branche $b$ d'une arête
$e$ aboutissant en un sommet $v$, un morphisme d'anabélioïdes (\ie un
foncteur exact) $b_*:\mcal G_e\to \mcal G_v$.\end{itemize}
Les semigraphes d'anabélio\"ides forment une 2-catégorie. Rappelons aussi que
travailler avec des catégories galoisiennes équivaut à travailler avec les
groupes profinis à automorphismes intérieurs près.\\

Si $\mcal C$ est une catégorie galoisienne de groupe fondamental $\Pi$,
alors $\Ind\text{-}\mcal C$ est équivalente au topos $\Pi\text{-}\Set$.\\
Un \emph{revêtement} $S$ d'un semigraphe d'anab\'elio\"ides $\mcal G$, de semigraphe sous-jacent $\mbb G$, consiste en~:
\begin{itemize}
\item pour tout sommet $v$ de $\mbb G$, un objet $S_v$ de $\Ind\text{-}\mcal G_v$,
\item pour toute arête $e$ de branches $b_1$ et $b_2$ aboutissant en $v_1$
  et $v_2$ respectivement, un isomorphisme $\phi_{e}$ entre $b_{1*}S_{v_1}$ et $b_{2*}S_{v_2}$.\end{itemize}
On a une définition naturelle de morphisme de revêtements, d'où une
catégorie $\Bcov(\mcal G)$. Mochizuki définit aussi un 2-foncteur de la
catégorie des revêtements de $\mcal G$ vers la 2-categorie des semigraphes
d'anabélioides au-dessus de $\mcal G$.\\
Un objet de $\Bcov(\mcal G)$ est fini si tout $S_v$ est dans $\mcal G_v$,
\emph{topologique}\index{Rev\^etement!topologique!d'un graphe d'anabélioïdes}  si tout $S_v$ est un objet constant de
$\Ind\text{-}\mcal G_v$, et \emph{tempéré} si il devient topologique apr\`es
changement de base par un revêtement fini. La sous-catégorie pleine des
revêtements tempérés (resp. finis, resp. topologiques) est notée $\Btemp(\mcal G)$ (resp. $\Ba(\mcal G)$, resp. $\Btop(\mcal G)$).\\
Si $\mbb G$ est connexe, $\Ba(\mcal G)$ est une catégorie galoisienne et
son groupe fondamental est noté $\ga(\mcal G)$.\\
Si $v$ est un sommet de $\mbb G$ et $F$ est un foncteur exact et
conservatif de $\mcal G_v$ vers la cat\'egorie $\fSet$ des ensembles finis (un tel foncteur est appelé foncteur
fondamental dans~\cite[section 5]{sga}; il s'étend en un point du topos
$\Ind\text{-}G_v$, aussi noté $F$), on peut définir un foncteur
$F_{(v,F)}:\Bcov(\mcal G)\to\Set$ qui envoie $S$ sur $F(S_v)$ (si on change
le point-base $(v,F)$, les foncteurs obtenus sont isomorphes), et soit
$F^{\temp}_{(v,F)}$ sa restriction à $\Btemp(\mcal G)$. On définit alors
\[\gtemp(\mcal G,(v,F)):=\Aut(F^{\temp}_{(v,F)}).\]\\

Illustrons ces définitions en associant à une courbe un semigraphe d'anabélioïdes.\\
Supposons $K$ de valuation discrète, et soit $\overline K$ le complété
d'une clôture s\'eparable de $K$.\\
Soit $(\overline X,D)$ une courbe hyperbolique lisse $n$-pointée de genre
$g$ sur $K$, soit $X=\overline X\backslash D$. Si $K'$ est une extension
finie de $K$, un modèle semistable de $X_{K'}$ sur $O_{K'}$ est donné par un
morphisme de schémas $\overline{\mcal X}_{O_K'}\to\Spec O_{K'}$ propre et
plat dont la fibre générique est $\overline X_{K'}$ et dont la fibre
spéciale géométrique est réduite et a seulement des points doubles
ordinaires comme singularités et par un diviseur $D_{O_{K'}}$ de
$\overline{\mcal X}_{O_K'}$ prolongeant $D_{K'}$ à support dans le lieu
lisse de $\overline{\mcal X}_{O_K'}$ et tel que $D_{O_{K'}}\to \Spec
O_{K'}$ soit étale. Il existe toujours une
extension finie de $K$ pour laquelle il existe un tel modèle
semistable~(cf.~\cite{DeligneMumford}).\\
Soit $(\overline{\mcal
  X_{O_K'}},\mcal D_{K'})$ un modèle semistable sur $O_{K'}$ où
$K'$ est une extension finie de $K$. Soit $(\overline{\mcal X},\mcal D)$
son pullback à $O_{\overline K}$ et soit $\mcal X=\overline{\mcal X}\backslash\mcal D$.\\
Le semigraphe de la fibre sp\'eciale $\mcal X_s$  de $\mcal X$ est défini ainsi~: les sommets sont les
composantes irréductibles de $\mcal X_s$, les arêtes sont les points
doubles et les points marqués. Un point double $e$ a deux branches qui
aboutissent aux composantes irréductibles contenant $e$; un point marqué
$e$ a une seule branche aboutissant à la composante irréductible contenant
le point marqué.\\
Quand $\mcal X$ est le modèle stable de $X$, ce semigraphe est noté $\mbb
G_{X}^c$ (ou $\mbb G^c$ s'il n'y a pas de risque de confusion).\\
On peut munir ce semigraphe $\mbb G^c$ d'une structure de semigraphe d'anabélioïdes $\mcal G^c$.
Pour un sommet $v_i$ correspondant à une composante irréductible $C_i$ de
$\mcal X_s$, considérons l'ouvert $U_i$ du normalisé 
$\overline C'_i$ de $\overline C_i$ qui est le complémentaire des points
marqués et des préimages des points doubles de $\mcal X_s$ (les points de
$\overline C'_i-U_i$ 
correspondent exactement au branches aboutissant en $v_i$). Alors le groupe
$\Pi_{(v_i)}$ est le groupe fondamental modéré $\gt(U_i)$ de $U_i$. Le groupe d'une arête
est $\widehat{\mathbf Z(1)}^{(p')}=\varprojlim_{(n,p)=1} \mu_n (\simeq
\widehat{\mathbf Z}^{(p')})$ (comme d'habitude, l'exposant $(p')$ indique
le quotient maximal premier à $p$, où $p$ est la caractéristique résiduelle
de $O_K$), 
qui est canoniquement isomorphe au sous-groupe de monodromie de $\gt(U_i)$
d'un point de $\overline C'_i -U_i$. Le morphisme correspondant à une
branche est le plongement du groupe de monodromie du point correspondant de
$\overline C'_i-U_i$ (qui est défini à conjugaison près), 
tandis que, pour une arête avec deux branches, on identifie les deux $\widehat{\mathbf Z(1)}^{(p')}$ par $x \mapsto x^{-1}$.\\
Si $\mcal G^{c,(p')}$ est le semigraphe d'anabélioïdes obtenu à partir de
$\mcal G^c$ en remplaçant chaque groupe profini par son complété pro-$(p')$,
\[\label{revkumm} \gtemp(\mathcal G_X^{c,(p')})=\gtemp(X_{\overline K})^{(p')}\quad(\cite[\text{ex. 3.10}]{mochi}).\]
Nous g\'en\'eraliserons cette description en dimension sup\'erieure dans le chapitre~\ref{chap4}.\\
Supposons maintenant de plus $K$ de caract\'eristique mixte.
Mochizuki montre alors\footnote{Dans~\cite{mochi}, Mochizuki suppose que $K$ est une extension finie de $\mbb Q_p$, mais la preuve s'adapte pour $K$  \`a valuation discr\`ete de caract\'eristique mixte (voir~\cite[rem. 2.11.1]{mochiabstopic2}).}~:
\begin{thm}[{\cite[cor. 3.11]{mochi}}]\label{mochi} Si $X_{\alpha}$ et
  $X_\beta$ sont deux courbes, tout isomorphisme
  $\gamma:\gtemp(X_{\alpha,\overline K})\simeq\gtemp(X_{\beta,\overline
    K})$ détermine, fonctoriellement en $\gamma$ à 2-isomorphisme près, un
  isomorphisme de semigraphe d'anabélioïdes $\gamma':\mcal G^c_{X_\alpha}\simeq\mcal G^c_{X_\beta}$.\end{thm}
Plus précisément, le diagramme induit de groupes topologiques est commutatif~:
\[\xymatrix{\gtemp(\mathcal G_{X_\alpha}^{c,(p')}) \ar@{=}[d] \ar[r]^{\gtemp(\gamma'^{(p')})} & \gtemp(\mathcal G_{X_\beta}^{c,(p')})\ar@{=}[d]\\ \gtemp(X_{\alpha,\overline K})^{(p')} \ar[r]^{\gamma^{(p')}} & \gtemp(X_{\beta,\overline K})^{(p')}}.\]
Les sommets du graphe correspondent aux classes de conjugaison de
sous-groupes compacts maximaux de $\gtemp(X)^{(p')}$ (appelés sous-groupes
\emph{verticiels}\index{Verticiel (sous-groupe)} de $\gtemp(X)^{(p')}$), et les arêtes ayant deux branches
correspondent aux classes de
conjugaison d'intersections non triviales de sous-groupes verticiels~(\cite[th. 3.7]{mochi}).

\section{Groupes fondamentaux logarithmiques}\label{loggeometry}

Cette section rappelle les notions de base et quelques r\'esultats de la théorie des log schémas~(voir~\cite{kato} et \cite{ogus}) et du groupe fondamental logarithmique~(voir~\cite{ill} or~\cite{stix}).

\subsection{Log schémas}

Tous les monoïdes\index{Monoïde} considérés ici sont commutatifs et unitaires, et un
morphisme de monoïdes est supposé préserver l'élément unité. Si $P$ est un
monoïde, $P^{\gp}$ désignera le groupe associ\'e (symétrisé), $P^*$ désignera le groupe des
éléments inversibles de $P$ et $\overline
P=P/P^*$.\\
Si $a,b\in P$, nous considérerons la relation de préordre suivante : $a|b$
si et seulement si il existe $c\in P$ tel que $b=ac$.\\
Un monoïde $P$ est \emph{aigu}\index{Monoïde!aigu} (sharp) si $P^*=\{1\}$.\\ 
Un monoïde $P$ est \emph{intègre}\index{Monoïde!intègre} si le morphisme $P\to P^{\gp}$ est
injectif. Un monoïde est \emph{fin}\index{Monoïde!fin} si il est intègre et de type fini.
Un monoïde intègre $P$ est \emph{saturé}\index{Monoïde!saturé} si $a\in P^{\gp}$ est dans $P$
dès qu'il existe $n\geq 1$ tel que $a^n\in P$.\\
Si $P$ est un monoïde fin et saturé (ce qu'on abrégera en fs),
$P^{\gp}/P^*=\overline P^{\gp}$ est un groupe
abélien libre de type fini, et donc il existe une section (non canonique) $\overline
P^{\gp}\to P^{\gp}$ qui induit une décomposition $P=P^*\times \overline P$.\\
Un morphisme $f:P\to Q$ de monoïdes est \emph{local}\index{Morphisme!de monoïdes!local} si $f^{-1}(Q^*)=P^*$.
Un morphisme $P\to Q$ de monoïdes intègres est
\emph{exact}\index{Morphisme!de monoïdes!exact} si $P$ est
l'image inverse de $Q$ dans $P^{\gp}$.\\

Un \emph{idéal premier}\index{Idéal premier(d'un monoïde)} $\fk p$ d'un monoïde intègre $P$ est un sous-ensemble de
$P$ tel que si $p\in\fk p$ et $p'\in P$ alors $p\cdot p'\in\fk p$, et si $p,p'\in P$ et
$p\cdot p'\in\fk p$ alors $p\in\fk p$ ou $p'\in \fk p$.\\
Un sous-ensemble $F$ d'un monoïde intègre $P$ est une
\emph{face}\index{Face (d'un monoïde)} si
$P\backslash F$ est un
idéal premier (en particulier $F$ est un sous-monoïde de $P$).\\
$\Spec P$ désigne l'espace topologique des idéaux premiers de $P$, où $(D(f)=\{\fk
p,f\notin\fk p\})_{f\in P}$ est une base de la topologie de $\Spec P$.\\
Si $f:P\to Q$ est un morphisme de monoïdes intègres et $\fk q$ est un idéal
premier de $Q$,
alors $f^{-1}(\fk q)$ est un idéal premier de $P$, d'où une application continue $\Spec Q\to\Spec P$.\\

Une \emph{pré-log structure} sur un schéma $X$ est un couple $(M,\alpha)$
où $M$ est un
faisceau en monoïdes sur $X_{\et}$, et $\alpha:M\to (O_X,.)$ est un
morphisme de faisceaux en monoïdes, où $O_X$ est le faisceau canonique de
$X$ et $.$
est la multiplication sur $O_X$. Une pré-log structure est une \emph{log
  structure}\index{Log structure} si le morphisme induit
$\alpha^{-1}(O_X^*)\to O_X^*$ est un isomorphisme. Un \emph{log
  schéma}\index{Log schéma} $X$
est un
schéma (le \emph{schéma sous-jacent}\index{Schéma sous-jacent à un log schéma} $\mring{X}$ au log-schéma $X$) muni
d'une log structure.\\
Le foncteur banal des log structures sur $X$ vers les pre-log structures sur $X$ admet un
adjoint à gauche $(M,\alpha)\mapsto (M^a,\alpha^a)$ où $M^a$ est la
somme amalgamée de $M$ et $O_X$ le long de $\alpha^{-1}(O^*_X)$ (cette log
structure est dite \emph{associée}\index{Log structure!associée} à $(M,\alpha)$).\\
Un \emph{morphisme de log schémas}\index{Morphisme!de log schémas} $f:(X,M,\alpha)\to (Y,N,\beta)$ est un
morphisme de schémas
$f:X\to Y$ muni d'un morphisme de faisceaux en monoïdes $f^{-1}N\to M$ compatible
avec $\alpha$ et $\beta$. Alors $f^{-1}N\to M$ est nécessairement un morphisme
local de faisceaux de monoïdes.\\
Un log schéma $X$ est \emph{intègre} si pour tout point géométrique $\bar x$
de $\mring{X}$, $M_{\bar x}$ est intègre.\\
Si $Y=(Y,N,\beta)$ est un log schéma et $X\to\mring{Y}$ est un morphisme de
schémas, la log structure sur $X$ associée à $(f^{-1}N,f^{-1}\beta)$ est appelée
\emph{log structure image inverse}\index{Log structure!image inverse} et est notée $f^*N$. Un morphisme de log
schémas $f:(X,M,\alpha)\to (Y,N,\beta)$ est \emph{strict}\index{Morphisme!de log schémas!strict} si le
 morphisme induit $f^*N\to M$ est un isomorphisme (si $X$ est intègre,
cela équivaut à dire que
$\overline N_{f(\bar x)}\to\overline M_{\bar x}$ est un isomorphisme pour
tout point géométrique $\bar x$ de $X$).\\

Si $P$ est un monoïde, la log structure sur $\Spec \mbf Z[P]$ associée à la
pré-log structure d\'efinie par $P\to\mbf Z[P]$ est appelée log structure
canonique. Il y a une application canonique $\Spec \mbf Z[P]\to\Spec P$ qui
envoie un idéal premier $I$ de $\mbf Z[P]$ sur $I\cap P$.\\
Une \emph{carte}\index{Carte!d'un log schéma} (globale) modelée sur un monoïde $P$ d'un log schéma $X$
est un morphisme de faisceaux en monoïdes
du faisceau constant $P_X$ vers $M$ induisant un isomorphisme sur les log
structures associées. Cela revient à se donner un morphisme strict
$X\to\Spec \mbf Z[P]$, où $\mbf Z[P]$ est muni de sa log
structure canonique.\\
Si $\bar x$ est un point géométrique, une carte $P_X\to M$ modelée sur un
monoïde intègre $P$ est \emph{exacte en $\bar
  x$}\index{Carte!exacte} (resp. \emph{bonne en $\bar
  x$})\index{Carte!bonne} si le morphisme induit $\overline P\to \overline
M_{\bar x}$ (resp. $P\to\overline M_{\bar x}$) est un isomorphisme.\\
Un log schéma est \emph{fin}\index{Log schéma!fin} (resp. \emph{fin et
  saturé}\index{Log schéma!fs} ou fs pour
abréger) si il est intègre et, localement pour la topologie étale, il admet
une carte modelée sur un monoïde fin (resp. fs).\\
Nous travaillerons principalement avec des log schémas fs. La catégorie des
log schémas fs admet des produits fibrés, mais les produits fibrés ne
commutent en général pas avec le foncteur d'oubli de la log structure, qui
va de la catégorie des log schémas fs vers celle des schémas (cf.~\cite[\S{}
II.2.4]{ogus} pour une discussion détaillée).\\
Si $X\to\Spec P$ est une carte fs et $\overline x$ est un point géométrique de
$X$ qui s'envoie sur $\fk p\in\Spec P$ et soit $F=P\backslash \fk p$, alors
$P\to M_{\bar x}$ induit un isomorphisme $\overline{F^{-1}P}\to\overline
M_{\bar x}$. De plus $\Spec \mbf Z[F^{-1}P]\to\Spec\mbf Z[P]$ est une
immersion ouverte, correspondant à la préimage de $D(\fk p)=\{\fk p'|\fk p\subset
\fk p'\}\subset \Spec P$. Ainsi $F^{-1}P$ induit une carte exacte d'un
voisinage ouvert pour la topologie de Zariski de $\overline x$. Si l'on
choisit une décomposition $F^{-1}P=\overline{F^{-1}P}\oplus(F^{-1}P)^*$, le
morphisme induit $\overline{F^{-1}P}\to F^{-1}P\to M_X$ est une bonne carte en
$\overline x$.\\

Parfois, nous aurons besoin de log structure sur le site de Zariski. Soit
$\epsilon:X_{\Zar}\to X_{\et}$ la projection naturelle. Nous dirons qu'une
log structure $M$ sur $X$ est \emph{zariskienne} (et le log schéma $X$ est
\emph{log zariskien}\index{Log schéma!log zariskien}) si $\epsilon^*\epsilon_*M\to M$ est un
isomorphisme. Si $X$ est un log schéma fs, la log structure est zariskienne si
et seulement si $X$ admet des cartes localement pour la topologie de
Zariski. En particulier tout log schéma fs est log zariskien localement pour
la topologie étale.\\

Si $f:X\to Y$ est un morphisme de log schémas fins, une
\emph{carte}\index{Carte!d'un morphisme de log schémas} de $f$ est
donnée par une carte $X\to\Spec \mbf Z[P]$, une carte $Y\to\Spec \mbf Z[Q]$
et un morphisme $Q\to P$ tel que le carré correspondant de log schémas soit
commutatif.
Tout morphisme de log schémas fins admet des cartes localement pour la
topologie étale.\\

Un morphisme de log schémas fins $f:X\to Y$ est \emph{log
  lisse}\index{Morphisme!de log schémas!log lisse}
(resp. \emph{log étale}\index{Morphisme!de log schémas!log étale}) étale localement sur $X$ et $Y$, $f$ admet une
carte $u:Q\to P$ telle que le noyau et la torsion du conoyau
(resp. le noyau et le conoyau) de $u^{\gp}$ soient des groupes finis d'ordre
inversible sur $X$ et telle que $X\to Y\times_{\Spec \mbf Z[Q]}\Spec \mbf Z[P]$ soit
étale.\\
Il existe une caractérisation valuative des morphismes log lisses et log
étales. Les morphismes log étales et log lisses sont stables par changement
de base et par composition.

Un morphisme $h:Q\to P$ de monoïdes fs est \emph{kumm\'erien}
\index{Morphisme!de monoïdes!kumm\'erien}(resp. $\mbb
L$-kumm\'erien)
si $h$ est injectif
et pour tout $a\in P$, il existe un entier (resp. un $\mbb L$-entier) $n$
tel que $na\in h(Q)$
(remarquons que si $Q\to P$ est kumm\'erien, $\Spec P\to\Spec Q$ est un homéomorphisme).\\
Un morphisme $f:X\to Y$ de log schémas fs est
\emph{kumm\'erien}\index{Morphisme!de log schémas!kumm\'erien}
(resp. \emph{exact}\index{Morphisme!de log schémas!exact})
si pour tout point géométrique $\bar x$ de $X$, $\overline M_{Y,f(\bar x)}\to \overline
M_{X,\bar x}$ est kumm\'erien (resp. exact).\\
Un morphisme kummérien $X\to Y$ de log schémas fs est un \emph{hom\'eomorphisme kumm\'erien
   universel}\index{Morphisme!de log schémas!kuh} (ou \emph{kuh} pour abréger)
si le morphisme de schémas sous-jacents reste un homéomorphisme après
n'importe quel changement de base de log schémas fs (\cite[def. 2.1]{vidal}).\\
Un morphisme Kummer $q:X\to Y$ est kuh si et seulement si $\mring q$ est un
homéomorphisme universel (\ie est entier, radiciel et surjectif ;
cf.~\cite[cor. 18.12.11]{ega4}) et pour
tout point géométrique $\bar x$ de $X$, $\overline M_{Y,q(\bar x)}=\overline
M_{X,\bar x}$ est $p$-Kummer, où $p$ est la caractéristique résiduelle de
$\bar x$ (\cite[thm. 2.7]{vidal}).\\
Par exemple, si $P\to Q$ est un morphisme $p$-Kummer de monoïdes fs, et $A$
est une $\mbf F_p$-algèbre,
alors $\Spec A[Q]\to\Spec A[P]$ est kuh.\\

Un log schéma fs $X$ est \emph{log régulier}\index{Morphisme!de log schémas!log régulier} si pour tout point géométrique 
$\bar x$, $O_{X,\bar x}/I_{\bar x}O_{X,\bar x}$ est régulier et
$\dim(O_{X,\bar x})=\dim(O_{X,\bar x}/I_{\bar x}O_{X,\bar x})+\rk(\overline
M_{\bar x}^{\gp})$ où $I_{\bar x}$ est l'idéal de $O_{X,\bar x}$
engendré par l'image de $M_{\bar x}\backslash O_{X,\bar x}^*$
(\cite[def. 2.2]{niziol}).\\
Les log schémas log Zariski log réguliers sont étudiés dans~\cite{kato2}.\\
Si $X$ est log régulier, $\mring X$ est normal.\\
 La partie triviale d'un log schéma fs
$X_{\tr}=\{x\in X,\overline M_{\bar x}=\{1\}\}$ (cette propriété ne dépend
pas du choix du point géométrique $\bar x$ au-dessus de $x$ car ils sont
tous isomorphes) de $X$ est un sous-sch\'ema ouvert de
$X$. Si $X$ est log régulier, on peut retrouver la log structure à partir
de la partie trivial grâce à l'isomorphisme suivant~:
 \[M=O_X\cap j_*O_{X_{\tr}}^*\] où $j$ d\'esigne l'immersion ouverte
 $X_{\tr}\to X$ (\cite[prop. 2.6]{niziol}).\\
Si $Y$ est log régulier et $X\to Y$ est log lisse, alors $X$ est log
régulier (\cite[thm. 8.2]{kato2}).

\begin{exs}
Soit $X$ un sch\'ema r\'egulier et $D$ un diviseur \`a croisement normaux, alors $M=O_X\cap j^*O_{X\backslash D}$ est une log structure fs qui fait de 
de $X$ un log sch\'ema log r\'egulier.\\
Si $K$ est un corps complet pour une valuation discr\`ete, on muni $\Spec O_K$ de la log structure associ\'ee \`a $O_K\backslash\{0\}\to O_K$. $\Spec O_K$ est alors un log sch\'ema log r\'egulier.\\
Si $X\to \Spec O_K$ est un  morphisme plurinodal de sch\'emas tel que $X_K\to K$ soit lisse, alors $M=O_X\cap j^*O_{X_{\tr}}$ d\'efinit une log structure sur $X$, qui fait de $X\to \Spec O_K$ un morphisme de log sch\'emas log lisse (cf. lemme~\ref{logstrpluri} pour un r\'esultat un peu plus g\'en\'eral et une d\'emonstration).
\end{exs}

\subsection{Revêtements két}\label{revket}

Un morphisme de log schémas fs est \emph{Kummer étale}\index{Morphisme!de log schémas!két}  (ou két pour
abréger) si il est Kummer et log étale.\\
Un morphisme $f$ est két si et seulement si, localement pour la
topologie étale, il est déduit par changement de base strict et par
localisation étale d'un morphisme $\Spec \mbf Z[P]\to\Spec \mbf
Z[Q]$ induit par un morphisme $\mbb L$-Kummer $Q\to P$ de monoïdes fs avec
$\mbb L$ inversible sur $X$.\\
En fait, si $f:Y\to X$ est két, $\bar y$ est un point géométrique de $Y$, et
$P\to M_X$ est une carte étale de $X$ en $f(\bar y)$, il y a un voisinage étale
$U$ de $\bar x$ et un voisinage ouvert de Zariski $V\subset
f^{-1}(U)$ de $\bar y$ tel que $V\to U$ est isomorphe à $U\times_{\Spec
  \mbf Z[P]}\Spec\mbf Z[Q]\to U$ avec $P\to Q$ un morphisme $\mbb L$-Kummer où
$\mbb L$ est inversible sur $U$ (\cite[Prop. 3.1.4]{stix}).\\
Les morphismes két sont ouverts. Les morphismes k\'et quasi-compacts sont quasi-finis.\\

N'importe quel
morphisme de log schéma au-dessus de $X$ entre deux log schémas két sur $X$
est két. Pour une famille de morphismes két sur $X$, le fait d'être une famille
ensemblistement couvrante est stable par changement de base fs.
Ainsi, la catégorie $X_{\ket}$ des log schémas fs két au-dessus de $X$, où les recouvrements
$(T_i\to T)$ de $T$ sont les familles ensemblistement couvrante, est un site. $\widetilde X_{\ket}$
désignera le topos correspondant.\\
Tout faisceau localement constant et fini sur $\widetilde X_{\ket}$ est
représentable. Un log schéma fs két au-dessus de $X$ qui représente un tel
faisceau localement constant et fini
est un \emph{revêtement két}\index{Revêtement!két} de $X$. Nous noterons $\KCov(X)$ la catégorie
des revêtements két de $X$.\\
Un \emph{point log géométrique}\index{Point!log géométrique} est un log schéma $s$ tel que~:
\begin{itemize}
\item $\mring{s}$ est le
spectre d'un corps séparablement clos $k$;
\item $M_s$ est saturé~;
\item la multiplication
par $n$ sur $\overline
M_s$ est un isomorphisme pour tout $n$ premier à la caractéristique de $k$.\end{itemize}
Un point log géométrique de $X$ est un morphisme $x:s\to X$ de log schémas où
$s$ est un point log géométrique. Un voisinage két $U$ de $s$ dans $X$ est
un morphisme $s\to
U$ de $X$-log schemes où $U\to X$ est két. Alors si $x$ est un log point
géométrique de $X$, le foncteur $F_x$ de $X_{\ket}$ vers $\Ens$ defini par $\mcal F\mapsto
\varinjlim_U \mcal F(U)$ où $U$ parcourt la catégorie filtrante des
voisinages két de $x$ est un point du topos $X_{\ket}$ et n'importe quel
point de ce topos est isomorphe à $F_x$ pour un certain point log
géométrique et ils forment un système conservatif.\\
 On définit aussi la log stricte localisation $X(x)$ comme la limite
 inverse des voisinages két de $X$ dans la catégorie des log schémas
 saturés. Si $x$ et $y$ sont des points log géométriques de $X$, une
 \emph{spécialisation (két)}\index{Spécialisation két} de points log géométriques $x\to y$
est un morphisme $X(x)\to X(y)$ au-dessus de $X$.\\
Une spécialisation $x\to y$ induit un morphisme canonique  $F_x\to F_y$ de foncteurs.\\
Si on a une spécialisation des points topologiques sous-jacents $x\to
y$, alors il existe au moins une spécialisation $x\to y$ de points log géométriques.\\
Si $X$ est connexe, pour tout point log géométrique $x$ de $X$, $F_x$ induit
un foncteur fondamental $\KCov(X)\to\fSet$ de la catégorie galoisienne $\KCov(X)$.\\
On note $\glog(X,x)$ le groupe profini des automorphismes de ce foncteur.\\

Les morphismes de log schémas fs $f$ exacts tels que $\mring f$ soit
propre, surjectif et de présentation finie, et les morphismes  de log
schémas fs $f$ exacts tels que $\mring f$ soit surjectif, de présentation finie et
universellement ouvert, sont de descente effective pour les revêtements két (\cite[th. 3.2.25, cor. 3.2.21]{stix}).\\

Si $X$ est un log schéma fs log régulier et si $\mbb L$ est
inversible sur $X$, alors $\KCov(X)^{\mbb L}\to\Covalg(X_{\tr})^{\mbb L}$
est une équivalence de categories (\cite[th. 7.6]{ill}).\\
\begin{thm}[{\cite[cor. 2.3]{org2}}]\label{orgsp} Soit $S$ un schéma
  strictement local de point fermé $s$, et soit $X$ un log schéma fs
  connexe tel que
  $\mring X$ soit propre sur $S$. Alors
\[\KCov(X)\to\KCov(X_s)\] est une équivalence de categories.\\
Si $x$ est un point log g\'eom\'etrique de $X_s$,
\[\glog(X_s,x)\to\glog(X,x)\] est un isomorphisme\end{thm}
Si $q:X\to Y$ est kuh, $q^*:Y_{\ket}\to X_{\ket}$ est une équivalence de
catégories (\cite[th. 0.1]{vidal})

\subsection{Mophismes saturés}
L'article de référence sur le sujet est~\cite{sattsuji}, malheureusement
non publié.\\

Un morphisme de monoïdes fs $P\to Q$ est \emph{intègre}\index{Morphisme!de monoïdes!intègre} si, pour tout morphisme
de monoïdes intègres $P\to Q'$, la somme amalgamée $Q\oplus_PQ'$ (dans la
catégorie des monoïdes) est encore
intègre.\\
Un morphisme intègre de monoides saturés $P\to Q$ est
\emph{saturé}\index{Morphisme!de monoïdes!saturé} si,
pour tout morphisme de monoïdes fs $P\to Q'$, la somme amalgamée
$Q\oplus_PQ'$ (dans la catégories des monoïdes) est encore
fs.\\
Si $f:P\to Q$ est intègre (resp. saturé) et $F'$ est une face de $Q$,
$\overline{F^{-1}P}\to\overline{{F'}^{-1}Q}$ est aussi intègre (resp. saturé), où $F=f^{-1}(F')$.\\
\begin{lem}\label{satface}Si $\phi:P\to Q$ est un morphisme intègre
  (resp. saturé) de monoïdes fs et $F'$ est une face de $Q$, soit
$F=\phi^{-1}(F')$. Alors $F\to F'$ est aussi intègre (resp. saturé).\end{lem}
\dem
Pour prouver que $F\to F'$ est intègre, grâce à~\cite[prop. I.4.3.11]{ogus},
il faut et il suffit de prouver que si $f'_1,f'_2\in F'$ et $f_1,f_2\in F$
vérifient $f'_1\phi(f_1)=f'_2\phi(f_2)$, alors il existe $g'\in F'$ et $g_1,g_2\in
F$ tels que $f'_1=g'\phi(g_1)$ et $f'_2=g'\phi(g_2)$.\\
Mais il existe $g'\in Q$ et $g_1,g_2\in P$ qui satisfont les propriétés
vooulues puisque $P\to Q$ est intègre. Comme $F'$ est une face de $Q$,
$g',\phi(g_1),\phi(g_2)$ doivent être dans $F'$, et donc $g_1$ et $g_2$
sont dans $F$.\\
Un critère de T. Tsuji~(\cite[prop. 4.1]{sattsuji}) dit qu'un morphisme intègre de monoides $f:P_0\to
Q_0$ est saturé si et seulement si pour tous $a\in P_0, b\in Q_0$ et tout
nombre premier $p$ tels que
$f(a)|b^p$, il existe $c\in P_0$ tel que $a|c^p$ et $f(c)|b$. Soient
$a\in F, b\in F'$ et $p$ un nombre premier tel que $\phi(a)|b^p$. Alors puisque
$\phi:P\to Q$ est saturé, il existe $c\in P$ tel que $a|c^p$ et
$f(c)|b$. Mais $f(c)|b$ implique que $f(c)\in F'$, donc $c\in F$.
\findem

Un morphisme $f:Y\to X$ de log schémas fs est
\emph{saturé}\index{Morphisme!de log schémas!saturé} si pour tout
point géométrique $\bar y$ de $Y$, $\bar M_{X,f(\bar y)}\to\bar M_{Y,\bar y}$ est
saturé.\\
Si $Y\to X$ est saturé et $Z\to X$ est un morphisme de log schémas fs, alors
le schéma sous-jacent à $Z\times_XY$ est $\mring Z\times_{\mring X}\mring
Y$.\\
 
Si $P\to Q$ est un morphisme local et intègre (resp. saturé) de monoïdes fs
et $P$ est aigu, le morphisme $\Spec \mbf
Z[Q]\to\Spec\mbf Z[P]$ est plat (resp. séparable\index{Morphisme!de schémas!séparables}, \emph{i.e.} plat
à fibres géométriquement réduites, cf.~\cite[cor. 4.3.16]{ogus} et \cite[rem. 6.3.3]{satmorph}).\\
Soit $f:X\to Y$ un morphisme log lisse et $\bar x$ un point géométrique de 
$X$. Localement sur $Y$ pour la topologie étale, il existe une bonne carte 
$Y\to\Spec P$ en $\bar y$. Alors, grâce à~\cite[prop. A.3.1.1]{satmorph},
étale localement en $x$, il existe une carte $P\to Q$ de $Y\to X$ telle que $Y\to
\Spec \mbf Z[Q]\times_{\mbf Z[P]}X$ soit étale et $X\to\Spec Q$ est exact à
$x$. Ainsi si $f$ est intègre (resp. saturé), $P\to Q$ est un morphisme
intègre (resp. saturé) et local de monoïdes fs et $P$ est aigu. Donc
$f$ est plat (resp. séparable).\\
 
Si $P\to Q$ est un morphisme intègre de monoïdes fs, il existe un entier
$n$ tel que le changement de base fs $P_n\to Q'$ de $P\to Q$ le long
 de $P\stackrel{n}{\to}P=P_n$ soit saturé (\cite[th. A.4.2]{satmorph}).\\
De plus, si $P\to Q$ se factorise à travers $Q_0$ de manière à ce que $P\to
Q_0$ soit saturé
et $Q_0\to Q$ soit $\mbb L$-Kummer, on peut choisir pour $n$ un $\mbb L$-entier.\\
Ainsi, si $Z''\to Z'$ est un revêtement két et $Z'\to Z$ est saturé log
lisse et propre, alors pour tout log point géométrique $z$ de $Z$, il
existe un voisinage két $U$ de $z$ tel
que $Z''_U\to U$ soit saturé (et l'hypothèse de propreté peut être remplacée
par la quasicompacité de $Z''$ si $Z$ est simplement un log point fs, \ie
son schéma sous-jacent est le spectre d'un corps).

\section{Squelette d'un espace de Berkovich à réduction pluristable}\label{berkspaces}

Soit $K$ un corps complet non archimédien et soit $O_{K}$ son anneau d'entiers.\\
Si $\fk X$ est un schéma formel localement de présentation finie sur $O_K$, $\fk
X_\eta$ sera la fibre générique de $\fk X$ au sens de Berkovich
(\cite[section 1]{berkvc1}).\\ 
Rappelons la définition d'un morphisme polystable de schémas formels~:
\begin{dfn}[{\cite[def. 1.2]{berk2}, \cite[section 4.1]{berk3}}]\label{defpolystable} Soit
  $\phi:\fk Y\to \fk X$  un morphisme de schémas formels localement de
  présentation finie. 
\begin{enumerate}[(i)]
\item $\phi$ est \emph{strictement polystable} si, pour tout point
  $y\in\fk Y$, il existe un voisinage ouvert affine $\fk X'=\Spf(A)$ de
  $x:=\phi(y)$ et un voisinage ouvert $\fk Y'\subset\phi^{-1}(\fk X')$ de
  $x$ tel que le morphisme induit $\fk Y'\to \fk X'$ se factorise
  à travers un morphisme étale $\fk Y'\to\Spf(B_0)\times_{\fk X'}\cdots\times_{\fk
    X'}\Spf(B_p)$ où chaque $B_i$ est de la forme
  $A\{T_0,\cdots,T_{n_i}\}/(T_0\cdots T_{n_i}-{a_i})$ avec $a\in A$ et $n\geq
  0$.\\
  Si $\fk X$ est de type fini sur $O_K$, on dit que $\phi$ est \emph{non dégénéré} si l'on peut choisir $X'$, $Y'$
  et $(B_i,a_i)$ tels que $\{x\in(\Spf(A)_\eta)|a_i(x)=0\}$ soit nulle part
  dense. 
\item $\phi$ est \emph{polystable}\index{Morphisme!polystable} si il existe un morphisme étale surjectif
  $\fk Y'\to\fk Y$ tel que $\fk Y'\to \fk X$ soit strictement
  polystable. Si $\fk X$ est de type fini sur $O_K$, $\phi$ est alors dit \emph{non dégénéré} si l'on peut choisir $\fk
  Y'$ tel que $\fk Y'\to \fk X$ soit non dégénéré. \end{enumerate}\end{dfn} 
Une \emph{fibration polystable}\index{Fibration polystable} (resp. \emph{non dégénérée}) de longueur $l$
au-dessus de $\fk S$ est une suite de $l$ morphismes polystables (resp. non dégénérés)
$\underline{\fk X}=(\fk X_l\to\cdots\to\fk X_1\to\fk S)$.\\ 
On note $K\text{-}\mcal Pstf_l^{\et}$ la catégorie des fibrations
polystables de longueur $l$ sur $O_K$, où un morphisme $\underline{\fk X}\to
\underline{\fk Y}$ est une collection de morphismes étales $(\fk X_i\to\fk
Y_i)_{1\leq i\leq l}$ qui satisfait aux hypothèses de commutation évidentes.\\ 
On note $\mcal Pstf_l^{\et}$ la catégorie des couples $(\underline{\fk
  X},K_1)$ où $K_1$ est un corps complet non archimédien et
$\underline{\fk X}$ est une fibration polystable sur $O_{K_1}$, et un morphisme
$(\underline{\fk X},K_1)\to (\underline{\fk Y},K_2)$ est un couple
$(\phi,\psi)$ où $\phi$ est une extension isométrique $K_2\to K_1$ et
$\psi$ est un morphisme $\underline{\fk X}\to\underline{\fk
  Y}\otimes_{O_{K_2}}O_{K_1}$ dans $K_1\text{-}\mcal Pstf_l^{\et}$.\\

Soit $X$ un schéma localement de type fini sur un corps $k$.\\
Le lieu normal d'un schéma réduit, localement de type fini sur un corps, en
est un ouvert dense. Définissons
par récurrence
$X^{(0)}=X^{\red}$, $X^{(i+1)}=X^{(i)}\backslash \Norm(X^{(i)})$. Les
composantes irréductibles de $X^{(i)}\backslash X^{(i+1)}$ sont appelées
les \emph{strates}\index{Strate} de
$X$ (de rang $i$). Cela fournit une partition de $X$. L'ensemble des points
génériques des strates de $X$ est noté $\Str(X)$ (il est en bijection
naturelle avec l'ensemble des strates de $X$).\\
Berkovich definit une autre filtration $X=X_{(0)}\subset X_{(1)}\subset\cdots$
de la mani\`ere suivante~: $X_{(i+1)}$ est le sous-ensemble fermé des points contenu dans
au moins deux composantes irréductibles de $X_{(i)}$.\\
On dit que $X$ est
\emph{quasinormal}\index{Quasinormal} si, pour tout $i$, toute composante connexe de $X_{(i)}$,
muni de la structure de sous-schéma réduit, est normale (cette propriété est
locale pour la topologie de Zariski et reste vraie après composition avec
un morphisme étale). Si
$X$ est quasinormal, alors $X_{(i)}=X^{(i)}$, et $X$ est quasinormal si et seulement si
l'adhérence de toute strate est normale.\\  
L'ensemble $\Str(X)$ est naturellement ordonné~: $x\leqslant y$ si et seulement si $y\in
\overline{\{x\}}$.\\

Berkovich définit les \emph{ensembles polysimpliciaux}\index{Ensemble polysimplicial} dans~\cite[section
3]{berk2} comme suit.\\ 
Si $n$ est un entier, posons $[n]:=\{0,1,\cdots,n\}$.\\ 
Pour un $p$-uplet $\mbf n=(n_0,\cdots,n_p)$ avec soit $p=n_0=0$ soit $n_i\geq 1$
pour tout $i$, posons $[\mbf n]:=[n_0]\times\cdots\times [n_p]$
et $w(\mbf n):=p$. $[\mbf n]$ est un espace métrique pour la distance
$d((i_0,\dots,i_p),(j_0,\dots,j_p))=\sharp\{k|i_k\neq j_k\}$. \\ 
Soit $\mbf \Lambda$ la cat\'egorie dont les objets sont $[\mbf n]$
et les morphismes sont les fonctions $[\mbf m]\to [\mbf n]$ associées à un triplet 
$(J,f,\alpha)$, où~: 
\begin{itemize}
\item $J$ est un sous-ensemble de $[w(\mbf m)]$ supposé vide si $[\mbf m]=[0]$, 
\item $f$ est une fonction injective $J\to [w(\mbf n)]$,
\item $\alpha$ est une famille $\{\alpha_l\}_{0\leq l\leq p}$, où
  $\alpha_l$ est une fonction injective $[m_{f^{-1}(l)}]\to [n_l]$  si
  $l\in\Image(f)$, et $\alpha_l$ est une fonction $[0]\to [n_l]$
  sinon.\end{itemize} 
La fonction $\gamma:[\mbf m]\to [\mbf n]$ associée à $(J,f,\alpha)$ envoie alors
$\mbf j=(j_0,\cdots,j_{w(m)})\in [m]$ en $\mbf i=(i_0,\cdots,i_{w(\mbf
  n)})$ avec $i_l=\alpha_l(j_{f^{-1}(l)})$ pour $l\in\Image(f)$, et
$i_l=\alpha_l(0)$ sinon.\\ 
Un ensemble polysimplicial $\C$ est un foncteur $\mbf
\Lambda^{\op}\to\Set$. On notera $\C_{\mbf n}$ l'image de $[\mbf n]$ par ce
foncteur. Les ensembles polysimpliciaux forment une catégorie notée $\mbf \Lambda^{\circ}\Set$.\\
$\mbf \Lambda$ est une sous-catégorie pleine de $\mbf
\Lambda^{\circ}\Set$ grâce au foncteur de Yoneda. Si $\C$ est un ensemble polysimplicial,
la catégorie $\mbf \Lambda/\C$ des \emph{polysimplexes}\index{Polysimplex} de $\C$ est la catégorie dont les objets sont les morphismes $[\mbf
n]\to\C$ dans $\mbf \Lambda^{\circ}\Set$ et les morphismes de $[\mbf n]\to \C$
vers $[\mbf m]\to \C$ sont des morphismes $[\mbf n]\to [\mbf m]$ qui rendent
le triangle commutatif.\\
Un polysimplexe $x$ d'un ensemble polysimplicial $\C$ est
\emph{dégénéré}\index{Polysimplex!dégénéré} si
il existe un morphisme surjectif $f$ de $\mbf \Lambda$, qui n'est pas un
isomorphisme, tel que $x$ soit l'image par $f$ d'un polysimplexe de
$\C$. Soit $\C_{\mbf n}^{\nd}$ le sous-ensembles des polysimplexes non
dégénérés de $\C_{\mbf n}$.\\
Grâce à un analogue du lemme d'Eilenberg-Zilber pour les ensembles polysimpliciaux
(\cite[lem. 3.2]{berk2}), un morphisme $\C'\to \C$ est bijectif si et
seulement si il envoie les polysimplexes non dégénérés sur les polysimplexes
non dégénérés et
$(\C')_{\mbf n}^{\nd}\to\C^{\nd}_{\mbf n}$ est bijectif pour tout $\mbf n$.\\
Il y a un foncteur \[O:\mbf \Lambda^{\circ}\Set\to\Poset\] où $O(\C)$ est
l'ensemble ordonné associé à $\Ob(\mbf\Lambda/\C)$ muni du préordre $x\leq
y$ si et seulement si il existe un morphisme $x\to y$ dans $\mbf
\Lambda/\C$. L'ensemble sous-jacent à $O(\C)$ coïncide avec l'ensemble des
classes d'équivalences de polysimplexes non dégénérés de $\C$.\\
Nous dirons qu'un ensemble polysimplicial $\C$ est \emph{intérieurement
  libre}\index{Ensemble polysimplicial!intérieurement libre} si $\Aut(\mbf n)$
agit librement sur $\C^{\nd}_{\mbf n}$. Si $\C_1\to\C_2$ est un morphisme d'ensembles
polysimpliciaux qui envoie les polysimplexes nondégénérés sur les polysimplexes
non dégénérés tels que $O(\C_1)\to O(\C_2)$ soit un ismorphisme et $\C_2$
soit intérieurement libre, alors $C_1\to\C_2$ est un isomorphisme.\\
Berkovich définit aussi une \emph{catégorie polysimpliciale stricte}
$\Lambda$ dont les objets sont ceux de $\mbf \Lambda$, mais avec seulement
les morphismes injectifs. Le foncteur $\Lambda\to\mbf \Lambda
\to\mbf \Lambda^{\circ}\Set$ s'étend en un foncteur
\[\Lambda^{\circ}\Ens\to\mbf \Lambda^{\circ}\Ens\]
qui commute aux limites
inductives (les objets de $\Lambda^{\circ}\Ens$ sont les \emph{ensembles
  strictement polysimpliciaux}\index{Ensemble strictement polysimplicial}).\\
Berkovich construit alors un foncteur $\Sigma:\mbf \Lambda\to\Ke$ vers
la catégorie des espaces de Kelley, \ie les espaces topologiques $X$ pour lesquels
un sous-ensemble est fermé dès que son intersection avec tout
compact est compacte. Ce foncteur envoie $[\mbf n]$ en $\Sigma_{\mbf
  n}=\{(u_{il})_{0\leq i\leq p,0\leq l\leq n_i}\in [0,1]^{[\mbf n]}|\sum_l
u_{il}=1\}$ (si $\mbf n=(n_0,\dots,n_p)$, $\Sigma_{\mbf n}$ est le produit
des simplexes standards de dimensions $n_0,\dots, n_p$), et si $\gamma$ est un morphisme associé à $(J,f,\alpha)$,
$\Sigma(\gamma)$ envoie $\mbf u=(u_{jk})$ vers $\mbf u'=(u'_{il})$ d\'efini
comme suit~: si $[\mbf m]\neq [0]$ et $i\notin\Image(f)$ ou $[\mbf m]=[0]$
alors $u'_{il}=1$ pour $l=\alpha_i(0)$ et $u'_{il}=0$ sinon; si $[\mbf
m]\neq [0]$ et $i\in\Image(f)$, alors 
$u'_{il}=u_{f^{-1}(i),\alpha_i^{-1}(l)}$ pour $l\in\Image(\alpha_i)$ et
$u'_{il}=0$ sinon.\\ 
En étendant $\Sigma$ pour que $|\ |$
commute aux limites inductives,
$\Sigma$ induit un foncteur, la \emph{réalisation géométrique},
\index{Réalisation géométrique} 
\[|\ |:\mbf\Lambda^{\circ}\Set\to\Ke.\]
On a aussi un bifoncteur
\[\sq:\mbf\Lambda^{\circ}\Set\times\mbf\Lambda^{\circ}\Set\to\mbf\Lambda^{\circ}\Set\]
qui commute aux limites inductives et qui v\'erifie
$[(n_0,\cdots,n_p)]\sq[(n'_0,\cdots,n'_{p'})]=[(n_0,\cdots,n_p,n'_0,\cdots,n'_{p'})]$. Ainsi
$|\C\sq\C'|=|\C|\times|\C'|$ (où le produit de droite est le produit
dans la catégorie des espaces de Kelley).\\

Si $X$ est strictement polystable sur $k$ et $x\in \Str(X)$, on note $\Irr(X,x)$ 
  l'espace métrique des composantes irréductibles de $X$ passant par $x$, muni de la distance 
  $d(X_1,X_2)=\codim_x(X_1\cap X_2)$. Il existe alors un n-uplet $\mbf n$
  tel que $\Irr(X,x)$ soit isométriquement en bijection avec $[\mbf n]$, et si
  $[\mbf m]\to[\mbf n]$ est isométrique, il existe un unique
  $y\in\Str(X)$ avec $y\leqslant x$ et une unique bijection isométrique
  $[\mbf m]\to\Irr(X,y)$ tels que
\[\begin{array}{ccc} {[\mbf n]} & \to & \Irr(X,x) \\ \uar & & \uar \\ {[\mbf m]}
  & \to & \Irr(X,y)\end{array}\] commute.\\
Le foncteur qui à $[\mbf n]$ associe l'ensemble des couples $(x,\mu)$, où
$x\in\Str(X)$ et $\mu$ est une bijection isométrique $[\mbf n]\to\Irr(X,x)$,
d\'efinit un ensemble polysimplicial strict $C(X)$ (et donc un ensemble polysimplicial
$\C(X)$).\\
On a un isomorphisme fonctoriel d'ensembles ordonnés $O(\C(X))\simeq\Str(X)$.
\begin{prop}[{\cite[prop. 3.14]{berk2}}]\label{berk314} On a un foncteur $\C:\mcal
  Pst^{\sm}\to\mbf\Lambda^{\circ}\Ens$ tel que $\C(X)$ est le complexe qu'on
  vient de définir si $X$ est strictement polystable et pour tout morphisme étale surjectif $X'\to X$~:
\[\C(X)=\Coker(\C(X'\times_XX')\rightrightarrows \C(X')).\]
\end{prop}
Ici $\Coker$ désigne simplement la limite inductive du diagramme (les
catégories $\Lambda^\circ\Ens$ et $\mbf\Lambda^\circ\Ens$ admettent toutes
les limites inductives en tant que catégories de préfaisceaux).\\
Ce foncteur $\C$ s'étend aux filtrations strictement polystables sur $K$
de longueur $l$.\\
Supposons que l'on a construit $\C$ pour les fibrations strictement
polystables de longueur $l-1$ de manière à ce que $O(\C(\underline X))=\Str(X_{l-1})$. Soit $\underline X:X_l\to
X_{l-1}\to\cdots\to\Spec k$ une fibration strictement polystable, et soit $\underline
X_{l-1}:X_{l-1}\to\cdots\to\Spec k$. Alors pour tout $x'\leqslant
x\in\Str(X_{l-1})$, on a le lemme suivant~:
\begin{lem}[{\cite[cor. 6.2]{berk2}}]\label{berk62} Il existe un foncteur
  canonique de cospécialisation $\C(X_{l,x})\to \C(X_{l,x'})$
et si $x''\leqslant x'\leqslant x$, le foncteur  $\C(X_{l,x})\to \C(X_{l,x''})$
coïncide avec la composition $\C(X_{l,x})\to \C(X_{l,x'})\to \C(X_{l,x''})$.\end{lem}
Ceci s'étend en un foncteur \[D:(\mbf
\Lambda/(\C(\underline X_{l-1})))^{\op}\to\Str(X_{l-1})^{\op}\to\mbf\Lambda^{\circ}\Ens.\]
Berkovich définit alors un ensemble polysimplicial (où l'on a posé
$\C=\C(\underline X_{l-1})$)~:
\[ \C(\underline X)=\C\sq D = \Coker(\coprod_{N_1(\mbf\Lambda/\C)} \mbf\Lambda [\mbf n_y]\sq
D_x\rightrightarrows\coprod_{N_0(\mbf\Lambda/\C)}\mbf\Lambda [\mbf n_x]
\sq D_x). \]
Cette construction s'étend aux fibrations polystables~:
\begin{prop}[{\cite[prop 6.9]{berk2}}]\label{berk69}
Il existe un foncteur $\C:\mcal Pst_l^{\tps}\to\mbf\Lambda^\circ\Ens$ tel que~:
\begin{enumerate}[(i)]
\item  pour tout morphisme étale surjectif de fibrations polystables $X'\to X$: 
\[\C(X)=\Coker(\C(X'\times_XX')\rightrightarrows \C(X')).\]
\item $O(\C(\underline X))\simeq\Str(X).$\end{enumerate}
\end{prop}

Si $k$ est algébriquement clos, le complexe polysimplicial d'une
fibration polystable est invariant par changement de corps de base:
\begin{prop}[{\cite[prop. 6.10]{berk2}}] \label{berk610}Si
  $\underline{\fk X}$ est une fibration polystable sur $k$, alors
  pour toute extension de corps $k\to k'$, $\C(X_{k'})\to\C(X)$ est un isomorphisme.\end{prop}

Berkovich associe à une fibration polystable $\underline{\fk X}=(\fk X=\fk
X_l\to\fk X_{l-1}\to\cdots\to\Spf(O_K))$ un sous-ensemble de la fibre
générique $\fk X_{\eta}$ de $\fk X$, le \emph{squelette}\index{Squelette} $S(\underline{\fk
  X})$ de $\underline{\fk X}$, qui est canoniquement homéomorphe à $|\C(\fk
X_s)|$ (voir~\cite[th. 8.2]{berk2}), et tel que $\fk X_{\eta}$
se rétracte par déformation forte et propre sur $S(\underline{\fk X})$.\\ 

En fait, quand $\underline{\fk X}$ est non dégénéré --- par exemple
g\'en\'eriquement lisse (nous utiliserons seulement les résultats de Berkovich
à de telles fibrations) --- le squelette de $\underline{\fk X}$ dépend seulement
de $\fk X$  d'après~\cite[prop. 4.3.1.(ii)]{berk3}. Un schéma formel $\fk
X$ sur
$O_K$ pour lequel il existe une telle fibration polystable est alors dit
\emph{pluristable}\index{Morphisme!pluristable}, et on notera simplement $S(\fk X)$ son squelette.\\ 
Dans ce cas, \cite[prop. 4.3.1.(ii)]{berk3} donne une description de $S(\fk
X)$, qui est indépendante de la fibration. Pour $x,y\in\fk
X_{\eta}$, on pose $x\preceq y$ si pour tout morphisme étale $\fk
X'\to\fk X$ et tout $x'$ au-dessus de $x$, il existe $y'$ au-dessus de
$y$ tel que pour tout $f\in O(\fk X_\eta)$, $|f(x')|\leq |f(y')|$
($\preceq$ est un ordre sur $\fk X_{\eta}$). Alors $S(\fk X)$ est
simplement l'ensemble des points maximaux de $\fk X_{\eta}$ pour $\preceq$.\\ 

La rétraction à $S(\underline{\fk X})$ commute aux morphismes étales~:
\begin{thm}[{\cite[th. 8.1]{berk2}}]\label{berk81}
On peut construire, pour toute fibration polystable $\underline{\fk X}=(\fk
X_l\stackrel{f_{l-1}}{\to}\cdots\stackrel{f_1}{\to}\fk X_1\to\Spf(O_K))$,
une rétraction par déformation forte et propre $\Phi^l:\fk X_{l,\eta}\times [0,l]\to \fk X_{l,\eta}$ de $\fk
X_{l,\eta}$ sur le squelette $S(\underline{\fk X})$ de $\underline{\fk X}$
tel que, en notant  $x_t=\Phi^l(x,t)$, on ait~:\begin{enumerate}[(i)]
\item $S(\underline{\fk X})=\bigcup_{x\in S(\underline{\fk X}_{l-1})}S(\fk
  X_{l,x})$ (union disjointe d'ensembles), où $\underline{\fk X}_{l-1}=(\fk
  X_{l-1}\to\cdots\to\Spf(O_K))$~;
\item $(x_t)_{t'}=x_{\max(t,t')}$
\item Si $\phi: \underline{\fk Y}\to \underline{\fk X}$ est un morphisme de fibrations dans $\mcal
  Pst f^{\et}_l$, on a $\phi_{l,\eta}(y_t)=\phi_{l,\eta}(y)_t$ pour tout
  $y\in\fk Y_{l,\eta}$.\end{enumerate}\end{thm}
Décrivons plus précisément la déformation.\\
Si $\fk X=\Spf O_K\{P\}/(p_i-z_i)$ où $P$ est isomorphe à
$\oplus_{0\leq i\leq p}\mbf N^{n_i+1}$,
$p_i=(1,\cdots,1)\in\mbf N^{n_i+1}$ et $z_i\in O_K$, soit $\fGm$ le groupe
multiplicatif $\Spf O_K\{T,\frac{1}{T}\}$ sur $O_K$, notons pour tout
$n$ par $\fGm^{(n)}$ le noyau de la multiplication $\fGm^{n+1}\to\fGm$ et
soit $\fk G$ le complété formel de l'identité dans $\prod_i\fGm^{(n_i)}$
(c'est un groupe formel). Alors $\fk G$ agit sur $\fk X$. $G=\fk
G_\eta$ agit alors sur $\fk X_\eta$. $G$ a un sous groupe canonique $G_t$ pour
$t\in [0,1]$ d\'efini par les inégalités $|T_{ij}-1|\leq t$ où $T_{ij}$
sont les coordonnées dans $G$, qui est un quotient de 
$\prod_i\Gm^{n_i+1}$. $G_t$ a un point maximal $g_t$.\\
Alors pour $x\in X$, $x_t=g_t*x$ d\'efinit la déformation forte (où $*$ est
la multiplication définie dans~\cite[§ 5.2]{berk}).\\
Si $\fk X$ est étale sur $\Spf O_K\{P\}/(p_i-z_i)$, l'action de $\fk G$
s'étend de façon unique sur $X$, et $x_t$ encore définit par
$g_t*x$. Pour tout $\fk X$ polystable sur $O_K$, on a ainsi défini la
déformation localement sur $\fk X_\eta^{\an}$ pour la topologie
quasi-étale, et Berkovich vérifie qu'elle se descend bien en une
déformation de $\fk X$.\\
Pour une fibration polystable $\fk X\to\fk X_{l-1}\to\cdots\to \Spf O_K$,
on suppose d'abord $\fk X=\Spf B\to\fk X_{l-1}=\Spf A$ avec
$B=A\{P\}/(p_i-a_i)$ (on appellera un tel morphisme polystable \emph{standard}),
on déforme d'abord sur $S(\fk X/\fk X_{l-1})$ fibre par fibre (les fibres
étant strictement polystables). L'image
obtenue s'identifie avec $S=\{(x,\mbf r_0,\cdots,\mbf r_p)\in \fk
X_{l-1,\eta}, r_{i0}\cdots r_{in_i}=|a_i(x)|\}$, on a alors une homotopie $\Psi:S\times [0,1]\to S$ d\'efinie par $\Psi(x,\mbf r_0,\cdots,\mbf
r_p,t)=(x_t,\psi_{n_0}(\mbf r_0,|a_0(x_t)|),\cdots,\psi_{n_p}(\mbf
r_p,|a_p(x_t)|))$, où $\psi_n$ est une déformation forte de
$[0,1]^{n+1}$ sur $(1,\cdots,1)\in [0,1]^{n+1}$ d\'efinie par Berkovich (nous
utiliserons seulement le fait que
$\psi_n(r_i,t)_k^\lambda=\psi_n(r_i^\lambda,t^\lambda)_k$ pour tout
$\lambda\in\mbf R^{*+}$ et tout $k\in [[0,n]]$)
, et $x_t$ est défini par la déformation forte de $\fk X_{l-1,\eta}$.\\
Si $\fk X\to\fk X'\to\fk X_{l-1}$ est la composée géométriquement
élémentaire d'un morphisme étale et d'un morphisme polystable standard, $S(\fk X/\fk
X_{l-1})\to S(\fk X'/\fk X_{l-1})$ est un isomorphisme, et donc on déforme $\fk X'$ fibre par
fibre sur $S(\fk X/\fk X_{l-1})$, puis on fait la même déformation que pour
$S(\fk X'/\fk X_{l-1})$. Pour une fibration polystable quelconque
$X\to\cdots\to O_K$, on a défini la déformation localement pour la topologie
quasi-étale de $\fk X_{\eta}$, et Berkovich vérifie que l'on peut la
descendre en une déformation de $X$.\\

Berkovich deduit de (\ref{berk81}.(ii)) le corollaire suivant~:
\begin{cor}[{\cite[cor. 8.5]{berk2}}]\label{berk85} Soit $K'$ une extension
  galoisienne finie de $K$ et soit $\underline{\fk X}$ une fibration
  polystable sur   
  $O_{K'}$ avec une fibre générique $\fk X_{l,\eta}$ normale. Supposons
  qu'on ait une action d'un groupe fini $G$ sur $\underline{\fk X}$ au
  dessus de $O_K$ et un ouvert dense de Zariski $U$ de $\fk X_{l,\eta}$
  stable par $G$. 
  Alors il existe une rétraction par déformation forte de 
  $G\backslash U$ sur un fermé homéomorphe à $G\backslash |\C(\underline{\fk X})|$.\end{cor}
Plus précisément, dans ce corollaire, le fermé en question est l'image de
$S(\underline{\fk X})$ (qui est $G$-équivariant et contenu dans $U$) par $U\to G\backslash U$.\\

Le théorème~\ref{berk81} implique aussi que le squelette est fonctoriel dans le cas pluristable~:
\begin{prop}[{\cite[prop. 4.3.2.(i)]{berk3}}] Si $\phi:\fk X\to\fk Y$ est
  un morphisme pluristable entre schémas formels pluristables non dégénérés
  de $O_K$, $\phi_\eta(S(\fk X))\subset S(\fk Y)$.\end{prop}
Plus précisément, d'après la construction de $S$, on a $S(\fk
Y)=\bigcup_{x\in S(\fk X)}S(\fk Y_x)$.\\

En vue d'utiliser la description précédente de l'espace de Berkovich d'un
schéma pluristable sur $O_K$ pour comprendre la topologie d'un schéma lisse
sur $k$, nous aurons besoin du résultat de de Jong sur l'existence
d'altérations\index{Altération} par de tels schémas pluristables sur $O_K$.\\
Plus précisément nous utiliserons la conséquence suivante des résultats de
de Jong donnée par Berkovich (comme nous supposerons $K$ alg\'ebriquement clos de caractéristique
0, nous simplifierons l'énoncé de Berkovich):
\begin{lem}\label{berk92} \emph{(\cite[lem. 9.2]{berk2})} Supposons $K$ alg\'ebriquement clos de
  caractéristique nulle, soit $X$ un schéma intègre propre plat et de
  présentation finie sur
  $O_K$, de fibre générique irréductible de dimension $l$. Alors il existe:\begin{enumerate}[(a)]
\item une fibration polystable $\underline X'=(X'_l\to\cdots\to X'_0=\Spec
  O_K)$, où tous les morphismes ont des fibres génériques lisses et
  géométriquement irréductibles ;
\item une action d'un groupe fini $G$ sur $\underline X'$ au-dessus de $O_K$,
\item un morphisme dominant $G$-équivariant $\phi:X'_l\to X$ au-dessus de $O_K$,
 dont la fibre générique est génériquement étale galoisien de groupe $G$.\end{enumerate}\end{lem}

\chapter{Quelques propriétés du groupe fondamental tempéré}
Dans ce chapitre, nous prouvons pour le groupe fondamental tempéré
certains résultats qui sont classiques pour le groupe fondamental
profini. Ainsi nous montrerons la localité des revêtements tempérés pour la
topologie étale (mais seulement dans le cas des courbes compactes),
l'invariance birationnelle du groupe fondamental tempéré des variétés
propres et lisses, l'invariance par extension isométrique algébriquement close du corps
de base, la formule de Künneth et l'isomorphisme entre l'abélianisé du
groupe fondamental tempéré d'une courbe et le groupe fondamental tempéré de
sa jacobienne.
\section{Localité des revêtements tempérés des courbes}
Dans ce paragraphe, nous nous intéresserons au caract\`ere local de la propri\'et\'e d'être un
revêtement tempéré. Il découle directement de la définition des revêtements
tempérés qu'un revêtement étale fini est un morphisme de descente effective
pour les revêtements tempérés. En revanche, la définition des revêtements
tempérés n'est pas locale pour la topologie de Berkovich. Par exemple, le
logarithme induit un revêtement étale de $\mbf A^1_{\mbf C_p}$ qui est
tempéré au-dessus de tous les ouverts pr\'ecompacts (\ie d'adh\'erence compacte), mais qui n'est
pas lui-même tempéré.\\
Plus précisément, pour tout entier $m\in \mbf N$, l'image inversedu disque
ouvert $\mbf D_m$ de centre $0$ et de rayon $|p|^{-m+\frac{1}{p-1}}$ est une
union disjointe
\[\coprod_{\zeta\in\mu_{p^{\infty}}}\{q\in \mbf D(1,1^{-})|\ |\zeta q^{p^m}-1|<1\} \]
et chaque composante est finie étale sur $\mbf D_m$
(\cite[ex. III.1.2.6.(i)]{andre1}). Par contre le revêtement n'est pas
tempéré puisque $\gtemp(\mbf A^1)$ est trivial.\\
Nous verrons ici qu'au-dessus d'un bon espace analytique compact de dimension
1 (par exemple une courbe projective ou une courbe affinoïde), la notion de
revêtement tempéré est locale pour la topologie étale de Berkovich. Ainsi
dans ce cas, les notions de revêtement tempéré et de faisceau localement
constant pour la topologie étale co\"incident.\\
Comme application, nous en déduirons qu'un revêtement tempéré d'une courbe
algébrique qui est scindé au-dessus de tous les sommets du squelette de la
courbe est un revêtement topologique.\\
 
Si $\overline X$ est un bon espace analytique de dimension 1 et $D$ est un
ensemble localement fini de points fermé pour la topologie de Zariski, un
revêtement ramifié de $(\overline X,D)$ sera un morphisme $Y\to X$ tel que il existe
un recouvrement ouvert $(U_i)$ de $X$ tel que $Y_{U_i}\to U_i$ soit une somme
directe de revêtements finis non ramifiés hors de $D\cap U_i$ et
kummérien au voisinage de tout point de $D$.

\begin{prop}\label{courbesgalrev}
Soit $\overline X$ un bon espace analytique compact de dimension 1, soit $D$ un ensemble
fini de points fermés pour la topologie de Zariski et soit $X=\overline X\backslash D$.\\
Soit $\overline Y\to \overline X$  un revêtement ramifié de $(\overline X,D)$  qui est un $G$-torseur (\ie $Y\to X$ est
galoisien mais pas nécessairement connexe) pour un certain groupe
$G$. Alors $Y=\overline Y_X\to X$ est tempéré.\\
Plus g\'en\'eralement, tout revêtement ramifié $\overline Y\to \overline X$
de $(\overline X,D)$ tel qu'il existe $Z\to \overline X$ galoisien
dominant toute composante connexe de $\overline Y$ est tempéré au-dessus de $X$.\end{prop}
\dem
Il existe un recouvrement ouvert $(V_i)_{i\in I}$ de $\overline X$ tel que
$Y_{V_i}=\coprod_j U_{ij}$ où $U_{ij}$ est un revêtement fini de $V_i$.\\
Comme $X$ est de dimension 1, quitte à raffiner le recouvrement $(V_i)$, on peut supposer que
$V_{ijk}:=V_i\cap V_j\cap V_k=\emptyset$ pour $i,j,k$ deux à deux
distincts~(\cite[cor. 3.2.8]{berk}).\\
Comme tout bon espace analytique, $V_i$ est localement
connexe~(\cite[cor. 2.2.8]{berk}), et quitte
à remplacer $V_i$ par l'ensemble de ses composantes connexes, on
peut supposer de plus les $V_i$ connexes (le nouveau recouvrement vérifie
bien encore la propriété précédente). Comme $\overline X$ est compact, on peut supposer de plus $I$ fini.\\

Pour $i\in I$, soit $U_{i0}$ une composante connexe de $Y_{V_i}$~: c'est un
revêtement fini galoisien de $V_i$, et toutes les composantes
connexes de $Y_{V_i}$ sont (non canoniquement) isomorphes à $U_{i0}$. Soit $n_i$ le cardinal
des fibres géométriques de $U_{i0}$, et posons $n:=\prod_{i\in I}n_i$.\\
Soit $S_i$ l'union disjointe de $n/n_i$ copies de $U_{i0}$. Montrons qu'il
existe un revêtement $S$ de $X$ dont la restriction à $V_i$ est isomorphe à
$S_i$. Pour cela il suffit, pour tout $i\neq j$, de construire un
$V_{ij}$-isomorphisme entre $S_{i|V_{ij}}$ et $S_{j|V_{ij}}$, où
$V_{ij}:=V_i\cap V_j$ (puisque
$V_{ijk}$ est vide, il n'y a pas de conditions de cocycle à vérifier).\\

Comme $V_{ij}$ est localement connexe, il suffit de construire
un tel isomorphisme au-dessus de toute composante connexe de $V_{ij}$.\\
Soit $Z$
une telle composante connexe. Les composantes connexes de $Y_Z$ sont toutes isomorphes entre elles (non canoniquement). Soit
$U_Z$ une telle composante connexe, et soit $n_Z$ le cardinal des fibres
géométriques de $U_Z\to Z$ (elles ont toutes même cardinal). Alors $U_{i0|Z}$ est union disjointe de composantes
connexes de $Y_Z$, et donc est isomorphe à une union disjointe de copies
de $U_Z$ (plus précisément $n_i/n_Z$ copies). Donc $S_{i|Z}$ est isomorphe
à une union disjointe de $n/n_Z$ copies de $U_Z$. On a le même résultat
pour $S_{j|Z}$. Il existe donc effectivement un isomorphisme entre
$S_{i|Z}$ et $S_{j|Z}$, ce qui permet de construire $S$ ayant la propriété
voulue. Mais $Y_{S_{i,X}}\to S_{i,X}$ est scindé, donc $Y_{S_X}\to {S_X}$ est un revêtement topologique.

\findem

Pour $D$ vide, on obtient que tout revêtement étale galoisien de $\overline
X$ est tempéré.

\begin{cor}\label{courbestempet}Soit $X$ un bon espace analytique compact de dimension 1. Alors,
  la catégorie des revêtements tempérés de $X$ est canoniquement équivalente à la catégorie des
  faisceaux localement constants sur le site étale $X_{\et}$.
\end{cor}
\dem
De manière générale (sans hypothèse sur $X$), l'objet de $\tilde X_{\et}$ associé à un revêtement tempéré
est clairement localement constant pour la topologie étale (il est trivialisé par un recouvrement
ouvert d'un revêtement fini). D'où un foncteur pleinement fidèle de la cat\'egorie $\Covtemp(X)$ vers la cat\'egorie $\lcs(X_{\et})$ des objets localement constants de $X_\et$.\\
En outre, un faisceau localement constant sur $X_{\et}$ est représentable par un
revêtement étale grâce à~\cite[lem. 2.3]{dJ1}, ce qui fait de la catégorie
des faisceaux étales localement constants de $X$ une sous-cat\'egorie pleine
de la catégorie des revêtements étales. 

Revenons à notre cas. On peut supposer $X$ connexe.\\
Comme $\lcs(X_{\et})$ est engendrée par les
objets galoisiens~(car le topos étale de $X$ est localement connexe),
tout faisceau localement constant est tempéré d'après~\ref{courbesgalrev}.
\findem
On en déduit en particulier des résultats de type Van Kampen~:
\begin{cor} Soit $X$ un bon espace analytique compact de dimension 1. Soit
  $U_i$ un recouvrement étale de $X$. Alors $\coprod U_i\to X$ est de
  descente effective pour les revêtements tempérés.\end{cor}

\begin{cor} Soit $X$ un bon espace analytique compact de dimension
  1. Alors $\gtemp(X)$ est le complété prodiscret $\get(X)^\pd$ de $\get(X)$.
\end{cor}
\dem
La catégorie des $\get(X)^\pd$-ensembles est équivalente à la
sous-catégorie pleine de la catégorie des $\get(X)$-ensembles dont toute
composante connexe est dominé par un revêtement
galoisien. D'après~\ref{courbesgalrev}, cette catégorie est équivalente à
celle des unions disjointes de revêtements tempérés.
\findem
Si $Y$ est un revêtement tempéré de $X$, soit $H$ le sous-groupe des éléments
de
$\gtemp(X,x)$ qui agissent trivialement sur $Y_x$. C'est un sous-groupe
distingué de $\gtemp(X,x)$, ouvert d'après la prodiscrétion de
$\gtemp(X,x)$. Le revêtement tempéré galoisien correspondant $Z\to X$ est appelé
clôture galoisienne de $Y\to X$ (il ne dépend pas de $x$ à isomorphisme
près). Si $Y\to X$ est scindé au-dessus d'un point $x$ de $X$, $Z\to X$ est
aussi scindé au-dessus de $x$.
 
\begin{cor}
Soit $X$ une courbe algébrique, $Y$ un revêtement étale de
$X^{\an}$ qui est tempéré au-dessus d'un ouvert de Zariski $U$ de $X$, alors $Y\to X$ est tempéré.
\end{cor}
\dem
Soit $\overline X$ une courbe compl\`ete contenant $X$.
Soit $Z_U\to U$ la clôture galoisienne de $Y_U\to U$. Alors $Z_U\to U$ se
prolonge en un revêtement galoisien $Z\to \overline X$ (\cite[th.2.1.11]{andre1}). Il est clair que $Z_X\to X$ est non ramifié. Donc $Z_X$, et a fortiori $Y$,
est tempéré.
\findem

\begin{prop}
Soit $X$ une courbe algébrique, de compactification $\overline X$. Soit
$V$ l'ensemble des sommets du squelette de $X$ associé à une réduction semistable. Si $S\to X$ est un
revêtement tempéré scindé en tout point $v$ de $V$ (\ie $S_{v}$ est un $\ga(\mcal H(v))$-ensemble trivial), alors $S\to X$ est un
revêtement topologique.
\end{prop}
\dem
Quitte à remplacer $S$ par sa clôture galoisienne, on peut supposer $S\to
X$ galoisien. Soit $G$ son groupe de Galois.\\
On peut aussi prolonger $S$ en un revêtement ramifié $\overline S\to\overline
X$ (\cite[th.2.1.11]{andre1}).\\
Supposons par l'absurde qu'il existe $x\in X$ tel que $S_x\to x$ soit non
trivial~(évidemment, $x\notin V$).\\
Soit $A$ la composante connexe de l'ouvert $\overline X\backslash V$ qui
contient $x$. Alors $A$ est ou bien isomorphe à la boule ouverte $B(0,1)$ ou bien isomorphe à une couronne ouverte.\\
De plus $\overline S_A\to A$ est ramifié au-dessus d'au plus un point de
$A$ (il y a au plus un point cuspidal de $X$ dans $A$).
Commençons par supposer $A\simeq B(0,1)$, et fixons un isomorphisme (ce qui
nous permettra d'identifier $A$ à $B(0,1)$).\\
Soit $v\in V$ le point frontière de $A$. Par définition d'un revêtement
ramifié, il existe un voisinage $V$ de $v$ tel que $\overline
S_V=\coprod_{i\in I} U_i $ avec $U_i$ un revêtement étale fini de $V$ d'ordre
$n_i$. Comme $S$ est supposé galoisien, on peut supposer que tous les $U_i$
sont isomorphes, et notons-le simplement $U$.\\
Comme $U\to V$ est scindé en $v$, quitte à réduire $V$, on peut supposer
que $U\to V$ est scindé. On peut donc supposer $U\simeq V$, \ie
$S_V=\coprod_{i\in I} V$.\\
$V\cap  A$ contient une couronne ouverte $A(]\epsilon,1[)$ pour $\epsilon<
1$ assez grand (quand on identifie $A$ avec $B(0,1)$).\\
Construisons alors $T$ par recollement de $\overline S_A$ et de
$\coprod_IB(\infty,\epsilon)$ (où $B(\infty,\epsilon)=\{z\in \mbf P_1|\
|z|>\epsilon\}$) le long des immersions ouvertes~:
\[\overline S_A\supset \coprod_I A(]\epsilon,1[)\subset \coprod_I
B(\infty,\epsilon).\]
On a un morphisme $T\to \mbf P^1$ dont la restriction à $\overline S_A$ est
$\overline S_A\to A\simeq B(0,1)\subset\mbf P^1$ et la restriction à $\coprod_I
B(\infty,\epsilon)$ est $\coprod_I B(\infty,\epsilon)\to
B(\infty,\epsilon)\subset \mbf P^1$.\\
La restriction de $T\to \mbf P^1$ au-dessus de $B(0,1)$ est isomorphe à
$\overline S_A\to A$, et la restriction à $B(\infty,\epsilon)$ est un
revêtement trivial, donc $T\to \mbf P^1$ est un revêtement, ramifié en au
plus un point $x_0$.\\
De plus l'action de $G$ sur $\overline S_A$ se prolonge à $T$, donc $T\to
\mbf P^1$ est un revêtement galoisien non ramifié $\mbf
P^1\backslash\{x_0\}$. Grâce à la remarque suivant la proposition~\ref{courbesgalrev}, la
restriction de $T$ à $\mbf P^1\backslash\{x_0\}$ est tempérée. Donc $T$ est
nécessairement trivial, ce qui contredit l'hypothèse faite sur $x$.\\
Si $A$ est isomorphe à une couronne ouverte la construction est
parfaitement similaire (on prolonge $\overline S_A$ en un revêtement
étale galoisien de $\mbf P^1$). 
\findem

Si $x\in X$ est un point, l'image de $G_{\mcal
  H(x)}=\ga(x)=\gtemp(x)\to\gtemp(X)$ est appelé sous-groupe de
décomposition de $x$ (il est défini à conjugaison près). C'est un sous-groupe compact de $\gtemp(X)$. 
\begin{cor} Soit $X$ une courbe algébrique munie d'une réduction semi-stable. Le noyau de
  $\gtemp(X)\to\gtop(X)$ est le sous-groupe distingué de $\gtemp(X)$ topologiquement engendré par les sous-groupes
  de décomposition en les sommets du squelette de $X$ (pour cette réduction semis-stable).
\end{cor}
\dem
Soit $H$ le noyau de $\gtemp(X)\to\gtop(X)$.
Puisque $\gtop(X)$ est discret et sans torsion (car isomorphe
au groupe fondamental de son squelette qui est un graphe), $H$ contient tous les sous-groupes compacts de $\gtemp(X)$, en particulier les sous-groupes de d\'ecomposition.\\
Soit $S$ un rev\^etement temp\'er\'e de $X$, et $T$ le $\gtemp(X)$-ensemble correspondant.\\
L'action de $H$ sur $T$ est triviale 
\begin{itemize}
\item
si et seulement si $S$ est un rev\^etement topologique,
\item si et seulement si $S$ est scind\'e en tout $v\in V$,
\item si et seulement si l'action de $G_{\mcal H(v)}$ sur
  $T$ est triviale pour tout sommet $v\in V$,
\item si et seulement si l'action du sous-groupe distingu\'e $H'$
  engendr\'e par les $G_{\mcal H(v)}$ sur $T$ est triviale,
\item si et seulement si l'action de l'adh\'erence $\overline{H'}$ de $H'$ sur $T$ est triviale.\end{itemize}
Mais $\overline{H'}$ est aussi l'intersection des sous-groupes ouverts contenant $H'$, donc $\overline{H'}=H$.
\findem
\section{Invariance birationnelle}
Soit $K$ un corps complet non archimédien.

\begin{prop}\label{birat} Soit $f:X\to Y$ un morphisme birationnel entre
  $K$-schémas propres et lisses. Alors
\[\Covtemp(Y)^{\mbb L}\to\Covtemp(X)^{\mbb L}\] est une équivalence de
catégories. En particulier,
\[\gtemp(X)^{\mbb L}\to\gtemp(Y)^{\mbb L}\]
est un isomorphisme.\end{prop}
\dem
On a un diagramme 2-commutatif de foncteurs~:
\[\begin{array}{ccc}\Dtop(Y) & \to & \Dtop(X) \\ \dar & & \dar \\ \Covalg(Y)^{\mbb L} & \simeq
& \Covalg(X)^{\mbb L}\end{array}\]
où la flèche du bas est une équivalence d'après~\cite[cor. X.3.4]{sga}.\\
Soit $S$ un revêtement étale fini de $Y$ et $T:=f^*S$ son pullback à $X$
(notons $g$ le morphisme $T\to S$ obtenu par changement de base de $f$).\\
Soit $i:U\to X$ l'immersion d'un ouvert dense, telle que $fi$ soit aussi
une immersion ouverte ($U$ est alors aussi un ouvert dense de $Y$).
$j:V:=i^*T\to T$ est aussi une immersion d'un ouvert dense et $gj$ aussi.\\
On a un diagramme 2-commutatif~:
\[\xymatrix{\Covtop(S) \ar[r]^{g^*}\ar[rd]^{(gj)^*} & \Covtop(T) \ar[d]^{j^*} \\ & \Covtop(V)}\] 
où les flèches verticales sont des équivalences d'après la proposition~\ref{andre114}. Donc
$\Dtop(Y)_S\to\Dtop(X)_T$ est une équivalence de catégories. $\Dtop(Y)\to
\Dtop(X)$ est donc une équivalence de catégories fibrées, et en passant aux
champs associés, on en déduit que $\Dtemp(Y)^{\mbb L}\to\Dtemp(X)^{\mbb L}$
est une équivalence de champs. En prenant les sections globales, on en
déduit le résultat voulu.
\findem

Dans le cas où $K$ est algébriquement clos de caractéristique nulle, on
peut également prouver l'invariance birationnelle directement en suivant la
preuve de~\cite[cor. X.3.4]{sga}.\\
En effet, soit $f:X\to Y$ une application rationnelle dominante de $K$-schemas lisses
et propres. $f$ est défini sur un ouvert de Zariski $U$ de $X$ (notons
$f_U$ le morphisme $U\to Y$ et $i_U$ l'immersion $U\to X$) dont le
complémentaire est de codimension $\geq 2$ dans $X$, on obtient un foncteur
de $\Covtemp(Y^{\mbb L})$ vers $\Covtemp(X)^{\mbb L}$ et, grâce à la
proposition~\ref{andre2111} on peut le
composer avec un quasi-inverse $\Covtemp(U)^{\mbb L}\to\Covtemp(X)^{\mbb
  L}$~: on obtient un foncteur $f^*_{(U)}:\Covtemp(Y)\to\Covtemp(X)$ tel
que $i_U^*f^*_{(U)}$ soit isomorphe à $f^*_U$. Si l'on choisit un autre
ouvert de Zariski $U'$ de $X$ vérifiant les mêmes propriétés, $i_{U\cap
  U'}^*f^*_{(U)}$ et $i_{U\cap U'}^*f^*_{(U')}$ sont tous les deux
isomorphes à $f^*_{U\cap U'}$, et donc $f^*_{(U)}$ et $f^*_{(U')}$ sont
isomorphes, puisque $X\backslash U\cap U'$ est également de codimension
$\geq 2$ dans $X$. On obtient donc un homomorphisme extérieur de groupes
topologiques $f_*:\gtemp(X)\to\gtemp(Y)$, qui ne dépend pas de $U$. En
particulier si $f$ est un morphisme de schémas, on peut choisir $U=X$ et
$f_*$ est l'homomorphisme extérieur usuel $\gtemp(X)\to\gtemp(Y)$.\\
Soit $g:Y\to Z$ une autre application rationnelle dominante entre
$K$-sch\'emas propres et lisses. Elle est définie sur un ouvert de Zariski $V$
de $Y$ dont le complémentaire est de codimension $\geq 2$ dans $Y$.
$gf:X\to Z$ est aussi une application rationnelle dominante de $K$-schémas
propres et lisses. $gf$ est donc définie sur un ouvert de Zariski $W$ de $X$
dont le complémentaire est de dimension $\geq 2$. Soit $U_0=U\cap
f_U^{-1}(V)\cap W$ (remarquons que $X\backslash U_0$ peut être de
codimension $< 2$). On a des morphismes $U_0\to V$ et $V\to Z$ représentant
$f$ et $g$ tel que le morphisme composé $(gf)_{U_0}:U_0\to Z$ représente
$gf$. Donc $i_{U_0}^*f^*_{(U)}g^*_{(V)}$ et $i_{U_0}^*(gf)^*_{(W)}$ sont
tous les deux isomorphes à $(gf)_{U_0}^*$. Puisque $i_{U_0}^*$ est
pleinement fidèle (prop. \ref{andre2111}), $f^*_{(U)}g^*_{(V)}$ et $(gf)^*_{(W)}$ sont isomorphes
(et donc $g_*f_*=(gf)_*$). On obtient ainsi un foncteur de la catégorie des
$K$-schémas propres et lisses avec pour morphismes les applications rationnelles
dominantes vers la catégorie des groupes topologiques avec pour morphismes
les morphismes extérieurs. En particulier, il envoie les isomorphismes (\ie
les applications birationnelles) sur
des isomorphismes.\\

\section{Altérations et groupe fondamental tempéré}
Dans cette section, nous allons décrire deux applications des théorèmes
d'existences d'altérations semi-stables de de Jong au groupe fondamental
tempéré. En effet, Berkovich a déjà montré dans~\cite[§ 9]{berk2} comment ces
altérations permettent de construire un squelette, homéomorphe à la
réalisation géométrique d'un ensemble polysimplicial, sur lequel se
rétracte un ouvert de Zariski de la variété.\\
Nous allons ici déduire des résultats de Berkovich que le groupe
fondamental tempéré d'une variété algébrique lisse est invariant par
changement de base algébriquement clos, et que le groupe fondamental
tempéré du produit de deux variétés lisses (sur un corps de base
algébriquement clos) est canoniquement isomorphe au produit des groupes
fondamentaux tempérés de chacune des variétés. 

\subsection{Invariance de $\gtemp$ par extension algébriquement close de corps de base}
Soit $X$ une variété algébrique lisse et connexe sur un corps complet non
archimédien algébriquement clos $K$ de caractéristique nulle.\\
Soit $K'/K$ une extension isométrique de corps valués complets.
\begin{lem}\label{leminvariance} Le foncteur $\Covtop(X)\to\Covtop(X_{K'})$
  est une équivalence de catégories.\\
Ainsi, si $x'$ est un point géométrique de $X^{\an}_{K'}$ d'image $x$ dans $X^{\an}$, alors
\[\gtop(X_{K'},x')\to\gtop(X,x)\]
est un isomorphisme.\end{lem}
\begin{proof}
Plongeons $X$ dans un schéma $\overline{X}$ propre,
plat et de présentation finie sur $O_K$.\\
Alors, d'après le lemme~\ref{berk92}, il existe une
fibration polystable $X'$ sur $O_K$ génériquement lisse muni d'une action
de groupe $G$ telle que
$(X',G)$ soit une altération galoisienne de $\overline{X}$.\\
Soit $U$ un ouvert de Zariski dense de $\overline{X}$ (et donc dans $X$) tel que $U'\to U$ soit fini (ou $U'$ est
l'image réciproque de $U$ dans $X'$).\\
Alors $U^{\an}$ se rétracte par déformation forte sur $G\backslash
S(X'_{s})$, d'après le corollaire~\ref{berk85}.\\
$X'_{O_{K'}}$ est aussi une fibration polystable
muni d'une action de $G$ et $X'_{O_{K'}}$ est aussi une altération
galoisienne de $\overline{X}_{K'}$, finie sur $U_{K'}$.\\
Ainsi comme dans le cas précédent, $U^{\an}_{K'}$ se rétracte par
déformation forte sur le sous-espace 
 $G\backslash S(X'_{K'})$ (et le morphisme naturel $U^{\an}_{K'}\to U^{\an}$ envoie $S(X'_{K',s})$ sur
$S(X'_s)$.\\
Mais $\C(X'_{K',s})\to \C(X'_{s})$ est un isomorphisme
d'après la proposition~\ref{berk610}. Le morphisme $U^{\an}_{K'}\to
U^{\an}$ est donc une équivalence d'homotopie.\\
On a le diagramme 2-commutatif suivant~:
\[\xymatrix{\Covtop(U_{K'}) & \Covtop(U) \ar[l] \\ \Covtop(X_{K'})\ar[u]
  & \Covtop(X)\ar[u]  \ar[l]}.\]
Les flèches verticales sont des équivalences d'après~\ref{andre114}, et
  nous venons de montrer que la flèche du haut est une équivalence.\\
La flèche du bas est donc également une équivalence.\end{proof}

Supposons maintenant $K'$ aussi algébriquement clos.
\begin{prop}\label{invariance} Soit $x'$ un point géométrique de $X_{K'}^{\an}$
  d'image $x$ dans $X^{\an}$, alors le morphisme $\gtemp(X_{K'},x')\to
  \gtemp(X,x)$ est un isomorphisme.\end{prop}
\begin{proof}
Considérons un système projectif cofinal dénombrable de revêtements
étales finis galoisiens géométriquement pointés $((X_i)_i,x_i)_{i\in \mbf N}$ de $X$.\\
Alors $(X_{i,K'},x'_i)_{i\in \mbf N}$, où $x'_i$ est un point de $X_{i,K'}$
au-dessus de $x_i$, est aussi un système projectif cofinal de revêtements étales finis
galoisiens géométriquement pointés de $X_{K'}$ d'après~\cite[lecture XIII]{sga}).\\
Si $X^{\infty}_i$ est le revêtement topologique universel de $X_i$,
  $X^{\infty}_{i,K'}:=(X^{\infty}_i)_{K'}$ est le revêtement topologique
  universel de $X_{i,K'}$ d'après~\ref{leminvariance}.\\
Comme $\Gal(X^{\infty}_{i,K'}/X_{K'})=\Gal(X^{\infty}_i/X)$,
en prenant la limite projective pour $i\in \mbf N$, on obtient le résultat désiré.\end{proof}

\subsection{Produits et groupe fondamental tempéré}
Soient $X,Y$ des variétés algébriques lisses et connexes sur un corps complet
non archimédien algébriquement clos $K$ de caractéristique nulle.\\
\begin{lem}\label{lemprod} Si $x$ et $y$ sont des points géométriques de
  $X^{\an}$ et $Y^{\an}$ respectivement (à valeur dans un même corps $\Omega$), alors
\[\gtop(X\times Y,(x,y))\to\gtop(X,x)\times\gtop(Y,y)\]
est un isomorphisme.\end{lem}
\begin{proof}
Soit $\overline{X}$ (resp. $\overline{Y}$) un schéma propre, plat et de
présentation finie sur $O_K$, dans lequel se plonge $X$ (resp. $Y$), et soit $(X',G)\to\overline{X}$ (resp. $(Y',H)\to \overline{Y}$) une
altération galoisienne où $X'\to O_K$ (resp. $Y'\to O_K$) est une fibration
polystable sur laquelle agit $G$ (resp. $H$).\\
Soit aussi $U\subset X$ (resp. $V\subset Y$) une immersion d'un ouvert de
Zariski dense
tel que $U'\to U$ (resp. $V'\to V$) soit fini.\\
Le fait que $\gtop(X\times Y,(x,y))\to\gtop(X,x)\times\gtop(Y,y)$ soit un
isomorphisme ne dépend pas de $x$ et de $y$, on peut donc supposer $x\in U$ et $y\in V$.\\
Alors, comme dans la preuve du lemme~\ref{leminvariance}, $U^{\an}$
(resp. $V^{\an}$) se rétracte par déformation forte sur $G\backslash S(X'_s)$ (resp. 
$H\backslash S(Y'_s)$).\\

On obtient le même résultat pour $X\times Y$~:\\
$X'\times Y'\to X\times Y$ est une altération galoisienne de groupe $G\times
H$, et $U\times V$ se rétracte par déformation forte sur $(G\times H)\backslash S(X\times Y)$.\\
Mais les morphismes pluristables $X\times Y\to X$ et $X\times Y\to Y$
envoient $S(X\times Y)$ sur $S(X)$ et $S(X)$ respectivement. D'où une
application continue $f:S(X\times Y)\to S(X)\times S(Y)$ (compatible avec
l'action de $G\times H$). Mais comme 
\[S(X\times Y)=\bigcup_{x\in S(X)}S((X\times Y)_x)=\bigcup_{x\in S(X)}S(Y_j\otimes \mcal H(x))\]
et puisque $S(Y\otimes\mcal H(x))\to S(Y)$ est un homéomorphisme d'après la
proposition~\ref{berk610}, $f$ est bijectif, et un homéomorphisme
puisque $S(X\times Y)$ est compact.\\
Ainsi $(G\times H)\backslash S(X\times Y)\to G\backslash S(X)\times
H\backslash S(Y)$ est un homéomorphisme.\\
Donc $(U\times V)^{\an}\to U^{\an}\times V^{\an}$ est une équivalence d'homotopie
(le produit à droite \'etant le produit usuel d'espaces topologiques), et
$\gtop(U\times V,(x,y))\to\gtop(U,x)\times\gtop(V,y)$ est un isomorphisme.\\
En appliquant la proposition~\ref{andre114} à $U\subset X$, $V\subset
Y$ et $U\times V\subset X\times Y$, on obtient que \[\gtop(X\times
Y,(x,y))\to \gtop(X,x)\times \gtop(Y,y)\] est un isomorphisme.\end{proof}

\begin{prop}\label{prod}Si $x$ et $y$ sont des points géométriques de $X$ et $Y$, alors
\[\gtemp(X\times Y,(x,y))\to\gtemp(X,x)\times\gtemp(Y,y)\] est un isomorphisme.\end{prop}
\begin{proof} Soient $(X_i,x_i)_i$ et $(Y_j,y_j)_j$ des systèmes projectifs
  dénombrables cofinaux de revêtements étales finis galoisiens connexes géométriquement pointés de $X$ et $Y$. Alors
$(X_i\times Y_j,(x_i,y_j))_{(i,j)}$ est un système projectif cofinal de
revêtementétales finis galoisiens  connexes géométriquement pointés de $X\times Y$
d'après~\cite[lecture XIII]{sga}.\\
D'après le lemme~\ref{lemprod}, $(X_i\times Y_j)^{\infty}=X^{\infty}_i\times Y^{\infty}_j$ et
\[\Gal((X_i\times Y_j)^{\infty}/(X\times Y))=\Gal(X_i^{\infty}/X)\times
\Gal(Y_j^{\infty}/Y).\]
Ainsi, en prenant la limite projective pour $(i,j)\in \mbf N^2$ des
isomorphismes précédents, on obtient le résultat voulu.\end{proof}

\section{Groupe fondamental tempéré abélianisé}
Dans ce paragraphe, nous nous intéressons à l'abélianisé du groupe
fondamental tempéré d'une courbe sur un corps non archimédien algébriquement clos de
caractéristique nulle. Nous montrerons qu'il est isomorphe au groupe
fondamental tempéré de la jacobienne de la courbe (l'isomorphisme étant
induit par le morphisme de la courbe dans sa jacobienne).\\
Nous commencerons par décrire le groupe fondamental tempéré d'une variété
abélienne en termes de l'uniformisation de cette variété abélienne.\\
Dans un second temps, nous nous restreindrons aux groupes
$(p')$-tempérés. La description du groupe $(p')$-tempéré d'une courbe en
terme d'un graphe de groupe d'une réduction semistable nous donne tout de
suite une description de son abélianisé. On en déduira une variante $(p')$
du résultat de comparaison voulu.\\
Enfin, nous montrerons directement que le morphisme du groupe fondamental
tempéré abélianisé de la courbe vers celui de sa jacobienne est un
isomorphisme (nous utiliserons l'invariance birationnelle et la formule de
Künneth).

\subsection{Groupe fondamental tempéré d'une variété abélienne}\label{varab}
Soit $K$ un corps valué non archimédien complet algébriquement clos de caractéristique $0$.\\
Soit $A$ une variété abélienne sur $K$, de dimension $g$. Décrivons son
groupe fondamental tempéré.\\

D'après~\cite{FvdP}, il existe un groupe algébrique commutatif $G$ (plus
précisément une variété semi-abélienne) et un morphisme
analytique surjectif $u:G^{\an}\to A^{\an}$ qui est le revêtement topologique
universel de $A$ et $\ker u$ est un sous-groupe discret $\Lambda$, libre de rang
$d$.\\
On sait de plus que $(A^{(n)}\to A)_{n\in \mathbf N}$ ($\mathbf N$ est ordonné par la divisibilité), où $A^{(n)}$ est une
copie de $A$ et où $A^{(n)}\to A$ est la multiplication par $n$, est une
famille cofinale de revêtements finis galoisiens pointés de $A$ (d'après par
exemple~\cite[exposé XI]{sga}).\\
Soit $G^{(n)}$ le revêtement universel de $A^{(n)}$ (qui est isomorphe à $G$
puisque $A^{(n)}$ est isomorphe à $A$) : on a
\[\gtemp(A)=\varprojlim_n \Gal(G^{(n)}/A).\]\\
$\gtop(A^{(n)})=\gtop(A)=\Lambda \simeq \mathbf Z^d$ est résiduellement
fini pour tout $n$. $\Gal(G^{(n)}/A)$ a un sous-groupe d'indice fini qui
est résiduellement fini donc est résiduellement fini. $\gtemp(A)$ est donc
résiduellement fini en tant que limite projective de groupes résiduellement
finis. Ainsi
$\gtemp(A)$ s'injecte dans son complété profini ,$\ga(A)$, qui est abélien, donc $\gtemp(A)$ est
lui-même abélien.\\

Si $n|m$, on a le diagramme commutatif suivant :
\[\xymatrix{G^{(m)} \ar[d] \ar[r] & G^{(n)} \ar[d] \ar[r] & G \ar[d]\\
  A^{(m)} \ar[r] & A^{(n)} \ar[r] & A},\]
où la flèche $G^{(m)}\to G^{(n)}$ n'est autre que la multiplication  par $m/n$
(c'est un revêtement étale fini).\\

D'où un diagramme commutatif de suites exactes :
\[\xymatrix{0 \ar[r] & \Gal(G^{(m)}/G) \ar[d] \ar[r] & \Gal(G^{(m)}/A)
  \ar[d] \ar[r] & \Gal(G/A)=\Lambda \ar@{=}[d] \ar[r] & 0\\
0 \ar[r] & \Gal(G^{(n)}/G) \ar[r] & \Gal(G^{(n)}/A) \ar[r] &
\Gal(G/A)=\Lambda \ar[r] & 0}.\]

En passant à la limite projective, on obtient la suite exacte suivante
(l'exactitude à droite résulte de la finitude des $\Gal(G^{(n)}/G)$) :
\begin{equation}\xymatrix{0 \ar[r] & \varprojlim \Gal(G^{(n)}/G) \ar[r] & \gtemp(A)
  \ar[r] & \Lambda \ar[r] & 0},\label{seA}\end{equation}
avec $\gtemp(A)$ abélien. C'est donc en fait une suite exacte de groupes
abéliens.\\
Or $\Lambda$ est un groupe abélien libre donc la suite exacte est scindée,
d'où un isomorphisme non canonique :
\[\gtemp(A)\simeq \Lambda \times \varprojlim \Gal(G^{(n)}/G).\]
Notons $T(G)= \varprojlim \Gal(G^{(n)}/G)$. C'est un groupe abélien
profini, donc il se décompose canoniquement comme produit de ses
pro-$l$-Sylows $T_l(G)$.\\
En prenant le complété pro-$l$ de l'isomorphisme ci-dessus, on obtient :
\[\mathbf Z_l^{2g}\simeq \mathbf Z_l^d \times T_l(G),\]
d'où $T_l(G)\simeq\mathbf Z_l^{2g-d}$.
En résumé, on obtient donc un isomorphisme non canonique :
\[\gtemp(A)\simeq \mathbf Z^d \times \widehat{\mathbf Z}^{2g-d}.\]
\begin{rem}
On a une description analogue pour le groupe fondamental $\mbb
L$-temp\'er\'e (qui d'ailleurs peut être reconstruit, comme pour les
courbes, à partir du groupe fondamental temp\'er\'e puisque le groupe
fondamental topologique de tout revêtement \'etale fini de $A$ est sans
torsion).\\

Soit $K_0$ un sous-corps complet de $K$ \`a valuation discr\`ete. Soit
$A_{K_0}$ une vari\'et\'e ab\'elienne sur $K_0$ telle que
$A=A_{K_0}\times_{K_0}K$ et supposons également $A_{K_0}$ \`a r\'eduction
semistable. Alors $T_l(G)$ est appel\'e la \emph{partie fixe} de $T_l(A)=\ga(A)$
et est not\'ee $T_l(A)^f$ dans~\cite[exposé IX]{sga7} (cf. \cite[\S
IX.7]{sga7} ; si $A_{K_0}$ n'est plus suppos\'ee \`a r\'eduction 
semistable, $T_l(G)$ est la partie essentiellement fixe
$T_l(A)^{ef}$). Cette partie fixe s'interpr\`ete aussi comme les
\'el\'ements fixes par l'action de l'inertie sur $T_l(A)$.\\
Remarquons qu'on peut reconstruire $T_l(A)^{ef}$ \`a partir de $\gtemp(A)^l$
puisque c'est son unique sous-groupe compact maximal. On peut \'egalement
caract\'eriser la partie essentiellement torique $T_l(A)^{et}$ en termes de
groupes fondamentaux temp\'er\'es et de l'accouplement de dualit\'e puisque c'est l'orthogonal de la partie
essentiellement fixe du module de Tate de la vari\'et\'e duale.
\end{rem}

\subsection{Abélianisation du groupe fondamental $(p')$-tempéré d'une
  courbe}
Soit $p$ la caractéristique résiduelle de $K$.\\
Soit $K$ un corps complet à valuation discrète, de complété de clôture algébrique
$\overline K$. Soit $C$ une courbe lisse projective géométriquement connexe
sur $K$ de genre $g$, et soit $C_{\bar s}$ une réduction
semi-stable de $C_{\bar K}$ de nombre de cycles $h$. On note $\mcal G$ le
graphe de pro-$(p')$-groupes associ\'e (\ref{mochipart}). On a alors un isomorphisme~:
\[\gtemp(C_{\bar K})^{(p')}\simeq \gtemp(\mathcal G).\]
Si $G$ est un groupe topologique prodiscret ayant une base dénombrable
de voisinage de $1$, on note $G^{\ab}$ le groupe
topologique $G/\overline{D(G)}$ (où $\overline D(G)$ est l'adhérence du
sous-groupe dérivé de $G$), c'est aussi un groupe prodiscret, et $G\to
G^{\ab}$ identifie $G^{\ab}-\Ens$ à une sous-catégorie pleine de $G-\Ens$.\\

Soit $\mathring{\mcal G}$ le graphe de groupes obtenus en remplaçant les
composantes $G_x$ de $\mathcal G$ par les groupes discrets
sous-jacents $\mathring{G}_x$. Soit $\get(\mring{\mcal G})$ son groupe fondamental. $\Btemp(\mathcal G)$ est une sous-catégorie
pleine de $\mathcal B(\mathring{\mathcal G})$, et $\gtemp(\mathcal
G)^{\ab}\tEns$ est une sous-catégorie pleine de $\get(\mathring{\cal G})^{\ab}\tEns$.\\
Alors, si $\mathbb T$ est un sous-arbre maximal de $\mathbb G$,
on dispose d'une présentation de $\get(\mathring{\mathcal G})$ de la forme suivante (voir par
exemple~\cite[2.2.3]{stix}) :
$$(\Astt_{e\notin \mathbb T}\mathbf Ze\ast\Astt_v \mathring G_v)/H,$$
où $H$ est le sous-groupe normal engendré par :
\begin{itemize}
\item $b_{1*}(a)b_{2*}(a^{-1})$ pour toute arête $e$ de $\mathbb T$ (il y
  en a $h$), de branches $b_1$ et $b_2$, et pour tout
$a\in G_e$,
\item $b_{1*}(a)eb_{2*}(a^{-1})e^{-1}$ pour chaque arête $e \notin \mathbb T$, de branches $b_1$ et $b_2$, et pour tout
$a\in G_e$ (on a choisit pour chacune de ces arêtes une orientation arbitraire).
\end{itemize}
On obtient, en envoyant les $\mathring G_v$ sur l'élément neutre un quotient
discret $\Astt_{e \notin \mathbb T}\mathbf Ze$ qui correspond au revêtement
topologique universel de $\mathbb G$ (c'est-à-dire $\gtop(\mathcal G)=\Astt_{e \notin \mathbb
  T}\mathbf Ze$).\\
En abélianisant cette présentation, on en déduit que $$\get(\mathring{\mathcal
G})^{\ab}=\prod_{e\notin T} \mathbf
Ze\times(\prod_v\mathring{G}^{\ab}_v/H),$$ où $H$ est
le sous-groupe engendré par les $b_{1*}(a)b_{2*}(a^{-1})$ et
$\prod_{e\notin T} \mathbf Ze=\gtop(\bb G)^{\ab}$.\\
Comme $G_v$ agit continûment sur les revêtements tempérés abéliens, le
morphisme \[\mathring G_v \to \gtemp(\cal G)^{\ab}\] se
factorise par $$G_v\to \gtemp(\cal G)^{\ab},$$ et donc $$\prod_{e\notin T}
\mathbf Ze\times(\prod_v\mathring{G}^{\ab}_v/H)\to \gtemp(\cal G)^{\ab}$$ se
factorise par \[\prod_{e\notin T} \mathbf
Ze\times(\prod_vG^{\ab}_v/\overline{H})\to \gtemp(\cal G)^{\ab}.\]\\
En notant $G_0=\prod_vG^{\ab}_v/\overline{H}\simeq (\mathbf Z^{(p')})^{2g-h}$ (où $\mathbf
Z^{(p')}=\prod_{q\neq p}\mathbf Z_q$), on a
donc un morphisme d'image dense (qui correspond à un foncteur pleinement
fidèle sur les catégories classifiantes correspondantes) $$G_1:=G_{0}\times
\mathbf Z^h \to \gtemp(\mathcal G)^{\ab},$$
(et le quotient $G_1\to \mathbf Z^h$ identifie $\mbf Z^h$ à $\gtop(\bb G)^{\ab})$.\\
Plus canoniquement, on a une suite exacte scindée de groupes abéliens
topologiques~:
\begin{equation}0\to G_0 \to G_1 \to \gtop (\mathcal G)^{\ab}\to 0.\label{seC}\end{equation}\\

Soit $S$ un  $G_1$-ensemble connexe (qui correspond à un revêtement de
$\mathring{\cal G})$ ; comme chaque composante $G_v$ agit continûment sur
$S$, il correspond en 
fait à un revêtement de $\cal G$
), correspondant à un
sous-groupe ouvert $U$ de $G_1$.\\
Soit $U_0=U\cap G_0$, c'est un sous-groupe ouvert de $G_0$, mais aussi de
$G_1$. Il correspond donc à un $G$-ensemble connexe $S_0$ qui domine $S$ car
$U \subset U_0$.\\
Soit $U_1=U_0\times \mathbf Z^h \subset G_1$, c'est un sous-groupe ouvert
d'indice fini de $G_1$, correspondant donc à un $G_1$-ensemble $S_1$
fini connexe. L'inclusion $U_0 \subset U_1$ induit un morphisme $S_0 \to
S_1$, qui est scindé en toute composante $G_v$ : $S_0\to S_1$ est donc un
revêtement topologique, donc $S$ est un revêtement tempéré.\\

Tout revêtement étale abélien de $\mathcal G$ est donc aussi tempéré (et il
est même dominé par un revêtement abélien qui est topologique sur
un revêtement fini), d'où un isomorphisme :
$$\gtemp(C_{\bar K})^{(p') \ab}=\gtemp(\mathcal G)^{\ab}\simeq (\mathbf Z^{(p')})^{2g-h}\times \mathbf Z^h.$$
\begin{rem} On aurait aussi pu utiliser
directement~\cite[Th. III.2.1.9]{andre1} pour montrer que tout
$G_1$-ensemble connexe est tempéré, car tout quotient discret de $G_1$ est
virtuellement sans torsion.
\end{rem}

\subsection{Jacobienne et $\gtemp^{(p'),\ab}$ d'une courbe}
Soit $K$ un corps de caract\'eristique 0 complet pour une valuation discrète, soit
$\bar K$ le complété de sa clôture algébrique. Soit $p$ la caractéristique
résiduelle. Soit $C$ une courbe projective lisse géométriquement connexe sur $K$,
et soit $A$ la jacobienne de $C_{\bar K}$.\\
Si $x$ est un point fermé de $C_{\bar K}$, on a un morphisme canonique
$C_{\bar K}\to A$ qui envoie $x$ sur l'élément neutre de $A$. Il induit un morphisme $\gtemp(C_{\bar
  K})^{(p'),\ab}\to \gtemp(A)^{(p')}$, qui ne dépend pas de $x$.

\begin{prop}Le morphisme $\gtemp(C_{\bar K})^{(p'),\ab}\to \gtemp(A)^{(p')}$
  est un isomorphisme.\end{prop}
\dem
Les deux suites exactes~\eqref{seA} et~\eqref{seC} fournissent un diagramme
commutatif (où l'on a conservé les notations précédentes) :
$$\xymatrix{0 \ar[r] & G_0\ar[r] \ar[d]^a & \gtemp(C_{\bar K})^{(p'),\ab} \ar[r] \ar[d]^b &
  \gtop(C_{\bar K}) \ar[r] \ar[d]^c & 0\\ 0 \ar[r] &  T_{(p')}(G) \ar[r] &
  \gtemp(A)^{(p')} \ar[r] & \gtop(A) \ar[r] &
  0},$$
où $b$ et $c$ se déduisent de la fonctorialité du $\gtemp^{(p')}$ et du
  $\gtop$, et $a$ s'obtient par passage au noyau.\\
En passant au complété pro-$(p')$, on obtient le diagramme suivant ($T_{(p')}(G)$ et $G_0$ sont déjà pro-$(p')$) :
$$\xymatrix{0 \ar[r] & G_0 \ar[r] \ar[d]^a & \ga(C_{\bar K})^{(p'),\ab} \ar[r]
  \ar[d]^{b'} &
  \gtop(C_{\bar K})^{\ab}_{(p')} \ar[r] \ar[d]^{c'} & 0\\ 0 \ar[r] & T_{(p')}(G) \ar[r] &
  \ga(A)^{(p')} \ar[r] & \gtop(A)_{(p')} \ar[r] &
  0},$$
où, si $\Pi$ est un groupe topologique, $\Pi_{(p')}$ désigne ici le
  complété pro-$(p')$.\\

Rappelons que $b'$ est un isomorphisme d'après la théorie de Kummer.\\
Montrons que $c'$ est injectif.\\
Soit $\bar y\in \Ker c'$, et $y$ un antécédent de $\bar y$ dans
$\ga(C_{\bar K})^{(p'),\ab}$. Alors $G_0$ et $<y>$ sont en somme directe dans
$\ga(C_{\bar K})^{(p'),\ab}$. Si $l$ est un nombre premier tel que $y_l$ soit non nul,
$(G_0\oplus <y>)_l$ est un $\mathbf Z_l$-module de rang $2g-h+1$, et
son image par l'isomorphisme $b'$ est incluse dans $T_l(G)$ qui est de
rang $2g-h$ sur $\mathbf Z_l$. D'où une contradiction, et $c'$ est bien
injectif.\\
On a le diagramme commutatif :
$$\xymatrix{\gtop(C_{\bar K})^{\ab} \ar@{^{(}->}[r] \ar[d]^c &
  \gtop(C_{\bar K})^{\ab}_{(p')} \ar@{^{(}->}[d]^{c'}\\ \gtop(A)
  \ar@{^{(}->}[r] & \gtop(A)_{(p')},}$$
d'o\`u on déduit immédiatement l'injectivité de $c$.\\

Si $c$ n'est pas surjectif, son image $\Image c$ dans $\gtop(A)$ est un
sous-groupe strict de $\gtop(A)\simeq \mathbf Z^h$ et est alors inclus dans
un sous-groupe strict d'indice fini. On en déduit que le complété profini
$\hat c$ de $c$ n'est pas non plus surjectif.\\
Or on a le diagramme commutatif :
$$\xymatrix{\ga(C_{\bar K})^{\ab} \ar[r] \ar@{=}[d] & \widehat{\gtop(C_{\bar
      K})^{\ab}} \ar[d]^{\hat c}\\ \ga(A) \ar@{->>}[r] &\widehat{\gtop(A)}},$$
qui montre la surjectivité de $\hat c$. Donc $c$ est bien surjectif.\\
$c$ est donc bijectif (et donc un isomorphisme puisque les groupes de
      départ et d'arrivée sont discrets).\\

Le morphisme $c'$ est donc aussi un isomorphisme et $b'$ aussi, donc $a=\Ker(b'\to c')$
est un isomorphisme. Comme
$a$ et $c$ sont des isomorphismes, $b$ est également un isomorphisme.
\findem

\subsection{Jacobienne et $\gtemp^{\ab}$ d'une courbe}
Soit $K$, $C$ et $A$ comme au paragraphe pr\'ec\'edent.\\

Soit $x_0$ est un point fermé de $C_{\overline K}$. Il induit un morphisme de $C_{\overline K}$ vers sa jacobienne $A$\\

\begin{thm}\label{abel} Le morphisme $\gtemp(C_{\bar K},x_0)^{\ab}\to\gtemp(A,0)$ est un isomorphisme.\end{thm}
\dem
On a un morphisme $C_{\bar K}^g\to A$ qui à $(x_1,\dots,x_g)$ associe le
diviseur de $C_{\bar K}$ $[x_1]+\cdots+[x_g]-g[x_0]$. Ce morphisme est
invariant par l'action de $\fk S_g$ sur $C_{\bar K}^g$ et se factorise
donc en un morphisme $C_{\bar K}^{(g)}:=C_{\bar K}^g/\fk S_g\to
A$. Rappelons que ce morphisme est birationnel et que $C_{\bar K}^{(g)}$
est lisse sur $\bar K$ (voir~\cite[Th. 5.1.(a), prop. 3.2]{jacmilne}).\\
On a donc une suite de morphismes~:
\[C_{\bar K}\to C_{\bar K}^g\to C_{\bar K}^{(g)}\to A,\]
où le morphisme de gauche envoie $x$ sur $(x,x_0,\dots, x_0)$ et où la
composée est le morphisme qui envoie $x$ sur $[x]-[x_0]$, comme à
la section précédente.\\

$C_{\bar K}^{(g)}\to A$ est un morphisme birationnel entre $\bar
K$-variétés lisses et propres, donc d'après~\ref{birat}, $\gtemp(C_{\bar
  K}^{(g)},(x_0,\dots,x_0))\to\gtemp(A,0)$ est un isomorphisme.\\
On en déduit en particulier que $\gtemp(C_{\bar K}^{(g)},(x_0,\dots,x_0))$ est
abélien et résiduellement fini.\\
Donc $\gtemp(C_{\bar K},x_0)\to\gtemp(C_{\bar K}^{(g)},(x_0,\dots,x_0))$ se
factorise à travers $\phi:\gtemp(C_{\bar K},x_0)^{\ab}\to\gtemp(C_{\bar
  K}^{(g)},(x_0,\dots,x_0))$.\\

Comme $\gtemp(C_{\bar K},x_0)^{\ab}$ est une limite projective de groupes
abéliens de type fini donc résiduellement finis, il est lui-même
résiduellement fini.\\
Le diagramme commutatif suivant~:
\[\xymatrix{\gtemp(C_{\bar K},x_0)^{\ab} \ar[r]^{\phi} \ar[d] & \gtemp(C_{\bar
    K}^{(g)},(x_0,\dots,x_0)) \ar[d]\\ \ga(C_{\bar K},x_0)^{\ab} \ar[r]^{\simeq} & \ga(C_{\bar
    K}^{(g)},(x_0,\dots,x_0))},\]
dont les flèches verticales sont injectives, puisque $\gtemp(C_{\bar
  K},x_0)$ et $\gtemp(C_{\bar
    K}^{(g)},(x_0,\dots,x_0))$ sont résiduellement finis, montre que $\phi$ est
injectif.\\

On peut munir également $C_{\bar K}^{(g)}$ d'une structure d'orbifold en
tant que $C_{\bar K}^g/\fk S_g$ (voir~\cite[\S{} III.4.4.1]{andre1}), et
définir un groupe fondamental tempéré d'orbifold correspondant (voir~\cite[\S{} III.4.5]{andre1}).\\
On a le diagramme dont la ligne est exacte
(voir~\cite[prop. III.4.5.8]{andre1} pour la ligne supérieure):
\[\xymatrix{
1 \ar[r] & \gtemp(C_{\bar K}^g,(x_0,\dots,x_0)) \ar[r]^i \ar[ddr]_\alpha &
\gorb(C_{\bar K}^{(g)},(x_0,\dots,x_0)) \ar[r]\ar[d]^{\pi_1} & \fk S_g \ar[r] & 1\\
& & \gorb(C_{\bar K}^{(g)},(x_0,\dots,x_0))^{\ab} \ar[d]^{\pi_2} 
& & \\  & & \gtemp(C_{\bar K}^{(g)},(x_0,\dots,x_0)) & &}\]
Les morphismes $i$, $\pi_1$ et $\pi_2$ sont ouverts, donc $\alpha=\pi_2 \pi_1 i$ est ouvert.\\
Or $\alpha$ est d'image dense car $\hat \alpha$ est surjectif (puisque $\ga(C_{\bar K})\to\ga(C_{\bar K}^{(g)})$, qui se
factorise par $\hat \alpha$ l'est) et $\gtemp(C_{\bar
  K}^{(g)},(x_0,\dots,x_0))$ est résiduellement fini, donc $\alpha$ est
également surjectif.\\

Considérons maintenant $\delta:\gtemp(C_{\bar K},x_0)\to\gtemp(C_{\bar
K}^{g},(x_0,\dots,x_0))$. En identifiant $\gtemp(C_{\bar
K}^{g},(x_0,\dots,x_0))$ à $\gtemp(C_{\bar
K},(x_0,\dots,x_0))^g$, $\delta$ peut être identifié à $\gtemp(C_{\bar
K},x_0)\to\gtemp(C_{\bar K},x_0)^g:g\mapsto (g,1,\dots,1)$.\\
Ainsi $\gtemp(C_{\bar
K}^{g},(x_0,\dots,x_0))$ est engendré par les images de $\sigma\circ\delta$
quand $\sigma$ décrit $\fk S_g$, et comme $\alpha$ est invariant par
$\sigma$, $\alpha\circ\delta$ est surjectif donc $\gtemp(C_{\bar
  K},x_0)^{\ab}\to\gtemp(C_{\bar K}^{(g)},(x_0,\dots,x_0))$ aussi et c'est donc
un isomorphisme.\\
De même si $U$ est un sous-groupe ouvert de $\gtemp(C_{\bar K},x_0)$,le
groupe engendré par les $\sigma(\delta(U))$ est un groupe ouvert de $\gtemp(C_{\bar
K}^{g},(x_0,\dots,x_0))$, donc, comme $\alpha$ est ouvert et $\fk
S_g$-invariant, $\alpha\circ\delta$ est ouvert, donc $\gtemp(C_{\bar
  K},x_0)^{\ab}\to\gtemp(C_{\bar K}^{(g)},(x_0,\dots,x_0))$ est aussi ouvert donc
un isomorphisme.\findem

\chapter{Métrique du graphe de la réduction stable d'une courbe de Mumford}
Le but principal de ce chapitre est de prouver que le graphe métrisé
associé au modèle stable d'une courbe de Mumford (ou, de manière équivalente, du graphe
sous-jacent au squelette de l'espace de Berkovich de la courbe,
cf.~\ref{berkspaces}) ne dépend que du groupe fondamental tempéré de la courbe.\\
Nous considérerons seulement le cas de caractéristique mixte (dans le cas
d'égale caractéristique nulle, c'est certainement faux puisque le groupe
fondamental tempéré dépend uniquement du graphe de groupe associé à la
courbe selon~\cite[ex. 3.10]{mochi}).
\begin{dfn} Une \emph{métrique} sur un graphe $\mbb G$ est une
  fonction \[d:\{\text{arêtes de } \mbb G\}\to \mbf R^{>0}. \]Pour toute arête $e$,
  $d(e)$ est appelé la \emph{longueur} de $e$ pour la métrique $d$.\\
Un graphe muni d'une métrique est appelé un \emph{graphe métrique}.\end{dfn}
Soit $X=\overline X\backslash D$ une courbe sur un corps complet non archimédien algébriquement
clos $K$, munie d'un modèle semistable $(\mcal X,\mcal D)$ (cf. \ref{mochipart}). Soit $e$ une arête du
graphe de la réduction semistable (\ie un point double de $\mcal
X_s$). Alors, localement pour la topologie étale en ce point double, $\mcal
X$ est étale sur $O_K[X_0,X_1]/(X_0X_1-a)$, avec $a\in
O_K$. D'après~\cite[cor. 2.2.18]{thuillier}, $|a|$ ne dépend d'aucun choix.
Nous noterons alors
\[d(e)=-\log_p(|a|),\]
ce qui définit une métrique naturelle sur le graphe de cette réduction
semistable de $X$.\\
Par exemple, si $\mcal X$ est le modèle stable de $\mbf
P^1\backslash\{0,1,\infty,\lambda\}$ avec $|\lambda|<1$, alors le graphe de
$\mcal X$ a une unique arête de longueur $-\log_p(|\lambda|)$.\\ 

Nous savons déjà, grâce à~\cite[ex. 3.10]{mochi}, que l'on peut
reconstruire le graphe de la réduction stable de la courbe à partir du
groupe fondamental tempéré. On peut déduire de l'étude de Mochizuki que
l'on peut déterminer, pour tout sous-groupe ouvert d'indice fini du groupe
fondamental tempéré et pour tout sommet du graphe de la réduction stable de
la courbe, si le revêtement correspondant de la courbe est totalement décomposé
en ce sommet (un revêtement $X'\to X$ de variétés est dit
\emph{totalement décomposé} en un
point $x\in X$ si $\mcal H(x)\to\mcal H(x')$ est un isomorphisme pour tout
$x'\in X'$ au-dessus de $x$; un revêtement $X'\to X$ d'ordre $n$ est
décomposé en $x$ si et seulement si la fibre de $x$ a pour cardinal $n$, ce
qui revient aussi à dire que localement au voisinage de $x$, $X'\to X$ est
un revêtement topologique~(\cite[III.1.2.1]{andre1})).\\
Ceci suggère de regarder quels revêtements étales finis de notre courbe de
Mumford sont décomposés en un sommet (ou un ensemble de sommets). En \'etudiant
des revêtements étales finis assez simples, on se convainc que la métrique
du graphe doit jouer un rôle dans le groupe fondamental tempéré.\\
Commençons par un cas élémentaire mais qui illustre le rôle que peut avoir la
distance dans le fait qu'un revêtement soit décomposé en un point.
\begin{lem}\label{decompositionpuissance}Le revêtement $\mbf
  G_m\stackrel{z\mapsto z^{p^h}}{\to}\mbf G_m$ est totalement décomposé au
  dessus du point de Berkovich $b(1,r)$ (correspondant à la boule ouverte $B(1,r)$ de
  centre 1 et de rayon $r$ avec $r<1$) si et seulement si $r<p^{-h-\frac{1}{p-1}}$.\\
Plus précisément, l'image réciproque de $b(1,r)$ a $p^i$ éléments~:
\begin{itemize}
\item avec $i=0$ quand $r\in ]p^{-\frac{p}{p-1}},1]$;
\item quand $r\in ]p^{-i-\frac{p}{p-1}},p^{-i-\frac{1}{p-1}}]$ et $1\leq i\leq h-1$;
\item avec $i=h$ quand $r\in [0,p^{-h-\frac{1}{p-1}}]$.
\end{itemize}\end{lem}
\dem Nous 
Posons $g:z\mapsto z^p$, et calculons $g(b(z_1,r))$ avec $|z_1|=1$ et $r<1$.
Posons $g_{1}:z\mapsto (z+z_1)^p-z_1^p=\sum_{i=1}^p a_iz^i$ avec $|a_p|=1$
et $|a_i|=p^{-1}$ si $1\leq i\leq p-1$.\\
Si $f\in\mbf C_p[X]$,
\[|f|_{g(b(z_1,r))}=|f\circ g|_{b(z_1,r)}=\sup_{x\in B(z_1,r)_{K'}}|f\circ
  g(x)|=\sup_{y\in g(B(z_1,r)_{K'})}|f(y)|,\]
où $K'$ est une extension isométrique algébriquement close sphériquement complète de $\mbf
C_p$.\\
Or $g(B(z_1,r)_{K'})=B(z_1^p,r')_{K'}$ avec 
\[r'=\sup_{|z|<r}|g_1(z)|=\max
|a_i|r^i=\max\{r^p,rp^{-1}\}=\left\{\begin{array}{ll} rp^{-1} & \text{if
    }r\leq p^{-\frac{1}{p-1}}\\ r^p& \text{if }r\geq
    p^{-\frac{1}{p-1}}\end{array}\right.\]
Donc \[|f|_{g(b(z_1,r))}=\sup_{y\in
  B(z_1^p,r')_{K'}}|f(y)|=|f|_{b(z_1^p,r')}\] et donc $g(b(z_1,r))=b(z_1^p,r')$.\\
De plus, soit $z_0$ de norme $1$, soit $z_0^{1/p}$ une racine $p$\ieme{} de
$z_0$ et soit $r'<1$. Posons
\[r=\left\{\begin{array}{ll} pr' & \text{si }r'\leq p^{-\frac{p}{p-1}}\\ r'^{1/p}
& \text{si }r'\geq p^{-\frac{p}{p-1}}\end{array}\right.\]
Alors $b(z_0^{1/p},r)\in g^{-1}(b(z_0,r'))$. Les autres  éléments de
$g^{-1}(b(z_0,r'))$ sont les conjugués par l'action du groupe de Galois
$\mu_p$ du revêtement $g:z\mapsto z^p$. Or, si $\zeta\in\mu_p$, alors $\zeta
b(z^{1/p}_0,r)=b(\zeta z^{1/p}_0,r)$. Donc 
\[g^{-1}(B(z_0,r'))=\left\{\begin{array}{ll} \{B(\zeta z_0^{1/p},pr')\}_{\zeta\in\mu_p} & \text{si }r'\leq p^{-\frac{p}{p-1}}\\ \{B(\zeta z_0^{1/p},r'^{1/p})\}_{\zeta\in\mu_p} & \text{si }r'\geq p^{-\frac{p}{p-1}}\end{array}\right.\]
On obtient que $g^{-1}(B(z_0,r'))$ a un unique élément si $r'\geq
p^{-\frac{p}{p-1}}$ et $p$ \'el\'ements sinon (car $|\zeta-\zeta'|=p^{-\frac{1}{p-1}}$
si $\zeta\neq\zeta'\in\mu_p$). On obtient donc le résultat voulu quand $h=1$.\\
Dans le cas général, on raisonne par récurrence sur $h$ en décomposant $z\mapsto z^{p^h}$ en $z\mapsto z^{p^{h-1}}\mapsto z^{p^h}$.\findem

Grâce à ce résultat, nous pourrons étudier le cas d'une droite projective
privée de quelques points et, par découpage et recollage de ce type de
revêtements, le cas d'une courbe elliptique épointée.\\
Dans le cas plus général d'une courbe de Mumford $X$ (th\'eor\`eme \ref{mumford}), nous étudierons aussi
la structure des $\mbf Z/p^h\mbf Z$-torseurs. La théorie des fonctions
thêta telle qu'on peut la trouver dans~\cite{vdp1} et~\cite{vdp2} nous dit
que le pullback d'un tel torseur au revêtement topologique universel
$\Omega$ est en fait le pullback de $\mbf G_m\stackrel{z\mapsto
  z^{p^h}}{\to}\mbf G_m$ le long d'une fonction thêta $\Omega\to\mbf G_m$,
qui à son tour correspond à un courant sur l'arbre $\mbb T(\Omega)$ de
$\Omega$ (proposition~\ref{torseurstheta}). Ainsi, nous commencerons notre étude en prouvant que si deux
courants coïncident sur une partie "assez grosse" de $\mbb T(\Omega)$, le
quotient des deux fonctions inversibles associées varie assez peu sur une
partie légèrement plus petite de $\mbb T(\Omega)$ et donc les deux $\mbf
Z/p^h\mbf Z$-torseurs seront décomposés en les mêmes points de cette partie
de $\mbb T(\Omega)$. Ainsi, nous considérerons des courants qui, sur une
"assez grosse" partie de $\Omega$, coïncident avec le courant correspondant à
une homographie inversible sur $\Omega$ et qui est équivariant pour un
sous-groupe d'indice fini de $\Gal(\Omega/X)$. Nous considérerons ensuite le
$\mbf Z/p^h\mbf Z$-torseur sur un certain revêtement topologique fini de
$X$ qui se comporte comme  $\mbf G_m\stackrel{z\mapsto z^{p^h}}{\to}\mbf
G_m$ sur une certaine partie de $\mbb T(\Omega)$ (d'apr\`es corollaire~\ref{totdecrev}). Nous en déduirons que la longueur de tout cycle de tout revêtement
topologique du graphe du modèle stable de $X$ peut-être retrouvé à partir
du groupe fondamental tempéré. Un dernier argument combinatoire 
donnera le résultat souhaité.

\section{Préliminaires}
Soient $K$ un corps valué complet non archimédien de caractéristique 0, $O_K$
son anneau des entiers, $k$ son corps résiduel, supposé de caractéristique
$p>0$, et soit $\overline K$ le complété d'une clôture algébrique de 
$K$ (d'anneau des entiers $O_{\overline K}$).\\
On fixe un système compatible de racines de 1
dans $\overline K$, ce qui nous permet d'identifier $\mu_n$ et $\mbf Z/n\mbf Z$. Nous parlerons donc de $\mbf Z/n\mbf Z$-torseurs sur une courbe sur
$\overline K$ quand il serait plus naturel de parler de $\mu_n$-torseurs.\\
Soient $X_{1,K}$, $X_{2,K}$ deux courbes hyperboliques de type $(g,n)$ sur
$K$.
Pour $i=1,2$,
nous noterons $\Pi_i:=\gtemp(X_{i,\overline K})$.
Soit $\phi$ un isomorphisme $\Pi_1\simeq\Pi_2$.\\
Soit $H_1$ un sous-groupe ouvert d'indice fini de $\gtemp(X_{1,\overline
  K})$ et soit $H_2:=\phi(H_1)$.\\

Soit $Y_{i,\overline K}\to X_{i,\overline K}$ le
revêtement étale fini connexe de $X_{i,\overline K}$ correspondant à $H_i$.\\
Soient $X_i$ et $Y_i$ les modèles stables de $X_{i,\overline K}$ and
$Y_{i,\overline K}$, et soit $\psi_i:Y_i\to X_i$ l'unique morphisme
étendant $Y_{i,\overline K}\to X_{i,\overline K}$.
On notera $X_{i,s}$ lafibre spéciale de $X_i$ et $X_{i,\eta}$ la fibre
générique de $X_i$ (\ie $X_{i,\overline K}$).\\

D'après~\ref{mochi}, $\phi$ induit un isomorphisme $\mcal G^c_1\simeq\mcal
G^c_2$ entre les semigraphes d'anabélioïdes de $X_1$ et $X_2$. De même,
$H_1\simeq H_2$ induit un isomorphisme $\mcal H^c_1\to\mcal H^c_2$ entre
les semigraphes d'anabélioïdes de $Y_1$ et $Y_2$.\\
 Un \emph{point cuspidal}\index{Point cuspidal} d'une courbe hyperbolique de type
$(g,n)$ est un des $n$ points du complété de la courbe qui n'est pas dans
la courbe elle-même.\\
Nous nous intéressons à savoir quelles données peuvent être récupérées
sur la préimage des sommets du squelette et des points cuspidaux de
$X_{i,\overline K}$ \`a partir de
$\phi:\gtemp(Y_{i,\overline K})\to\gtemp(X_{i,\overline K})$.\\

\begin{lem}\label{cusp} Soit $x_1$ un point cuspidal de $X_{1,s}$, et $x_2$ le point
  cuspidal de $X_{2,s}$ correspondant à $x_1$ par l'isomorphisme de
  semigraphes $\mbb G^c_1\simeq\mbb
  G^c_2$. Soit $x'_i$ le point cuspidal correspondant de
  $X_{i,\eta}$. Alors $\psi_{1,\eta}^{-1}(x'_1)$ et
  $\psi_{2,\eta}^{-1}(x'_2)$ ont le même nombre d'éléments.\end{lem}
\begin{proof}
Soit $y_i$ un point cuspidal de $Y_{i,s}$ (correspondant à un point cuspidal
 $y'_i$ de $Y_{i,\eta}$) et soit $z_i$ son image dans
$X_{i,s}$ (correspondant à un point cuspidal $z_i'$ de $X_{i,\eta}$).\\
Supposons que $y_1$ et $y_2$ se correspondent l'un l'autre par l'isomorphisme $\mbb G^c_1\simeq\mbb G^c_2$.\\
Soit $I_i\subset H_i^{(p')}$ un sous-groupe d'inertie de $y_i$. L'image de $I_i$ dans
$\Pi_i^{(p')}$ est un sous-groupe ouvert (et donc non trivial) d'un
groupe d'inertie de $z_i$. Comme l'intersection de deux groupes d'inertie
différents est réduite à  $\{1\}$, l'image de $I_i$
n'est contenu dans aucun autre groupe d'inertie de
$X_{i,s}$, donc $z_i$ est caractérisé par le morphisme $H_i\to \Pi_i$
comme étant l'unique point cuspidal de $X_{i,s}$ tel que les groupes
d'inerties de $y_i$ s'envoient par $H_i\to \Pi_i$ sur les groupes
d'inerties de $z_i$. Or $\Pi_1\to\Pi_2$ (resp. $H_1\simeq H_2$) envoie
les groupes d'inertie d'un point cuspidal de $X_{1,s}$ (resp. $Y_{1,s}$)
sur les groupes d'inerties du point cuspidal correspondant de $X_{2,s}$
(resp. $Y_{2,s}$). Ainsi $z_2$ est le point cuspidal correspondant à $z_1$ par $H_1\to H_2$.\\
Ainsi le diagramme suivant commute~:
\[\begin{array}{ccc}\cusp(Y_{1,s}) & \simeq & \cusp(Y_{2,s})\\
\dar & & \dar\\
\cusp(X_{1,s}) & \simeq & \cusp(X_{2,s})\end{array},\]
ce qui donne le résultat voulu.\end{proof}
En particulier, le morphisme $Y_{1,\eta}\to X_{1,\eta}$ est ramifié en
$x'_1$ si et seulement si $Y_{2,\eta}\to X_{2,\eta}$ est ramifié en $x'_2$.\\
\begin{lem}\label{reductioncusps} Soit $\widetilde X_{i,\eta}$ un
  ouvert de $\overline X_{i,\eta}$ contenant $X_{i,\eta}$. Supposons
  que les points cuspidaux de $\widetilde X_{1,\eta}$ correspondent aux
  points cuspidaux de $\widetilde X_{2,\eta}$ par la bijection $\cusp(X_{1,\eta})\simeq\cusp(X_{2,\eta})$ induite par $\phi$.\\
Alors $\phi$ induit un isomorphisme $\gtemp(\widetilde{X}_{1,\eta})\to \gtemp(\widetilde X_{2,\eta})$.\end{lem}
\begin{proof}
Soit $S^{\infty}_i$ le revêtement topologique universel d'un revêtement
galoisien connexe $S_i$ de $X_{i,\eta}$. Alors le quotient discret
correspondant de $\Pi_i$ est une extension $G_i$ d'un groupe libre par un
quotient fini $G_{1,i}$ de $\Pi_i$, qui correspond au revêtement $S_i$ de $X_{i,\eta}$.\\
En outre, si $G'_i$ est un quotient discret de $\Pi_i$ (correspondant à un
revêtement tempéré de $S''_i$) c'est une extension $G'_i$ d'un groupe libre
par un quotient fini $G'_{1,i}$ (qui correspond à un revêtement fini
$S'_i$), alors $S''_i\to S'_i$ est un revêtement topologique~(\cite[th. III.2.1.9.a]{andre1}).\\
Alors, d'après la proposition~\ref{limproj},
$\gtemp(\widetilde{X}_{i,\eta})=\varprojlim_j \Pi_i/H_{i,j}$, où
$(\Pi_i/H_{i,j})_j$ sont les quotients discrets de $\gtemp(X_{i,\eta})$ qui
sont des extensions d'un groupe libre par un quotient fini de
$\gtemp(X_{i,\eta})$ correspondant à un revêtement fini qui est non ramifié
au-dessus de $\widetilde X_{i,\eta}$. D'après le lemme~\ref{cusp}, la
famille $(H_{2,j})_j$ est juste la famille $(\phi(H_{1,j}))_j)$. Cela donne
l'isomorphisme voulu.\end{proof}

\begin{lem}\label{vertex} Soit $x_1$ le point générique d'une composante
  irréductible de $X_{1,s}$, et $x_2$ le point générique de la composante
  irréductible correspondante de $X_{2,s}$ par $\mbb G^c_1\simeq\mbb
  G^c_2$. Soit $x'_i$  le point correspondant du squelette de la fibre
  générique $X_{i,\eta}$. Alors $\psi_{1,\eta}^{-1}(x'_1)$ et
  $\psi_{2,\eta}^{-1}(x'_2)$ ont le même cardinal.\end{lem}
\begin{proof}
Si $X_{i,0}$ est une composante irréductible de $X_{i,s}$, soit $Y_{i,0}$
une composante irréductible de $Y_{i,s}$ qui s'envoie surjectivement sur 
$X_{i,0}$. Alors le morphisme entre les composantes des graphes de groupes
$\Pi^{(p')}_{Y_{i,0}}\to\Pi^{(p')}_{X_{i,0}}$ est ouvert (en particulier
son image est non commutative) puisqu'il se plonge dans le diagramme
commutatif suivant
\[\xymatrix{\Gal\big(\overline{K(Y_{i,0})}/K(Y_{i,0})\big) \ar@{->>}[d] \ar@{^{(}->}[r] &
  \Gal\big(\overline{K(X_{i,0})}/K(X_{i,0})\big) \ar@{->>}[d] \\ \Pi^{(p')}_{Y_{i,0}} \ar[r] &
  \Pi^{(p')}_{X_{i,0}}},\]
où la flèche du haut est une immersion ouverte et les flèches verticales
sont des projections.\\
Puisque $\Pi^{(p')}_{X_{i,0}}\to \Pi_i^{(p')}$ (définie à conjugaison près)
est injectif, l'image de
$\Pi^{(p')}_{Y_{i,0}}$ dans $\Pi_i^{(p')}$ (défini à conjugaison près) est
non commutative, et donc
$\Pi_{X_{i,0}}^{(p')}$ est le seul sous-groupe verticiel de $\Pi_i^{(p')}$
qui contient l'image de $\Pi^{(p')}_{Y_{i,0}}$ (\`a conjugaison pr\`es).\\
De plus, si $Y_{i,0}$ est une composante irréductible de $Y_{i,s}$ qui
ne s'envoie pas surjectivement sur une composante irréductible de $X_{i,s}$,
l'image de $\Pi^{(p')}_{Y_{i,0}}$ dans $\Pi_i^{(p')}$ est commutative, donc
le morphisme $H_i\to\Pi_i$ détermine quelles composantes de $Y_{i,s}$
s'envoient surjectivement sur une composante de
$X_{i,s}$.\\
En particulier, si $x_1$ et $x_2$ sont les points génériques de composantes
irréductibles de $X_{1,s}$ et $X_{2,s}$ se correspondant l'un l'autre,
le nombre d'antécédents par $Y_{i,s}\to X_{i,s}$ de $x_i$ est indépendant de $i$.\\

Notons maintenant $\pi$ l'application de spécialisation de la fibre
générique (analytique) vers la fibre spéciale. $x_i=\pi(x'_i)$ est le point
générique d'une composante irréductible de $X_{\eta}$, toute préimage $y_i$ de $x_i$
par $\psi_i$ est un point générique d'une composante irréductible de
$Y_{\eta}$, $\pi^{-1}(y_i)$ est réduit à un unique élément
par~\cite[th. 4.3.1.(i)]{berk}, qui doit s'envoyer sur $x'_i$ puisque
$\pi^{-1}(x_i)=\{x'_i\}$ d'après~\cite[th. 4.3.1.(i)]{berk}. Ainsi
$\psi_i^{\an,-1}(x'_i)$ est en bijection naturelle avec $\psi_i^{-1}(x_i)$. \end{proof}

\section{Cas de $\mbf{P}^1\backslash\{z_1,\dots,z_n\}$}
Soit $i=1$ ou $2$.\\
Soit $z_{i,1},\cdots,z_{i,n}\in\mbf Q_p^{\text{nr}}$, avec
$n\geqslant 4$.\\
Notons $X_{i}:=\mbf{P}_{\mbf C_p}^1\backslash\{z_{i, 1},\dots,z_{i, n}\}$.\\
Soit $\Pi_{i}=\gtemp(X_i)$. Nous savons déjà qu'un isomorphisme $\phi:\Pi_{1}\simeq\Pi_{2}$ induit un isomorphisme
entre les semi-graphes des réductions stables de
$\mbf{P}^1\backslash\{z_{1,1},\dots,z_{1,n}\}$ et de $\mbf{P}^1\backslash\{z_{2,1},\dots,z_{2,n}\}$. Quitte
à réordonner les $z_{2, j}$, on peut supposer que ce morphisme
de semigraphes identifie le sous-groupe d'inertie (défini à conjugaison
près) de $z_{1,j}$ avec le sous-groupe d'inertie de $z_{2,j}$.\\

\begin{thm}\label{thmp1} L'isomorphisme de graphes déterminé par $\phi$
  entre les squelettes de
$(\mbf{P}^1\backslash\{z_{1, 1},\dots,z_{1, n}\})^{\an}$ et de
$(\mbf{P}^1\backslash\{z_{2, 1},\dots,z_{2, n}\})^{\an}$ préserve
les longueurs des arêtes (\ie induit un isomorphisme de graphes métriques). \end{thm}
De façon équivalente, pour tout $(j_1,j_2,j_3,j_4)$, on obtient une égalité
entre les birapports suivants~: 
\[|(z_{1,j_1},z_{1,j_2},z_{1,j_3},z_{1,j_4})|=|(z_{2,j_1},z_{2,j_2},z_{2,j_3},z_{2,j_4})|.\]
En fait, nous prouverons ultérieurement ce résultat sans supposer
$z_{i,1},\cdots,z_{i,n}\in\mbf Q_p^{\text{nr}}$, après l'étude du cas des
courbes elliptiques dans le cas $p\neq 2$, et sans hypothèses sur $p$
après l'étude du cas des courbes de Mumford.\\

\begin{proof}

D'après le lemme~\ref{reductioncusps}, on peut supposer $n=4$.
On peut supposer, quitte à faire une homographie, que les points cuspidaux sont $0,1,\infty$ et
$\lambda_{i}$. De plus, on peut supposer que $X_{i}$ n'a pas bonne
réduction (sinon il n'y a rien à prouver). Le cas de bonne réduction
correspond à $v_p(\lambda_i)=v_p(\lambda_i-1)=0$. Quitte à permuter
$0,1$ et $\infty$ par une autre homographie, on peut supposer $v_p(\lambda_{i}-1)>0$
(c'est un entier puisque $\lambda\in\mbf Q_p^{\text{nr}}$). Il nous faut donc
prouver que
$v_p(\lambda_{1}-1)=v_p(\lambda_{2}-1)$.\\
Supposons, par l'absurde, que
$v_p(\lambda_{1}-1)<v_p(\lambda_{2}-1)$.
Posons $h:=v_p(\lambda_{2}-1)-1\geqslant v_p(\lambda_{1}-1)$ et soit
$H_{i}$ le sous-groupe de
$\Pi_{i}$ d'indice $p^h$
correspondant à l'unique revêtement connexe $Y_{i}\to X_{i}$
de degré $p^h$ ramifié seulement en $0$ et $\infty$ (c'est le morphisme $z\mapsto z^{p^h}$
de $\mbf P^1$ vers $\mbf P^1$). On a $\phi(H_1)=H_2$ grâce au lemme~\ref{cusp}.\\
Mais d'après le lemme~\ref{decompositionpuissance}, le point $B(1,r)$ de
l'espace de Berkovich $\mbf P^1$ correspondant à la boule de rayon $r$ et
de centre $1$ a $p^{h-1}$ préimages si
$p^{-(h+\frac{1}{p-1})}\leqslant r <p^{-(h-1+\frac{1}{p-1})}$ et a $p^{h}$
préimages si $r<p^{-(h+\frac{1}{p-1})}$.\\
Ainsi, si $p\neq 2$, $B(1,|\lambda_{1}-1|)$ a $p^{h-1}$ préimages dans
$Y_1$ et $B(1,|\lambda_{2}-1|)$ a $p^h$ préimages dans $Y_2$,
ce qui contredit le lemme~\ref{vertex}.\\

Dans le cas $p=2$ et $h\geqslant 2$,  $B(1,|\lambda_{1}-1|)$ a $2^{h-2}$ préimages dans
$Y_1$ et $B(1,|\lambda_{2}-1|)$ a $2^{h-1}$ préimages dans $Y_2$,
ce qui contredit le lemme~\ref{vertex}.\\
Si $h=1$, et donc $v_2(\lambda_{1})=1$ et $v_2(\lambda_{2})=2$, le
semi-graphe de la réduction de $Y_{1}\to X_{1}$ est~:  
\[\xymatrix{ \infty & & & \sqrt{\lambda_{1}} \\  & & \sommet \ar@{-}[ur] \ar@{-}[r] & -\sqrt{\lambda_{1}}\\ \sommet
  \ar@{-}[uu] \ar@{-}[dd] \ar@{-}[r] & \sommet \ar@{-}[ur] \ar@{-}[dr]&  & \\ & & \sommet
  \ar@{-}[r] \ar@{-}[dr] & 1\\ 0 & & & -1}\]
Le semi-graphe de la réduction stable de $Y_{2}$ est 
\[\xymatrix{ \infty & & & \sqrt{\lambda_{2}} \\ & & & -\sqrt{\lambda_{2}}\\ \sommet \ar@{-}[uu] \ar@{-}[dd] \ar@{-}[rr] & & \sommet \ar@{-}[uur]
  \ar@{-}[ur] \ar@{-}[dr] \ar@{-}[ddr] & \\ & & & 1\\  0 & & & -1}\]
Les deux graphes ne sont pas isomorphes, et l'on obtient donc une contradiction.\end{proof}

\section{Cas d'une courbe elliptique épointée}
Pour $i=1$ ou $2$,
soit $X_{i}$ une courbe de Tate épointée sur $\mbf C_p$ qui est définie sur
une extension finie de $\mbf Q_p$. On a un isomorphisme $X_i=\mbf C_p^*/q_{i}^{\mbf
  Z}-\{1\}$ avec $q_i\in\overline{\mbf Q}_p^*$ et $|q_{i}|<1$.\\
Posons $\Pi_{i}=\gtemp(X_{i})$.
\begin{thm} Si il existe un isomorphisme $\phi:\Pi_1\simeq\Pi_2$, alors
$|q_{1}|=|q_{2}|$\label{ellcurve}
\end{thm}
L'idée principale de la preuve est de considérer un $\mbf Z/p\mbf
Z$-torseur sur $\Gm=\mbf
P^1\backslash\{0,\infty\}$ ($\Gm$ est le rev\^etement topologique universel de la courbe elliptique) qui est ramifié en au
plus deux points (et qui n'est pas ramifié ni en $0$ ni en $\infty$). Comme
le torseur est décomposé au voisinage de $0$ et de $\infty$, on peut
recoller des "morceaux" de ce torseur pour le rendre périodique (ainsi on
pourra le descendre en un $\mbf Z/p\mbf Z$-torseur d'un revêtement
topologique fini). Nous utiliserons alors le
lemme~\ref{decompositionpuissance} comme dans le cas de la droite épointée.\\
Pour vérifier que les torseurs que nous aurons construits se
correspondent bien l'un l'autre par $\phi$ (à un torseur topologique et
multiplication par une constante près), nous aurons à étudier le
$\mbf F_p$-espace vectoriel des torseurs de la courbe elliptique ramifiés
en au plus deux points (fixés)~: c'est un espace vectoriel de dimension
3. Une base sera obtenue en considérant un torseur topologique, le torseur
que nous venons de décrire brièvement et un autre construit de façon
similaire. En termes de courants (dont
nous n'aurons pas besoin pour la preuve), on pourra dire que, si le torseur que nous aurons
construit correspond au courant qui suit un chemin reliant les deux
points cuspidaux, l'autre torseur correspond au courant suivant l'autre
chemin reliant les deux points cuspidaux.\\
Nous n'aurons pas ici besoin de supposer $q_i\in\mbf Q_p^{\text{nr}}$ comme
dans le cas de la droite épointée parce qu'en prenant un revêtement non
ramifié de la courbe elliptique, on obtient autant de sommets sur le
squelette de la courbe que l'on veut.\\
\begin{proof}[D\'emonstration de \ref{ellcurve}]
Choisissons des entiers $n$, $l$ et $m$ tels que:
\begin{itemize}
\item $n$ soit premier à $p$ et $n\geqslant \frac{p-1}{p}\times\frac{v_p(q_2)v_p(q_1)}{|v_p(q_2)-v_p(q_1)|}$,
\item $l\geqslant 1+\frac{2np}{(p-1)\cdot v_p(q_{i})}$,
\item $m\geqslant \frac{2l}{n}$.
\end{itemize}
Soit  $H_{i,0}=[\Pi_{i},\Pi_{i}]\Pi_{i}^n$ la préimage dans
$\Pi_{i}$ de l'image par la multiplication par $n$ dans l'abélianisé de
$\Pi_{i}$. $\phi$ induit un isomorphisme
$H_{1,0}\to H_{2,0}$.\\
Le revêtement $Y_{i,0}$ de $X_{i}$ correspondant à $H_{i}$ est la
multiplication $\overline X_{i}\stackrel{\times n}{\to}\overline X_{i}$
par $n$ sur la courbe elliptique $\overline{X}_{i}$.\\

Soit maintenant $H_{i,1}$ le sous-groupe de $H_{i,0}$ correspondant à
l'unique revêtement topologique connexe $Y_{i,1}$ de degré $n$ de
$Y_{i,0}$. Alors $Y_{i,1}=\mbf C_p^*/q_i^{m\mbf
  Z}-\{q_i^{\frac{a}{n}}\zeta^b\}_{(a,b)\in\mbf Z^2}$ où $q_i^{\frac{1}{n}}$
est une racine $n$\ieme{} de $q_i$ et $\zeta$ est une racine $n$\ieme{} de $1$.\\ 
Le semi-graphe de la réduction stable de $Y_{i,1}$ a $mn$ sommets
qui sont dispos\'es sur le squelette de $\overline Y_{i,1}=\mbf C_p^*/q_i^{m\mbf Z}$ qui est hom\'eomorphe \`a un cercle (la distance entre deux sommets consécutifs est
$v_p(q_{i})/n$) et  $n$ points cuspidaux aboutissent à chaque sommet. L'isomorphisme
$\phi$ induit un isomorphisme $H_{1,1}\simeq H_{2,1}$ qui lui-même
induit un isomorphisme entre les semi-graphes de la réduction stable de
$Y_{i,1}$. Numérotons de $0$ à $mn-1$ les sommets de ce graphe en suivant
le cercle, de façon compatible avec l'isomorphisme induit par $\phi$.
Notons $x_{i,0},\cdots, x_{i,mn-1}$ les sommets correspondant du squelette
de $Y_{i,1}$.
Soient $z_{1,0}$ et $z_{1,l}$ deux points cuspidaux de
$Y_{1,1}$ aoutissant aux sommets du graphe numérotés $0$ et $l$
respectivement. Soient $z_{2,0}$ et $z_{2,1}$ les points cuspidaux
correspondant de $Y_{2,1}$.\\

Concentrons-nous maintenant sur les $\mbf Z/p\mbf Z$-torseurs sur $\overline
Y_{i,1}$ qui sont ramifiés seulement en $z_{i,0}$ et $z_{i,l}$. Ce sont les \'el\'ements de $V_i:=\Hom(\ga(\overline Y_{i,1}\backslash\{z_{i,0},z_{i,l}\}),\mbf Z/p\mbf Z)$. $V_i$ est un un groupe ab\'elien sur lequel la multiplication par $p$ est triviale~: c'est un $\mbf F_p$-espace vectoriel. Comme $\ga(\overline Y_{i,1}\backslash\{z_{i,0},z_{i,l}\})$ est un groupe libre \`a 3 g\'en\'erateurs, $V_i$ est de dimension 3 sur $\mbf F_p$.\\
Rappelons que l'on a choisi $m,n$ et $l$ de manière à ce que $\frac{l-1}{n}v_p(q_{i})>\frac{2p}{p-1}$ et
$\frac{mn-l}{n}v_p(q_{i})>\frac{l}{n}v_p(q_{i})$.\\
Décrivons maintenant une base de ce $\mbf F_p$-espace vectoriel.\\
\begin{center}
\includegraphics[scale=0.5]{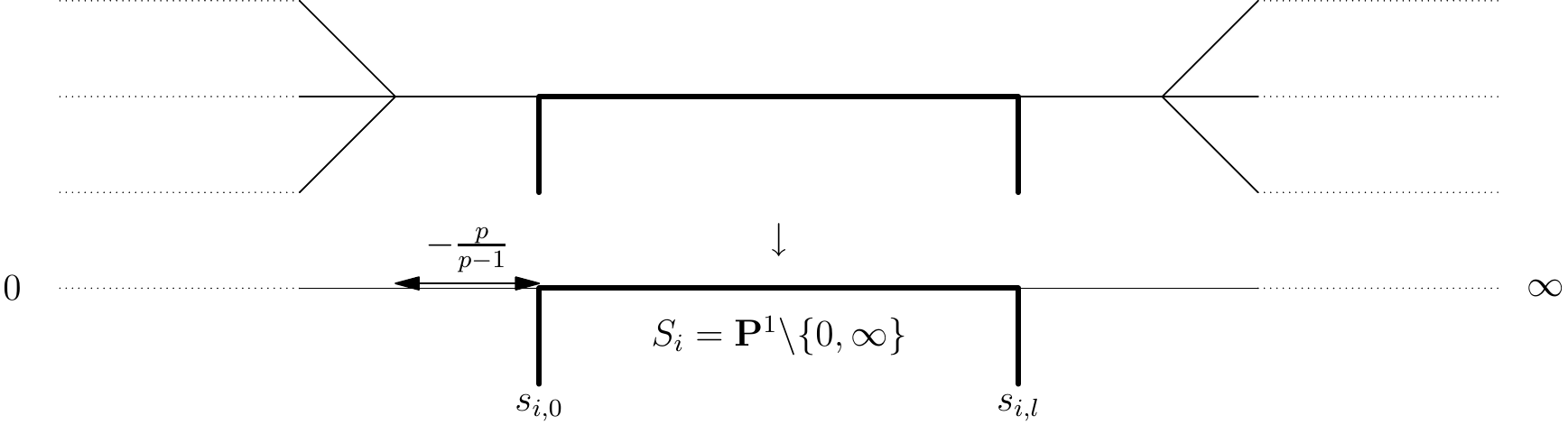}
\end{center}
Soit $S_{i}$ le revêtement topologique universel de $\overline
Y_{i,1}$, que nous identifions à $\mbf P^1\backslash\{ 0,\infty
\}\subset\mbf P^1$. Soient
$s_{i,0}$ et $s_{i,l}$ les seuls antécédents dans $S_{i}$ de $z_{i,0}$ et
$z_{i,l}$ de norme $1$ et $|q|^{\frac{l}{mn}}$ et soient $U_{i,1}\subset
S_{i}$ la couronne ouverte $\{|q_{i}|^{\frac{l-mn}{2n}}>|z|>|q_{i}|^{\frac{l+mn}{2n}}\}$
et $U_{2,i}\subset S_{i}$ la couronne ouverte
$\{|q_{i}|^{\frac{l}{n}}p^{-\frac{p}{p-1}}>|z|>|q_{i}|^mp^{\frac{p}{p-1}}\}$. Les
applications de $U_{i,1}$ et de $U_{i,2}$ vers $\overline Y_{i,1}$  sont
encore des immersions ouvertes, et ensemble recouvrent $\overline Y_{i,1}$.\\
Soit $T_{i,1}$ la restriction à $U_{i,1}$ du revêtement ramifié
(seulement au-dessus de $s_{i,0}$ et de $s_{i,l}$) $\mbf
P^1\to\mbf P^1:z\mapsto \frac{s_{i,0}z^p+s_{i,l}}{z^p+1}$, qui est
galoisien de groupe isomorphe à $\mbf Z/p\mbf Z$, et choisissons un tel
isomorphisme pour obtenir un $\mbf Z/p\mbf Z$-torseur. Soit $T_{i,2}$ le
$\mbf Z/p\mbf Z$-torseur trivial sur $U_{i,2}$  et soit $T_{i,3}=T_{i,1}\coprod
T_{i,2}\to U_{i,3}=U_{i,1}\coprod U_{i,2}$.\\
Sur $U_{i,1}\times_{\overline Y_{i,1}} U_{i,2}$ (qui a deux composantes
connexes), $T_{i,1}$ est trivial. Choisissons 
une trivialisation~: on peut maintenant descendre $T_{i,3}\to U_{i,3}$ en un
 $\mbf Z/p\mbf Z$-torseur $T_{i}\to \overline Y_{i,1}$, qui est ramifié
seulement au-dessus de $z_{i,0}$ et de $z_{i,l}$.\\
D'après le lemme~\ref{decompositionpuissance}, $T_{i}\to \overline Y_{i,1}$
est décomposé au-dessus de $x_j$ (avec $j\in [0,mn-1]$)
si et seulement si
$|q_{i}|^{\frac{j}{n}}\in
[|q_{i}|^mp^{\frac{p}{p-1}},|q_{i}|^{\frac{l}{n}}p^{-\frac{p}{p-1}}]$,
\ie \[j\in I_{i,1}:=\big[l+\frac{np}{v_p(q_{i})(p-1)},mn-\frac{np}{v_p(q_{i})(p-1)}\big].\]
Il existe un tel entier grâce à l'hypothèse faite sur $l,m$ et $n$ car
\[\lg(I_{i,1})=mn-l-2\frac{np}{v_p(q_{i})(p-1)}\geqslant 1\] (où $\lg$
est la longueur d'un intervalle).\\

De même, soit $s'_{i,l}$ une préimage de $z_{i,l}$ de norme
$|q|^{\frac{l-mn}{n}}$. Soient $U'_{1,i}$ la couronne
$\{|q_{i}|^{\frac{l-2mn}{2n}}>|z|>|q_{i}|^{\frac{l}{2n}}\}$ et
$U'_{2,i}$ la couronne $\{p^{-\frac{p}{p-1}}>|z|>|q_{i}|^{\frac{mn-l}{n}}p^{\frac{p}{p-1}}\}$.
Ce sont des ouverts de $\overline Y_{i,1}$, qu'ils recouvrent.\\
Soit $T'_{i,1}$ un $\mbf Z/p\mbf Z$-torseur sur $U'_1$ obtenu comme la
restriction d'un $\mbf Z/p\mbf Z$-torseur sur $\mbf P^1$ ramifié seulement
au-dessus de
$z_{i,0}$ et de $z'_{i,l}$. Il est trivial au-dessus de $U'_{i,1}\cap
U'_{i,2}$, et, en choisissant une telle trivialisation, on obtient par
descente un
$\mbf Z/p\mbf Z$-torseur $T'_{i}$ sur $\overline Y_{i,1}$,
ramifié seulement au-dessus de $z_{i,0}$ et de $z_{i,l}$.\\
$T_{i}\to \overline Y_{i,1}$ est décomposé au-dessus de $x_j$ (avec $j\in [0,mn-1]$)
si et seulement si
$|q_{i}|^{\frac{j}{n}}\in
[|q_{i}|^{\frac{j}{n}}p^{\frac{p}{p-1}},p^{-\frac{p}{p-1}}]$, \ie
\[j\in
I_{i,2}:=\big[\frac{np}{v_p(q_{i})(p-1)},l-\frac{np}{v_p(q_{i})(p-1)}\big].\]
Il existe un tel entier grâce à l'hypothèse faite sur $l,m$ et $n$
car \[\lg(I_{i,2})=l-2\frac{np}{v_p(q_{i})(p-1)}\geqslant 1.\]
Soit finalement $T''_{i}$ le revêtement topologique connexe
(unique \`a isomorphisme pr\`es) de degré $p$ de $\overline Y_{i}$, et fixons un
isomorphisme de $\mbf
Z/p\mbf Z$ sur son groupe de Galois, pour le consid\'erer comme un $\mbf Z/p\mbf
Z$-torseur. Montrons qu'alors $T_{i},T'_{i},$ et $T''_{i}$
constituent une base de $V_{i}$.\\
Soit $i$ un entier dans $I_{i,1}$. Comme $I_{i,1}\cap
I_{i,2}=\emptyset$ et comme $T''$
est partout décomposé, si $aT_{i}+bT'_{i}+cT''_{i}$ (la combinaison
linéaire est à prendre au sens de l'espace vectoriel $V_{i}$) est
décomposé en $x_j$, $c=0$. Par un même argument avec $I_{i,2}$, si
$aT_{i}+bT'_{i}+cT''_{i}=0$, on obtient $b=0$, mais comme $T_{i}$
n'est pas trivial non plus, on obtient en fait que $T_{i},T'_{i}$ et
$T''_{i}$ constituent une base de $V_{i}$.\\

Supposons maintenant, par l'absurde, que $|q_1|>|q_2|$, alors $I_{1,1}\subset
I_{2,1}$.\\
Soit $T_{2,0}=aT_{2}+bT'_{2}+cT''_2$ l'image de $T_{1}$ par $(\phi^{-1})^*:V_{1}\to
V_{2}$ (en effet, $\phi$ induit un tel $(\phi^{-1})^*$ grâce au lemme~\ref{cusp}).\\
Soit $j\in I_{1,1}$, $T_{1}$ est décomposé au-dessus de $x_{1,j}$,
donc d'après le lemme~\ref{vertex}, $T_{2,0}$ est décomposé au-dessus de $x_{2,j}$,
et donc $c=0$. Alors, $T_{2,0}$ est décomposé au-dessus de chaque $x_j$ si
$b=0$ ou précisément au-dessus des $x_j$ tels que $j\in I_{2,1}$ sinon.\\
Or, d'après le lemme~\ref{vertex}, $T_{2,0}$ doit être décomposé
précisément au-dessus $x_{2,j}$ tel que $T_{1}$ est décomposé au-dessus de $x_{1,j}$,
\ie $j\in I_{1,1}$.\\
Ainsi le cas $b=0$ est impossible et $I_{2,1}$ et $I_{1,1}$
doivent contenir exactement les mêmes entiers.\\
Mais ce n'est pas possible si  $n\geqslant
\frac{v_p(q_2)v_p(q_{1})(p-1)}{(v_p(q_2)-v_p(q_1))p}$, car
\[\lg(I_{2,1})-\lg(I_{1,1})=2\frac{np}{v_p(q_{1})(p-1)}-2\frac{np}{v_p(q_{2})(p-1)}\geqslant
2.\]   
\end{proof}

\begin{rem}
Supposons $p\neq 2$. Soit $\{z_{i,1},\cdots, z_{i,4}\}$ quatre éléments de
$\overline{\mbf Q}_p$, et soit $\phi$ un isomorphisme entre les
  $\gtemp(\mbf P^1\backslash\{z_{i,1},\cdots, z_{i,4}\})$ pour $i=1,2$, qui
  identifie $z_{1,j}$ à $z_{2,j}$ par l'identification des
  points cuspidaux des graphes des réductions stables (nous supposerons les
  courbes ayant mauvaise réduction, et la longueur des arêtes du graphe \'egale \`a
  $l_{i}$). Soit $E_{i}$ l'unique
  $\mbf Z/2\mbf Z$-torseur de $\mbf P^1$ ramifié au-dessus de
  $\{z_{i,1},\cdots, z_{i,4}\}$ et seulement au-dessus de ces points (les sous-groupes
  d'indice 2 correspondant aux $E_{i}$ s'envoient donc l'un sur l'autre). $E_{i}$
  est une courbe de Tate avec $|q_{i}|=2l_{i}$. D'après le théorème~\ref{ellcurve},
  $|q_{1}|=|q_{2}|$ donc $l_{1}=l_{2}$, ce qui prouve encore
  le théorème~\ref{thmp1} sans supposer que les points soient dans
  $\mbf Q_p^{\text{nr}}$.\\
Si $p=2$ et $l_{i}>4$, alors $E_{i}$ est aussi une courbe de Tate,
$|q_{i}|=2l_{i}-8$ et l'argument précédent fonctionne également.
\end{rem}

\section{Cas d'une courbe de Mumford}
\subsection{Rappels sur les courbes de Mumford et les courants}
Soit $X$ une courbe de Mumford de genre $g\geq 2$ sur $\overline K$, soit  $\Omega\subset \mbf P^1$
son revêtement topologique universel, et $\Gamma=\Gal(\Omega/X)$. Ainsi $X=\Omega/\Gamma$.\\
Soit $\Phi$ la rétraction de $\Omega$, comme espace de Berkovich, sur son squelette $\mbb T=\mbb T(\Omega)$.\\
Pour $z,z'\in \Omega$, posons \[d_{\Omega}(z,z')=\sup_{x_1,x_2\in\mbf
  P^1\backslash\Omega}|v_p(\frac{z'-x_1}{z-x_1}\frac{z-x_2}{z'-x_2})|,\]
distance qui est invariante par homographie. De plus $d(z,z')$ dépend uniquement de
$\Phi(z)$ et $\Phi(z')$, et ce n'est rien d'autre que la distance entre $\Phi(z)$
et $\Phi(z')$ pour la métrique usuelle $d$ sur l'arbre de $\Omega$. Ainsi $d_\Omega(z,z')=d(\Phi(z),\Phi(z'))$.\\
Si $z\in\Omega$ et $\lambda>0$, notons
$U_{z,\lambda}:=\{z\in\Omega|d_{\Omega}(z,z')\leqslant \lambda\}$. C'est un
sous-espace affinoïde de $\Omega$ (c'est une intersection de disques fermés).\\

Soit $\mcal L=\mbf P^1\backslash \Omega$, c'est un compact de $\mbf
P^1$. Il est constitué de points de $\mbf P^1(\mbf C_p)$ qui est totalement
discontinu, donc $\mcal L$ est aussi totalement discontinu, donc profini.
Dans~\cite[1.8.9]{FvdP2}, Fresnel et van der Put definissent une
\emph{mesure} sur un espace topologique profini $Z$ comme étant une
fonction $\mu:\{\text{ouverts compacts de }Z\}\to\mbf Z$ telle que
$\mu(U_1\cup U_2)+\mu(U_1\cap U_2)=\mu(U_1)+\mu(U_2)$ pour tous ouverts
compacts $U_1,U_2$ de $Z$. Le groupe des mesures sur $Z$ telles que
$\mu(Z)=0$ est noté $M_0(Z)$.\\
On peut alors associer à $f\in
O(\Omega)^*$ une mesure $\mu_f$ sur $\mcal L$. En suivant~\cite{FvdP2}, on
appellera aussi trou d'un sous-espace affinoïde $U$ de $\mbf P^1$ toute
composante connexe du complémentaire de cet affinoïde, et nous noterons
$t(U)$ l'ensemble des trous de $U$.\\
On a une suite exacte~(\cite[Prop. 1.8.9]{FvdP2}):
\[1\to \overline K^*\to O(\Omega)^*\to M_0(\mcal L)\to 0\]
Plus précisément, d'après~\cite[1.8.10 ex. $\beta$]{FvdP2}, si $a,b\in\mcal L$ et $f:z\mapsto\frac{z-a}{z-b}$, alors on a $\mu_f=\delta_a-\delta_b$.
En général, si $f\in O(\Omega)^*$ et $\mu=\mu_{f}$, $\mu$ est une limite
faible de $\mbf Z$-combinaisons linéaires $(\mu_k)_{k\in\mbf N}$ de mesures
de Dirac. Si $\mu_k=\sum n_i(k)\delta_{a_i(k)}$, soit $f_k=\prod
(1-\frac{a_i(k)}{z})^{n_i(k)}$. Alors $\mu_{f_k}=\mu_k$, et $(f_k)$ tend
uniformément sur tout affinoïde vers $g$. Alors $\mu_g=\mu$ et donc $g=\lambda f$ avec $\lambda \in \overline K^*$.\\

Un courant $C$ sur un graphe $\mbb G$ à coefficient dans un groupe abélien $A$
est une fonction de l'ensemble des branches de $\mbb G$, telle que~:
\begin{itemize}
\item pour toute arête $e$, $C(b_1)=-C(b_2)$, où $b_1$ et $b_2$ sont les
  branches de $e$~;
\item pour tout sommet $v$, $\sum_{b}C(b)=0$ où la somme porte sur toutes
  les branches $b$ aboutissant en $v$.
\end{itemize}
Notons $C\mbb G)$ est le groupe des courants à coefficients entiers sur
$\mbb G$. 
D'après~\cite[prop. 1.1]{vdp1}, on a la suite exacte suivante~:
\[1\to \overline K^*\to O(\Omega)^*\to C(\mbb T)\to 0\]
qui donne, par comparaison avec la suite précédente, un isomorphisme $M_0(\mcal L)\to C(\mbb T)$.\\
On peut décrire cet isomorphisme de la façon suivante~: si $e$ est une
arête orientée de $\mbb T$, $\mbf P^1\backslash\Phi^{-1}(e)$ a deux
composantes connexes, et $\Phi^{-1}(e)\cap\mcal L=\emptyset$, d'où une
partition de $\mcal L$ par deux ouverts, $\mcal L_1(e)$ au départ de $e$ et
$\mcal L_2(e)$ à l'arrivée de $e$ (on a pour $\mu\in M_0(\mcal L)$, $\mu(\mcal L_1(e))=-\mu(\mcal L_2(e))$). Alors $\mcal C(e)=\mu(\mcal L_2(e))$.\\
Plus généralement, si $\mbb K$ est un sous-graphe fini connexe de $\mbb T$
(contenant au moins une arête), $\Phi^{-1}(\mbb K)$ est un affinoïde
contenu dans $\Omega$. Il y a une bijection naturelle entre
$t(\Phi^{-1}(\mbb K))$ et l'ensemble des arêtes de $\mbb T$ ayant
exactement une extrémité dans $\mbb K$. Si
$\mu\in M_0(\mcal L)$, $C\in t(\Phi^{-1}(\mbb K))$ et $e$ est l'unique
arête orientée de $\mbb T$ partant de $\mbb K$ correspondant
à $C$, alors $C\cap \mcal L=\mcal L_2(e)$ et donc $\mcal C(e)=\mu(C\cap \mcal
L)$. En particulier $\mu(C\cap \mcal L)=0$ pour tout $C\in t(\Phi^{-1}(\mbb
K))$ si et seulement si $\mcal C(e)=0$ pour toute arête terminale partant
de $\mbb K$, si et seulement si $\mcal C$ est nul l'étoile $\widetilde{\mbb K}$
de $\mbb K$ (l'\emph{étoile} de $\mbb K$ est par définition l'ensemble des
arêtes ayant une extrémité dans $\mbb K$).\\

Soit $\Theta$ le groupe des fonctions th\^eta de $X$, \ie le sous-groupe des
fonctions $f\in O(\Omega)^*$ telles que pour tout $\gamma\in \Gamma$,
$z\mapsto f(\gamma z)/f(z)$ est une fonction constante (cela signifie que
le courant correspondant, ou de manière équivalente la mesure correspondante, est $\Gamma$-équivariant). Alors on a la suite
exacte suivante~:
\[1\to \overline K^*\to \Theta \to C(\mbb G)\to 0\]
où $\mbb G=\mbb T/\Gamma$ est le graphe de la réduction stable de $X$.\\
et donc $\Theta/\overline K^*$ est un $\mbf Z$-module libre de rang $g$.\\
Le r\'esultat suivant est inspir\'e de~\cite[prop. 2.1]{vdp2}~:
\begin{prop}\label{torseurstheta} Pour tout $n\geq 2$ et tout $\mbf Z/n\mbf Z$-torseur $Y$ sur $X$,
il existe un élément $\theta$ dans $\Theta$, unique modulo $\overline
K^*\Theta^n$, tel que $Y\times_X \Omega=\Omega[f]/(f^n=\theta)$ où
$\Omega[f]/(f^n=\theta)$ est le pullback du $\mbf Z/n\mbf Z$-torseur $\mbf
G_m \stackrel{z\mapsto z^n}{\to}\mbf G_m$ le long de $\Omega\stackrel{\theta}{\to}\mbf G_m$.\\
Réciproquement, pour tout $\theta$ dans $\Theta$, il existe un $\mbf
Z/n\mbf Z$-torseur $Y$ sur $X$ tel que $Y\times\Omega=\Omega[f]/(f^n=\theta)$.\end{prop}
\begin{proof} Nous suivons les notations de~\cite[section 2]{vdp2}. Soit $\Omega_*$
  une composante connexe de $Y\times_X\Omega$.\\
Supposons d'abord $Y\times_X\Omega$ connexe. Ainsi
$\Omega_*=Y\times_X\Omega$. Alors, d'après~\cite[prop. 2.1]{vdp2}, il
existe un unique réseau $T$ dans $(\Theta/\overline K^*)\otimes \mbf Q$
contenant $\Theta/\overline K^*$ tel que, comme revêtement de $X$,
$\Omega_*=\Omega(T):=\Omega\times_{\Spec \overline K[\Theta/\overline
  K^*]}\Spec \overline K[T]$. Alors $T/(\Theta/\overline K^*)$ est
isomorphe à $\mbf Z/n\mbf Z$, et choisir un générateur $\bar f$ de
$T/(\Theta/\overline K^*)$) revient à choisir un structure de $\mbf Z/n\mbf
Z$-torseur sur $\Omega_*$.\\
Si l'on prend $\bar f^i$ (avec $i\in (\mbf Z/n\mbf Z)^*$) comme autre
générateur, le torseur correspondant est $i\cdot\Omega_*$ (pour la
structure de $\mbf Z/n\mbf Z$-module sur l'ensemble des $\mbf Z/n\mbf
Z$-torseurs). Ainsi en changeant de générateur, on obtient les
$\phi(n)$ structures distinctes de $\mbf Z/n\mbf Z$-torseur sur le revêtement $\Omega_*$ de $\Omega$.\\
Si $\bar f$ est le générateur correspondant à la structure de $\mbf Z/n\mbf
Z$-torseur sur $\Omega_*=Y\times_X\Omega$ h\'erit\'ee de celle de $Y$, alors $\theta$ est n'importe quel
relèvement de $\bar f^n\in T^n/(\Theta/\overline K^*)^n$.\\
Dans le cas général, $\Omega_*$ acquiert la structure de $\mbf Z/m\mbf
Z$-torseur sur $\Omega$ pour $m|n$ convenable, et comme précédemment, on obtient
un unique $\theta_0$ modulo $\overline K^*\Theta^m$ tel que
$\Omega_*=\Omega[f]/(f^m=\theta_0)$. Alors $Y\times_X\Omega=\Ind_{\mbf
  Z/m\mbf Z}^{\mbf Z/n\mbf Z}\Omega_*$, et ainsi
$Y\times_X\Omega=\Omega[f]/(f^n=\theta_0^{n/m})$ et $\theta=\theta_0^{n/m}$.\\
La seconde assertion découle du fait que si $\Omega_*$ est une composante
connexe de $\Omega[f]/(f^n=\theta)$, alors $\Gal(\Omega_*/X)$ est (non
canoniquement) isomorphe au produit direct de $\Gal(\Omega_*/\Omega)$ et
$\Gamma$ (d'après~\cite[section 2, intro.]{vdp2}). Ainsi $\Omega_*$ se
descend (non canoniquement) à $X$ en considérant $Y_0=\Omega_*/N$ où $N$
est un complément de $\Gal(\Omega_*/\Omega)$ dans $\Gal(\Omega_*/X)$ (et
alors $Y\times \Omega=\Omega[f]/(f^n=\theta)$ où $Y$ est une somme
disjointe de copies de $Y_0$).\end{proof}

\begin{rem}
On peut aussi montrer cela en considérant l'homomorphisme $\widetilde J=\mcal
Hom(\Theta/\overline K^*,\mbf G_m)\to J$ où $J$ est la jacobienne de $X$ et
$\widetilde J$ est son revêtement topologique universel.\\
Plus pr\'ecis\'ement, le diagramme commutatif suivant~:
\[\xymatrix{\Omega \ar[r] \ar[d] & \widetilde J \ar[d] \\ X \ar[r] & J}\]
induit un diagramme commutatif~:
\[\xymatrix{\Tors(\Omega,\mbf Z/n\mbf Z) & \Tors(\widetilde J,\mbf Z/n\mbf Z)=\Theta/\overline K^*\Theta^n \ar[l]\\
\Tors(X,\mbf Z/n\mbf Z) \ar[u] & \Tors(J,\mbf Z/n\mbf Z) \ar[u]\ar[l]}.\]
La premi\`ere assertion de la proposition d\'ecoule alors du fait que la fl\`eche du bas soit un isomorphisme et la seconde assertion du fait que la fl\`eche de droite soit surjective (ce qui se d\'eduit du fait que $\ga(\widetilde J)$ soit un facteur direct de $\ga(J)$).
\end{rem}

\subsection{Résultats préliminaires sur la décompositions des torseurs
  correspondant aux courants}
Rappelons que $U_{z,\lambda}$ d\'esigne le sous-espace affinoïde $\{z'|d(z,z')\leq\lambda\}$ de $\Omega$.
\begin{prop} Soit $z\in\Omega,\lambda>0$. Soit $f\in O(\Omega)^*$ tel que $f(z)=1$.
  Soit $\mu$ la mesure sur $\mcal L$ correspondant à $f$
  et supposons $\mu(C\cap \mcal L)=0$ pour tout trou $C$ de $U_{z,\lambda}$.\\
Alors $\forall z'\in U_{z,\lambda}$, $|f(z')-1|\leqslant
p^{d(z,z')-\lambda}$\end{prop}
\begin{proof} Pour simplifier, supposons $z=\infty$.\\
D'après~\cite[1.8.10 ex. $\beta$]{FvdP2} et~\cite[prop. 1.8.9.(i)]{FvdP2},
on a $f=\lim f_k$ uniformément sur tout affinoïde de $\Omega$ (en particulier
sur $U_{z,\lambda}$) où $f_k$ est de la forme suivante
\[f_k(z')=\prod^{s_k}_{i=1} (1-\frac{x_{i,k}}{z'})^{n_{i,k}},\]
et $\mu_k=\sum_i n_{i,k}\delta_{a_i(k)}\to \mu$.\\
Pour $k$ assez grand, on a $\mu_k(C\cap \mcal L)=0$ pour tout trou $C$ de
$U_{z,\lambda}$ et $|f-f_k|_{U_{z,\lambda}}\leqslant p^{-\lambda}$.\\
Il suffit donc de prouver le résultat pour $f_k$, qui est un produit de
fonctions de la forme $g:z'\mapsto
z'=(1-\frac{x_1}{z'})(1-\frac{x_2}{z'})^{-1}=\frac{z'-x_1}{z'-x_2}$ avec $x_1$
et $x_2$ dans le même trou $C$ de $U_{z,\lambda}$.
Nous n'avons qu'à montrer le résultat pour $g$, ce qui se voit
aisément.\end{proof}
Si $x$ et $x'$ sont deux points de la réalisation géométrique $|\mbb T_0|$ d'un arbre $\mbb T_0$,l'ensemble des
parties connexes de $|\mbb T_0|$ qui contienne $x$ et $x'$ admet un minimum
pour l'inclusion. Il est noté $[x,x']$.
\begin{cor} Soit $f\in O(\Omega)^*$ telle que $f(z)=1$. Soit $U$ un affinoïde
  de $\Omega$ tel que pour tout trou $C$ de $U$, $\mu(C\cap\mcal L)=0$.\\
Supposons que pour tout $z''\in\Phi^{-1}([\Phi(z),\Phi(z')]),
U_{z'',\lambda}\subset U$.\\
Alors $|f(z')-1|\leqslant p^{-\lambda}$.\end{cor}
\begin{proof} Soit $\epsilon>0$. Soit $(z_i)_{i=0\cdots n}$ une suite telle que $z_0=z$, $z_n=z'$, $\Phi(z_i)\in
[\Phi(z),\Phi(z')]$ et $d(\Phi(z_i),\Phi(z_{i+1}))\leqslant \epsilon$. Par
hypothèse, on a $U_{z_i,\lambda}\subset U$. Donc, pour tout trou de
$U_{z_i,\lambda}$, $C\cap\mcal L=\coprod_i C_i\cap\mcal L$ où l'union est
sur tous les trous $C_i$ de $U$ qui sont contenu dans $C$. Ainsi,
$\mu(C\cap\mcal L)=0$ pour tout trou de $U_{z_i,\lambda}$. Comme
$z_{i+1}\in U_{z_i,\lambda}$, la proposition précédente nous donne
$|\frac{f(z_{i+1})}{f(z_i)}-1|\leqslant p^{\epsilon-\lambda}$.
Ainsi $|f(z'')-1|\leqslant \sup|f(z_{i+1})-f(z_i)|\leqslant
p^{\epsilon-\lambda}$. On obtient le résultat en faisant tendre $\epsilon$
vers 0.\end{proof}
\begin{cor}\label{totdecrev} Soit $f$ comme précédemment, $U$ tel que
  $\mu(C\cap\mcal L)=0$ pour tout trou $C$ de $U$. Soit $e$ un entier positif. Soit $\lambda>e+\frac{1}{p-1}$, soit $Y\to
  \Omega$ un revêtement fini obtenu par pullback de $\mbf
  G_m\stackrel{z\mapsto z^{p^e}}{\to}\mbf G_m$ le long de $f:\Omega\to \mbf G_m$. Soit
  $V\subset U$ tel que $\forall z\in V, U_{z,\lambda}\subset U$.\\
Alors $Y$ est scindé sur $V$.\end{cor}
\begin{proof} On peut supposer $V$ connexe, car il suffit de prouver le
  résultat pour toute composante connexe de $V$. Soit $z\in V$. Quitte à
  multiplier $f$ par une constante, ce qui ne modifie pas $Y$, on peut supposer $f(z)=1$.\\
D'après le corollaire précédent, $f(V)\subset
D(1,p^{-\lambda})$. Mais selon~\ref{decompositionpuissance}, $\mbf G_m\stackrel{z\mapsto z^{p^e}}{\to}\mbf G_m$
est décomposé sur $D(1,p^{-\lambda})$, ce qui prouve le
résultat.\end{proof}
En particulier, supposons que $\mbb K$ et $\mbb K'$ sont des sous-graphes de $\mbb T$
tels que $\{z\in\mbb T|d(z,K)\leqslant \lambda\}\subset\mbb K'$ pour un
certain $\lambda>e+\frac{1}{p-1}$, et que
$\mcal C$ est un courant sur $\mbb T$ nul sur l'étoile de $\mbb
K'$. Alors $V=\phi^{-1}(|\mbb K|)$ et $U=\phi^{-1}(|\mbb K'|)$ sont affinoïdes~;
$\mu_{\mcal C}(C\cap L)=0$ pour tout trou de  $U$, et donc le corollaire précédent
implique que le $\mbf Z/p^e\mbf Z$-torseur associé à $\mcal C$ est scindé
sur $V$ (et donc a fortiori totalement décomposé sur $\mbb K$).

\begin{prop}\label{nondecrev} Soit $\mcal C$ un courant sur le squelette $\mbb T:=\mbb T(\Omega)$, correspondant
  à une fonction inversible $f$ sur $\Omega$. Soit $a$ un sommet de $\mbb
  T$ telle que la restriction $\mcal C_a$ de $\mcal C$ à l'étoile de $a$
  soit non nulle modulo $n$. Alors, si $Y\to
  \Omega$ est le pullback de $\mbf G_m\stackrel{z\mapsto z^n}{\to}\mbf G_m$
  le long de $f:\Omega\to  \mbf G_m$ , $Y\to
  \Omega$ n'est pas totalement décomposé en $a$.\end{prop}
\begin{proof} Le rev\^etement $Y\to \Omega$ est décomposé en $a$ si et seulement si il
  existe $f_1\in \mcal O_{\Omega,a}$ telle que $f_1^n=f_{|\mcal  O_{\Omega,a}}$.
Quitte à multiplier $f$ par une constante, on peut supposer $|f|_a=1$.\\
Si $f_1$ existe, en regardant le corps résiduel $\tilde
k_{\Omega,a}\simeq\tilde k(X)$, $\tilde f_1^n=\tilde f$. Si $\tilde
f_1(z)=\lambda\prod(z-a_i)$, $\tilde f=\lambda^n\prod(z-a_i)^n$, et donc
tous les zéros et pôles de $\tilde f$ sont d'ordre multiple de $n$, ce qui
conclut la preuve.\end{proof}

\subsection{Graphe métrique de la réduction stable et groupe fondamental
tempéré}

Supposons maintenant que $X_1$ et $X_2$ soient deux
$\overline K$-courbes de Mumford définies sur $K$ (pour que
l'on puisse utiliser le théorème~\ref{mochi}, qui suppose que les courbes
sont déjà définies sur un corps à valuation discrète), et que l'on ait un isomorphisme \[\phi:\gtemp(X_1)=\gtemp(X_2),\]
qui induit donc un isomorphisme de graphes \[\mbb G_1\to \mbb
G_2,\] d'où un isomorphisme $\mbb T(\Omega_{1})\simeq
\mbb T(\Omega_{2})$.\\
\begin{thm}\label{mumford}L'isomorphisme $\mbb G_{1}\to\mbb G_{2}$ de
  graphes est en fait un isomorphisme de graphes métriques.\end{thm}

\begin{rem} Supposons que $X_1$ et $X_2$ soient des droites projectives
  privées de quatre points $a_{i},b_{i},c_{i},d_{i}$ avec
  $h_i:=v_p(a_{i},b_{i},c_{i},d_{i})>0$ et que l'on ait un isomorphisme
  $\phi$ entre leurs groupes fondamentaux tempérés. Soit $l>2$ un nombre
  premier autre que $p$. On peut considérer un revêtement galoisien $X'_1$
  d'ordre $l$ de $X_1$ tel que la restriction à chaque composante
  irréductible de la réduction stable soit connexe mais tel que $X'_1$ soit
  totalement décomposé au-dessus du point double de la réduction (et soit
  $X'_2$ le revêtement de $X_2$ correspondans à $X'_1$ par $\phi$, il
  vérifie les mêmes propriétés). Alors $\overline X'_{i}$ est une courbe de
  Mumford (le revêtement correspondant de chaque composante irréductible de
  la réduction stable de $X_{i}$ est ramifié en au plus deux points, donc
  c'est un revêtement par une droite projective) dont l'arbre a $l$ arêtes,
  chacune de longueur $h_{i}$. Ainsi d'après le théorème~\ref{mumford}
  $h_1=h_2$, ce qui cl\^ot le cas de la droite épointée.
\end{rem}

Commençons par donner une esquisse de la preuve de~\ref{mumford}. Fixons deux points
terminaux du squelette $\mbb T:=\mbb T(\Omega_1)\simeq\mbb T(\Omega_2)$. Les points correspondants de $\mbf P^1\backslash
\Omega_i$ jouerons le rôle de $0$ et $\infty$. Etant donné
un sous-arbre fini $\mbf K_0$ de $\mbb T$, Nous construirons
un $\mbf Z/p^h\mbf Z$-torseur $\widetilde X''_1$ sur $\Omega_1$ se
descendant à un revêtement topologique fini de $X_1$ qui, restreint à
$\Phi_1^{-1}(\mbf K_0)$, soit isomorphe à la restriction à $\Omega_1$ de $\mbf
G_m\stackrel{z\mapsto z^{p^h}}{\to}\mbf G_m $ et tel que le torseur $\widetilde
X''_2$ sur $\Omega_2$ correspondant à $\widetilde X''_1$ par $\phi$ (ce qui
a un sens car $\widetilde X''_1$ est défini sur un revêtement fini de
$X_1$) soit aussi isomorphe à un multiple (pour la structure de $\mbf Z/p^h\mbf Z$-modules sur l'ensemble des $\mbf Z/p^h\mbf Z$-torseurs) de $\mbf
G_m\stackrel{z\mapsto z^{p^h}}{\to}\mbf G_m $. Nous obtiendrons alors l'\'egalit\'e des deux distances sur $\mbb T(\Omega_1)\simeq\mbb T(\Omega_2)$
en appliquant le lemme~\ref{decompositionpuissance} et le résultat
combinatoire~\ref{combi}.\\
Pour ce faire, considérons le courant $\mcal C_0$ sur $\mbb T$ longeant la ligne liant $0$ à
$\infty$ et rendons le invariant par un sous-groupe d'indice fini $\Gamma'$
de
$\Gamma:=\Gal(\mbb T/\mbb G)$ (pour cela, il faudra rajouter comme hypothèse sur la
ligne $0$ à $\infty$ d'être stabilisée par un sous-groupe non
trivial de $\Gal(\mbb T/\mbb G)$), pour qu'il induise un $\mbf Z/p^h\mbf
Z$-torseur sur un revêtement topologique fini de $X_1$~; $\widetilde X''_1$
sera alors son pullback à $\Omega_1$. Nous choisirons
$\Gamma'$ pour que le courant ainsi défini $\mbf C_1$ coïncide avec $\mbf
C_0$ sur un sous-graphe $K'$ de $\mbb T$ suffisamment grand par rapport
$K_0$.\\
Alors, en utilisant les lemmes du paragraphe précédent, $\widetilde X''_1$
sera isomorphe à $\mbf
G_m\stackrel{z\mapsto z^{p^h}}{\to}\mbf G_m $ sur un sous-graphe $K$ de
$K'$ mais encore assez grand par rapport à $K_0$. Les
lemmes~\ref{decompositionpuissance} et~\ref{vertex} permettent alors de déterminer en quels
sommets de $K$ $\widetilde X''_1$ et $\widetilde X''_2$ seront décomposés. En particulier, en choisissant $K'$ de
façon appropriée, nous pourrons nous arranger pour que $\widetilde X''_i$ soit
totalement décomposé en tous les sommets de la frontière de $K$ autre que
ceux sur la ligne reliant $0$ à $\infty$. Le courant $\mbf C_2$ correspondant
à $\widetilde X''_2$ sera alors nul en tous ces points frontières et sera
donc égal sur $K$ à un multiple de $\mbf C_0$. Si  $K$ était assez grand
par rapport à $K_0$, $\widetilde X''_2$ sera donc isomorphe à un multiple (pour la structure de $\mbf Z/p^h\mbf Z$-module de l'ensemble des $\mbf Z/p^h\mbf Z$-torseurs)
de la restriction de $\mbf
G_m\stackrel{z\mapsto z^{p^h}}{\to}\mbf G_m $, comme souhait\'e.
\begin{figure}[ht]
\begin{center}
\includegraphics[scale=0.6]{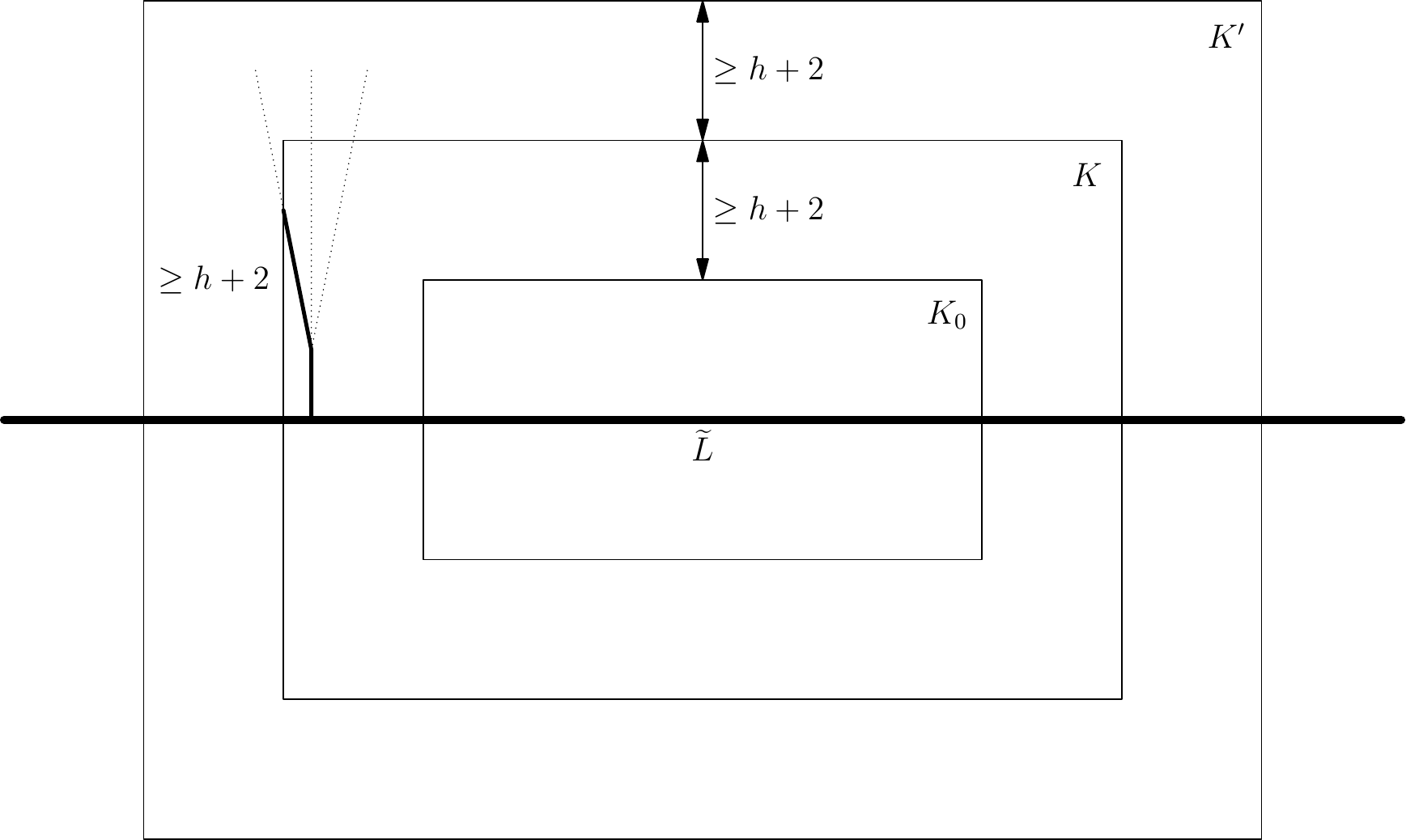}
\end{center}
\end{figure}

\begin{proof}[D\'emonstration de \ref{mumford}]
Un \emph{cycle} d'un graphe est une suite cyclique d'arêtes orientées du
graphe telle que le but d'une arête coïncide avec la source de l'arête
suivante, et qui ne repasse jamais par le même sommet ou par la même arête
(non orientée).\\
Soit $C$ un cycle du graphe $\mbb G$, et notons $\lg_{i}(C)$ la longueur de $C$
pour la métrique $d_{i}$ sur $\mbb G$. Soit $\widetilde C$ le revêtement
universel de $C$, soit $\widetilde C\to\mbb T$ un relèvement de
 $\widetilde C\to \mbb G$  et soit
$z_0$ un sommet de $\mbb T$ qui appartient à
$\widetilde C$, numérotons alors les sommets $(z_j)_{j\in \mbf Z}$ de
$\widetilde C$ ayant la même image que $z_0$ dans $\mbb G$. Soit $L$ un
autre cycle (muni d'une orientation) de $\mbb G$ (il doit exister une autre
boucle puisque $g>1$), soit $\widetilde L$ un relèvement du revêtement
universel de $L$ à
$\mbb T$, soit $r_{i}=d_{i}(\widetilde L,z_0)$ (on peut supposer, quitte à
renuméroter les $z_j$ que
pour $n\geqslant 0$, $d_{i}(\widetilde L,z_n)=r_{i}+n\lg_{i}(C)$) et
soit $z'_{0}$ le point de $\widetilde L$ le plus proche de $z_0$ (cela ne
dépend pas de $i$).\\
Si $z$ est un sommet de $\mbb T$, notons $F_z$ la composante connexe de
$\Omega\backslash\{\text{arêtes ouvertes de } \widetilde L\}$ qui contient $z$.\\

Soit $e\geqslant 1$ un entier.\\
Soit $K_0$  un sous-graphe connexe fini de $\mbb T$ contenant $z_1$ (ainsi
$\Phi_i^{-1}(|K_0|)$ est compact par propreté de $\Phi_i$).\\
Soit $K$ un sous-graphe connexe fini de $\mbb T$  tel que, pour tout $z\in
\Phi_i^{-1}(|K_0|)$,
$U_{z_,h+2}\subset \Phi_i^{-1}(|K|)$ et tel que, pour tout sommet $z$ de
$\widetilde L\cap K_{i}$, $\{z'\in 
F_z|d(z',z)\leqslant h+2\}\subset \Phi_i^{-1}(|K|)$ (en particulier, si
$z'$ est un sommet de la frontière de $K$ dans $\mbb T$ qui n'est pas une
extrémité du segment
$\widetilde L\cap K$, $d(z',\widetilde L)\geqslant h+2$). Soit $K'_{i}$
un sous-graphe compact de $\mbb T$ coïncidant pour $i=1$ et $i=2$ tel que
$\Phi_{i}^{-1}(|K'_{i}|)$ contient $U_{z,h+2}$ pour tout $z$ dans $\Phi^{-1}(K_{i})$.\\

Soit $\Gamma=\Gal(\mbb T/\mbb G)$, soit $H=\Stab(\widetilde L) (\simeq \mbf
Z)$ et soit $\Gamma'$ un sous-groupe d'indice fini de $\Gamma$ tel que,
pour tout $g\neq 1 \in \Gamma'$, $\Phi_{i}^{-1}(|K'|)\cap
g\cdot\Phi_{i}^{-1}(|K'|)=\emptyset$ et pour tout $g\in \Gamma\backslash H$,
$d_{i}(g\cdot\widetilde L,\widetilde L)>\diametre_{i}(|K'|)$.\\
Un tel $\Gamma'$ existe. En effet, $A:=\{g\neq
1\in \Gamma|\Phi^{-1}(|K'|)\cap g\cdot\Phi^{-1}(|K'|)\neq \emptyset\}$ est fini
par compacité de $\Phi^{-1}(|K'|)$. Comme $\Gamma$ est résiduellement fini, il
existe $\Gamma'_1$ d'indice fini dans $\Gamma$ qui n'intersecte pas $A$. $B:=\{g\in \Gamma/H-H |d_{i}(g\cdot \widetilde
L,\widetilde L)\leqslant \diametre_{i}(|K'|)\}$ est aussi fini, et comme
 $\overline H\cap \Gamma=H$ dans le complété profini de $\Gamma$, il
existe $\Gamma'_2$ d'indice fini dans $\Gamma$ contenant $H$ tel que
$\Gamma'_2\cap B\cdot H=\emptyset$. On peut alors choisir $\Gamma'=\Gamma'_1\cap\Gamma'_2$.\\
Soit $H'=H\cap \Gamma'$. Soit
$X'_{i}=\Omega_{i}/\Gamma'$ : c'est un revêtement topologique fini de
$X_{i}$, et l'isomorphisme $\phi:\gtemp(X_1)\simeq \gtemp(X_2)$
induit un isomorphisme
$\phi':\gtemp(X'_1)\simeq\gtemp(X'_2)$.\\

Soit $\mcal C_0$ le courant sur $\mbb T$ tel que  $\mcal C_0(e)=+1$ si $e$
est une arête de
$\widetilde L$ (et $e$ a la bonne orientation) et $0$ sinon (sauf si $e$
est une arête de $\widetilde L$ avec la mauvaise orientation, auquel cas
$\mcal C_0(e)=-1$)~: ce courant est invariant par $H$. Soit $\mcal C_{1}=\sum_{g\in
  \Gamma'/H'}g\cdot \mcal C_0$. C'est un courant sur $\mbb T$, invariant par
$\Gamma'$ et qui coïncide avec $\mcal C_0$ sur $K'$.\\
Soit $f_1\in O(\Omega_1)^*$ "la" fonction inversible correspondante, et soit
$X''_1$ un $(\mbf Z/p^h\mbf Z)$-torseur de $X'_1$ correspondant à ce
courant, c'est-à-dire tel que son pullback $\widetilde X''_1$ à
$\Omega_1$ soit isomorphe à $\Omega_1\times_{\mbf G_m}\mbf G_m\to
\Omega_1$ où le produit fibré est pris à gauche le long de $f_1$ et à
droite le long de $z\mapsto z^{p^h}$. Soit
$X''_{2}=\phi'^*X''_1$ ($X''_2$ n'a a priori aucune raison de correspondre
à $\mcal C_1$).\\
Soit également $f_{0,i}\in O(\Omega_{i})^*$ "la" fonction inversible
correspondant au courant $\mcal C_0$ et soit $\widetilde X_{0,i}$ le $(\mbf
Z/p^h\mbf Z)$-torseur correspondant sur $\Omega_{i}$. Rappelons que ce
torseur est totalement décomposé au-dessus du point $z\in \mbb T(\Omega_{i})$
si et seulement si $d_{i}(z,\widetilde L)>h+\frac{1}{p-1}$.\\

D'après le corollaire~\ref{totdecrev} appliqué à $U=\Phi^{-1}(|K'|)$ et à $V=\Phi^{-1}(|K|)$,
le torseur $\widetilde X''_1-\widetilde X_{0,1}$ sur
$\Omega_1$, qui correspond au courant $\mcal C_{1}-\mcal C_0$ qui
est nul au voisinage de $K'$, est totalement décomposé sur $\Phi^{-1}(|K|)$ puisque
pour tout $z\in \Phi^{-1}(|K|)$, $U_{z,h+2}\in \Phi^{-1}(|K'|)$. Ainsi, pour
$z\in K$, $\widetilde X''_1$ est décomposé si et seulement si $\widetilde
X_{0,1}$ l'est, si et seulement si $d_{1}(z,\widetilde L)>h+\frac{1}{p-1}$.\\
En particulier $\widetilde X''_1$ est décomposé en les sommets de la
frontière de $\mbb K$ qui ne sont pas des points terminaux de $K\cap
\widetilde L$. Ainsi, d'après le lemme~\ref{vertex}
appliqué à $X''_1$ et $X''_2$ (comme $\Omega_{i}\to
X'_{i}$ est un revêtement topologique, $\widetilde X''_{i}\to
\Omega_{i}$ est décomposé au-dessus d'un point si et seulement si
$X''_{i}\to X'_{i}$ l'est au-dessus de l'image de ce point), $\widetilde
X''_{0,2}$ est aussi décomposé au-dessus des sommets de la frontière de $K$
qui ne sont pas les points terminaux de $K\cap \widetilde L$.\\

Soit $\mcal C_{2}$ un courant sur $\mbb T(\Omega_2)$
correspondant au $\mbf Z/p^h\mbf Z$-torseur $\widetilde X''_2$ (le torseur
correspondant à
$\widetilde X''_2$  est bien défini seulement modulo $p^h$). D'après~\ref{nondecrev},
la restriction de $\mcal C_{2}$ à
l'étoile d'un sommet de la frontière de $K$ qui n'est pas un point terminal de $K\cap\widetilde
L$ est nul modulo $p^h$. On en déduit que, modulo $p^h$, la
restriction de $\mathcal C_{2}$ à l'étoile de $K$ doit être congrue à la restriction de
$a\mcal C_0$ pour un certain $a\in \mbf Z$. Quitte à ajouter à $\mcal
C_{2}$ un courant qui est multiple de $p^h$, on peut supposer que $\mcal
C_{2}-a\mcal C_0$ est nul sur l'étoile de $K$ (car tout courant avec bord
sur l'étoile de $K$, \ie
qui vérifie la loi de Kirchoff en tout sommet de $K$ mais
sans condition sur la frontière de l'étoile de $K$, peut s'étendre en un courant sur tout $\mbb T$).\\
Ainsi en appliquant le corollaire~\ref{totdecrev} à $U=\Phi^{-1}(|K|)$ et à
$V=\Phi^{-1}(K_0)$, on en déduit que $\widetilde X''_{2}-a\widetilde
X_{0,2}$ est décomposé au-dessus de $|K_0|$, donc si $z$ est un sommet de
$K_0$, $\widetilde X''_{2}$ est décomposé au-dessus de
$z$ si et seulement si $\widetilde X_{0,2}$ l'est ($a$ est nécessairement
non nul modulo $p^h$ car $\widetilde X''_{2}$ ne peut pas être décomposé en
$z'_0$), si et seulement si
$d_{2}(z,\widetilde L)>h+\frac{1}{p-1}$, d'après le lemme~\ref{decompositionpuissance}.\\

Par conséquent, $d_{2}(z,\widetilde L)>h+\frac{1}{p-1}$ si et seulement si
$d_{2}(z,\widetilde L)>h+\frac{1}{p-1}$ si et seulement si
$d_{1}(z,\widetilde L)>h+\frac{1}{p-1}$ pour tout sommet $z\in K_0$.\\
Comme on peut choisir $K_0$ aussi gros que voulu, on en déduit que pour tout $z$ de $\mbb T$, $d_{1}(z,\widetilde L)>h+\frac{1}{p-1}$
si et seulement si $d_{2}(z, \widetilde L)>h+\frac{1}{p-1}$, et ceci pour
tout entier $h\geqslant 1$.\\
Donc $\max(1,\lceil d_{2}(z,\widetilde
L)-\frac{1}{p-1}\rceil)=\max(1,\lceil d_{1}(z,\widetilde
L)-\frac{1}{p-1}\rceil)$.\\
En l'appliquant à $(z_j)_{j\geqslant 0}$, on obtient que pour tout
$j\geqslant 0$, \[\max(1,\lceil
j\lg_{1}+r_{1}-\frac{1}{p-1}\rceil)=\max(1,\lceil
j\lg_{2}+r_{2}-\frac{1}{p-1}\rceil).\]
On en déduit que pour tout cycle $C$ de $\mbb G$,
\[\lg_{1}(C)=\lg_{2}(C).\]
Si l'on considère un revêtement topologique fini $\mbb G'$ de $\mbb G$,
l'isomorphisme $\phi$ de groupes fondamentaux tempérés induit un
isomorphisme entre les groupes fondamentaux tempérés des revêtements
topologiques finis $X'_1$ et $X'_2$ de $X_1$ et $X_2$ qui correspondent au
revêtement $\mbb G'\to\mbb G$. $X'_1$ et $X'_2$ sont aussi des courbes de
Mumford, et on peut donc appliquer le résultat qu'on vient d'obtenir à
$\mbb G'$. On obtient donc que $\lg_{1}(C)=\lg_2(C)$ pour tout cycle de
$\mbb G'$.\\
Par définition de la réduction stable, toutes les composantes de la
réduction stable de $X_i$, qui sont toutes rationnels car $X_i$ est une
courbe de Mumford, ont au moins trois points marqués. Donc tous les sommets
de $\mbb G$ sont d'arité au moins trois.
On conclut grâce à la proposition~\ref{combi}.\end{proof} 
On peut généraliser le résultat précédent au cas de courbes de Mumford \'epoint\'ees. Soit $X_1$,
$X_2$ sur $\overline K$ (mais déjà définies sur $K$) deux courbes de
Mumford hyperboliques connexes de type $(g,n)$. Soit
$\phi:\gtemp(X_1)\simeq\gtemp(X_2)$ un isomorphisme. Il induit un
isomorphisme entre les graphes de leurs réductions stables $\mbb
G_1\simeq\mbb G_2$, et même entre les semigraphes d'anabélioïdes
$\mcal G_1^c\simeq\mcal G_2^c$.
\begin{cor}
L'isomorphisme de graphes $\mbb G_1\simeq \mbb G_2$ induit par $\phi$ est
un isomorphisme de graphes métriques.
\end{cor}
\dem
On a déjà traité le cas $g=0$. Supposons donc $g>0$.\\
Soit $\overline X_i$ la compactification de $X_i$. $\phi$ induit un
isomorphisme $\bar \phi:\gtemp(\overline X_1)\simeq\gtemp(\overline X_2)$
d'après~\ref{reductioncusps}
Soit $\mbb G'_i$ le graphe de la réduction stable de la compactification de
$\overline X_i$. On a un plongement isométrique $|\mbb G'_i|\to |\mbb G_i|$ compatible avec
les isomorphismes induits par $\phi$ et $\bar \phi$, c'est une équivalence d'homotopie.\\
Soit $e$ une arête de $\mbb G$. Distinguons trois cas~:
\begin{itemize}
\item $e$ est déjà une arête de $\mbb G'$. On sait déjà grâce
  au cas propre que $d_1(e)=d_2(e)$.
\item $e$ est un morceau d'une arête $e'$ de $\mbb G'$ qui a été
  divisée. $e$ a donc deux sommets distincts (le cas $g=1,n=1$ ayant déjà
  été traité). $e$ a nécessairement un sommet $v$ qui n'est pas un sommet
  de $\mbb G'$, et donc il existe un point cuspidal $x$ du semigraphe
  $\mbb G^c$ tel que $v$ soit le point de $|\mbb G'|$ le plus
  près de $x$. Soit $\mbb H'$ un revêtement (topologique) connexe d'ordre 2
  de $\mbb G'$, $\mbb H^c$ le revêtement correspondant de
  $\mbb G^c$ et $Y_i$ le revêtement topologique correspondant de
  $X_i$. Soit $x_1$ et $x_2$ (resp. $v_1$ et $v_2$, $e_1$ et $e_2$) les
  préimages de $x$ 
  (resp. $v$, $e$) et $L$ un
  chemin injectif reliant $x_1$ à $x_2$.  Soit $\mcal C$ le courant sur
  $\mbb H^c$ qui vaut
  $1$ le long de $L$ et $0$ ailleurs. Soit $l$ premier à $p$, et soit $Z_i$
  le $\mbf Z/l\mbf Z$-torseur sur $Y_i$ correspondant à $\mcal C$. Comme en
  tout sommet le
  courant est non nul en au plus deux branches, $Z_i$ est encore une courbe
  de Mumford. Soit
  $\mbb I$ le semigraphe de la réduction stable de $Z_i$~: une arête $a$ de
  $\mbb H$ a $1$ préimage de longueur $d_i(a)/l$ si $a$ est sur $L$ et $l$
  préimages de longueurs $d_i(a)$ sinon. La réalisation
  géométrique du graphe $\mbb I'$ de la
  réduction stable de la compactification $\overline Z_i$ de $Z_i$ contient la
  préimage $I''$ de $|\mbb H'|$ dans $|\mbb I|$ ($\mbb I'$ ne dépend pas de
  $i$ grâce à~\ref{reductioncusps}). Mais l'unique
  préimage $v'_1$ (resp. $v'_2$) de $v_1$ (resp. $v_2$) dans $I''$
  est d'arité $l+1$ donc d'arité au moins $l+1>2$ dans $|\mbb I'|$, donc $v$ est un sommet de $\mbb I'$.\\
Si l'autre sommet de $e$ était un sommet de $\mbb G'$, toute préimage $e''$
de $e$ dans $\mbb I$ serait une arête de $\mbb I'$, donc
$d_1(e'')=d_2(e'')$. Mais comme $d_i(e)=jd_i(e'')$ où $j=1$ ou $l$ mais ne
dépend pas de $i$, on déduit du premier cas appliqué à $Z_1$ et $Z_2$ que $d_1(e)=d_2(e)$.\\
Si aucun des deux sommets de $e$ n'était un sommet de $\mbb G'$,  une
préimage $e''$ de $e$ dans $\mbb I$ aurait maintenant au moins un de ces sommet
qui est dans $\mbb I'$, donc $d_1(e'')=d_2(e'')$ d'après le cas qui
précède appliqué à $Z_1$ et $Z_2$. On conclut alors comme dans le cas précédent.
\item $e$ n'est pas dans $|\mbb G'|$. Soit $v$ le sommet de $e$ le plus
  éloigné de $|\mbb G|$. Soit $e_1$ et $e_2$ deux autres arêtes distinctes
  ayant $v$ pour sommet. Soit $v_1$ (resp. $v_2$) un point cuspidal du
  semigraphe $\mbb G^c$ plus proche de $e_1$ (resp. $e_2$) que de
  $v$. Soit alors $L$ l'unique chemin injectif reliant $v_1$ à $v_2$. $L$
  passe par $v$ mais ne passe pas par l'intérieur de $e$. 
  soit $\mcal C$ le courant sur $\mbb G^c$ qui suit $L$. Soit $l>1$
  premier à $p$ et soit $Y_i$ le $\mbf Z/l\mbf Z$-torseur sur $X_i$ défini
  par $\mcal C$. C'est encore une courbe de Mumford puisqu'en tout sommet
  le courant est non nul sur au plus deux branches. Soit $\mbb H^c$ le semigraphe correspondant et
  $\mbb H$ le graphe de la réduction stable de la compactification
  $\overline Y_i$ de $Y_i$ ($\mbb H$ ne dépend pas de $i$). $e$
  a $l$ préimages, toutes de longueur $d_i(e)$. De plus une telle préimage
  est dans $|\mbb H|$. Donc $d_1(e)=d_2(e)$ d'après le cas précédent
  appliqué à $Y_1$ et $Y_2$. 
\end{itemize} 
\findem
\subsection{Un résultat combinatoire}

\begin{prop}\label{combi} Soit $\mbb G$ un graphe fini dont tous les sommets sont
  d'arité au moins 3. Soit $f:\{\text{arêtes de } \mbb G\}\to \mbf R$. On
  notera encore $f$ la fonction induite sur l'ensemble des arêtes d'un
  revêtement de $\mbb G$. Posons, pour $C$ un cycle d'un revêtement de
  $\mbb G$, $f(C)=\sum_{x\in\{\text{arêtes de } C\}}f(x)$.\\
Si $f(C)$ est nul pour tout cycle $C$ de tout revêtement fini de $\mbb G$, alors
$f$ est nulle.\end{prop}
\dem
Commençons par remarquer que si $\mbb G$ est un graphe fini dont toutes les arêtes
sont d'arité au moins 3, et si $\mbb H$ est un sous-graphe connexe tel que
le nombre de demi-arêtes de $\mbb G\backslash \mbb H$ dont l'extrémité est dans
$\mbb H$ est strictement inférieur à 3, alors $\mbb H$ n'est pas un
arbre (si $\mbb H$ est un arbre avec au moins une arête, $\mbb H$ a au
moins deux sommets d'arité 1, et donc déjà 4 demi-arêtes de $\mbb
G\backslash \mbb H$ doivent terminer en l'un de ces deux sommets ; si $\mbb
H$ est réduit à un point, c'est également évident).\\

Revenons à l'énoncé. Nous allons procéder par récurrence sur le nombre d'arêtes de $\mbb G$.\\
Soit donc $(\mbb G,f)$ un graphe à $n\geqslant 1$ arêtes muni d'une fonction $f$ sur
l'ensemble des  arêtes de $\mbb G$ vérifiant les hypothèses de l'énoncé, et
supposons l'énoncé vrai dès que $\mbb G$ a strictement moins de $n$
arêtes.\\
On peut supposer $\mbb G$ connexe (sinon on applique l'hypothèse de
récurrence à chaque composante connexe).\\
Soit $e$ une arête de $\mbb G$, et commençons par supposer que :
\begin{enumerate}
\item les deux sommets de $e$ sont distincts. Soient $(m,n)$ les arités des
  sommets de $e$ dans $\mbb G\backslash\{e\}$. On a, par hypothèse sur $\mbb
  G$, $m\geqslant   2$ et $n\geqslant 2$.
  \begin{enumerate} 
    \item Si $m\geqslant 3$ et $n\geqslant 3$, $\mbb G'=\mbb
    G\backslash\{e\}$ est encore un graphe dont tous les sommets sont d'arité au
    moins $3$. De plus, tout revêtement de $\mbb G'$ peut se prolonger en
    un revêtement de $\mbb G$ donc si $C$ est un cycle d'un revêtement
    fini de
    $\mbb G'$, $C$ est aussi un cycle d'un revêtement fini de $\mbb G$ et
    donc $f(C)=0$. On peut ainsi appliquer l'hypothèse de récurrence à $\mbb G'$
    et à $f$ : $f(x)=0$, pour toute les arêtes de $\mbb G$ autre que $e$.
    \begin{enumerate}
      \item Si $\mbb G'$ est connexe, on peut trouver un cycle $C$ de $\mbb
        G$ passant par $e$, et alors $f(e)=f(C)=0$. D'où le résultat.
      \item Si $\mbb G'$ a deux composantes connexes $A$ et $B$, qui ne
        peuvent pas être des arbres d'après la remarque faite en début de démonstration, on considère $A'$ et $B'$
        des revêtements connexes d'ordre 2 de $A$ et $B$ respectivement,
        qu'on recolle en un revêtement $\mbb G'$ de $\mbb G$ d'ordre 2.\\
        Il existe alors un cycle $C$ de $\mbb G'$ passant par les 2
        préimages de $e$. Alors $2f(e)=f(C)=0$. D'où le résultat.
    \end{enumerate}
    \item si $m\geqslant 3$ et $n= 2$ (ou l'inverse), et soit $a$
      et $b$ les deux arêtes partant du second sommet de $e$ (si $a$ et $b$
      ne sont en fait que les deux demi-arêtes d'une seule arête, le graphe
      a alors la même structure que pour le cas 2.(c).i, où il sera traité
      ; nous supposerons donc ici que $a$ et $b$ sont deux arêtes distinctes). Soit $\mbb
      G'$ le graphe obtenu à partir de $\mbb G$ en enlevant $e$ et en
      concaténant $a$ et $b$ en une seule arête notée $a+b$ (et posons
      $f(a+b)=f(a)+f(b)$). Alors $(\mbb G',f)$ vérifie les hypothèses de
      l'énoncé, donc $f(x)=0$ pour toute arête de $\mbb G$ autre que $a,b$
      et $e$. En fonction du nombre de composantes connexes de $\mbb G''=\mbb
      G\backslash\{a,b,e\}$, distinguons plusieurs cas :
     \begin{enumerate}
       \item 
           Dans le cas où $\mbb G''$ n'a qu'une composante connexe, on peut contracter $\mbb
           G''$ en un point pour obtenir un graphe $\mbb G_1$ à trois arêtes
           (en effet toute cycle $C_1$ de $\mbb G_1$ admet un
           relèvement $C$ à $\mbb G$ par connexité de $\mbb G''$, et
           $f(C_1)=f(C)$ puisque $f$ est nulle sur $\mbb G''$), et l'on en déduit que $f(a)+f(b)=0$,
           $f(a)+f(e)=0$ et $f(b)+f(e)=0$. D'où le résultat.
\begin{figure}[ht]
\begin{center}
\includegraphics{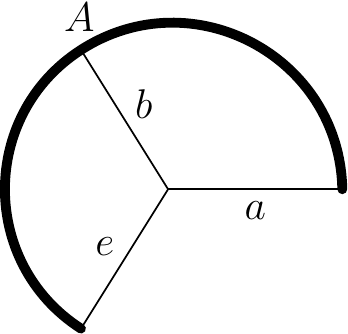}
\includegraphics{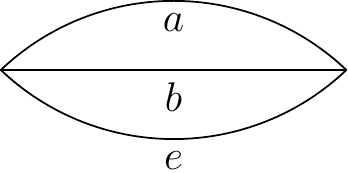}
\end{center}
\end{figure}
       \item
           Dans le cas où $\mbb G''$ a deux composantes $A$ et $B$ comme
           sur le dessin (maintenant $a$, $b$ et $e$ joueront le même rôle
           et le raisonnement sera le même si on les échange), on commence
           par considérer un revêtement connexe $\mbb G_1$ d'ordre 2 de $\mbb G$
           dont les restrictions $A'$ et $B'$ à $A$ et à $B$ sont connexes (il en existe
           car $A$ et $B$ ne sont pas des arbres d'après la remarque du
           début de la démonstration), puis on
           contracte $A'$ et $B'$ en un graphe $\mbb G_2$. On obtient
           $f(b)+f(e)=0$, $2f(a)+2f(b)=0$, $2f(a)+2f(e)=0$, d'où le
           résultat.
 \begin{figure}[ht]
 \begin{center}
 \includegraphics{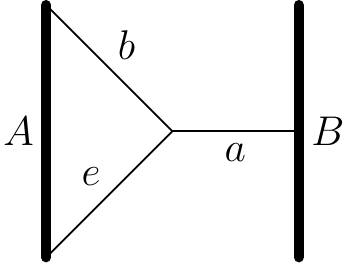}\hfill
 \includegraphics{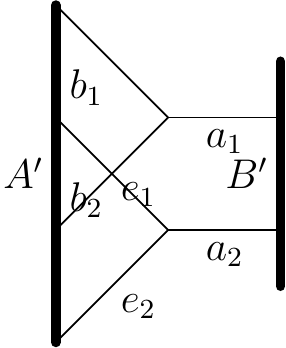}\hfill
 \includegraphics{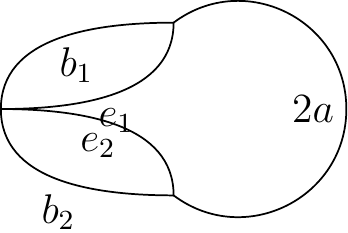}
 \end{center}
 \end{figure}
         \item
           Dans le cas où $\mbb G''$ a trois composantes connexes $A$, $B$ et
           $C$, on commence par considérer un revêtement connexe $\mbb G_1$
           d'ordre 2 de $\mbb G$ dont les restrictions $A'$, $B'$ et $C'$ à
           $A$, $B$ et $C$ sont connexes, puis on contracte $A'$, $B'$ et
           $C'$. On en déduit
           $2f(a)+2f(b)=0,2f(b)+2f(e)=0,2f(e)+2f(a)=0$. D'où le résultat.
\begin{figure}[ht]
\begin{center}
\includegraphics{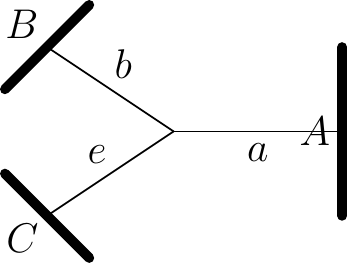}
\includegraphics{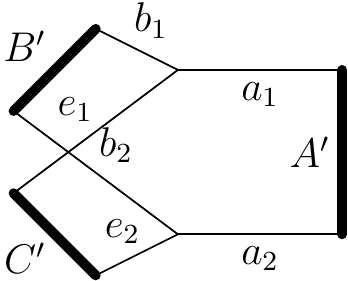}
\includegraphics{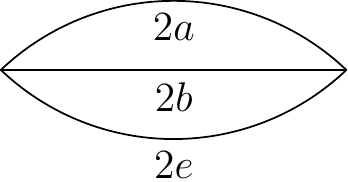}
\end{center}
\end{figure}
       \end{enumerate}
       \item Si $m=n=2$, soient $a,b$ les (demi-)arêtes partant d'un sommet et $c,d$ les
       (demi-)arêtes partant de l'autre sommet de $e$. Soit $\mbb G'$ le graphe obtenu
       à partir de $\mbb G$ en retirant $e$ puis en concaténant $a$ et $b$ en
       $a+b$ d'un
       coté, et $c$ et $d$ en $c+d$ de l'autre. Si l'on pose
       $f(a+b)=f(a)+f(b)$ et $f(c+d)=f(c)+f(d)$, $(\mbb G',f)$ vérifie les
       condition de l'énoncé, et donc, par hypothèse de récurrence,
       $f(x)=0$ pour toute arête $x$ de $\mbb G$ différente de $a,b,c,d$ et
       $e$. Soit $\mbb G''=\mbb G\backslash\{a,b,c,d,e\}$. En fonction des
       composantes connexes de $\mbb G''$, distinguons plusieurs cas :
       \begin{enumerate}
         \item
           Supposons que $\mbb G$ ait deux composantes connexes, $A$
           contenant une extrémité de $a$ et de $b$, et $B$ contenant une
           extrémité de $c$ et de $d$ (si $(a,b)$ ou $(c,d)$ sont les
           demi-arêtes d'une même arête, le graphe à la même structure que
           dans le cas 2.(c).ii.B ; si $(a,b)$ et $(c,d)$ ne forment que
           deux arêtes, la structure est celle du cas dégénéré de 2.(c).(ii)). On Commence par considérer un
           revêtement connexe $\mbb G_1$ d'ordre 2 de $\mbb G$ dont les
           restrictions $A'$ et $B'$ à $A$ et $B$ sont connexes, puis
           on contracte $A'$ et $B'$. On obtient par exemple que
           $f(c)+f(d),f(a)+f(b),2(f(a)+f(e)+f(d)),2(f(a)+f(e)+f(c)$ et $2(f(b)+f(e)+f(d))$
           sont nuls. D'où le résultat.
\begin{figure}[ht]
\begin{center}
\includegraphics{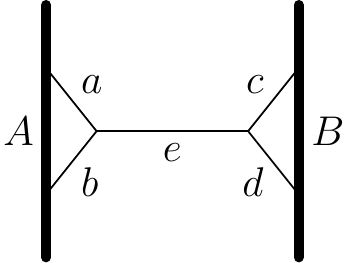}
\includegraphics{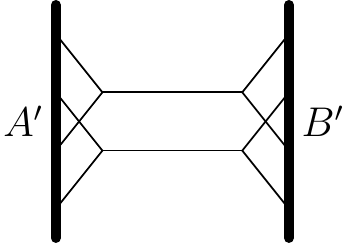}
\end{center}
\end{figure}
           \item
             Supposons que $\mbb G$ ait deux composantes connexes, $A$
             contenant une extrémité de $a$ et $c$, et $B$ contenant une
             extrémité de $b$ et $d$. On commence, comme d'habitude par
             considérer un revêtement connexe $\mbb G_1$ d'ordre 2 de $\mbb
             G$ dont les restrictions $A'$ et $B'$ à $A$ et $B$ sont
             connexes puis on contracte $A'$ et $B'$. On en déduit par
             exemple que
             $f(a)+f(b)+f(c)+f(d)=2(f(a)+f(e)+f(d))=2(f(b)+f(e)+f(c))=2(f(a)+f(c))=f(b)+f(c)+f(e)=0$. D'où
             le résultat.
\begin{center}
\includegraphics{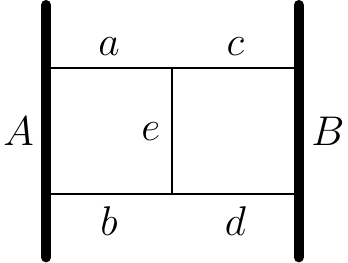}
\includegraphics{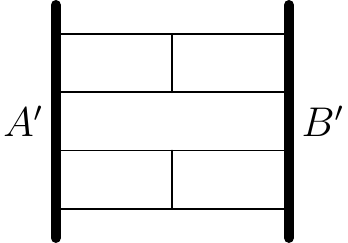}
\end{center}
             Si $(a,c)$ (ou symétriquement $(b,d)$) ne sont que les deux
             demi-arêtes d'une seule arête que l'on notera $a$, on
             considère un
             revêtement $\mbb G_1$ d'ordre 2 de $\mbb G$ dont les
             restrictions $B'$ à $B$ et $(a\cup e)'$ à $a\cup e$ sont
             connexes. On en déduit $f(a)+f(b)+f(c)= f(e)+f(b)+f(c)=
             2(f(a)+f(e)) = f(a)+f(e)+2f(b)=0$. D'où le résultat.
             \begin{center}
               \includegraphics{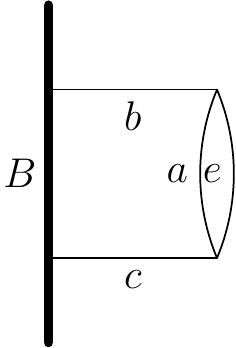}
               \includegraphics{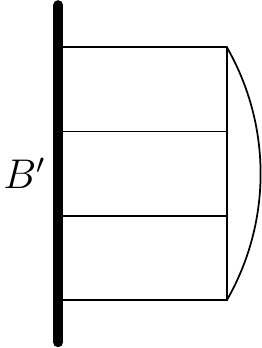}
             \end{center}
             Si $(a,c)$ et $(b,c)$ sont dégénérés en deux arêtes, alors
             $\mbb G$ n'a que deux sommets et trois arêtes reliant ces deux
             sommet, et le résultat est évident. 
             \item
             Supposons que $\mbb G$ ait deux composantes connexes, $A$
             contenant une extrémité de $a$,$b$ et $c$, et $B$ contenant une
             extrémité de $d$. On sait déjà par hypothèse de récurrence
             appliqué à $\mbb G'$ que $f(c)+f(d)=0$. Maintenant, en
             contractant $A$, on en déduit que
             $f(a)+f(b)=f(b)+f(c)+f(e)=f(a)+f(c)+f(e)=0$. D'où le résultat.
\begin{center}
\includegraphics{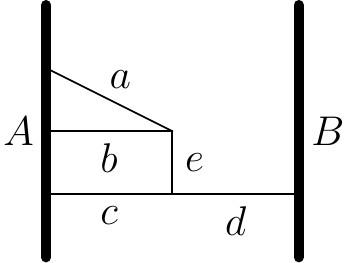}
\end{center}
             \item
             Supposons que $\mbb G$ ait trois composantes connexes, $A$
             contenant une extrémité de $a$, $B$ contenant une extrémité de
             $b$, et $C$ contenant une
             extrémité de $c$ et $d$. On commence, comme d'habitude par
             considérer un revêtement connexe $\mbb G_1$ d'ordre 2 de $\mbb
             G$ dont les restrictions $A'$, $B'$ et $C'$ à $A$, $B$ et $C$ sont
             connexes, puis on contracte $A'$, $B'$ et $C'$. On en déduit par
             exemple que $f(c)+f(d) = 2(f(a)+f(b))=2(f(a)+f(e)+f(c)) =
             2(f(b)+f(e)+f(c)) = 2f(a)+2f(e)+f(c)+f(d)=0$. D'où le
             résultat.
\begin{center}
\includegraphics{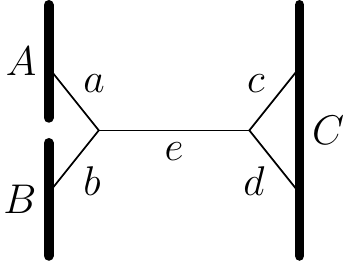}
\includegraphics{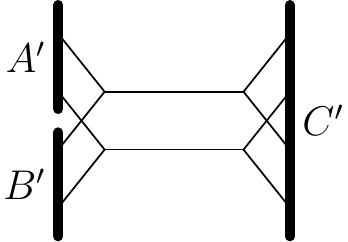}
\end{center}
             Si $(c,d)$ est dégénéré en une seule arête $c$, le graphe a la
             même structure qu'en 2.(c).ii.A.
             \item
             Supposons que $\mbb G$ ait trois composantes connexes, $A$
             contenant une extrémité de $a$, $B$ contenant une extrémité de
             $c$, et $C$ contenant une
             extrémité de $b$ et $d$. On commence, comme d'habitude par
             considérer un revêtement connexe $\mbb G_1$ d'ordre 2 de $\mbb
             G$ dont les restrictions $A'$, $B'$ et $C'$ à $A$, $B$ et $C$ sont
             connexes, puis on contracte $A'$, $B'$ et $C'$. On en déduit par
             exemple que $f(b)+f(d)+f(e)=f(b)+f(d)+f(e)+2f(a)=2(f(a)+f(b))
             = 2(f(d)+f(c))=2(f(c)+f(e)+f(a))=0$.
             \begin{center}
                 \includegraphics{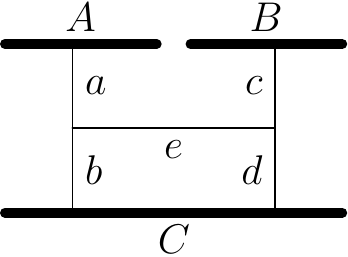}
                 \includegraphics{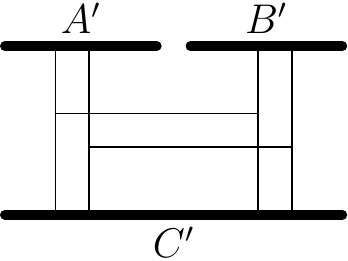}
               \end{center}
             Si $(b,d)$ est dégénéré en une seule arête $b$, considérons un
             revêtement $\mbb G_1$ d'ordre 2 de $\mbb G$ dont les restrictions à $A$,
             $B$ et $b\cap e$ sont connexes, puis on contracte $A'$ et $B'$. On en déduit que
             $2f(a)+f(e)+f(b) = 2(f(b)+f(e)) =2f(c)+f(b)+f(e) =
             2(f(a)+f(b)+f(c))=0$.
             \begin{center}
               \includegraphics{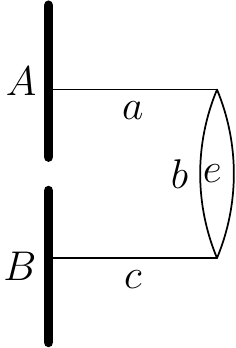}
               \includegraphics{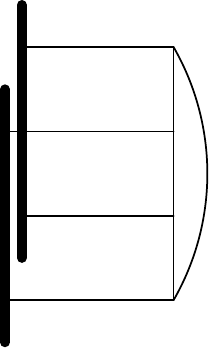}
             \end{center}
             \item
             Supposons que $\mbb G$ ait quatre composantes connexes, $A$
             contenant une extrémité de $a$, $B$ contenant une extrémité de
             $c$, $C$ contenant une
             extrémité de $b$ et $D$ contenant une extrémité de $d$. On commence, comme d'habitude par
             considérer un revêtement connexe $\mbb G_1$ d'ordre 2 de $\mbb
             G$ dont les restrictions $A'$, $B'$, $C'$ et $D'$ à $A$, $B$,
             $C$ et $D$ sont
             connexes, puis on contracte $A'$, $B'$, $C'$ et $D'$. On en déduit par
             exemple que $2(f(b)+f(d)+f(e))= f(b)+f(d)+f(e)+2f(a)=
             f(b)+f(d)+f(e)+2f(c)=2(f(b)+f(a))= 2(f(c)+f(d))=0$.
             \begin{center}
               \includegraphics{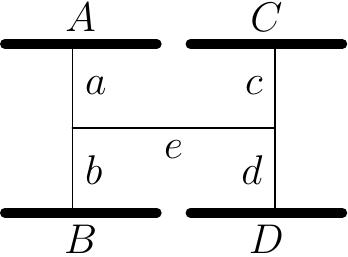}
               \includegraphics{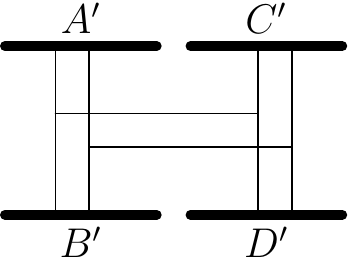}
             \end{center}
             \item
             Supposons que $\mbb G$ ait une seule composante connexe $A$,
             et contractons-la. On en déduit que $f(a)+f(b)=f(c)+f(d)=
             f(a)+f(c)+f(e)=f(b)+f(c)+f(e)=f(a)+f(d)+f(e)=0$. D'où le
             résultat.
             \begin{center}
               \includegraphics{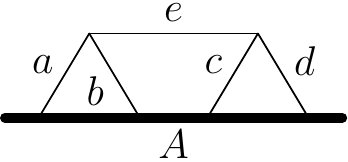}
             \end{center}
         \end{enumerate}
  \end{enumerate}
\item
Supposons que les deux extrémités de $e$ soit un même sommet, et soit $m$
l'arité de ce sommet dans $\mbb G\backslash\{e\}$. On a $m\geqslant 1$.
  \begin{enumerate}
    \item
      Supposons $m\geqslant 3$. Alors $\mbb G\backslash\{e\}$ vérifie les
      conditions de l'énoncé et donc $f(x)=0$ pour toute arête de $\mbb G$
      autre que $e$, et comme $e$ est une boucle, $f(e)=0$ aussi.
    \item
      Supposons $m=2$, et soit $a$ et $b$ les deux arêtes partant du sommet
      de $e$ (si $a$ et $b$ sont en fait les deux demi-arêtes d'une même
      arête, $\mbb G$ est réduit à deux arêtes formant chacune une boucle,
      et le résultat est évident). Soit $\mbb G'$ le graphe obtenu à partir
      de $\mbb G$ en supprimant $e$ et en concaténant $a$ et $b$ en une
      arête $a+b$ et posons $f(a+b)=f(a)+f(b)$. Alors $\mbb G'$ vérifie les
      hypothèses de l'énoncé et donc, par hypothèse de récurrence, $f$ est
      nulle en toute arête autre que $e$, $a$ et $b$. Distinguons selon les
      composantes connexes de $\mbb G''=\mbb G\backslash\{a,b,e\}$.
      \begin{enumerate}
        \item
          Supposons que $\mbb G''$ ait une unique composante connexe $A$. On consid\`ere un revêtement $\mbb
          G_1$ d'ordre 2 de $\mbb G$ dont la restriction $A'$ à $A$ est
          connexe et dont la restriction à $e$ est également connexe, puis
          on contracte $A'$ pour obtenir un graphe $\mbb G_2$. On en déduit
          $2f(a)+f(e)= 2f(b)+f(e)=2f(e)=0$. D'où le résultat.
\begin{center}
\includegraphics{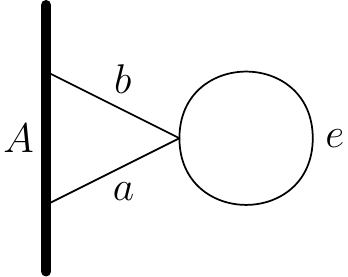}
\includegraphics{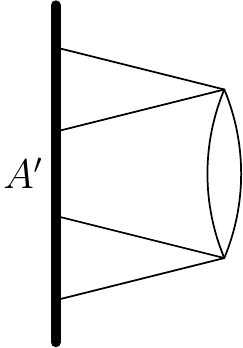}
\end{center}
        \item
          Supposons que $\mbb G''$ ait deux composantes connexes, $A$ ayant
          l'extrémité de $a$, et $B$ ayant l'extrémité de $b$. On consid\`ere un revêtement $\mbb
          G_1$ d'ordre 2 de $\mbb G$ dont les restrictions $A'$ et $B'$ à
          $A$ et $B$ sont
          connexes et dont la restriction à $e$ est également connexe, puis
          on contracte $A'$ et $B'$ pour obtenir un graphe $\mbb G_2$. On en
          déduit
          $2f(a)+f(e)= 2f(b)+f(e)=2f(e)=0$. D'où le résultat.
\begin{center}
\includegraphics{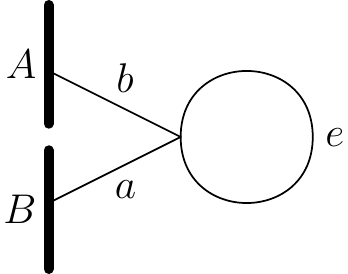}
\includegraphics{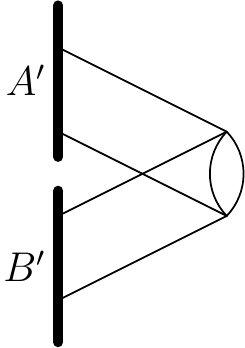}
\end{center}
      \end{enumerate}
   \item
     Supposons $m=1$, et soit $a$ l'arête partant l'extrémité de $e$. Soit
     $n$ l'arité de l'autre extrémité de $a$ dans $\mbb G'=\mbb G\backslash\{a,e\}$. On a nécessairement
     $n\geqslant 2$.
     \begin{enumerate}
       \item
         Supposons $n\geqslant 3$. Alors $\mbb G'$
         vérifie encore les conditions de l'énoncé, donc, par hypothèse de
         récurrence, $f$ est nulle sur $\mbb G'$. On consid\`ere maintenant un
         revêtement $\mbb G_1$ d'ordre 2 de $\mbb G$ dont les restrictions
         à $\mbb G'$ et à $e$ sont connexes, puis on contracte l'image
         réciproque de $\mbb G'$. On obtient $2f(e)=2f(a)+f(e)=0$. D'où le
         résultat.
\begin{center}
\includegraphics{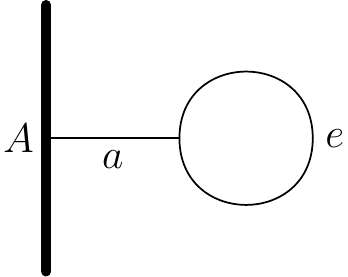}
\includegraphics{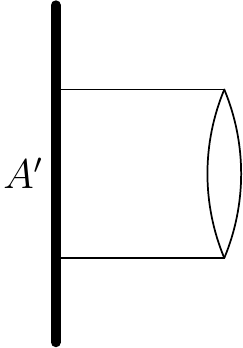}
\end{center}
       \item
         Supposons $n\geqslant 2$, et soit $c$ et $d$ les deux arêtes en
         question (si ce sont les deux demi-arêtes d'une unique arête,
         $\mbb G$ est réduit à deux boucles reliées entre elles par une
         arête ; graphe pour lequel on montre le résultat en considérant un
         graphe d'ordre 2 dont la restriction à chacune des boucles est
         connexe). Soit $\mbb G''$ le graphe obtenu à partir de $\mbb G'$
         en concaténant $c$ et $d$ en $c+d$ et en posant
         $f(c+d)=f(c)+f(d)$. $\mbb G''$ vérifie les conditions de l'énoncé,
         et donc $f$ est nulle en toute arête autre que $a,e,c$ et
         $d$ (et en fait $f(e)=0$ aussi). Soit $\mbb G'''=\mbb G\backslash\{a,c,d,e\}$, et distinguons
         suivant le nombre de composantes connexes de $\mbb G'''$.
         \begin{enumerate}
           \item
             Si $\mbb G'''$ a deux composantes connexes $C$ contenant
             l'extrémité de $c$ et $D$ contenant l'extrémité de $d$, on
             consid\`ere un revêtement $\mbb G_1$ de $\mbb G$ dont les
             restrictions $C'$, $D'$ et $e'$ à $C$, $D$ et $e$ sont connexes, puis on contracte
             $C'$, $D'$ et $e'$, on obtient
             $2f(a)+2f(c)=2f(a)+2f(d)=2f(c)+2f(d)=0$. D'où le résultat.
             \begin{center}
               \includegraphics{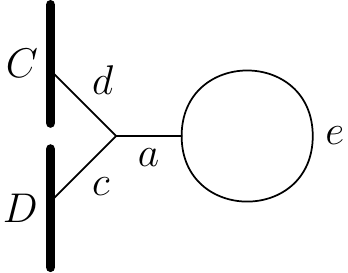}
               \includegraphics{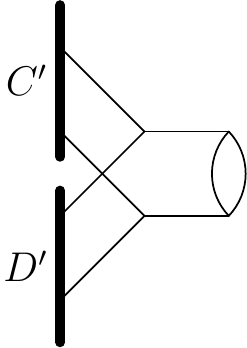}
             \end{center}
           \item
             Si $\mbb G'''$ a une composante connexe $A$, on consid\`ere un
             revêtement $\mbb G_1$ de $\mbb G$ dont les restrictions $A'$
             et $e'$ à $A$ et $e$ sont connexes, puis on contracte $A'$ et
             $e'$, on obtient
             $2f(a)+2f(c)=2f(a)+2f(d)=2f(c)+2f(d)=0$. D'où le résultat.
             \begin{center}
               \includegraphics{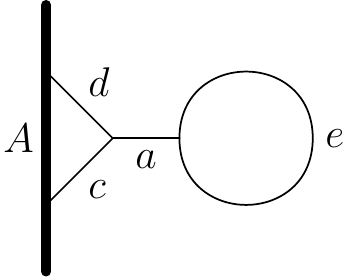}
               \includegraphics{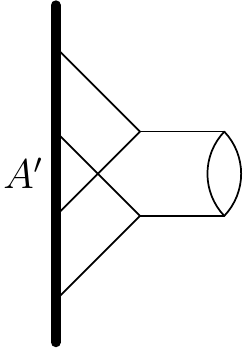}
               \end{center}
         \end{enumerate}
      \end{enumerate}
  \end{enumerate}
\end{enumerate} 
\findem
\chapter{Complexes classifiants de groupes}

Nous introduisons ici des objets combinatoires, que nous appellerons
complexes classifiants, généralisant la notion de
complexes de groupe de Haefliger (\cite{haef}). L'espace combinatoire de
base n'est pas nécessairement un ensemble simplicial. En fait, nous nous en
servirons plutôt avec des ensembles polysimpliciaux qu'avec des ensembles
simpliciaux. Les groupes servant de constituants aux complexes ne seront
pas non plus supposés discrets. En fait nous utiliserons des
complexes dont les composantes sont des groupes profinis.\\
A un tel complexe classifiant nous associerons une catégorie de revêtements
et un groupe fondamental, qui classifie ces revêtements. Ce groupe
fondamental sera encore un groupe topologique.\\
Nous définirons également une notion de revêtement tempérés dans ce contexte,et donc un groupe fondamental temp\'er\'e associ\'e.

\section{Groupes quasiprodiscrets}
Soit $G$ un groupe topologique. Soit $G\tEns$ la catégorie des
$G$-ensembles discrets, muni du foncteur canonique
$F:G\tEns \to \Ens$.\\
$G\tEns$ est un topos~(\cite[ex. A.2.1.6]{johnstone}). Il est clair
qu'il est connexe (son objet final n'a pas de sous-objet non trivial),
localement connexe et
atomique. $G\tEns \to \Ens$ est le pullback
d'un point du topos $G\tEns$ (et ce point est conservatif~;~\cite[ex. A.4.1.7]{johnstone}).\\

Notons $\hat G=\Aut F$, on a un morphisme évident $G \to \hat G$. $\hat G$ est un groupe topologique pour lequel
les $(\Stab^{\hat G}_{F(S),s})_{S\in \Ob(G\tEns),s\in F(S)}$ forment une base
de voisinages de $1$.\\
\begin{prop}
On a alors une immersion topologique $\hat G\to \varprojlim \hat
G/\Stab^{\hat G}_{F(S),s}= \varprojlim_H G/H$ o\`u $H$ parcours les
sous-groupes ouverts de $G$. Plus pr\'ecis\'ement, on a $\hat G=\Aut
F\subset \End F=\varprojlim_H G/H$.
\end{prop}
\dem
En effet, soit $\gamma\in \End
F$. Etant donné un sous-groupe ouvert $H$ de $\mbb G$, on considère le
$G$-ensemble $S=G/H$. Notons $\gamma_H$ l'image du point $H$ de $F(S)=G/H$
par $\gamma$, c'est un élément de $G/H$. La fonctorialité de $\gamma$
assure que l'on obtient un système compatible d'objets de $G/H$. D'où une
fonction $\End F\to\varprojlim_H G/H$.\\
 Réciproquement étant donné
$(\gamma_H)_H\in\varprojlim_H G/H$, définissons un endomorphisme $\gamma$
de $F$. Soit $S$ un $G$-ensemble et soit $s\in F(S)$. Soit
$H=\Stab^G_{F(S),s}$. Alors $\gamma(s)=g\cdot s$ où $g\in \gamma_H$. Si
$\phi:S'\to S$ est un morphisme de $G$-ensembles qui envoie $s'$ sur $s$,
on a $H'\subset H$ où $H'=\Stab^G_{s'}$ et $H=\Stab^G_s$. Soit
$g\in\gamma_{H'}$ (et donc $g\in\gamma_H$ par compatibilité de
$(\gamma_H)$). Alors $\gamma(s')=g\cdot s'=g\cdot\phi(s)=phi(g\cdot
s)=\phi(\gamma(s))$. Donc $\gamma$ définit bien un endomorphisme de $F$.
\findem

On d\'eduit en particulier de cette proposition que $G$ est d'image dense
dans $\hat G$.\\
On peut décrire plus directement la structure de monoïde sur $\varprojlim_H
G/H$ induite par celle de $\End F$. Si
$(\bar \alpha_H)_H,(\bar \beta_H)_H\in\varprojlim_H G/H$, leur produit est
donné par $\gamma_H=\alpha_{\beta_HH\beta_H^{-1}}\beta_H$. En effet
l'image$\beta\in \End(F)$ de $(\bar \beta_H)_H$ envoie $s=H\in G/H$ sur
$s'=\beta_HH$. Le stabilisateur de $s'$ est $\beta_HH\beta_H^{-1}$, donc
l'endomorphisme $\alpha$ envoie $s'$ sur
$s''=\alpha_{\beta_HH\beta_H^{-1}}\beta_HH$.\\
Si l'on note $\widetilde G$ le groupe $G$ mais muni de la topologie où les
sous-groupes ouverts de $G$ forment une base d'ouverts. Alors $\varprojlim
G/H$ n'est autre que le complété de $\widetilde G$ pour la structure
uniforme gauche de~\cite[\S{} III.3.1]{bourbakitopgen}. La multiplication est alors
celle donnée par~\cite[prop. III.3.6]{bourbakitopgen}.  
\begin{rem}
L'immersion $\Aut F\subset \End F$ n'est en g\'en\'eral pas surjective. Par exemple, soit $I$ un ensemble et soit $G$ le groupe des permutations de $I$ qui fixent tous les \'el\'ements sauf un nombre fini, muni de la topologie pour laquelle les fixateurs des parties finies de $I$ forment une base de voisinage ouverts de $1$, alors $\End F$ est l'ensemble de toutes les applications de $I$ dans $I$. Ainsi, $\End F$ contient des \'el\'ements non inversibles d\`es que $I$ est infini.
\end{rem}
\begin{dfn} On appelera groupe \emph{quasiprodiscret}\index{Groupe quasiprodiscret} (ou groupe qpd) un groupe
  topologique tel que $G\to\hat G$ soit un isomorphisme de groupes topologiques.\end{dfn}
En g\'en\'eral, l'homomorphisme $f:G\to \hat G$ induit un foncteur canonique $f^*:\hat G\tEns\to G\tEns$, et
  comme $\hat G=\Aut F$, $F(S)$ est naturellement muni d'une structure de
  $\hat G$-ensemble, ce qui donne un foncteur $g:G\tEns\to \hat
  G\tEns$ : ce sont en fait des foncteurs quasi-inverses (en particulier
  $f^*$ et $g$ sont des équivalences de catégories). En effet, $f^*g$
  est clairement isomorphe à l'identité, et $f^*$ est pleinement fidèle car
  $G$ est d'image dense dans $\hat G$.\\
Donc $\hat{\hat G}=\hat G$ ($\hat G$ est qpd). On appellera $\hat G$ le
complété qpd de $G$.\\ 

\emph{Remarque :} Un objet connexe $G/H$ de $G\tEns$ est localement
constant en tant qu'objet de $G\tEns$ si et seulement si $H$ contient un
sous-groupe ouvert distingué dans $G$. $G\tEns$ est donc un topos galoisien si et seulement si $\hat G$ est prodiscret.\\

Soit $G$ un groupe qpd, et soit $H$ un sous-groupe de $G$. On notera
$\overline H$ l'adhérence de $H$ et $\overline{\overline H}$ l'intersection
des sous-groupes ouverts contenant $H$.\\
On a $\overline{\overline H}=\overline H$ dans les deux cas suivants :
\begin{itemize}
\item $G$ est un groupe prodiscret.
En effet, comme la famille $(gO)$, $O$ d\'ecrivant les sous-groupes ouverts distingu\'es de $G$ et $g$ d\'ecrivant les \'el\'ements de $G$, forme une base d'ouverts de $G$, \[\begin{array}{rcl} G\backslash \overline H=\mring{G\backslash H} & = & \bigcup_{(O,g)|H\cap gO=\emptyset}gO \\ 
& = & \bigcup_O\bigcup_{g\notin HO} gO.\end{array}\] o\`u $HO$ est un sous-groupe ouvert de $G$ car $O$ est distingu\'e. Donc
\[\overline H=\bigcap_O\bigcap_{g\in HO} gO=\bigcap_{O}HO.\]
\item $H$ est distingué. En effet, on a $\overline H=\bigcap_{O}HO$ o\`u $O$ d\'ecrit les sous-groupes ouverts de $G$. Comme $H$ est maintenant suppos\'e distingu\'e, $HO$ est encore un sous-groupe ouvert de $G$.
\end{itemize}
Soient $G_1$, $G_2$ deux groupes qpd et $f:G_1\to G_2$ un morphisme,
$f^*:G_2\tEns \to G_1\tEns$ le foncteur induit.
\begin{prop}
\begin{itemize}
\item $f^*$ est fidèle ;
\item si $f^*$ est pleinement fidèle, $\overline{\overline{\Image f}}=G_2$
  ;
\item si $\overline{\Image f}=G_2$, $f^*$ est pleinement fidèle ;
\item $f^*$ est une équivalence de catégorie si et seulement si $f$ est un isomorphisme.
\end{itemize}\end{prop}
\dem
Soit $\psi$ le foncteur oubli $G_1\tEns\to\Ens$.\\
Le foncteur $\psi f^*:G_2\tEns\to\Ens$ est fid\`ele, donc $f^*$ aussi.\\
Supposons $\overline{\overline{\Image f}}\neq G_2$. Soit $H$ un sous-groupe ouvert strict de $G_2$ contenant $\Image f$. Alors $G_2/H$ est un $G_2$-ensemble connexe donc $|\Aut(G_2/H)|\leqslant [G_2:H]$, alors que $f^*(G_2/H)$ est trivial, et donc $|\Aut(f^*(G_2/H))|=[G_2:H]!$.\\
Supposons $\overline{\Image f}=G_2$. Soit $\phi:f^*S\to f^*S'$. $\psi(\phi)$ est $G_1$-\'equivariante, donc $\Image f$-\'equivariante. Mais $\{g\in G_2|\ \forall x\in f^*(S),\phi(gx)=g\phi(x)\}$ est un sous-groupe ferm\'e de $G_2$. Donc $\psi(\phi)$ est bien $G_2$-\'equivariante, et d\'efinit donc un morphisme $S\to S'$.\\
Si $f^*$ est une \'equivalence, $\Aut f^*\to\Aut \psi f^*$ est un isomorphisme. Or ce n'est autre que $f$.
\findem

Soient $\mcal B$ une catégorie et $F:\mcal B\to\Ens$ un foncteur. Munissons
$\Aut F$ de la topologie pour laquelle les $\Stab^{\Aut F}_{F(S),s}$ forment une
base de voisinages de 1. Alors $F$ s'enrichit en un foncteur $\mcal B\to
\Aut F\tEns$. $(\mcal B,F)$ est une catégorie classifiante pointée si le
foncteur $\mcal B\to\Aut F\tEns$ est un isomorphisme. Si $(\mcal B,F)$ est
une catégorie classifiante pointée, et $F'$ est isomorphe à $F$, alors
$(\mcal B,F')$ est aussi une catégorie classifiante pointée.
On appellera \emph{catégorie classifiante}\index{Cat\'egorie classifiante} un couple $(\mathcal B,\mcal F)$ où $\mathcal
B$ est une catégorie, $\mcal F$ une classe d'isomorphisme de foncteurs $F:\mathcal
B\to \Ens$ tels que $(\mcal B,F)$ soit une catégorie classifiante pointée (un élément $F$ de $\mcal F$ sera alors appelé un foncteur
fibre de la catégorie classifiante). On notera aussi $\gf(\mcal B,F)=\Aut
F$, muni de sa structure topologique de groupe qpd, qu'on appellera
\emph{groupe fondamental}\index{Groupe fondamental!d'une cat\'egorie classifiante} de $\mcal B$.\\
Une catégorie classifiante admet des limites directes, des produits fibrés
et tout morphisme se décompose en la composé d'un épimorphisme strict et
de l'immersion d'une composante connexe.\\
Un \emph{1-morphisme} de catégories classifiantes $f:(\mathcal B_1,\mcal F_1)\to(\mathcal
B_2,\mcal F_2)$ est un foncteur $f^*:\mcal B_2\to \mcal B_1$ tel qu'il
existe un foncteur fibre $F_1\in\mcal F_1$ tel que
$F_1f^*$ soit dans $\mcal F_2$. Ceci vaut alors pour tout $F'_1\in \mcal
F_1)$~: en effet si $\alpha:F_1\to F'_1$ est un isomorphisme de foncteurs,
$\alpha\circ f^*:F_1 f^*\to F'_1 f^*$ est aussi un isomorphisme.\\
Un \emph{2-morphisme} de catégories classifiantes $\phi:f\to g$ avec $f,g:\mathcal
B_1\to \mathcal B_2$ est un isomorphisme de foncteur $f^*\to g^*$.\\
Avec ces d\'efinitions, les catégories classifiantes forment ainsi une 2-catégorie.\\

Si $G$ est un groupe topologique, $G\tEns$ est naturellement munie d'une
structure de catégorie classifiante, de groupe fondamental $\hat G$.\\
Si $\mathcal C$ est une catégorie galoisienne, $\Indob\mcal C$ est aussi naturellement
muni d'une structure de catégorie classifiante en prenant pour $\mcal F$ la
classe des foncteurs fondamentaux. Tout foncteur exact de catégories
galoisiennes définit un (unique) 1-morphisme de catégories classifiantes.\\

On a également, plus généralement, une notion de catégorie
multiclassifiante\index{Cat\'egorie multiclassifiante}. C'est une catégorie $\mcal B$ munie d'une famille de
classes d'isomorphismes $(\mcal F_i)_{i\in I} $ de foncteurs $F_i:\mcal B\to \Ens$,
tels que, si $F_i\in \mcal F_i$, le foncteur induit $\mcal C\to \prod_{i\in I}\Aut(F_i)\tEns$ soit
une équivalence de catégories (en particulier, les foncteurs $F_i$ sont deux à deux
non isomorphes, et donc les classes $\mcal F_i$ sont disjointes). On notera alors $\g0(\mcal B):=I$ et $\gf(\mcal
B,F_i):=\Aut F_i$. Le facteur direct de $\mcal B$ qui est la
sous-catégorie pleine de $\mcal B$ des $x$ tel que $F_j(x)=\emptyset$ pour tout
$j\neq i$ est noté $\mcal B_i$ (il ne dépend pas d'un choix de foncteurs
fibres).\\
Un 1-morphisme de catégories multiclassifiantes $F:(\mcal B,(F_i)_{i\in
  I})\to (\mcal B',(F'_j)_{j\in J})$ est un foncteur $F^*:\mcal B'\to \mcal 
B$ tel que, pour tout $i\in I$, $F_iF^*$ soit isomorphe à l'un des $F_j$ ($j$
est alors nécessairement unique, d'où une fonction $\g0(F):\g0(\mcal B')\to
\g0(\mcal B)$, compatible avec la composition de 1-morphismes).\\
Un 2-morphisme $u:F\to G$ de 1-morphismes est un isomorphisme de
foncteurs $F^*\to G^*$ (alors $F_i(u):F_iF^*\to F_iG^*$ est un
isomorphisme, donc $\g0(F)=\g0(G)$).\\

Soit $\mcal B$ une catégorie classifiante, et $S$ un objet de $\mcal B$. On
notera $\mcal B_{/S}$ la catégorie des morphismes $S'\to S$.\\
Si $S$ est connexe, $F$ est un foncteur fibre de $\mcal B$ et $s\in F(S)$,
on a un foncteur $F_s:\mcal B_{/S}\to \Ens$ défini par $F_s(\phi:S'\to
S)=F(\phi)^{-1}(\{s\})\subset F(S')$. Soit $H=\Stab_{S,s}\subset G=\gf(\mcal
B,F)$, l'action de $H$
sur $F(S')$ laisse stable $F_s(S'\to S)$ et définit ainsi un morphisme
$H\to \Aut F_s$ (d'où un foncteur $\mcal B_{/S}\to H\tEns$ noté également
par abus $F_s$). En identifiant $\mcal B$ à $G\tEns$ grâce à $F$, le
foncteur $\Ind_G^H:H\tEns\to (G\tEns)_{/(G/H)}$ fournit un quasi-inverse à
$F_s$. Ainsi $(\mcal B_{/S},F_s)$ définit une catégorie classifiante (et le
foncteur $\mcal B\to \mcal B_{/S}:S'\mapsto (S\times S'\to S)$ définit un
1-morphisme $\mcal B_{/S}\to\mathcal B$ de catégories classifiantes).\\
Si $s'$ est un autre élément de $F(S)$, comme $G$ agit transitivement sur
$F(S)$, soit $g\in G$ tel que $g\cdot s=s'$. Alors $g_{F(S')}$ envoie
$F_s(S'\to S)$ sur $F_{s'}(S'\to S)$ et définit ainsi un isomorphisme
$F_s\to F_{s'}$. De plus, si $\alpha : F\to F'$ est un isomorphisme de
foncteurs fibres, $\alpha$ induit un isomorphisme $F_s\to
F'_{\alpha(s)}$. Ainsi la classe d'isomorphisme $F_s$ ne dépend ni de $F$,
ni de $s$, donc la catégorie classifiante $\mcal B_{/S}$ est bien définie.\\
Si $\psi:\mcal B_1\to \mcal B_2$ est un 1-morphisme de catégories
classifiantes, $S_1$ un objet connexe de $\mcal B_1$, $S_2$ un objet
connexe de $\mcal B_2$ et $\phi:\psi^*(S_2)\to S_1$ un morphisme de $\mcal
B_2$, le foncteur $\mcal B_{1/S_1}\to\mcal B_{2/S_2}$ est un morphisme de
catégories classifiantes.\\
Pour $S$ quelconque, on peut décomposer $S$ comme somme directe de ses
composantes connexes $S=\coprod S_i$, et on obtient ainsi que $\mcal
B_{/S}=\prod \mcal B_{/S_i}$ est une catégorie multiclassifiante. 

\section{Complexes classifiants de groupes et groupes fondamentaux}
\subsection{Complexes classifiants}\label{cplxclassifiants}
Soit $\mathcal C$ une petite catégorie. Un \emph{complexe
classifiant}\index{Complexe classifiant} sur $\mathcal C$ est une catégorie fibrée $\mathcal G\to
\mathcal C$ dont les fibres $\mathcal G_x$ ($x\in \Ob(\mathcal C)$) sont
munies d'une structure de catégorie classifiante, telle que pour toute
flèche $f:x\to x'$ de $\mathcal C$, $f^*$ soit un morphisme de catégories
classifiantes.\\
Un 1-morphisme $\psi:\mcal G'\to \mcal G$ de complexes classifiants sur
$\mcal C$ est un 1-morphisme cartésien de catégories fibrées $\psi^*:\mcal
G\to \mcal G'$sur $\mcal C$ tel que pour tout $x\in \Ob(\mcal C)$ le foncteur
$\psi^*_x:\mcal G_x\to \mcal G'_x$ définisse un 1-morphisme $\psi_x:\mcal
G'_x\to \mcal G_x$ de catégories classifiantes.\\
Un 2-morphisme $u:\psi\to \psi'$ de 1-morphismes de complexes classifiants sur $\mcal C$ est
un 2-morphisme $\psi\to\psi'$ de catégories fibrées sur $\mcal C$.\\
Un 1-morphisme de complexes classifiants $(\mcal G'\to
\mcal C')\to(\mcal G\to\mcal C)$ est donné par un foncteur $\mcal
C'\to\mcal C$ et un 1-morphisme $\mcal G'\to \mcal
G\times_{\mcal C}\mcal C'$ de complexes classifiants sur $C'$. Un
2-morphisme de complexes classifiants est un 2-morphisme de catégories
fibrées.\\
On peut de même définir une notion de complexe multiclassifiant comme
catégorie fibrée dont les fibres sont des catégories multiclassifiantes.

\begin{rem} La terminologie de complexe classifiant est choisie pour insister sur le fait
que la définition ici présent\'ee est une généralisation de la notion de
complexe de groupe\index{Complexe de groupes} tel que l'a définie Haefliger dans~\cite{haef}. En effet,
la donnée d'un complexe de groupes au sens d'Haefliger équivaut à la donnée
d'un complexe classifiant $\mathcal G\to \mathcal C$ où $\mathcal
C$ est la catégorie des objets d'un ensemble 2-simplicial et où le groupe
fondamental $\gf(\mcal G_x)$ de chaque fibre est discret (la définition que
l'on donnera du groupe fondamental d'un complexe classifiant correspondera
dans le cas d'un tel complexe de groupes à la définition du groupe
fondamental d'Haefliger).\\
On peut trouver aussi une variation du concept d'Haefliger dans~\cite[§
  2.2]{stix} pour des complexes (où $\mathcal C$ est toujours la catégorie
des objets d'un ensemble 2-simplicial) dont les fibres sont en groupes
profinis et les morphismes continus. Le groupe fondamental que Stix définit alors correspondra à
notre $\ga(\mcal G)$.\\
On pourra aussi regarder~\cite{mochi} où Mochizuki étudie les complexes
en groupes profinis dans le cas où la base est un (semi)graphe. Il y définit dans
ce cas la notion de groupe fondamental tempéré.\end{rem}

La sous-2-catégorie pleine des complexes classifiants avec $\mcal C$
discrète est naturellement équivalente à la 2-catégorie des catégories multiclassifiantes.

La donnée d'un complexe classifiant sur $\mathcal C$ muni d'un clivage
équivaut,  compte tenu de~\cite[VI.8]{sga}, à :
\begin{itemize}
\item pour tout objet $x$ de $\mcal C$, une catégorie classifiante $\Pi_x$,
\item pour toute flèche $f:x\to x'$, un 1-morphisme de catégories classifiantes $f^*:\Pi_{x'}\to \Pi_{x}$ (qui est l'identité si $f$ est $Id_x$),
\item pour tous $f,g$ avec $f:x\to x'$ et $g:x'\to x''$, un 2-isomorphisme \[c_{f,g}: (gf)^* \stackrel{\sim}{\to} f^*\! g^*\] v\'erifiant
\[(f^*\circ c_{h,g}) \cdot c_{f,hg}=(c_{f,g}\circ h^*) \cdot c_{gf,h},\]
ce qui permet d'identifier $g^*\!f^*$ et $(fg)^*$.\\
\end{itemize}
Ces données, vues sur la catégorie de base $\mcal C^{\op}$, donnent aussi
naissance à une  catégorie cofibrée $\mcal G^{\circ}$ sur $\mcal
C^{\op}$. $\mcal G^{\circ}$ est alors un topos sur $\mcal C^{\op}$.\\

Si $\mcal G$ et $\mcal G'$ sont munis de clivages sur $\mcal C$, un
$1$-morphisme de catégories classifiante $\mcal G\to\mcal G'$ équivaut, compte tenu de~\cite[VI.12]{sga}, à :
\begin{itemize}
\item pour tout $x\in\Ob(\mcal C)$, un 1-morphisme de catégories
  classifiantes $\mcal G_x\to\mcal G'_x$ donné par un foncteur $\psi_x:\mcal
  G'_x\to\mcal G_x$,
\item pour tout $f\in\Hom_{\mcal C}(x,x')$, un isomorphisme de foncteurs
  $\phi_f:\psi_xf^*_{\mcal G'}\stackrel{\sim}{\to}f^*_{\mcal G}\psi_{x'}$, tels
  que\begin{itemize}
\item $\phi_{\id_x}=id_{\psi_x}$,
\item pour tous $x\stackrel{f}{\to}x'\stackrel{g}{\to}x''$, $(c^{\mcal
  G}_{f,g}\circ\psi_{x''})\cdot \phi_{gf}=(f^*_{\mcal G}\circ\phi_g) \cdot \phi_fg^*_{\mcal G'}\cdot (\psi_x\circ c^{\mcal
  G'}_{f,g})$\\\end{itemize}\end{itemize}
De même, et plus généralement, un 1-morphisme de
  catégories classifiantes $(\mcal G\to\mcal C)\to(\mcal G'\to\mcal C')$ est donné par :
\begin{itemize}
\item un foncteur $F:\mcal C\to\mcal C'$,
\item pour tout objet $x$ de $\mcal C$, un foncteur $\psi_x:\mcal
  G'_{F(x)}\to\mcal G_x$ correspondant à un 1-morphisme de catégories classifiantes, 
\item pour tout morphisme $f:x\to x'$ dans $\mcal C$, un isomorphisme de
  foncteur $\phi_f:\psi_xF(f)^*\stackrel{\sim}{\to} f^*\psi_{x'}$ vérifiant
  $\phi_{\id_x}=\id_{\psi_x}$ et la condition de compatibilité avec la composition.\\\end{itemize}

On note \[\mathcal{B(G)}:=\Cart_{\mathcal C}(\mathcal C,\mathcal G)=\projLim \mcal G/\mcal C\]
 la catégorie des sections cartésiennes de
$\mathcal G$, qu'on appellera aussi catégorie des revêtements de $\mcal G$.\\
Plus explicitement, étant donné un clivage de $\mcal G\to \mcal C$, un
objet $S$ est donné par :
\begin{itemize}
\item pour tout objet $x$ de $\mcal C$, un objet $S_x$ de $\mcal G_x$,
\item pour toute flèche $f:x\to x'$, un isomorphisme $\alpha_f:f^*S_{x'}\to S_x$,
tels que pour tous $g:x\to x'$, $f:x'\to x''$, $\alpha_g\cdot g^*\alpha_f=\alpha_{fg}$ (après avoir identifié $g^*\! f^*S_x$ et $(fg)^*S_x$ grâce à $c_{f,g}$).
\end{itemize}
Un morphisme $\phi:S\to S'$ est donné par un morphisme $\phi_x:S_x\to
S'_x$, pour tout $x$, de manière à  ce que, pour tout $f:x\to x'$,
$\alpha'_f\cdot f^*\phi_{x'}=\phi_x\cdot\alpha_f$. \\

$\mathcal{B(G)}$ est un topos localement connexe~(compte tenu de \cite[4.1.1,4.1.2]{leroy}). En, particulier, il admet des limites directes et des produits fibrés (il
suffit de les prendre composante par composante). Tout morphisme $S'\to S$ se
décompose de façon unique comme composé d'un épimorphisme strict $S'\to
S_0$ (\ie $S'\to S_0$ est le
coégalisateur de $S'\times_{S_0}S'\rightrightarrows S'$) et de
l'immersion d'un facteur direct $S_0\to S$ (\ie il existe un sous-objet
$S_1$ de $S$ tel que $S_0\coprod S_1\to S$ soit un isomorphisme). En particuliertout sous-objet est un facteur direct (comme $\mcal{B(G)}$ est
localement connexe, tout sous-objet de $S$ est donc somme direct de
composantes de $S$) et tout épimorphisme est strict (on parlera alors
simplement de quotient).\\
Comme exemple de limite directe qui servira par la suite, si $R$ est une relation d'équivalence sur $S$ (c'est-à-dire un sous-objet
de $S\times S$ qui vérifie les propriétés d'une relation d'équivalence), il
existe alors un quotient $S/R$ (c'est-à-dire un épimorphisme strict $S\to
S/R$ tel que $S\times_{S/R}S=R$), défini composante par composante par
$(S/R)_x=S_x/R_x$ et en recollant grâce à l'isomorphisme naturel
$f^*(S_{x'}/R_{x'})\simeq f^*(S_{x'})/f^*(R_{x'})$.\\

Si $S$ est un objet de $\mathcal{B(G)}$, on peut lui associer la catégorie
fibrée $\mathcal G_S$ sur $\mathcal C_S$, où $\mathcal C_S$ est la
catégorie cofibrée sur $\mcal C$ dont la fibre en $x$ est la catégorie discrète
$\pi_0(S_x)$, et $\mathcal G_S$ est la catégorie fibrée sur $\mathcal C_S$
dont la fibre en $(i,x)$ (où $S_i$ est une composante connexe de $S_x$) est
$\mathcal G_{x/S_i}$.\\
$\mathcal G_{x/S_i}$ étant muni de sa structure de catégorie classifiante,
$\mathcal G_S$ est un complexe classifiant sur $\mathcal C_S$, et on a un
morphisme de complexe classifiant $\mcal G_S\to \mcal G$.\\
$\mathcal{B(G}_S)$ est alors naturellement équivalente à
$\mathcal{B(G)}_{/S}$.\\
On peut reformuler cela en définissant sur $\mcal C$ le complexe
multiclassifiant $\mcal G/S$ dont la fibre en $x$ est la catégorie multiclassifiante
$\mcal G_{x/S_x}$. A un tel complexe multiclassifiant $\mcal G_0$ sur une catégorie
$\mcal C$, on peut alors associer un complexe classifiant $\mcal G_1$ sur une catégorie
$\mcal C'$ cofibrée sur $\mcal C$ dont la fibre en $x$ est $\g0(\mcal
G_{0,x})$, et si $i\in\g0(\mcal G_{0,x})$, $\mcal G_{1,(x,i)}$ est le facteur
direct $\mcal G_{0,x,i}$ de $\mcal G_{0,x}$.\\

On appellera \emph{point géométrique}
\index{Point!g\'eom\'etrique!d'un complexe classifiant} $\bar x=(x,F)$ de $\mathcal G$ la donnée d'un
objet $x$ de $\mcal C$ et d'un foncteur fibre $F$ de $\mathcal G_x$.\\
Si $\bar x=(x,F)$ est un point géométrique de $\mathcal G$, on a un
foncteur $F_{\bar x}:\mathcal{B(G)}\to Ens$ qui à $S\in\mathcal{B(G)}:\mathcal C\to
\mathcal G$ associe $F(S_x)$. Ce foncteur est appelé foncteur fibre en
$\bar x$ . Il commute aux limites directes et aux produits fibrés.\\
\begin{prop}\label{chemins} Si $\mathcal C$ est une catégorie connexe et $\bar x$ et
  $\bar{x'}$ sont deux points géométriques de $\mathcal G$, alors $F_{\bar
  x}$ et $F_{\bar{x'}}$ sont isomorphes.\end{prop}
\dem Posons $(x,F)=\bar x$ et $(x',F')=\bar x'$. Supposons qu'il existe un morphisme $f:x\to x'$ dans $\mathcal
  C$. Comme $f^*:\mathcal G_{x'}\to \mathcal G_x$ est un morphisme de
  catégories classifiantes, il existe un isomorphisme $\alpha:Ff^*\to F'$.\\
Or le morphisme naturel $\beta(S):S(x)\to f^*(S(x'))$ est un isomorphisme
  car $S$ est une section cartésienne.\\
$(F\circ\beta)(S):F(S(x))=F_{\bar x}(S)\to Ff^*(S(x'))$ est un
  isomorphisme, donc $\alpha\cdot(F\circ\beta)(S):F_{\bar x}(S)\to
  F'(S(x'))=F_{\bar{x'}}(S)$ est un isomorphisme, et donc
  $\alpha\cdot(F\circ\beta)$ fournit un isomorphisme $F_{\bar x}\to F_{\bar x'}$.\\
Dans le cas général, comme $\mathcal C$ est connexe, on peut trouver une
  chaîne d'objets $x, x_1,\dots,x_n,x'$ reliés par des morphismes (dans un
  sens ou dans l'autre), on obtient par ce qui précède
  que $F_{\overline{x_i}}$ et $F_{\overline{x_{i+1}}}$ sont isomorphes. On
  compose des isomorphismes pour obtenir un isomorphisme entre $F_{\bar x}$ et $F_{\bar x'}$.\findem

On appelle \emph{chemin} de $\bar x$ à $\bar{x'}$ un tel isomorphisme $F_{\bar
  x}\to F_{\bar{x'}}$.\\

Si $\psi:\mcal G'\to\mcal G$ est un 1-morphisme de complexes classifiants sur
$\mcal C$, $S\mapsto S\circ\psi^*$ définit un foncteur également noté
$\psi^*:\mcal{B(G)}\to \mcal{B(G}')$.\\
Si $\psi:(\mcal G'\to\mcal C')\to(\mcal G\to \mcal C)$ est un 1-morphisme
de complexes classifiants donné par un morphisme $\psi_0:\mcal G'\to \mcal
G\times_{\mcal C}\mcal C'$ de complexes classifiants sur $\mcal C'$,
$S\mapsto \psi_0^*(S\times_{\mcal C}\mcal C')$ définit un foncteur $\psi^*:\mcal{B(G)}\to\mcal{B(G}')$.\\
Si $u:\psi_1\to\psi_2$ est un
2-morphisme, il définit un isomorphisme de foncteurs
$u:\psi_1^*\to\psi_2^*$.\\
Si $\bar x'$ est un point géométrique $(x',F')$ de $\mcal G'$, soit $\psi(x)$ l'image
dans $\mcal C$ de $x\in\Ob(\mcal C')$ et soit $F$ la composée $\mcal
G_x\to\mcal G'_{x'}\stackrel{F}{\to} \Ens$. Alors $\psi_*(\bar x'):=\bar x=(x,F)$ est un
point géométrique de $\mcal G$, et $F_{\bar x}=F_{\bar x'}\circ \psi^*$.\\
Si $\alpha$ est un chemin de $\bar x'_1$ à $\bar x'_2$ (c'est-à-dire un
isomorphisme $\alpha:F_{\bar x'_1}\to F_{\bar x'_2}$), on a un chemin
$\psi_*(\alpha)$ de $\psi_*(\bar x'_1)$ à $\psi_*(\bar x'_2)$ défini par
$\psi_*(\alpha)(S)=\alpha(\psi^*(S))$ et $\psi_*$ est bien évidemment
compatible à la composition des chemins.\\
D'où en particulier un homomorphisme $\psi_*:\Aut F_{\bar x'}\to \Aut
F_{\psi_*(\bar x')}$.

\subsection{Groupe fondamental}
Soit $\mcal G$ un complexe classifiant sur $\mcal C$ et soit $\bar x=(x,F)$ un point géométrique de $\mathcal G$.\\
On appelle \emph{groupe fondamental}
\index{Groupe fondamental!d'un complexe classifiant} de $\mathcal G$ en $\bar x$ le groupe :
$$\gf(\mathcal G,\bar x)=\Aut F_{\bar x}.$$
C'est un groupe topologique, pour lequel les sous-groupes $\Stab_{S,s}$ ($S\in \mathcal{B(G)},
s\in F_{\bar x}(S)$) forment une base de voisinages de $1$ (les
$\Stab_{S,s}$ sont stables par intersection grâce à l'existence de produits
fibrés, et par conjugaison en remplaçant $s$ par $s'$).\\
L'existence de chemins montre que, à isomorphisme près, $\gf(\mathcal
G,\bar x)$ ne dépend pas de $\bar x$.\\

\begin{thm}\label{gfeqcplx}Le foncteur $F_{\bar x}:\mathcal{B(G)}\to \gf(\mathcal G,\bar x)\tEns$ est
  une équivalence de catégories.\end{thm}
\dem Si $S\in\mathcal{B(G)}$ est connexe, alors $\mathcal C_S$ est une
catégorie connexe. notons $f:\mathcal G_S\to\mathcal G$ le morphisme de
catégories fibrées en catégories classifiantes. Notons $S_x=\coprod S_i$ la décomposition en composantes
connexes de $S_x$. Si $s_1,s_2\in F_{\bar x}(S)$ (avec $s_1\in F_{\bar x}(S_1), s_2\in
F_{\bar x}(S_2)$), alors $(S_1,F_{s_1})$ et $(S_2,F_{s_2})$ sont des foncteurs fibres de 
$\mathcal G_S$. D'après la proposition~\ref{chemins}, il existe un chemin $\alpha$
de $(S_1,F_{s_1})$ vers $(S_2,F_{s_2})$. $f_*(\alpha)\in\gf(\mathcal G,\bar
x)$ est tel que $f_*(\alpha)s_1=s_2$. En résumé on a donc obtenu que si $S$
est connexe, $F_{\bar x}(S)$ est un $\gf(\mathcal G,\bar x)$-ensemble
connexe. Donc, pour tout $S$, $F_{\bar x}$ induit une bijection entre
l'ensemble des composantes connexes de $S$ et l'ensemble des composantes
connexes de $F_{\bar x}(S)$, ainsi qu'une bijection entre les sous-objets
de $S$ et les sous-objets de $F_{\bar x}(S)$ (car tout sous-objet est union
de composantes connexes).\\
On en déduit aisément que $F_{\bar x}:\mathcal{B(G)}\to \gf(\mathcal G,\bar
x)\tEns$ est pleinement fidèle, car
\begin{equation*} \Hom(S,S')=\left\{ \begin{array}{l} \text{sous-objets } U \text{ de } S\times S'\\ \text{tel que }U\to S \text{ soit un
  isomorphisme}\end{array}\right\},\end{equation*}
$F_{\bar x}(S\times S')=F_{\bar x}(S)\times F_{\bar x}(S')$ et $F_{\bar
  x}(U\to S)$ est un isomorphisme si et seulement si $U\to S$ est un
isomorphisme.\\
Comme les sommes directes existent dans $\mcal{B(G)}$, il ne reste plus qu'à
montrer qu'un $\gf(\mathcal G,\bar x)$-ensemble connexe $\gf(\mathcal G,\bar
x)/H$, où $H$ est un sous-groupe ouvert de $\gf(\mathcal G,\bar x)$, est
dans l'image essentielle de $F_{\bar x}$. Or puisque $H$ est ouvert, $H$
contient un sous-groupe de la forme $\Stab_{S',s'}$ pour un certain
$S'\in\Ob(\mcal B(\mcal G))$ et $s'\in F_{\bar x}(s')$ (comme
$\Stab_{S',s'}\subset H$, on a un unique morphisme $G/H\to F_{\bar x}(S')$
qui envoie $H$ en $s'$). Soit $R$ la r\'eunion
des composantes connexes de $S'\times S'$ correspondant à l'union des
composantes connexes $F_{\bar x}(R)=F_{\bar x}(S')\times_{G/H}F_{\bar
  x}(S')\subset F_{\bar x}(S')\times F_{\bar x}(S')$. $R$ est une relation
d'équivalence sur $S'$ et le quotient $S'/R$ existe dans $\mathcal B(\mcal
G)$. On a bien $F_{\bar x}(S'/R)=\gf(\mathcal G,\bar X)/H$.
\findem
Ainsi $(\mcal B(\mcal G),F_x)$ est une catégorie classifiante. De plus, si
$\bar x$ et $\bar x'$ sont deux points géométriques de $\mcal G$, $F_{\bar
  x}$ et $F_{\bar x'}$ sont isomorphes, donc la structure de catégorie
classifiante sur $\mcal B(\mcal G)$ est indépendante du point
géométrique.\\

\section{Sous-catégories d'une catégorie classifiante et applications}
Soient $\mcal B$ une catégorie classifiante et $\mcal B'$ une sous-catégorie
pleine qui contient l'objet final de $\mcal B$, stable par sous-objets, par produits fibrés, et
par quotients.\\

Donnons quelques exemples :
\begin{itemize}
\item la sous-catégorie pleine $\Bpd$ des objets $S$ localement constants
  du topos $\mcal B$. Si $\mcal B=G\tEns$, un objet $S$ est localement
  constant si et seulement si il existe un objet $S'$ tel que $S\times
  S'\to S'$ soit constant, \ie tel que $S\times S'$
  soit isomorphe à une somme disjointe de copies de $S'$, en tant qu'objets
  au-dessus de $S'$. $S'$ peut alors être choisi connexe, donc isomorphe à
  $G/H$ où $H$ est un sous-groupe ouvert de $G$. Un objet $S$ est constant
  au-dessus de $G/H$ si et seulement si, pour tout $s\in S$, $\Stab_s\subset
  H$ (ceci implique en particulier que, pour tout $s\in S$, $\Stab_s\subset
  \bigcap_{g\in G} gHg^{-1}$, et donc $H$ peut être choisi distingué).
\item la sous-catégorie pleine $\Ba$ des objets $S$ de $\mcal B$ tel que $F(S)$
  est fini (où $F$ est un foncteur fibre de $\mcal B$) ; on dira alors que
  $S$ est fini. C'est une catégorie galoisienne.\\
Dans le cas où $\mcal B$ est la catégorie $\mcal B(\mcal G)$ des
  revêtements d'un complexe classifiant, $\Ba$ peut aussi être vue comme la
  catégorie des sections cartésiennes de la catégorie fibrée $\mcal
  G^{\alg}$ dont la fibre en $x$ est la catégorie $\mcal G^{\alg}_x$ des
  objets finis de $\mcal G_x$.
\item Plus généralement, soit $\bb L$ un ensemble de nombres premiers, la
  catégorie $\mcal B^{\bb L}$ la sous-catégorie pleine de $\Ba$
  composée des objets finis dont le cardinal de la clôture galoisienne 
  soit un produit de nombre premiers dans $\bb L$. Un
  objet $\mcal B^{\bb L}$ sera dit $\bb L$-fini. $\mcal B^{\bb L}$ est
  également une catégorie galoisienne.\\
Dans le cas où $\mcal B$ est la catégorie $\mcal B(\mcal G)$ des
  revêtements d'un complexe classifiant, $\mcal B^{\bb L}$ peut aussi être vue comme la
  catégorie des sections cartésiennes de la catégorie fibrée $\mcal
  G^{\bb L}$ dont la fibre en $x$ est la catégorie $\mcal G^{\bb L}_x$ des
  objets $\bb L$-finis de $\mcal G_x$.
\item Soit $\mcal B(G)$ un complexe classifiant, la sous-catégorie pleine
  $\Bctop(\mcal G)$ de $\mcal B(\mcal G)$ dont les objets $S$ sont les
  revêtements catégoriquement topologiques de $\mcal G$, c'est-à-dire ceux tels
  que, pour tout $x\in \Ob(\mcal C)$, $S_x$ soit scindé (c'est-à-dire somme
  directe de copies de l'objet final de $\mcal G_x$).
\item La sous-catégorie pleine $\Bctemp(\mcal G)^{\bb L}$ 
  de $\mcal B(\mcal G)$ constituée des objets $S$ de $\mcal B(\mcal G)$ tel
  qu'il existe un objet $S'$ $\bb L$-fini non vide tel que $S'\times S\to S'$,
  considéré comme objet de $\mcal B(\mcal G_{S'})$ grâce à l'équivalence de
  catégorie $\mcal B(\mcal G)_{/S'}\simeq \mcal B(\mcal G_{S'})$, soit un
  revêtement catégoriquement topologique de $\mcal G_{S'}$ (si $\mbb L$ est l'ensemble de tous
  les nombres premiers, on omettra l'exposant $\bb L$ de la notation). Les
  objets de $\Bctemp(\mcal G)^{\bb L}$ sont appelés \emph{revêtements
  catégoriquement tempérés} pro-$\bb L$.\\
\end{itemize}
Si $\mcal B''$ désigne la sous-catégorie pleine de $\mcal B$ dont les
objets sont des sommes directes d'objets de $\mcal B'$, alors $\mcal B''$
est encore stable par sous-objets, par produits fibrés et par quotients, et
est de plus stable par somme directe.\\
Soit $F$ un foncteur fibre de $\mcal B$, soient $F'$ et $F''$ ses restrictions à
$\mcal B'$ et à $\mcal B''$.\\
Alors $\Aut F'=\Aut F''$ et on a un morphisme continu évident $\Aut F\to \Aut F'$.
\begin{prop}
$F'':\mcal B''\to (\Aut F'')\tEns$ est une équivalence de catégories.
\end{prop}
\dem Soit $S$ un objet connexe de $\mcal B''$. Il est connexe en tant
qu'objet de $\mcal B$, donc $\Aut F$ agit déjà transitivement sur
l'ensemble $F(S)=F''(S)$, donc
$\Aut F''$ aussi. $F''(S)$ est donc un $\Aut F''$-ensemble connexe. Donc,
pour tout $S$, $F''$ induit une bijection entre l'ensemble des composantes
connexes de $S$ et $F''(S)$.  On en
déduit que $F''$ est pleinement fidèle (par le même argument que dans la
démonstration du théorème \ref{gfeqcplx}).\\
Soit $T$ un $(\Aut F'')$-ensemble. Il faut montrer qu'il existe $S$ tel que
$F''(S)$ soit isomorphe à $T$. Comme $F''$ commute aux sommes directes et
que $\mcal B''$ admet des sommes directes, on peut supposer $T$ connexe :
$T=(\Aut F'')/H$ où $H$ est un sous-groupe ouvert de $\Aut F''$ et contient
donc un $\Stab_{S',s'}$ avec $S'\in\Ob(\mcal B'')$ et $s'\in F(S')$. On
peut définir une relation d'équivalence sur $S'$ comme le sous-objet $R$ de
$S'\times S'$ tel que $F''(R)=F''(S')\times_T F''(S')\subset F''(S'\times
S')$ (il existe bien un tel $R$ car $F''(S')\times_T F''(S')$ est un
sous-objet de $F''(S'\times
S')$ et $F''$ induit une bijection entre les sous-objets de $S'\times S'$
et ceux de $F''(S'\times S')$). Comme $\mcal B''$
est stable par quotients, $S=S'/R$ est un objet de $\mcal B''$ et
$F''(S)\simeq F''(S')/F''(R)\simeq T$.\findem
De plus, si $\alpha F_1\simeq F_2$ est un isomorphisme de foncteurs fibres de
$\mcal B$, la restriction de $\alpha$ à $\mcal B''$ définit un isomorphisme
$F''_1\simeq F''_2$.\\

Ainsi $\mcal B''$ est muni d'une structure naturelle de catégorie
classifiante.\\
Dans les exemples ci-dessus, le groupe $\Aut F'=\gf(\mcal B'',F'')$ est
noté  de la mani\`ere suivante (en supposant dans les cas issus d'un complexe classifiant, $F=F_{\bar
  x}$ pour $x$ un point géométrique $\bar x$ de $\mcal G$) :
\begin{itemize}
\item $\gpd(\mcal B,F)$ si $\mcal B'$ est la sous-catégorie pleines des
  objets localement constants de $\mcal B$. $\gpd(\mcal B,F)$ est alors le
  complété prodiscret de $\gf(\mcal B,F)$, c'est-à-dire $\gpd(\mcal
  B,F)=\varprojlim_H \gf(\mcal B,F)/H$ où $H$ décrit les sous-groupes
  ouverts distingués de $\gf(\mcal B,F)$.
\item $\ga(\mcal B,F)$ si $\mcal B'$ est la sous-catégorie pleine des objets finis
  de $\mcal B$. $\ga(\mcal B,F)$ est alors le complété profini de $\gf(
  \mcal B,F)$.
\item $\gf(\mcal B,F)^{\bb L}$ si $\mcal B'$ est la sous-catégorie pleine des objets
  $\bb L$-finis de $\mcal B$. $\gf(\mcal B,F)^{\bb L}$ est alors le complété
  pro-$\bb L$ de $\gf(\mcal B,F)$.
\item $\gctop(\mcal G,\bar x)$ si $\mcal B'$ est la sous-catégorie pleine des
  revêtements catégoriquement topologiques de $\mcal G$.
\item $\gctemp(\mcal G,\bar x)^{\bb L}$ si $\mcal B'$ est la sous-catégorie
  pleine des revêtements catégoriquement tempérés pro-$\bb L$ de $\mcal G$.\\
\end{itemize}
Si $\mcal C$ est une petite catégorie, on peut définir un complexe classifiant
sur $\mcal C$ (qui sera noté $\underline{\mcal C}$) en considérant la catégorie
fibrée $\Ens\times \mcal C\to \mcal C$ : les catégories fibres sont
canoniquement isomorphes à $\Ens$ qu'on munit de sa structure de catégorie
classifiante évidente. Tout revêtement de $\underline{\mcal C}$ est
catégoriquement topologique.\\
Si $\mcal G\to\mcal C$ est un complexe classifiant sur $\mcal C$, on a un
unique morphisme $\mcal G\to
\underline{\mcal C}$ de complexes classifiants sur $\mcal C$  défini sur une fibre par le foncteur $\Ens\to \mcal
G_x$ qui à un ensemble $S$ associe $\coprod_S e_x$ où $e_x$ est l'objet
final de $\mcal G_x$.\\
Le foncteur induit $\mcal B(\underline{\mcal C})=\Bctop(\underline{\mcal C})\to
\Bctop(\mcal G)$ est alors une équivalence de catégories (en fait
$\Bctop(\mcal G)$ a été définie comme l'image essentielle de ce foncteur). Le
groupe fondamental $\gctop(\mcal G,{\bar x})$ ne dépend en particulier que
de $\mcal C$ et de $x$, et est noté $\gctop(\mcal C,x)$ pour cette
raison.\\
Si $X$ est un espace topologique localement contractile, on notera
$\Cov(X)$ la catégorie de ces revêtements.\\ 
Soit $|N(\mcal C)|$ la réalisation géométrique du nerf de $\mcal C$, le
foncteur $\Cov(|N(\mcal C)|)\to \Bctop(\underline{\mcal C})$ qui à
$S\in\Ob(\Cov(|N(\mcal C)|))$ associe en $x\in\Ob(\mcal C)$ l'ensemble
$S_x$ et à $f:x\to x'$ l'isomorphisme $S_{x'}\to S_x$ obtenu par monodromie
le long de l'arête de $N(\mcal C)$ correspondant à $f$, est une équivalence
de catégories.\\
Ceci justifie la dénomination de revêtement (catégoriquement) topologique.\\
On en déduit également que $\gctop(\mcal C,x)$ est discret et l'existence
d'un revêtement catégoriquement topologique universel.\\
\begin{exs}
\begin{itemize}
\item Si $\mcal C$ n'a qu'un objet $x$, alors $\gctop(\mcal
C,x)=(\End(x))^{\gp}$.\\
En effet, un rev\^etement topologique de $\mcal C$, est la donn\'ee d'un ensemble $S$, et pour tout morphisme $f\in \End(x)$ d'une permutation $f^*$ de $S$ telle que $f^*g^*=g^*f^*$. En rempla\c cant $f^*$ par $f_*=(f^*)^{-1}$, ceci revient donc \`a se donner un ensemble $S$ avec un morphisme de mono\"ides $\End(x)\to\Aut(S)$, ce qui revient encore \`a ce donner un morphisme de groupe $\End(x)^{\gp}\to\Aut(S)$. Donc un rev\^etement de $\mcal C$ n'est autre qu'un $\End(x)^{\gp}$-ensemble (et la correspondance est clairement compatible aux morphismes), ce qui implique $\gctop(\mcal
C,x)=(\End(x))^{\gp}$
\item Plus généralement, si $\mcal C$ est une catégorie connexe quelconque,
  soit $\mcal C^{\pm}$ la catégorie $\mcal C[\Sigma^{-1}]$ des fractions de
  $\mcal C$ obtenu en inversant la classe $\Sigma$ de tous les morphismes
  de $C$ (c'est un groupoïde).\\
Alors, si $x\in\Ob(\mcal C)$, $\gctop(\mcal C,x)=\Aut_{\mcal C^{\pm}}(x)$
  (où $x$ est vu à droite comme objet de $\mcal C^{\pm}$).
En effet si $S$ est un rev\^etement topologique de $\mcal C$, pour tout $f\in\mcal C$, on peut d\'efinir $(f^{-1})^*:=(f^*)^{-1}$, ce qui permet de prolonger $S$ de fa\c con unique en un rev\^etement topologique de $\mcal C^{\pm}$. On en d\'eduit que $\Bctop(\mcal C^{\pm})\to\Bctop{\mcal C}$ est une \'equivalence. Comme $\mcal C^{\pm}$ est \'equivalente \`a une cat\'egorie ayant un seul \'el\'ement on se ram\`ene au cas pr\'ec\'edent.
\item Si $\mcal C$ est la catégorie $\Delta/C$ des simplexes d'un complexe
  simplicial $C$, $\gctop(\mcal C)$ est le groupe fondamental de la
  réalisation géométrique de $C$.
  En effet $N(\mcal C)$ est la subdivision
  barycentrique de $C$, d'o\`u un hom\'eomorphisme $|N(\mcal C)|\to |C|$. On conclut alors avec les r\'esultats de la discussion pr\'ec\'edant ces exemples.\end{itemize}\end{exs}

Pour les revêtements tempérés pro-$\bb L$, on a par définition que pour tout
revêtement tempéré pro-$\bb L$ connexe $S$, il existe un revêtement $\bb
L$-fini $S'$ non vide tel que $S'\times S\to S'$ soit un revêtement catégoriquement
topologique (on peut de plus supposer $S'$ connexe et même galoisien car
tout objet est dominé par un objet galoisien dans la catégorie classifiante
$\mcal B(\mcal G)^{\bb L}$ par finitude, puisque c'est une catégorie
galoisienne).\\ 
Soient $F$ un foncteur fibre de $\mcal B(\mcal G)$. $s\in F(S)$ et $s'\in F(S')$. Soit $(S'^{\infty},s'^{\infty})$ le
revêtement catégoriquement topologique universel de $(S',s')$. Alors comme $(S'\times S,(s',s))\to
(S',s')$ est un revêtement catégoriquement topologique, le revêtement
catégoriquement topologique universel de $S'$ se
factorise en $(S'^{\infty},s'^{\infty})\to (S'\times S,(s',s))\to
(S',s')$, d'où un morphisme $(S'^{\infty},s'^{\infty})\to (S,s)$ dans
$\mcal B(\mcal G)^{\bb L}$. Or le revêtement catégoriquement topologique universel d'un
  revêtement galoisien est galoisien, donc tout revêtement tempéré pro-$\bb
  L$ est dominé par un revêtement galoisien : $\gctemp(\mcal G)^{\bb L}$ est
  pro-discret.\\
De plus, si $\{(S_i,s_i)\}$ est une famille cofinale de revêtements $\bb
L$-finis galoisiens, et en notant $\{(S_i^{\infty},s_i^{\infty})\}$ la
famille de leurs revêtements universels, on en déduit un isomorphisme :
$$\gctemp(\mcal G,F)=\varprojlim \Gal(S_i^{\infty}).$$
Il est à remarquer que, en général, $\gctemp(\mcal G)^{\bb L}$ ne peut pas
être reconstruit à partir de $\gctemp(\mcal G)$. Ainsi, soit $\mbb L$ qui
ne contienne pas $2$. Si $\mcal C_1$ est
la catégorie triviale à un élément et $\mcal G_1=\mbf Z/2\mbf Z\tEns$,
$\gtemp(\mcal G_1)=\mbf Z/2\mbf Z$ et $\gtemp(\mcal
G_1)^{\mbb L}=\{1\}$,
alors que si $\mcal C_2$ est la catégorie à un élément ayant $\mbf
Z/2\mbf Z$ comme monoïde d'endomorphismes et $\mcal G_2$ est le complexe
classifiant trivial sur $\mcal C_2$, $\gtemp(\mcal G_2)=\mbf Z/2\mbf Z$ et
$\gtemp(\mcal G_2)^{\mbb L}=\mbf Z/2\mbf Z$.

\section{Complexes $l$-polysimpliciaux classifiants et groupe fondamental
  tempéré}
Nous renvoyons au paragraphe~\ref{berkspaces} pour la définition de la catégorie
polysimpliciale $\mbf\Lambda$.\\

Soit $l$ un entier naturel.
Soit $C$ un foncteur $(\mbf\Lambda^l)^{\op}\to\Ens$. On dira que $C$ est un
\emph{ensemble $l$-polysimplicial}\index{Ensemble $l$-polysimplicial}. Soit $\mbf\Lambda^{l\circ}\Ens$ la catégorie
des ensembles $l$-polysimpliciaux.\\
Le foncteur $\sq^{l-1}:\mbf\Lambda^l \to\mbf\Lambda$ qui envoie $(x_1,\cdots,x_l)$
sur $x_1\sq\cdots\sq x_l$ s'étend en un foncteur
$\sq^{l-1}_!:\mbf\Lambda^{l\circ}\Ens\to\mbf\Lambda^\circ\Ens$ qui commute
aux limites inductives. En le composant avec le foncteur réalisation
géométrique, on obtient un foncteur $|\ |:\mbf\Lambda^{l\circ}\Ens\to\Ke$
qui commute également aux limites inductives (on l'appellera aussi
réalisation géométrique).

Notons $\mbf \Lambda^l/C$ la cat\'egorie dont les objets sont les couples $(x,y)$ o\`u $x$ est un objet de $\mbf \Lambda^l$ et $y$ appartient \`a $C_x$.
Un morphisme $(x,y)\to (x',y')$ est un morphisme $f\in\Hom_{\mbf \Lambda^l}(x,x')$ tel que $f^*(y')=y$.\\
Un complexe classifiant $\mcal G$ sur $(\mbf\Lambda^l/C)^{\op}$ sera alors
appelé un complexe classifiant $l$-polysimplicial
\index{Complexe classifiant!$l$-polysimplicial} sur $\mcal C$.\\
Un 1-morphisme de complexes $l$-polysimpliciaux classifiants $(\mcal
G,C)\to(\mcal G',C')$ est donné par un morphisme $C\to C'$ d'ensembles
$l$-polysimpliciaux et un morphisme de complexes classifiants $\mcal
G\to\mcal G'\times_{(\mbf\Lambda^l/C')^{\op}}(\mbf\Lambda^l/C)^{\op}$.\\

Si $S$ est un objet de $\mcal B(\mcal G)$, on a un foncteur
$D_S:(\mbf\Lambda^l/C)^{\op}\to\Ens$ qui à $x$ associe $\g0(S_x)$. Cette donnée
permet de construire un ensemble $l$-polysimplicial $C_S=C\sq D_S$ définit par
$(C\sq D_S)_{\mbf n}=\coprod_{x\in C_{\mbf n}}D_S(x)$. On a alors $(\mbf\Lambda^l/(C\sq
D_S))^{\op}=((\mbf\Lambda^l/C)^{\op})_S$, ce qui permet de donner à $\mcal
G_S$ une structure naturelle de complexe $l$-polysimplicial. Cette
construction est fonctorielle en $S$.\\

On a un foncteur $F:\Cov(|C|)\to\Bctop(\mbf \Lambda^l/C)$, qui à un revêtement
topologique $\underline S$ de $|C|$ associe le revêtement topologique de 
$\mbf\Lambda^l/C$ défini de la façon suivante~:\begin{itemize}
\item pour $\alpha_x:x\to
C\in\Ob(\mbf\Lambda^l/C)$, $F(\underline S)_{\alpha_x}=\g0(|\alpha_x|^{*}\underline S)$ (il est à remarquer que comme $|x|$ est
contractile, $|\alpha_x|^*\underline S$ est un ensemble de copies de
$|x|$)~;
\item si $f:x'\to
x$ est un morphisme de $\mbf\Lambda^l/C$, $F(f)$ est la bijection $|f|^*:|\alpha_x|^*\underline S \to |\alpha_{x'}|^*\underline S$.\end{itemize}
Si $x$ est un objet de $\mbf\Lambda^l/C$, $F_xF$ est naturellement isomorphe
au foncteur fibre en $y$ pour tout point $y\in|C|$ dans l'image de $|x|$.\\
\begin{lem}
Le foncteur $F$ est pleinement fidèle.\end{lem}
\dem
On a un foncteur $G:\Bctop(\Lambda/C)\to \Top/|C|$, commutant aux limites
inductives, qui à $S$ associe
$|D_S|$. On a, fonctoriellement en $S$, \[\begin{array}{rcl} GF(S) & = & \Coker(\coprod_{y\in
  N_1(\mbf\Lambda^l/F(S))} |y|\rightrightarrows \coprod_{x\in
    N_0(\mbf\Lambda^l/F(S)} |x|)\\ & = & \Coker(S\times_{|C|}\coprod_{y\in
    N_{1}(\mbf\Lambda^l/C)} |y|\rightrightarrows S\times_{|C|}\coprod_{x\in
    N_{0}(\mbf\Lambda^l/C)}|x|)\\ & = & S.\end{array}\]
Donc $GF$ est isomorphe au plongement usuel $\Covtop(|C|)\to
\Top/|C|$, qui est pleinement fidèle.
Donc $F$ est bien pleinement fidèle.
\findem
Ainsi $\Covtop(|C|)$ s'identifie à une sous-catégorie pleine de $\mcal
B(\mcal G)$, stable par sous-objets, produit fibrés et quotients.\\
Un revêtement de $\mcal G$ qui est dans l'image effective de ce foncteur
sera dit \emph{topologique}\index{Rev\^etement!topologique!d'un complexe polysimplicial classifiant}. On notera aussi $\Btop(\mcal G)=\Btop(C)$ la sous-catégorie pleine
en question et $\gtop(\mcal G)=\gtop(C)=\gtop(|C|)$ son groupe fondamental
topologique.\\

Un objet $S$ de $\mcal B(\mcal G)$ est appelé revêtement tempéré pro-$\mbb L$ 
  s'il existe un revêtement $S'$ $\mbb L$-fini non vide tel que $S'\times S\to S'$,
  considéré comme objet de $\mcal B(\mcal G_{S'})$ grâce à l'équivalence de
  catégories $\mcal B(\mcal G)_{/S'}\simeq \mcal B(\mcal G_{S'})$, soit un
  revêtement topologique de $\mcal G_{S'}$ (si $\bb L$ est l'ensemble de tous
  les nombres premiers, on omettra l'exposant $\bb L$ de la notation). La
  sous-catégorie pleine de $\mcal B(\mcal G)$ constituée des revêtements
  tempérés pro-$\mbb L$ est notée $\Btemp(\mcal G)^{\bb L}$. Le groupe fondamental associé est noté $\gtemp(\mcal G)^{\bb
  L}$.\\

On a alors une description analogue à celle de $\gctemp(\mcal G)^{\bb L}$.\\ 
Pour les revêtements tempérés pro-$\bb L$, on a par définition que pour tout
revêtement tempéré pro-$\bb L$ connexe $S$, il existe un revêtement $\mbb
L$-fini $S'$ non vide tel que $S'\times S\to S'$ soit un revêtement
topologique (on peut de plus supposer $S'$ connexe et même galoisien car
tout objet est dominé par un objet galoisien dans la catégorie classifiante
$\mcal B(\mcal G)^{\bb L}$ par finitude, puisque c'est une catégorie
galoisienne).\\ 
Soient $F$ un foncteur fibre de $\mcal B(\mcal G)$, $s\in F(S)$ et $s'\in F(S')$. Soit $(S'^{\infty},s'^{\infty})$ le
revêtement topologique universel de $(S',s')$. Alors comme $(S'\times S,(s',s))\to
(S',s')$ est un revêtement topologique, le revêtement universel de $S'$ se
factorise en $(S'^{\infty},s'^{\infty})\to (S'\times S,(s',s))\to
(S',s')$, d'où un morphisme $(S'^{\infty},s'^{\infty})\to (S,s)$ dans
$\mcal B(\mcal G)^{\bb L}$. Or le revêtement topologique universel d'un
  revêtement galoisien est galoisien. En effet, soient $s'^{\infty}_1,
  s'^{\infty}_2\in F(S'^{\infty})$, soient $s'_1,s'_2$ leurs images dans
  $F(S')$. Il existe un automorphisme $\phi$ de $S'$ qui envoie $s'_1$ en
  $s'_2$. $C_\phi$ est un automorphisme de $C_{S'}$, qui envoie la strate
  correspondant à $s'_1$ sur la strate correspondant à $s'_2$. La
  fonctorialité du rev\^etement topologique universel fournit des
  automorphismes de $S'^{\infty}$ compatibles avec $C_{\phi}$. Le caractère
  galoisie du revêtement universel nous dit qu'il en existe un qui envoie
  $s'_1$ en $s'_2$. Ainsi $S'^{\infty}$ est bien galoisien. Donc tout revêtement tempéré pro-$\bb
  L$ est dominé par un revêtement galoisien : $\gtemp(\mcal G)^{\bb L}$ est
  pro-discret.\\
De plus, si $\{(S_i,s_i)\}$ est une famille cofinale de revêtements $\bb
L$-finis galoisiens, et en notant $\{(S_i^{\infty},s_i^{\infty})\}$ la
famille de leurs revêtements universels topologiques, on en déduit un isomorphisme :
$$\gtemp(\mcal G,F)=\varprojlim \Gal(S_i^{\infty}).$$\\

Dans le cas où le complexe $l$-polyclassifiant $C$ est en fait
intérieurement libre (ce sera le cas pour le complexe associé à une
fibration strictement polystable), le foncteur $\Covtop(|C|)\to\Bctop(\mbf\Lambda^l/C)$
admet un quasi-inverse.\\
En effet, soit $S$ un revêtement topologique de $\mbf\Lambda^l/C$,
soit $S_0$ le $l$-ensemble polysimplicial $S_0$ dont la composante en $x$
est $\coprod_{x\stackrel{\alpha_x}{\to}C}S_x$. On a un morphisme $S_0\to C$
et $|S_0|\to |C|$ est un revêtement topologique. En effet, l'image réciproque
d'une cellule $\mathring\Sigma^{\mbf n}$ (les cellules sont de cette forme
car $C$ est intérieurement libre) de $|S_0|$ correspondant à un
polysimplexe non dégénéré $x$ est $\coprod_{S_x}\mathring\Sigma^{\mbf n}$
grâce à~\cite[lem 3.4]{berk2}. La restriction de $|S_0|\to|C|$ à l'étoile
d'une cellule ouverte est alors un revêtement trivial.\\
Ainsi $\Btop(\mcal G)=\Bctop(\mcal G)$ et $\gtop(\mcal G)=\gctop(\mcal
G)$.\\
Comme le complexe $l$-polysimplicial d'un revêtement d'un complexe
$l$-polysimplicial classifiant intérieurement libre est encore
intérieurement libre, on a aussi :
\[\Btemp(\mcal G)^{\bb L}=\Bctemp(\mcal G)^{\bb L}\]
et
\[\gtemp(\mcal G)^{\bb L}=\gctemp(\mcal G)^{\bb L}.\]

\section{Changement de catégorie de base}
\subsection{Cas d'invariance par changement de base}
Soit $F : \mathcal C_1 \to \mathcal C_2$ un foncteur de petites cat\'egories, soit $\mathcal D_2$ une
cat\'egorie fibr\'ee sur $\mathcal C_2$ et soit $\mathcal D_1=\mathcal
D_2\times_{\mathcal C_2}\mathcal C_1$.\\
On a un foncteur $F^* : \Cart_{\mathcal C_2}(\mathcal C_2,\mathcal D_2) \to
\Cart_{\mathcal C_1}(\mathcal C_1,\mathcal D_1)$.\\ 
Soit $x$ un objet de $\mathcal C_2$ et soit $x\!\!\dar\!\!F$ la catégorie telle
que $\Ob(x\!\!\dar\!\!F)=\{(f,y), y\in \Ob(\mathcal
C_1), f\in \Hom_{\mathcal C_2}(x,F(y))\}$ et $\Hom_{x\!\dar\!
  F}((f_1,y_1),(f_2,y_2))=\{\psi\in\Hom_{\mcal
  C_2}(y_1,y_2)|f_2=F(\psi)f_1\}$.\\
Rappelons que pour $\mcal C$ une petite catégorie, on note $\mcal C^{\pm}$ le
groupoïde obtenu par localisation de la classe de tous les morphismes de
$\mcal C$. On dira qu'un groupoïde est trivial s'il est connexe non vide et
tout automorphisme est une identité.
\begin{prop}\label{propcart}\begin{enumerate}[(i)]
\item Supposons que pour tout $x\in \Ob(\mathcal C_2), x\!\!\dar\!\!F$ ne
  soit pas la catégorie vide.
  Alors $F^*$ est fid\`ele.
\item Si pour tout $x$, $x\!\!\dar\!\!F$ est connexe non vide (\ie $F$ est final), alors $F^*$
  est pleinement fid\`ele.
\item Si pour tout $x$, $(x\!\!\dar\!\!F)^{\pm}$  est un groupoïde 
  trivial (\ie $F$ est 1-final), alors $F^*$ est
  une équivalence de catégorie.
\end{enumerate}
\end{prop}  
\dem 
\begin{enumerate}[(i)]
\item
Supposons $\forall x \in \Ob(\mathcal C_2), x\!\!\dar\!\!F \neq \varnothing$.\\
Soient $\phi_1,\phi_2 : h_1 \rightrightarrows h_2$ des morphismes de
$\Cart(\mathcal C_1,\mathcal D_1)$ tels que $F^*(\phi_1)=F^*(\phi_2)$.Montrons que $\phi_1=\phi_2$.\\
Soit $x \in \Ob(\mathcal C_2)$.\\
Comme $x\!\!\dar\!\!F\neq \varnothing$, on peut choisir $y \in \Ob(\mathcal C_1)$ et
$f :x \to F(y)$.\\
Comme $F^*(\phi_1)=F^*(\phi_2)$, $\phi_1(F(y))=\phi_2(F(y))$.\\
Or, pour $i=1$ ou $2$, on a les diagrammes commutatif suivants :
\[\xymatrix{h_1(x) \ar[r]^{\phi_i(x)} \ar[d]^{h_1(f)} & h_2(x) \ar[d]^{h_2(f)}\\
 h_1(F(y)) \ar[r]^{\phi_i(F(y))} & h_2(F(y))}\]
Donc $h_2(f)\phi_1(x)=h_2(f)\phi_2(x)$.\\
Comme $h_2(f)$ est un morphisme cart\'esien, $\phi_1(x)=\phi_2(x)$, d'o\`u la
fid\'elit\'e.\\

\item
Pla\c cons-nous maintenant dans le cas o\`u $x\!\!\dar\!\!F$ est connexe et
non vide.\\
Soient $h_1,h_2\in \Ob(\Cart(\mathcal C_2,\mathcal D_2))$, et $\bar{\phi} :
F^*(h_1) \to F^*(h_2)$.\\
Soit $x\in \Ob(\mcal C_2)$. Soit $(f,y)\in\Ob(x\!\!\dar\!\!F)$.\\
Comme $h_2(f)$ est un morphisme cart\'esien, il existe un unique morphisme
$\phi_{(f,y)}(x) : h_1(x) \to h_2(x)$ tel que :
\[h_2(f)\phi_{(f,y)}(x)=\bar{\phi}(y)h_1(f).\]
Montrons que $\phi_{(f,y)}$ ne dépend pas de $(f,y)$.\\
Si maintenant, on a $f_1=F(g)f_2$ avec $g\in \Hom_{\mathcal C_1}(y_2,y_1)$,
consid\'erons
\[\xymatrix{h_1(x) \ar[r]^{\phi_{(f_2,y_2)}} \ar[d] & h_2(x) \ar[d]\\
 F^*(h_1)(y_2) \ar[r]^{\bar{\phi}(y_2)} \ar[d] &
F^*(h_2)(y_2) \ar[d]\\ F^*(h_1)(y_1) \ar[r]^{\bar{\phi}(y_1)} &
F^*(h_2)(y_1),}\]
o\`u le rectangle du haut commute par d\'efinition de $\phi_{(f_2,y_2)}(x)$
et o\`u le rectangle du bas commute car $\bar{\phi}$ est un morphisme de
foncteurs.\\
Par unicit\'e de $\phi_{(f_1,y_1)}(x)$ faisant commuter le grand rectangle,
on en d\'eduit que $\phi_{(f_2,y_2)}(x)=\phi_{(f_1,y_1)}(x)$, donc
  $\phi_{(f,y)}(x)$ ne d\'epend que de la composante connexe de $(f,y)$.\\
Comme, par hypoth\`ese, il n'y a qu'une seule composante connexe,
$\phi_{(f,y)}(x)$ ne dépend pas de $(f,y)$. On le note $\phi(x)$.\\
Soit $g : x' \to x$ un morphisme de $\mcal C_2$.  Soit $(f,y)\in\Ob(x\!\!\dar\!\!F)$. Alors
$(fg,y)\in\Ob(x'\!\!\dar\!\!F)$. On a $h_2(f)\phi(x)h_1(g)=\bar
\phi(y)h_1(f)h_1(g)=\bar
\phi(y)h_1(fg)=h_2(fg)\phi(x')=h_2(f)h_2(g)\phi(x')$. Comme $h_2(f)$ est
cartésien $\phi(x)h_1(g)=h_2(g)\phi(x')$ (cf.~\cite[prop. 6.11]{sga}). Donc $\phi$ d\'efinit bien un
morphisme $h_1\to h_2$ d'image $\bar{\phi}$. D'o\`u la pl\'enitude.
\item Choisissons un clivage de $\mcal C_2$ et de $\mcal C_1$.
Soit $\overline S\in\Ob(\Cart(\mcal C_1,\mcal D_1))$. Soit
  $x\in\Ob(\mcal C_2)$.\\
Si $(f,y)\in\Ob(x\!\!\dar\!\!F)$, on peut définir
  $S_{x,(f,y)}=f^*(\overline S_{y})$ ($\overline S_{y}$ est un objet
  de $\mcal D_{2,y}=\mcal D_{1,F(y)}$, et donc $f^*\overline S_{y}$ est
  bien un objet de $\mcal D_{1,x}$).\\
Soit $g\in\Hom_{x\!\dar\!F}((f_2,y_2),(f_1,y_1))$ (c'est-à-dire
  $g\in\Hom_{\mcal C_1}(y_2,y_1)$ tel que $f_1=F(g)f_2$). Alors $g$ définit
  un isomorphisme $g_*:S_{x,(f_2,y_2)}\to S_{x,(f_1,y_1)}$ par $f_2^*(\overline
  S_{y_2})\simeq f_2^*(g^*(\overline S_{y_1}))\simeq (F(g)f_2)^*(\overline
  S_{y_1})=S_{x,(f_1,y_1)}$.
On obtient ainsi un foncteur $x\!\!\dar\!\!F\to \mcal D_{1,x}$ dont l'image
de tout morphisme est un isomorphisme. Ce foncteur se prolonge donc en un
foncteur $\alpha_x:(x\!\!\dar\!\!F)^{\pm}\to \mcal D_{1,x}$. Comme nous supposons
ici que $(x\!\!\dar\!\!F)^{\pm}$ est un groupoïde trivial, cela montre que
l'on a des isomorphismes canoniques $S_{x,(f_1,y_1)}\simeq S_{x,(f_2,y_2)}$ pour
tous les couples d'objets de $x\!\!\dar\!\!F$. Fixons pour tout $x$ un
objet arbitraire $(f,y)$ de $x\!\!\dar\!\!F$ et définissons
$S_x=S_{x,(f,y)}$. On notera $\theta_{(f',y')}$ l'unique morphisme
$(f',y')\to (f,y)$ dans $(x\!\!\dar\!\!F)^{\pm}$.\\
Soient $\psi:x\to x'$ et
$g\in\Hom_{x'\!\dar\!F}((f'_2,y'_2),(f'_1,y'_1))$. $g$ définit également un
objet $\psi^*g$ de $\Hom_{x\!\dar\!F}((f'_2\psi,y'_2),(f'_1\psi,y'_1))$.
On a $\psi^*\circ g_*=(\psi^*g)_*$ (aux identifications
$(f'_i\psi)^*\simeq\psi^*f'^*_i$ près), et donc, par localisation, le diagramme
\[\xymatrix{(x'\!\dar\!F)^\pm \ar[r]^{\alpha_{x'}} \ar[d]^{\psi^*} & \mcal
  D_{1,x'} \ar[d]^{\psi^*}\\ (x\!\dar\!F)^{\pm}\ar[r]^{\alpha_x} & \mcal D_{1,x}
}\]
est 2-commutatif (où le 2-morphisme est
$\alpha_x\psi^*(f',y')=(f'\psi)^*S_{y'}\simeq \psi^*f'^*S_{y'}=\psi^*\alpha_{x'}(f',y')$).\\
Si $\psi:x\to x'$ est un morphisme dans $\mcal C_2$, $(f'\psi,y')$ est
un objet de $x\!\!\dar\!\!F$, d'où l'isomorphisme canonique
$\alpha_\psi:\psi^*S_{x'}=\psi^* f'^*\overline
S_{y}=S_{x,(f'\psi,y')}\simeq S_x$. Si, $\psi:x\to x'$ et $\psi':x'\to x''$
sont deux morphismes dans $\mcal C_2$, on a~:
\[\psi^*(\alpha_{\psi'})\alpha_{\psi}=\psi^*(\alpha_{x'}(\theta_{(f''\psi',y'')}))\alpha_x(\theta_{(f'\psi,y')})\]\[\simeq\alpha_x(\psi^*(\theta_{(f''\psi',y'')})\theta{(f'\psi,y')})=\alpha_x(\theta_{(f''\psi'\psi,y'')})=\alpha_{\psi'\psi}\]
o\`u l'isomorphisme est donné par l'identification $(f'\psi)\simeq \psi^*f'^*$
(la dernière égalité vient du fait qu'on a forcément
$\psi^*(\theta_{(f''\psi',y'')})\theta_{(f'\psi,y')}=\theta_{(f''\psi'\psi,y'')}$
par unicité des morphismes dans $(x\!\!\dar\!\!F)^{\pm}$).\\
On obtient ainsi un objet
$S\in\Ob(\Cart(\mcal C_2,\mcal D_2))$ d'image par $F^*$ isomorphe à $\overline S$.
\end{enumerate}
\findem 

\begin{exs}
\begin{itemize}
\item Si $\mcal C$ a un objet final $e$, considérant le foncteur
  $F:1\to\mcal C$ de la catégorie triviale $1$ dans $\mcal C$ d'image
  $e$. Soit $x$ un élément de $\mcal C$, alors $x\!\!\dar\!\!F$ est la
  catégorie discrète $\Hom(x,e)$, qui, comme $e$ est un objet final est
  bien connexe (donc triviale). On peut donc appliquer~\ref{propcart}.(iii)
  à $F$.
\item Si $\mcal C$ est une catégorie, on a un foncteur $F:(\Delta/N(\mcal
  C))^{\op}\to \mcal C$ (où $N(\mcal C)$ est le nerf de $\mcal C$, et
  $\Delta$ est la catégorie simpliciale) qui à $[x_1\stackrel{f_1}{\to}x_2\cdots
  x_{n-1}\stackrel{f_{n-1}}{\to}x_n]$ associe $x_1$. Soit $x\in\Ob(\mcal
  C)$, alors $x\!\!\dar\!\!F=(\Delta/N(_{x\backslash}\mcal C))^{\op}$, avec
  $_{x\backslash}\mcal C$ admettant un objet initial. On vérifie alors que si
  $\mcal C'$ admet un objet initial (ou un objet final), $(\Delta/N(\mcal
  C'))^{\pm}$ est un groupoïde trivial. On peut donc
  appliquer~\ref{propcart}.(iii) à $F$.
\end{itemize}
\end{exs}

\subsection{Descente de complexes classifiants}
Etant donné un foncteur cofibrant $\mcal C_2\to\mcal C_1$ (cf. \cite[\S{}
VI.10]{sga} pour la définition d'une catégorie cofibrée), on cherche ici
à construire, pour une catégorie fibrée $\mcal D_1/\mcal C_1$, une catégorie
fibrée sur $\mcal C_2$ dont les fibres correspondent à la catégorie des
sections cartésiennes de la restriction de $\mcal D_1$ à la fibre de $\mcal
C_1\to\mcal C_2$.\\
On trouvera des résultats plus détaillés dans~\cite[section 4]{giraud2}
(notamment~\cite[th. 4.7]{giraud2}.\\

Soit $G:\mcal C_2\to\mcal C_1$ un foncteur cofibrant (qu'on supposera muni
d'un coclivage normalisé)
et soit $\mcal D_2$ une catégorie fibrée sur $\mcal C_2$ (qu'on supposera
munie d'un clivage normalisé).\\
Soit $x\in\Ob(\mcal C_1)$, on a alors par restriction une catégorie $\mcal
D_{2,x}$ fibrée sur $\mcal C_{2,x}$ (où $\mcal C_{2,x}$
est la fibre de $x$ pour $G$). Appelons $\mcal D_{1,x}$ la catégorie
des sections cartésiennes $\projLim \mcal D_{2,x}/\mcal C_{2,x}$.\\

Si $f\in\Hom_{\mcal C_1}(x,x')$, on a un foncteur $f_*:\mcal
C_{2,x}\to\mcal C_{2,x'}$.\\
Soit $y$ un objet de $\mcal C_{2,x}$ (notons $\alpha_y$ le morphisme $y\to f_*(y)$).\\
On a $(\mcal D_{2,x})_y=\mcal D_{2,y}$ et
$(\mcal D_{2,x'}\times_{\mcal C_{2,x'}}\mcal C_{2,x})_y=\mcal D_{2,f_*(y)}$.\\
$\alpha^*$ définit donc un foncteur $F_y:(\mcal D_{2,x})_y\to(\mcal
D_{2,x'}\times_{\mcal C_{2,x'}} \mcal C_{2,x})_y$.\\
Si $\psi:y\to z$ est un morphisme de $\mcal C_{2,x}$, on a
$\alpha_z\psi=f_*(\psi)\alpha_y$, donc l'isomorphisme
$\psi^*\alpha_z^*\simeq (\alpha_z\psi)^*=(f_*(\psi)\alpha_y)^*\simeq
\alpha_y^*(f_*(\psi))^*$ induit un isomorphisme $F_y\psi_{\mcal
D_{2,x'}\times_{\mcal C_{2,x'}}\mcal C_{2,x}}^*\to \psi^*_{D_2,x}F_z$ qui vérifie la
condition de compatibilité à la composition.\\
Ceci définit donc un foncteur cartésien $\mcal D_{2,x}\to \mcal D_{2,x'}\times_{\mcal
  C_{2,x'}}\mcal C_{2,x}$ sur $\mcal C_{2,x}$, d'où un foncteur
\[f^*:\projLim \mcal D_{2,x'}/\mcal C_{2,x'}\to\projLim \mcal
D_{2,x'}\times_{\mcal C_{2,x'}}\mcal C_{2,x}/\mcal C_{2,x}\to \projLim \mcal
  D_{2,x}/\mcal C_{2,x}.\]
Les données $D_{1,x}:=\projLim \mcal D_{2,x}/\mcal C_{2,x}$ permettent
ainsi de construire sur $\mcal C_1$ une catégorie fibrée $\mcal D_1$.\\

On peu aussi définir $\mcal D_1$ beaucoup plus canoniquement :
\begin{itemize}
\item $\Ob(\mcal D_1)=\{(x,s)|x\in\Ob(\mcal C_1),s\in\projLim \mcal
  D_{2,x}/\mcal C_{2,x}\}$
\item
  $\Hom_{\mcal D_1}((x,s),(x',s'))$ est l'ensemble des
  $(f\in\Hom_{\mcal C_1}(x,x'),(\phi_{\alpha:y\to y'}\in\Hom_{\mcal D_2}(s(y),s'(y'))))$, où
  $\alpha:y\to y'$ décrit l'ensemble des morphismes cocartésiens de $\mcal C_2$ (en
  tant que $\mcal C_1$-catégorie) qui s'envoient dans $\mcal C_1$ sur $f$
  et qui vérifient :
\begin{itemize}
\item $\phi_\alpha$ s'envoie sur $\alpha$ par le foncteur $\mcal
  D_2\to\mcal C_2$,
\item pour tout $\psi : y\to z$ au-dessus de $\id_x$, pour tous $\alpha : y\to y',\beta : z\to z'$
  cocartésiens aux dessus de $f$,
  $\phi_{\beta}s(\psi)=s'(\psi')\phi_{\alpha}$ (où $\psi':y'\to z'$ est
  l'unique morphisme de $\mcal C_2$ au-dessus de $\id_x'$ tel que
  $\psi'\alpha=\beta\psi$ (l'unicité provient de la cocartésianité de
  $\alpha$).
\end{itemize}
\item Si $(f_1,\phi_{1,\alpha}):(x,s)\to(x',s')$ et si
$(f_2,\phi_{2,\beta}):(x',s')\to(x'',s'')$, la composée est définie par
$(f_2f_1,\phi_{\gamma})$, où, si $\gamma$ est un morphisme cocartésien de
$\mcal C_2$ au-dessus de $f_2f_1$, choisissons $\alpha$ cocartésien
au-dessus de $f_1$ de même source que $\gamma$ et soit $\beta$ l'unique
morphisme au-dessus de $f_2$ tel que $\beta\alpha=\gamma$ ($\beta$ est bien
cocartésien aussi d'après~\cite[cor. VI.6.13]{sga}), alors
$\phi_{\gamma}=\phi_{2,\beta}\phi_{1,\alpha}$.\\
Cela ne dépend pas du choix
de $\alpha$ car si $\alpha'$ en est un autre (de $\beta'$ correspondant),
alors il existe $u$, nécessairement cocartésien, au-dessus de $id_{x'}$ tel
que $\alpha'=u\alpha$, donc $(\beta'u)\alpha=\beta\alpha$, donc
$\beta'u=\beta$ par unicité de $\beta$, et donc $\phi'_\gamma=\phi_{2,\beta'}\phi_{1,u\alpha}=\phi_{2,\beta'}s'(u)\phi_{1,\alpha}=\phi_{2,\beta'u}\phi_{1,\alpha}=\phi_\gamma$.
\end{itemize}
Le foncteur $F:\mcal D_1\to \mcal C_1$ est alors défini par $F((x,s))=x$
et $F((f,(\phi_\alpha)))=f$.\\
Vérifions que $\mcal D_1\to\mcal C_1$ est bien fibré.\\
Soit $f:x\to x'$ et soit $(x',s')\in\Ob(\mcal D_1)$.\\
Pour tout $y\in\Ob(\mcal C_2)$ au-dessus de $x$, choisissons $\alpha_y:y\to y'$ cocartésien
au-dessus de $f$ et choisissons un morphisme cartésien dans $\mcal D_2$
$\phi_y:\eta\to s'(y')$ au-dessus de $\alpha_y$. Appelons $s(y):=\eta$\\
Si $\psi:y\to z$, par cartésianité de $\alpha_y$, il existe un unique
$\psi'$ au-dessus de $\id_{x'}$ tel que $\psi'\alpha_y=\alpha_z\psi$. Par
cartésianité de $\phi_z$, il existe un unique $s(\psi):s(y)\to s(z)$ au
dessus de $\psi$ tel que $\phi_zs(\psi)=s'(\psi')\phi_y)$.\\
Cela définit ainsi un objet $(x,s)$  de $\mcal D_1$ au-dessus de $x$.\\
Pour tout $\alpha'':y\to y''$ cocartésien dans $\mcal C_2$, il existe un
unique isomorphisme $u:y'\to y''$ au-dessus de $\id_{x'}$ tel que
$\alpha''=u\alpha_y$. Définissons $\phi_\alpha:=s'(u)\phi_y:s(y)\to
s'(y'')$.\\
Cela définit un morphisme $(x,s)\to (x',s')$ dans $\mcal D_1$ au-dessus de
$f$, qui a la propriété que $\phi_\alpha$ est cartésien (sur $\mcal D_2$)
pour tout $\alpha$.\\

Soit $(f,(\phi_\alpha)):(x,s)\to (x',s')$ un tel morphisme dans $\mcal D_1$, montrons qu'il
est cartésien (la composée $(f'',(\phi''_\gamma))$ de deux tels morphismes
$(f,(\phi_\alpha))$ et $(f',(\phi'_\beta))$ vérifie encore la même
propriété parce que $\phi''_{\alpha\beta}=\phi_\alpha\phi'_\beta$, donc sera aussi cartésienne).\\
Soit donc $(f,(\phi_{0\alpha})):(x,s_0)\to (x',s')$ un autre morphisme de
même but que $(f,(\phi_\alpha))$ et au-dessus du même morphisme $f$.\\
Soit $u:y_0\to y$ un morphisme cocartésien de $\mcal C_2$ au-dessus de $\id_x$
(c'est donc un isomorphisme). Choisissons $\alpha:y\to y'$ un morphisme cocartésien
d'origine $y$ au-dessus de $f$ et soit $\alpha'=\alpha u$. Alors, comme
$\phi_\alpha$ est cartésien, il existe un unique $\phi_{1u}:s_0(y_0)\to s(y)$ au-dessus de
$u$ tel que $\phi_\alpha\phi_{1u}=\phi_{0\alpha'}$ ($\phi_{1u}$ doit
nécessairement vérifier cette propriété pour que
$(f,(\phi_\alpha))(\id_x,(\phi_{1u}))=(f,(\phi_{0\alpha'}))$ ; l'unicité
d'un tel $\phi_{1u}$ montre déjà l'unicité recherchée.\\
Vérifions que $\phi_{1u}$ ne dépend pas du choix de $\alpha$ mais
uniquement de $u$. Soit donc $\beta:y\to y''$ un autre morphisme
cocartésien au-dessus de $f$ et appelons $\beta'=\beta u$. Alors, comme $\alpha$ est cartésien, il
existe un unique morphisme $w:y'\to y''$ au-dessus $\id_{x'}$ tel que
$\beta=w\alpha$. On a aussi $\beta'=w\alpha'$. Donc
$\phi_\beta=s'(w)\phi_\alpha$ et $\phi_{\beta'} =s'(w)\phi_{\alpha'}$, donc
on a également $\phi_\beta\phi_{1u}=\phi_{0\beta'}$, donc $\phi_{1u}$ ne
dépend pas du choix de $\alpha$.\\
Si $\psi_0:y_0\to z_0$ est un morphisme au-dessus de $\id_x$, $u:y_0\to y$,
$v:z_0\to z$ des isomorphismes au-dessus de $\id_x$, posons $\psi=v\psi
u^{-1}:y\to z$. On veut montrer que
$s(\psi)\phi_{1u}=\phi_{1v}s_0(\psi_0)$. Soient donc $\alpha:y\to y'$ et
$\beta:z:\to z'$ des morphismes cocartésiens au-dessus de $f$, soit
$\psi':y'\to z'$ l'unique morphisme au-dessus de $\id_{x'}$ tel que
$\psi'\alpha=\beta\psi$ et soient $\alpha'=\alpha u$ et $\beta'=\beta
v$. Alors $\phi_{0\beta'}s_0(\psi_0)=s'(\psi')\phi_{0\alpha'}$, et donc
$\phi_\beta\phi_{1v}s_0(\psi_0)=s'(\psi')\phi_\alpha\phi_{1u}$. Or
$\phi_\beta s(\psi)=s'(\psi')\phi_\alpha$, donc
$\phi_\beta(\phi_{1v}s_0(\psi_0))=\phi_\beta(s(\psi)\phi_{1u}$. Comme
$\phi_\beta$ est cartésien, on a donc aussi
$\phi_{1v}s_0(\psi_0)=s(\psi)\phi_{1u}$.\\
Donc $(\id_x,(\phi_{1u}))$ définit bien un morphisme $(x,s_0)\to (x,s)$ qui
vérifie \[(f,(\phi_\alpha))(\id_x,(\phi_{1u}))=(f,(\phi_{0\alpha'})),\] donc
$(f,(\phi_\alpha))$ est bien cartésien, donc $\mcal D_1\to\mcal C_1$ est
bien fibrée.\\
De plus la fibre $\mcal D_{1,x}$ est bien naturellement équivalente à
$\projLim D_{2,x}/C_{2,x}$ (car si $(\id_x,\phi):(x,s)\to (x,s')$ est un
morphisme, il suffit de se donner les $\phi_{\id_y}$ pour tout $y$
au-dessus de $x$, car on peut reconstruire, pour $u:y\to y'$,
$\phi_u=\phi_{\id_{y'}}s(u) (=s'(u)\phi_{\id_{y}})$).

Quand on veut garder en tête l'importance du foncteur $F:\mcal C_2\to \mcal
C_1$, on notera aussi $F_*(\mcal D_2/\mcal C_2)$ la catégorie fibrée $\mcal
D_1/\mcal C_1$.

\begin{prop}[{\cite[th. 4.7]{giraud2}}]\begin{enumerate}[(i)]
\label{propdesc} 
\item On a une équivalence de catégories canonique 
\[\projLim \mcal D_2/\mcal C_2\to \projLim \mcal D_1/\mcal C_1\]
\item 
Plus généralement, si on a $\mcal C_2\stackrel{G}{\to} \mcal C_1 \stackrel{F}{\to} \mcal C_0$ des
foncteurs fibrants et si $\mcal D_2$ est une catégorie fibrée sur $\mcal
C_2$, on a une équivalence naturelle de catégories fibrées
$c_{F,G}:(FG)_*(\mcal D_2/\mcal C_2)\to F_*G_*(\mcal D_2/\mcal C_2)$.
\end{enumerate}
\end{prop}

\dem
\begin{enumerate}[$(i)$]
\item
A une section cartésienne $s$ de $\mcal D_2/\mcal C_2$, on associe à tout
$x$ sa restriction $(x,s_x)$ à $\mcal C_{2,x}$ (c'est un objet de $\mcal C_1$
au-dessus de $x$), et à tout $f:x\to x'$, on associe $(f,(\phi_{\alpha:y\to
  y'}=s(\alpha):s_x(y)\to s_{x'}(y')))$. Ceci définit une section $F(s)$ de $\mcal
D_1/\mcal C_1$.\\
Si $\psi: s\to s'$ est un morphisme dans la catégorie $\projLim \mcal
D_2/\mcal C_2$, on lui associe un morphisme $F(\psi):F(s)\to F(s')$ donné pour
tout $x$ par $F(\psi)_x=(\id_x,\phi_{u:y\to y'}=s'(u)\psi_y
(=\psi_{y'}s(u))$ (par l'identification précédente de $\mcal D_{1,x}$ avec
$\projLim \mcal D_{2,x}/\mcal C_{2,x}$, on peut voir $F(\psi)_x$ dans $\projLim
\mcal D_{2,x}/\mcal C_{2,x}$ comme la restriction de $\psi$ à $\mcal
D_{2,x}/\mcal C_{2,x}$).\\
On obtient ainsi un foncteur $F$ de $\projLim \mcal D_2/\mcal C_2\to
\projLim \mcal D_1/\mcal C_1$.\\

La donnée d'une section cartésienne équivaut à la donnée, pour tout $x$,
d'une section cartésienne $s_{x}$ de $\mcal D_{2,x}$, et pour tout $f$ et tout $\alpha$ cocartésien au-dessus de $f$ de $\phi_\alpha^f$ 
qui vérifient $\phi_\alpha^f\phi_\beta^{f'}=\phi_{\alpha\beta}^{ff'}$ pour
tous $\alpha$ et $\beta$ composables, et si $\beta \psi=\psi' \alpha$ au
dessus de $f \id_x=\id_{x'} f$ avec $\alpha$ et $\beta$ cocartésiens, alors
$\phi_\beta^f s_x(\psi)=s_{x'}(\psi')\phi_\alpha^f$.\\
On définit alors pour tout objet $y$ au-dessus de $x$~: $s(y):=s_x(y)$. Si
$\phi:y\to y'$ est au-dessus de $f:x\to x'$, alors, pour tout $\alpha$
cocartésien au-dessus de $f$ de source $y$, on peut factoriser de façon unique
$\phi=\psi\alpha$, et on pose $s(\phi)=s_{x'}(\psi)\phi_\alpha^f$ (si $\alpha'$ est
un autre morphisme cocartésien au-dessus de $f$, $\alpha'=u\alpha$, la
décomposition de $\phi$ est alors $\phi=(\psi u)\alpha'$ et $s_{x'}(\psi
u)\phi_{\alpha'}^f=s_{x'}(\psi)s_{x'}(u)\phi_{\alpha'}=s_{x'}(\psi)\phi_{\alpha}$
donc la construction de $s(f)$ ne dépend pas du choix de $\alpha$).\\
Si $\phi:y\to y'$ et $\phi':y'\to y''$, choisissons $\psi, \alpha,\psi'$ et
$\alpha'$ comme ci-dessus. Alors
$s(\phi')s(\phi)=s_{x''}(\psi')\phi^{f'}_{\alpha'}s_{x'}(\psi)\phi^f_\alpha$.
Soit $\beta$ cocartésien de source l'image de $\alpha$ au-dessus de $f'$,
alors il existe un unique $\psi_0$ tel que $\psi_0\beta=\alpha'\psi$, donc
$s_{x''}(\psi_0)\phi^{f'}_\beta=\phi^{f'}_{\alpha'}s_{x'}(\psi)$. Donc
$s(\phi')s(\phi)=s_{x''}(\psi')s_{x''}(\psi_0)\phi^{f'}_\beta\phi^f_\alpha
= s_{x''}(\psi'\psi_0)\phi^{f'f}_{\beta\alpha}=s(\phi\phi')$.\\

Un morphisme entre sections cartésiennes de $\mcal D_{2,x}$
$((s_x)_x,(\phi_\alpha)_{\alpha})\to ((s'_x)_x,(\phi'_{\alpha})_\alpha))$
est la donnée pour tout $x$ d'un morphisme $(\id_x,\psi_u):s_x\to s'_x$
dans $\mcal D_{1,x}$ tel que $\phi'_{\alpha}\psi_u=\psi_u\phi_{\alpha}$.\\
On définit alors pour tout $y$ $\psi_y:s(y)\to s'(y)$ par
$\psi_{\id_y}$. Pour que cela définisse bien un morphisme $s\to s'$, il
suffit de vérifier que pour tout $\phi:y\to y'$ au-dessus de $f:x\to x'$,
$s'(\phi)\psi_y=\psi_{y'}s(\phi)$, c'est-à-dire, en considérant une
décomposition $\phi=\psi\alpha$, que
$s'_{x'}(\psi)\phi'_\alpha\psi_{\id_y}=\psi_{\id_{y'}}s_{x'}(\psi)\phi_{\alpha}$.\\
Or
$s'_{x'}(\psi)\phi'_\alpha\psi_{\id_y}=s'_{x'}(\psi)\psi_{\id_{y}}\phi_{\alpha}$
car $(\id_x,\psi_u)_x$ est bien un morphisme de sections
$((s_x)_x,(\phi_\alpha)_{\alpha})\to ((s'_x)_x,(\phi'_{\alpha})_\alpha))$,
et
$s'_{x'}(\psi)\psi_{\id_{y}}\phi_{\alpha}=\psi_{\id_{y'}}s_{x'}(\psi)\phi_\alpha$
  car $(\id_x,\psi_u)$ est bien un morphisme $s_x\to s'_x$.\\
On obtient ainsi un foncteur $\projLim \mcal D_1/\mcal C_1\to
\projLim \mcal D_2/\mcal C_2$.\\

Les deux foncteurs ainsi obtenus sont inverses l'un de l'autre, d'où $(i)$.\\

 \item La preuve est similaire à celle de $(i)$ ($(i)$ nous donne déjà les
   équivalences fibres à fibres).
\end{enumerate}
\findem

On en déduit
que si on a un morphisme cocart\'esien de cat\'egories cofibr\'ees
\[\begin{array}{ccc}\mcal C_{2} & \to & \mcal C'_{2}\\ \dar & & \dar \\
\mcal C_{1} & \to & \mcal C'_{1}\end{array}\]
et qu'on a un morphisme cart\'esien de catégories fibrées $\mcal D_{2}/\mcal
C_{2}\to \mcal D'_{2}/\mcal C'_{2}$,
alors on a un morphisme naturel de catégories fibrée $\mcal D_{1}/\mcal
C_{1}\to \mcal D'_{1}/\mcal C'_{1}$ (la descente d'une catégorie fibrée
le long d'un foncteur $F$ est 2-fonctorielle en la catégorie fibrée).\\

$c_{F,G}$ est naturel au sens où si on a $\mcal
C_3\stackrel{H}{\to} \mcal C_2\stackrel{G}{\to} \mcal C_1 \stackrel{F}{\to}
\mcal C_0$ et $\mcal D_3/\mcal C_3$ une catégorie fibrée, on a
$(F\circ c_{G,H})\cdot c_{F,GH}=(c_{F,G}\circ H)\cdot c_{FG,H}$.\\

Si $\mcal D_2$ est un complexe classifiant sur $\mcal C_2$ (ou plus
généralement multiclassifiant), alors $\mcal D_1$ hérite aussi d'une
structure de complexe multiclassifiant (classifiant si les fibres de $\mcal
C_2\to \mcal C_1$ sont connexes). On a alors une équivalence $\mcal B(\mcal
C_2)\to\mcal B(\mcal C_1)$.\\

\begin{ex}
Soit $\mcal C$ une catégorie et $X$ un préfaisceau sur $\mcal C$. On a
alors un foncteur cofibrant $(X/\mcal C)^{\op}\to \mcal C^{\op}$ (nous
l'utiliserons principalement dans les cas où $\mcal C$ est la catégorie simpliciale, la catégorie
simpliciale stricte, la catégorie polysimpliciale, la catégorie bisimpliciale). Si donc
$\mcal D_2$ est un complexe classifiant sur $(X/\mcal C)^{\op}$, on peut se
ramener à l'étude du complexe multiclassifiant $\mcal D_1$ sur $\mcal C$ obtenu
en descendant $\mcal D_2$ (les fibres sont $(\mcal D_1)_x=\prod_{y\in X(x)}(\mcal
D_2)_y$).\end{ex}

Supposons que l'on ait un morphisme cocart\'esien de cat\'egories cofibr\'ees:
\[\xymatrix{\mcal C_2 \ar[d]^G \ar[r]^{F'} & \mcal C'_2 \ar[d]^{G'}\\ \mcal C_1 \ar[r]^F
  &\mcal C'_1.}\]
\begin{prop}Supposons de plus que $\mcal C_1\to \mcal C'_1$ vérifie la conclusion
de~\ref{propcart}.(iii) et que, pour tout objet $x$ de $\mcal C_1$, $\mcal
C_{2,x}\to\mcal C'_{2,F(x)}$ vérifie aussi la conlusion
de~\ref{propcart}.(iii), alors $F':\mcal C_2\to\mcal C'_2$ vérifie aussi la
conclusion de~\ref{propcart}.(iii).\end{prop}
\dem Soit $\mcal D'_2$ une catégorie fibrée sur $\mcal C'_2$ et soit
$\mcal D_2=\mcal D'_2\times_{\mcal C'_2}\mcal C_2$. On a, par fonctorialité de la descente, un
morphisme de catégories fibrées $G_*(\mcal D_2)/\mcal C_1\to G'_*(\mcal
D'_2)/\mcal C'_1$, d'où un morphisme $G_*(\mcal
D_2)\to G'_*(\mcal D'_2)\times_{\mcal C'_1}\mcal C_1$ de catégories fibrées
sur $\mcal C_1$. Mais si $x$ est un objet de $\mcal C_1$, la fibre en $x$
de ce morphisme est $\projLim \mcal D'_{2,F(x)}\times_{\mcal
  C'_{2,F(x)}}\mcal C_{2,x}/\mcal C_{2,x}\to \projLim \mcal
  D'_{2,F(x)}/\mcal C'_{2,F(x)}$, qui est une équivalence de catégories par
  hypothèse.\\
Donc $G_*(\mcal D_2)\to G'_*(\mcal D'_2)\times_{\mcal C'_1}\mcal C_1$ est une
équivalence de catégories, donc en appliquant l'hypothèse faite à  $\mcal
C_1\to \mcal C'_1$,
\[\projLim G_*(\mcal D_2)/\mcal C_1\to \projLim G'_*(\mcal D'_2)/\mcal
C'_1\]
est une équivalence de catégories.\\
Or on a le diagramme commutatif suivant :
\[\xymatrix{\projLim \mcal D_2/\mcal C_2 \ar[d] \ar[r] & \projLim \mcal
D'_2/\mcal C'_2 \ar[d]\\ \projLim G_*(\mcal D_2)/\mcal C_1\ar[r] & \projLim G'_*(\mcal D'_2)/\mcal
C'_1}\]
dont les flèches verticales sont des équivalences d'après la proposition~\ref{propdesc},
donc le foncteur du haut est aussi une équivalence.\findem

 \chapter{Groupe fondamental tempéré d'un log schéma pluristable}
\label{chap4}
Dans ce chapitre, nous d\'efinissons et \'etudions le groupe fondamental temp\'er\'e (et des variantes $\mbb L$) d'une fibration polystable $X\to\cdots\to k$ munie d'une log structure fs compatible, par analogie avec la construction sur les espaces de Berkovich. Les rev\^etements k\'et de $X$ joueront le r\^ole des rev\^etements finis. Berkovich a d\'ej\`a associ\'e  \`a une telle fibration polystable un espace polysimplicial (proposition \ref{berk69}) et sa r\'ealisation g\'eom\'etrique jouera le r\^ole de structure topologique de $X$. Cependant pour d\'efinir le groupe fondamental temp\'er\'e, il faut une structure topologique sur tous les rev\^etments finis. Nous serons donc amen\'es dans un premier temps \`a associer \`a un rev\^etement k\'et de $X$ un ensemble polysimplicial (\S \ref{logpolysimplicial}), dont la r\'ealisation g\'eom\'etrique jouera le r\^ole de structure topologique de $X$, ce qui nous permettra de d\'efinir le groupe fondamental temp\'er\'e de la fibration polystable.\\

Soit $K$ un corps complet pour une valuation discr\`ete. Soit $p$ sa
caract\'eristique r\'esiduelle. Pour une fibration polystable propre
$X\to\cdots\to O_K$  dont la fibre g\'en\'erique est lisse, nous
comparerons les variantes $(p')$ du groupe fondamental temp\'er\'e de la
fibre sp\'eciale de $X$ muni de sa structure logarithmique et le groupe
fondamental temp\'er\'e de l'espace analytique associ\'e \`a la fibre
g\'en\'erique de $X$ (\ref{isomfondtemp}). On dispose d\'ej\`a d'une \'equivalence entre les deux cat\'egories de rev\^etements $(p')$-finis. Il faudra donc principalement comparer la topologie de la fibre g\'en\'erique \`a celle de l'ensemble polysimplicial de la fibre sp\'eciale pour tous les rev\^etements $(p')$-finis de $X$. Pour ce faire, nous \'etendrons le th\'eor\`eme~\ref{berk81} aux rev\^etements k\'et de $X$.\\

Enfin nous donnerons une description plus combinatoire du groupe fondamental temp\'er\'e d'une fibration polystable logarithmique en terme d'un complexe polysimplicial classifiant. L'ensemble polysimplicial sous-jacent est celui de la proposition \ref{berk69}. Dans le cas d'une fibration strictement polystable, le groupe associ\'e \`a une composante de cet ensemble polysimplical est le groupe fondamental logarithmique de la strate correspondante. Dans le cas g\'en\'eral, le complexe classifiant sera construit par descente.

\section{Groupe fondamental tempéré d'un log-schéma pluristable}
Dans ce paragraphe nous définissons le groupe fondamental tempéré d'une
fibration polystable sur un corps $k$, muni d'une log structure compatible (ce
que nous appellerons une log fibration polystable). Pour définir notre
groupe fondamental tempéré, nous aurons besoin d'une notion de ``revêtement
topologique'' d'un revêtement két $Z$ de notre log fibration polystable
$X\to\cdots \to k$. Pour cela nous définirons pour tout tel $Z$ un ensemble polysimplicial
$\C(Z)$ au-dessus de $\C(X)$,
fonctoriellement en $Z$. Si $Z$ est un revêtement két galoisien 
connexe de $X$ de groupe $G$, l'action de $G$ sur $\C(Z)$ 
définira une extension de groupes~:
\[1\to\gtop(|\C(Z)|)\to \Pi_Z\to G\to 1.\]
Notre groupe fondamental tempéré sera la limite projective des $\Pi_Z$
quand $Z$ parcourt les revêtements két galoisiens pointés de $X$.

\subsection{Log schémas polystables}
Coomençons par définir un analogue logarithmique des notions de morphismes
polystables de schémas formels introduites dans la définition
\ref{defpolystable}.\\ 

Soit $S$ un log schéma fs.
\begin{dfn} Un morphisme $\phi:Y\to X$ de log schémas fs sera appelé~:
\begin{itemize}
\item \emph{nodal standard} si $X$
a une carte fs $X\to\Spec P$ et $Y$ est isomorphe à $X\times_{\Spec \mbf
  Z[P]}\mbf Z[Q]$ avec $Q=(P\oplus u\mbf N\oplus v\mbf N)/(u\cdot v=a)$ où $a\in
P$, et la log structure de $Y$ est celle associée à $Q$.
\item un \emph{morphisme strictement plurinodal de log schémas} si pour
  tout point $y\in Y$, il existe un voisinage ouvert de Zariski $X'$ de
  $\phi(y)$ et un voisinage ouvert de Zariski $Y'$ de $y$ dans
$Y\times_XX'$ tels que $Y'\to X'$ soit la composée de morphismes étales
stricts et de morphismes plurinodaux standard.
\item un \emph{morphisme plurinodal de log schémas} si, localement pour la topologie étale
  de $X$ et de $Y$, il est strictement plurinodal.
\item
 un \emph{morphisme strictement polystable de log schémas} si pour tout
 point $y\in Y$, il existe un voisinage ouvert de Zariski $X'$ de
 $\phi(y)$, une carte fs $P\to O(X')$ de la log structure de $X'$, un
 voisinage ouvert de Zariski $Y'$ de $y$ dans
$Y\times_XX'$ tel que $Y'\to X'$ se factorise à travers un morphisme strict étale
$Y'\to X'\times_{\mbf Z[P]}\mbf Z[Q]$ où
\[Q=(P\oplus\bigoplus_{i=0}^{p}<T_{i0},\cdots,T_{in_i}>)/(T_{i0}\cdots T_{in_i}=a_i)\]
avec $a_i\in P$.
\item un \emph{morphisme polystable de log schémas} si, localement pour la
  topologie étale de $Y$ et de $X$, c'est un morphisme strictement
  polystable de log schémas.
\end{itemize}
Une \emph{log fibration polystable} (resp. \emph{log
  fibration strictement polystable}) $\underline X$ sur $S$ de longueur $l$
est une suite de morphismes polystables (resp. strictement polystables) de
log schémas
$X_l\to\cdots\to X_1\to X_0=S$.\\
Un morphisme $\underline f:\underline
X\to\underline Y$ de log fibrations polystables de longueur $l$  est donné par des morphismes $f_i:X_i\to Y_i$ de log
schémas fs pour tout $i$ tel que le diagramme évident commute.\\
Un morphisme $\underline f$ de log fibrations polystables sera dit két
(resp. étale strict) si $f_i$ est két (resp. étale strict) pour tout $i$.\end{dfn}

Une composée de morphismes plurinodaux est plurinodale. Par contre une
composée de morphismes polystables n'est pas nécessairement polystable.
Un morphisme polystable (resp. strictement polystable) de log schémas est
plurinodal (resp. strictement plurinodal) car le morphisme \[P\to
Q=(P\oplus<T_{0},\cdots,T_{n}>/(T_{0}\cdots
T_{n}=a)\]
peut se décomposer en une suite de morphismes
\[P\to P_1=P\oplus<T_0,T'_0>/T_0T'_0=a\to
P_2=P_1\oplus<T_1,T'_1>/T_1T'_1=T'_0\]\[\to\cdots\to Q=P_{n-1}\oplus<T_{n-1},T_n>/T_{n-1}T_n=T'_{n-2}.\]

Un morphisme plurinodal de log schémas est log lisse et saturé.\\
\begin{rem}
Dans la définition de morphismes strictement polystables de log schémas,
si l'on choisit une autre carte $P'\to
O(X')$, il existe un voisinage ouvert de Zariski $X''$ de $x$ dans $X'$ tel que
$Y''=Y'\times_{X'}X''\to X''$ se factorise à travers un morphisme étale
strict $Y''\to X''\times_{\mbf Z[P']}\mbf Z[Q']$ où
\[Q'=(P'\oplus\bigoplus_{i=0}^{p}<T_{i0},\dots,T_{in_i}>)/(T_{i0}\cdot\dots\cdot T_{in_i}=a'_i)\]
avec $a'_i\in P'$. En effet, on peut supposer $X'$ affine. Ainsi
$X'=\Spec(A)$. Quite à remplacer $X'$ par un voisinage ouvert de Zariski de
$x$, on peut choisir $a'_i$ qui a la même image dans
$\overline M(X')$ que $a_i$ (rappelons que comme $X'$ admet une carte
globale, la log structure de $X'$ est
zariskienne, \ie $\epsilon^*\epsilon_*M=M$ où $\epsilon: X_{\et}\to X_{\Zar}$). Ainsi, dans $M(X')$, $a'_i=a_iu_i$ avec $u_i\in
A^*$. On remplace alors simplement $T_{i0}$ par $T_{i0}u_i$ pour tout $i$.
\end{rem}

\begin{lem}\label{logstrpluri} Soit $\phi:Y\to X$ un morphisme plurinodal (resp. polystable)
  de schémas (au sens de la définition~\ref{defpolystable}), tel que $X$
  ait une log structure log régulière $M_X$ et que $\phi$ soit lisse au-dessus de
  $X_{\text{tr}}$. Alors
  $(Y,O_Y\cap j_*O^*_{Y_{X_{\text{tr}}}})\to (X,M_X)$ est un morphisme 
  plurinodal (resp. polystable) de log schémas.\\
Par le même argument, on en déduit le résultat pour un morphisme
plurinodal standard. Le cas général s'en déduit par dévissage.\end{lem}
\dem
Prouvons-le dans le cas d'un morphisme polystable.\\
On peut supposer que $X=\Spec(A)$ a une carte $\psi:P\to A$ pour la log
structure $M_X$ et que
$Y=B_0\times_X\cdots\times_XB_p$ avec $B_i=\Spec
A[T_{i0},\cdots,T_{in_i}]/(T_{i0}\cdots T_{in_i}-a_i)$ avec $a_i\in A$. Puisque
$\phi$ est lisse sur $X_{\text{tr}}$, $a_i$ est
inversible sur $X_{\text{tr}}$, alors quitte à localiser pour la topologie
de Zariski et à multiplier $a_i$ par un
élément de $A^*$ (ce qu'on peut faire quitte \`a multiplier $T_{i0}$ par le même
élément), on peut supposer que $a_i=\psi(p_i)$ pour un certain $p_i\in P$. Ainsi on a  $Y=X\times_{\mbf Z[P]}\mbf Z[Q]$ où
$Q=(P\oplus\bigoplus_{i=0}^{p}<T_{i0},\cdots,T_{in_i}>)/(T_{i0}\cdot\dots\cdot+T_{in_i}=p_i)$
avec $p_i\in P$. Si l'on munit $Y$ de la log structure $M_Y$ associée à
$Q$, $Y\to X$ devient un morphisme polystable de log schémas. En
particulier $Y$ est log régulier~(\cite[th. 8.2]{kato2}). Puisque
l'ensemble des points de $Y$ où $M_Y$
est trivial est $Y_{X_{\text{tr}}}$, $M_Y=O_Y\cap
j_*O^*_{Y_{X_{\text{tr}}}}$ d'après~\cite[prop. 2.6]{niziol}. 
\findem
Dans le lemme précédent, si $(X,M_X)$ est supposé de plus Zariski log
régulier, et $\phi$ est strictement plurinodal (resp. strictement
polystable), alors pour la log structure du lemme sur $Y$, $\phi$ est un
morphisme strictement plurinodal (resp. strictement polystable) de log schémas.

\subsection{Ensemble polysimplicial d'un log schéma két au-dessus d'une log
  fibration polystable}\label{logpolysimplicial}
Soit $s$ un log point.\\
Nous construisons ici l'ensemble polysimplicial associé à un log schéma két
$Z$ au-dessus d'une log fibration polystable $X\to\cdots\to s$. Il
faut essentiellement montrer que, localement pour la topologie étale de
$X$, si $x\leq y$ sont deux strates de $X$, $Z\to X$ induit naturellement
une application $Z_y\to Z_x$, où $Z_y$ (resp. $Z_x$) est l'ensemble des
composantes connexes de la préimage de $Y$ (resp. $x$).\\ 
Pour
cela nous étudierons la stratification de $Z$ définie par
$\rk(z)=\rk(\overline M^{\gp}_{\bar z})$ (o\`u $\bar z$ est un point g\'eom\'etrique au-dessus de $z$). Cette stratification coïncide avec celle
de Berkovich pour les schémas plurinodaux, et nous montrerons que
localement pour la topologie étale un morphisme két $X\to Y$ induit un
isomorphisme entre les ensembles ordonnés des strates de $X$ et de
$Y$. Cela nous permettra de définir un ensemble polysimplicial de $Z$
localement pour la topologie étale. Nous construirons alors l'ensemble
polysimplicial associ\'e \`a $Z$ par descente (pour qu'il vérifie la même propriété  que dans la
proposition~\ref{berk69}).\\

Pour une (log) fibration polystable $\underline X:X\to\cdots\to \Spec k$,
Berkovich définit un ensemble polysimplicial $\C(\underline X)$. Dans cette
partie nous voulons généraliser cette construction à n'importe quel $Z$ két
au-dessus de $X$.\\
Quand $\underline X$ est strictement polystable, $\C(Z)$ sera défini de
manière à ce que pour toute strate $x$ de $X$ de point générique $\tilde x$
et pour tout objet $x'$ de
$\mbf \Lambda /\C(X)$ au-dessus de $x$, les objets de $\mbf \Lambda/\C(Z)$
au-dessus de
$x'$ seront en bijection naturelle avec $\tilde x$ dans $Z$.\\
Quand $\underline X$ n'est plus supposé strictement polystable, nous
définissons $\C(Z)$ par descente étale.\\ 

Soit $Z$ un log schéma fs, on obtient une stratification sur $Z$ en disant qu'un
point $z$ de $Z$ est de rang $r$ si $\rklog(z)=\rk(M^{\gp}_{\bar z}/\mathcal
O_{\bar z}^*)=r$ (où $\bar z$ est un point géométrique au-dessus de $z$ et
$\rk$ est le rang d'un groupe abélien de type fini).\\
Le sous-ensemble des points de $Z$ de rang $\leqslant r$ est ouvert
dans $Z$ (\cite[cor 2.3.5]{ogus}). Ainsi cela définit bien une stratification localement ferm\'ee.\\
Les strates de rang $r$ de $Z$ sont les composantes connexes de l'ensemble
des points $z$ d'ordre $r$. Les strates forment une partition de $Z$, et
une strate de rang $r$ est ouverte dans le fermé des points $x$ de rang $\geqslant
r$. Il est muni de la structure de sous-schéma réduit de $Z$.\\
L'ensemble des strates est ordonné par~: $x\leqslant y$ si et seulement si
$y\subset \bar x$. On note $\Str_x(Z)$ l'ensemble ordonné des strates
en dessous de $x$ (\ie contenant l'adh\'erence de $x$).\\

Si $f:Z' \to Z$ est un morphisme két, alors $\rklog(x)=\rklog(f(x))$. Ainsi
les strates de $Z'$ sont les composantes connexes des préimages des strates
de $Z$.\\

Soit $Z$ un log schéma plurinodal sur un log point
$s=(k,M_k)$ de caractéristique $p$ et de rang $r_0$, et soit $z$ un point de
$Z$.\\

On a $\rklog(z)=r_0+\rk(z)$ où $\rk(z)$ est la codimension de la strate
contenant $z$ dans $Z$ pour la stratification de Berkovich des schémas
plurinodaux. Ainsi les strates de la stratification logarithmique et de la
stratification de Berkovich sont les mêmes.\\

\begin{lem}
Soit $Z\to\Spec k$ un morphisme plurinodal de log schémas fs sur un log point
et soit $Z'\to Z$ un morphisme két. Alors les strates de $Z'$ (munies de la
structure de sous-schémas réduits) sont normales et donc irréductibles.\\
\end{lem}
Nous noterons parfois abusivement de la même façon une strate et son point générique.\\
\dem
Il suffit de le prouver localement pour la topologie étale.\\
Supposons donc que $Z\to \Spec k$ a une carte exacte \[\begin{array}{ccc} Z &
  \to & \Spec P \\ \dar & & \dar \\ \Spec k & \to & \Spec M\end{array}\]
où $P$ est une carte de $Z$, $M$ est une carte de $\Spec k$ et $Z\to\Spec k[P]\times_{\Spec k[M]}\Spec k$ est étale.\\
Les strates de $Z$ correspondent
aux composantes connexes des préimages au-dessus des différents points de 
$\Spec P$  (et donc des faces de $P$) qui sont au-dessus du point minimal
de $\Spec M$.\\
Soit $P\to Q$ un morphisme $(p')$-kumm\'erien de monoïdes (où $p$ est la
caractéristique de $k$, qui peut éventuellement être nulle). L'application
$\Spec Q \to \Spec P$ est alors bijective (à une face $F$ de $P$ correspond
réciproquement le
saturé $F_Q$ de $F$ dans $Q$).\\
Soit donc $F$ une face de $P$ et $\fk p$ l'idéal premier
$P\backslash F$ (et $F_Q$ et $\fk p_Q$
leurs saturés dans $Q$).\\
Alors la préimage de la cl\^oture du point $\fk p$ dans $\Spec
P$ par le morphisme $\Spec k[P]\to \Spec P$ correspond au sous-schéma $\Spec
k[P]/k[\fk p]$ de $\Spec k[P]$ d'après~\cite[I.3.2]{ogus} (et
c'est la clôture de la strate $\Str F$ de $\Spec k[P]$ correspondant à $F$). Mais
\[\Spec(k[Q]/k[\fk p_Q]) \to \Spec(k[Q]/(k[\fk p]k[Q]))=\Spec k[Q] \times_{k[P]}\Spec(k[P]/k[\fk p])\]
est juste le plongement du sous-schéma réduit (car $k[\fk p_Q]/(k[\fk
p]k[Q])$ est un idéal nilpotent $k[Q]/(k[\fk p]k[Q])$). Ainsi la
préimage de la clôture de la strate correspondant à $F$ est le support du
sous-schéma fermé $\Spec(k[Q]/k[\fk p_Q])$.\\
De plus, il existe d'après~\cite[I.3.2]{ogus} des isomorphismes canoniques
de $k$-algèbres $k[P]/k[\fk p] \to k[F]$ et
$k[Q]/k[\fk p_Q] \to k[F_Q]$ (mais la log structure sur $\Spec k[F]$
obtenue par l'isomorphisme de schéma affine correspondant
n'est
pas celle induite par $F$) et le diagramme de $k$-algèbres suivant commute~:
\[\begin{array}{ccc} k[P]/k[\fk p] & \simeq & k[F] \\ \dar & & \dar
  \\ k[Q]/k[\fk p_Q] & \simeq & k[F_Q] \end{array}\]
($k[F] \to k[F_Q]$ \'etant induit par le plongement de monoïdes $F \to F_Q$).\\
La préimage de $F$ par $\Spec k[P]\to \Spec P$ (\ie la strate
$\Str F$ de $\Spec(k[P])$) correspond alors à l'ouvert $\Spec k[F^{\gp}]$ de $\Spec k[F] \simeq \Spec k[P]/k[\fk p]$.\\
Le carré suivant
\[\begin{array}{ccc} k[F] & \to & k[F^{\gp}]\\ \dar & & \dar \\ k[F_Q] & \to & k[F^{\gp}_Q] \end{array}\]
est cocartésien. Mais $\Spec k[F^{\gp}_Q] \to \Spec k[F^{\gp}]$ est
étale. Ainsi, $(\Str F\times_{\Spec k[P]} \Spec k[Q])^{\red} \to \Str F$ est étale.\\
En faisant le changement de base $Z\to \Spec P$, on obtient que le
morphisme d'une strate de $Z_Q$ à la strate correspondante de $Z$ est étale. Comme les strates de $Z$ sont normales, celles de $Z_Q$ aussi (et donc celles de $Z'$ aussi).\findem

\begin{lem} Soit $Z'\to Z$ un morphisme két et $Z\to \Spec k$ un morphisme
  strictement plurinodal de log schémas, alors $Z'$ est quasinormal (cf. \ref{berkspaces}).\end{lem}
\dem

Nous montrerons que la clôture de la préimage dans $Z'$ d'une strate de $x$
de $Z$ est normale. On peut le faire localement pour la topologie étale de
$Z'$ (mais pas de $Z$).\\

Soit $z$ un point de $Z'$.\\
On peut supposer que $Z'$ est connexe et $Z'\to\Spec k$ a une carte 
fs exacte en $z$ \[\begin{array}{ccc} Z' &
  \to & \Spec P \\ \dar & & \dar \\ \Spec k & \to & \Spec M\end{array},\]
telle que $M$ soit aigu et $Z'\to \Spec k\times_{\Spec\mbf Z[M]}\Spec\mbf Z[P]$ soit étale.\\
Alors les préimages des différents points de $\Spec P$ (\ie les idéaux
premiers de $P$) s'envoient toutes sur des strates différentes de $Z$, puisque $Z$ est strictement
plurinodal.\\
Ainsi, l'image réciproque de $x$ est soit vide soit l'image réciproque de
$\fk p$ pour un certain idéal premier $\fk p$ de $\Spec P$.\\

$\Spec k\to \Spec k[M]$ est l'immersion fermée correspondant à la face 
$M^*=\{0\}$ de $M$ (et à l'idéal premier $M\backslash \{0\}$). Ainsi il s'identifie
à $\Spec k[M]/k[M\backslash \{0\}] \to \Spec k[M]$. On a alors
le diagramme commutatif suivant:
\[\begin{array}{ccc} \Spec k[P]/k[P(M\backslash \{0\})] & \to & \Spec k[P] \\
  \dar & & \dar \\ \Spec k[M]/k[M\backslash \{0\}] & \to & \Spec
  k[M] \end{array}\]
Le strates de $\Spec k[P]/k[P(M\backslash \{0\})]$ correspondent bijectivement
aux idéaux premiers $\fk p$ de $P$ qui contiennent $P(M\backslash \{0\})$
et l'adhérence de la strate correspondante est $\Spec k[P]/k[\fk
p]\simeq \Spec k[F]$ (où $F$ est la face $P\backslash \fk
p$). $\Spec k[F]$ est connexe, et grâce à~\cite[prop I.3.3.1
(2)]{ogus}, puisque $F$ est un monoïde fs (car $P$ est fs
et $F$ est une face de $P$), $\Spec k[F]$ est normal (donc irréductible~:
il existe une unique strate au-dessus de $\fk p$. Les clôtures des
strates de $\Spec k\times_{k[M]}\Spec k[P]$ sont donc également normales.\findem

On peut alors prouver l'analogue de~\cite[lem 2.10]{berk2} (dont la preuve
donne le résultat pour n'importe quel morphisme quasinormal) dans le cas où $Z$
est két au-dessus d'un log schéma strictement plurinodal~:
\begin{lem}\label{lemetale} Soit $\phi:Z'\to Z$ un morphisme étale avec $Z$
  két
  au-dessus de $Z_0$ et soit $Z_0\to k$ un morphisme strictement
  plurinodal. Soit $z'$ une strate de $Z'$ d'image $z$. Alors l'application
  $\Str_{z'}(Z')\to \Str_z(Z)$ est un isomorphisme d'ensembles
  ordonnés.\end{lem}
\dem
Les éléments minimaux des deux ensembles sont les points génériques des
composantes irréductibles de $Z'$ et $Z$ passant par $z'$ et $z$
respectivement.Il s'en suit que
$\Str_{z'}(Z')\cap\Norm(Z')\to\Str_{z}(Z)\cap\Norm(Z)$ est bijectif. On
s'est ramené à prouver le résultat pour $Z'^{(1)}\to Z^{(1)}$, ce qui
prouve par récurrence que $\Str_{z'}(Z')\to\Str_{z}(Z)$ est bijectif. Il
reste à prouver que l'application réciproque est croissante, \ie que si $\phi(z'_2)\leq\phi(z'_1)$, alors $z'_2\leq
z'_1$. Comme $\overline{\phi(z'_2)}$ est normal,
$\phi^{-1}(\overline{\phi(z'_2)})$ est l'union disjointe de ces
composantes irréductibles et donc $\overline z'_2$ est la seule contenant
$z'$. Il s'en suit que $z'_2\leq z'_1$.
\findem
On peut obtenir un raffinement de~\ref{lemetale} en prouvant le résultat pour $Z'\to Z$ két:
\begin{lem}\label{lemkum} Soit $Z\to Z'$ un morphisme két avec $Z'$
  két sur $Z''$ et soit $Z''\to k$ un morphisme strictement
  plurinodal. Soit $z$ une strate de $Z$ d'image $z'$ dans $Z'$. Alors
  l'application $\Str_{z}(Z)\to \Str_{z'}(Z')$ est un isomorphisme
  d'ensembles ordonnés.\end{lem}
\dem 
Il suffit de le prouver quand $Z'=Z''$, puisque si l'on connaît le résultat pour
$Z'\to Z''$ et $Z\to Z''$, cela implique le résultat pour $Z\to Z'$.\\

Si $Z'_0=\Spec k\times_{k[M]}\Spec k[P]$ et $Z_0=\Spec
k\times_{k[M]}\Spec k[Q]$ où $P\to Q$ est két, alors l'ensemble ordonné des strates de
$Z'_0$ (resp. $Z_0$) est isomorphe à l'ensemble ordonné des
strates de $P$
(resp. $Q$) qui s'envoient sur la face $M^*$ de $M$. Mais comme
l'application de l'ensemble des strates de $Q$ dans l'ensemble des strates
de $P$ est un isomorphisme, $\Str Z'_0 \to \Str Z_0$ est un isomorphisme,
et a fortiori $\Str_{z_0}(Z_0)\to \Str_{z'_0}(Z'_0)$.\\
Dans le cas général, localement dans un voisinage de Zariski de $z$ et de $z'$,
il existe un diagramme commutatif comme suit~:
\[\begin{array}{ccccc}Z & \leftarrow & U & \to & Z_0 \\ \dar & & \dar &
  & \dar \\ Z' & = & Z' & \to & Z'_0\end{array}\]
où les flèches horizontales sont étales: elles vérifient le
lemme~\ref{lemetale}. Grâce au cas particulier précédent, on obtient le résultat.
\findem

Considérons maintenant une log fibration strictement polystable $\underline
X:X\to X_{l-1}\to\cdots\to s$ où $s$ est un log point fs.\\
Si $f:Z \to X$ est két, on a un foncteur 
\[D_Z:(\mbf\Lambda/\C(\underline X))^{\op} \Ens\]
obtenu comme la composée des deux foncteurs
\[(\mbf\Lambda/\C(\underline X))^{\op} \to \Str(X_s)\to\Ens\]
où le foncteur de droite associe
à une strate de $X$ l'ensemble des composantes connexes de la préimage de
cette strate (la fonctorialité vient du lemme~\ref{lemkum}).\\
Si $\C$ est un ensemble polysimplicial et $D:(\Lambda/\C)^{\op}\to\Ens$ on
peut d\'efinir un ensemble polysimplicial $\C\sq D$ d\'efini par $(\C\sq
D)_{\mbf n}=\coprod_{x\in\C_{\mbf n}}D(x)$ et, pour $f:\mbf m\to\mbf n$ et
$d\in D(x)$ avec $x\in\C_{\mbf n}$, $f^*d=D(f_x)(d)\in
D(f^*x)\subset\coprod_{y\in\C_{\mbf m}}D(y)$ (c'est un cas particulier
trivial de la construction faite par Berkovich et rappelé dans~\ref{berkspaces}).\\
Cette construction définit donc un
ensemble polysimplicial 
\[\C_{\underline X}(Z)=\C(\underline X)\sq D_Z\]
(nous écrirons souvent $\C(Z)$
à la place de $\C_{\underline X}(Z)$). Cet ensemble polysimplicial est
encore intérieurement libre. $O(\C(Z))$ est isomorphe à $\Str(Z)$ fonctoriellement en $Z$.
\begin{rem}
Soit $\C\to \C'$ un morphisme d'ensembles polysimpliciaux. Soit $\alpha:S\to
O(\C)$ (resp. $\alpha':S'\to O(\C')$) un morphisme d'ensembles ordonnés tel que
$S_{\leq x}\stackrel{\simeq}{\to} O(\C)_{\leq\alpha(x)}$ (resp. $S'_{\leq
    x}\stackrel{\simeq}{\to} O(\C')_{\leq\alpha'(x)}$) pour tout $x$. On a
 des foncteurs \[D:\mbf \Lambda /\C\to O(\C)\to \Set \text{ (resp. } D':\mbf \Lambda
  /\C'\to O(\C')\to \Set),\] qui définissent un ensemble polysimplicial $\C_1=\C\sq D$ (resp.
$\C'_1=\C\sq D'$). On a $O(\C_1)=S$ et $O(\C'_1)=S'$.\\
Soit $f:S\to S'$ morphisme d'ensembles ordonnés  tel que
\[\begin{array}{ccc} S & \to & S' \\ \dar & & \dar \\ O(\C) & \to & O(\C')
\end{array}\]
commute. Alors on a une application naturelle $\C_{1,\mbf n}=\coprod_{x\in \C_{\mbf n}}S_x\to\coprod_{x'\in\C_{\mbf n}}S'_{x'}=\C'_{1,\mbf n}$ qui
d\'efinit un morphisme d'ensembles polysimpliciaux $\underline f:\C_1\to \C'_1$
au-dessus de $\C\to \C'$ tel que $O(\underline f)=f$ (il est facile de voir que $\underline f$ est en fait le seul morphisme $\tilde f:\C_1\to\C'_1$ au-dessus de $\C'\to C$ tel que $O(\tilde f)=f$).\\
Grâce à cela, pour construire des morphismes entre les ensembles
polysimpliciaux de log schémas két au-dessus de log schémas strictement
plurinodaux, il nous suffira souvent de construire un morphisme entre les ensembles
ordonnés des strates.
\end{rem}

Si $\underline X'\to\underline X$ est un morphisme két de log
fibrations polystables, alors $\C_{\underline X}(X'_l)$ est canoniquement
isomorphe à
$\C(\underline X')$.\\
Si l'on a un diagramme commutatif
\[\begin{array}{ccc} Z' & \to & Z' \\ \dar & & \dar \\ \underline X' & \to &
  \underline X\end{array}\]
où $\underline X' \to\underline X$ est un morphisme két de log
fibrations strictement polystables, il existe un morphisme induit $\C_{\underline
  X'}(Z')\to\C_{\underline X}(Z)$. Si
$\Str(Z')\to\Str(Z)$ est un isomorphisme d'ensembles ordonnés, alors $\C_{\underline
  X'}(Z')\to\C_{\underline X}(Z)$ est un isomorphisme (car il envoie les polysimplexes non d\'eg\'en\'er\'es
sur les polysimplexes non dégénérés et les complexes polysimpliciaux sont
intérieurement libres).\\
Soit $Z'\to Z$ un recouvrement két, soit $Z''=Z'\times_ZZ'$ et
soit $x$ une strate de $X_s$, alors $D_Z(x)=\Coker(D_{Z''}(x)\rightrightarrows
D_{Z'}(x))$. On en déduit que
\[\C(Z)\to\Coker(\C(Z'')\rightrightarrows \C(Z'))\]
est un isomorphisme.\\

On peut aussi définir $\C_{\underline X}(Z)$ pour $\underline X$ une
fibration polystable générale. Soit $\underline X'\to\underline X$ un
recouvrement étale avec $\underline X'$ strictement polystable, soit $\underline X''=\underline X'\times_{\underline X}\underline
X'$ et soit $Z'$ et $Z''$ les pullbacks de $Z$ à $X'$ et
$X''$. On définit alors \[\C_{\underline X}(Z)=\Coker(\C_{\underline
  X''}(Z'')\rightrightarrows \C_{\underline X'}(Z')).\] Par un argument
souvent utilis\'e dans~\cite{berk2}, on d\'eduit de l'\'equation
pr\'ec\'edente que cela ne dépend pas 
du choix de $\underline X'$. En effet soit $\underline X_1$ et $\underline
X_2$ deux recouvrement \'etales par une log fibration strictement polystable de $\underline
X$. Posons ${\underline X}_{11}:= {\underline X}_1\times_{\underline
  X}{\underline X}_1$, ${\underline X}_{12}:= {\underline
  X}_1\times_{\underline X}{\underline X}_2$, ${\underline
  X}_{22}:={\underline X}_2\times_{\underline X}{\underline X}_2$,
${\underline X}_{112}:= {\underline X}_{11}\times_{\underline X}{\underline
  X}_2$, ${\underline X}_{122}:={\underline X}_1\times_{\underline
  X}{\underline X}_{22}$, ${\underline X}_{1122}={\underline
  X}_{11}\times_{\underline X}{\underline X}_{22}$ et nous aurons des
notations similaires en rempla\c cant $\underline X$ par $Z$ quand on fait le
changement de base $Z\to X$. On obtient alors 
\[\begin{array}{rcll} \Coker(\C(Z_{11})\rightrightarrows\C(Z_1)) & = &
  \Coker\big( &\!\!\!\!\!\! \Coker(\C(Z_{1122})
  \rightrightarrows\C(Z_{112}))\\ & & & \rightrightarrows
  \Coker(\C(Z_{122})\rightrightarrows\C(Z_{12}))\big)\\ 
& = &
\Coker\big( &\!\!\!\!\!\!\Coker( \C(Z_{1122})\rightrightarrows\C(Z_{122}))\\ & & &
\rightrightarrows\Coker(\C(Z_{112})\rightrightarrows\C(Z_{12}))\big)\\  
& & & \!\!\!\!\text{par intervention de colimites}\\
& = & \Coker(&\!\!\!\!\!\!\!\C(Z_{22})\rightrightarrows\C(Z_2)).\end{array}\]
$\C_{\underline X}(Z)$ est donc bien d\'efini (qu'on notera souvent
simplement 
$\C(Z)$ par abus de notation).\\

Si $Z'\to Z$ est un morphisme két surjectif au-dessus de $\underline X$ et
$Z''=Z'\times_ZZ'$, alors
$\C(Z)=\Coker(\C(Z'')\rightrightarrows \C(Z'))$.\\
On obtient donc ($\ket(X)$ est la catégorie des log schémas két sur $X$):
\begin{prop}\label{ketpolysimpcplx} Soit $\underline X$ une log fibration
  polystable, on a un
  foncteur $\C_{\underline X}:\ket(X)\to\mbf\Lambda^{\circ}\Ens$ tel que:
\begin{itemize}
\item si $Z'\to Z$ est un recouvrement dans $\ket(X_s)$,
\[\C_{\underline X}(Z)=\Coker(\C_{\underline
  X}(Z'\times_ZZ')\rightrightarrows \C_{\underline X}(Z')).\]
\item $O(\C_{\underline X}(Z))$ est fonctoriellement isomorphe à $\Str(Z)$.\end{itemize}\end{prop}
\begin{rem} Si l'on a un morphisme két $\underline Y\to\underline X$ de
  fibrations polystables de longueur $l$, le complexe polysimplicial
  $\C(Y_l)$ que l'on vient de définir en considérant $Y_l$ comme étant dans
  $\ket(X_l)$ est canoniquement isomorphe au complexe polysimplicial de la
  fibration polystable $\C(\underline
  Y)$ défini par Berkovich.
\end{rem}

Si $s$ est un log point et $s'$ est un revêtement két connexe, on note
$\tilde s'$ le sous-log-schéma strict réduit de $s'$. C'est un log point.
\begin{prop}Supposons $Z$ quasicompact. Il existe un revêtement két connexe $s'\to s$ tel que pour tout
  morphisme de log points fs $s''\to\tilde s'$, $\C(Z_{s''})\to\C(Z_{s'})$
  soit un isomorphisme.
\end{prop}
\dem
Si $Z$ est quasicompact et $\underline X$ est strictement polystable, il
existe un revêtement két connexe $s'\to s$
tel que toutes les strates de $Z_{\tilde s'}$ soient géométriquement irréductibles
et que $Z_{s'}\to s'$ soit saturé. En particulier, pour tout morphisme de log
points $s''\to \tilde s'$, $\C(Z_{s''})\to\C(Z_{\tilde s'})$ est un
isomorphisme (puisque les complexes sont intérieurement libres).\\
Si l'on ne suppose plus $\underline X$ strictement polystable, soit
$\underline X'\to\underline X$ un recouvrement étale de type fini avec
$\underline X'$ strictement polystable. Alors, il existe un revêtement két
connexe $s'$ pour lequel $\C(Z''_{s''})\to\C(Z''_{\tilde s'})$ et
$\C(Z'_{s''})\to\C(Z'_{\tilde s'})$ sont des isomorphismes pour toute
extension $s''$ de $\tilde s'$. $\C(Z_{s''})\to\C(Z_{\tilde s'})$ est donc
également un isomorphisme.
\findem
Le complexe polysimplicial $\C(Z_{\tilde s'})$ pour un tel $s'$ est noté
$\Cgeom(Z/s)$.\\
Si $Z'\to Z$ est un recouvrement két dont le morphisme de schéma
sous-jacent est de type fini, notons $Z''=Z'\times
_ZZ'$. On a encore $\Cgeom(Z/s)=\Coker(\Cgeom(Z''/s)\rightrightarrows\Cgeom(Z/s))$. 

\subsection{Groupe fondamental tempéré d'une log fibration polystable}\label{tfgsp}
Nous définissons ici le groupe fondamental tempéré d'une log fibration polystable $\underline
X:X\to\cdots\to\Spec k$ sur un log point fs. Si $T$ est un revétement két
de $X$, les revêtements topologiques de
$|\C(T)|$ joueront le rôle de revêtements topologiques de $T$.\\
 
Commençons par une construction catégorique de groupes fondamentaux
tempérés que nous utiliserons plus tard dans notre situation logarithmique.\\
Considérons une catégorie fibrée
$ \mcal D\to \mcal C$ vérifiant les propriétés suivantes~:
\begin{itemize}
\item $\mcal C$ est une catégorie galoisienne,
\item pour tout objet connexe $U$ de $\mcal C$, $\mcal D_U$ est équivalente
  à la catégorie $\Pi_U\tSet$
pour un certain groupe discret $\Pi_U$,
\item Si $U$ et $V$ sont deux objets de $\mcal C$, le foncteur $\mcal
  D_{U\coprod V}\to \mcal D_U\times\mcal D_V$ est une équivalence,
\item si $f:U\to V$ est un morphisme dans $\mcal C$, $f^*:\mcal D_V\to \mcal
  D_U$ est exact.
\end{itemize}
Alors, on peut définir une catégorie fibrée $\mcal D'\to\mcal C$ telle que
la fibre en $U$ soit la catégorie des données de descente de $\mcal
D\to\mcal C$ pour le morphisme $U\to e$ (où $e$ est l'élément final de $\mcal
C$).\\
Soit $U$ un objet galoisien connexe de $\mcal C$ et soit $H$ le groupe de
Galois de $U/e$. On peut voir $H$ comme une cat\'egorie ayant un unique
objet $u$, tel que $\End(u)=H$. On a un foncteur $H\to\mcal C$, qui envoie
$u$ sur $U$. La cat\'egorie $\mcal D'_U$ est naturellement \'equivalente
\`a la cat\'egorie $\mcal B(\mcal D/H)$ des sections cart\'esiennes de la
cat\'egorie fibr\'ee $\mcal D\times_{\mcal C}H/H$. $\mcal D\times_{\mcal
  C}H/H$ est un complexe classifiant.\\
Supposons que l'on ait un scindage de $\mcal D/\mcal C$.\\
Alors $\mcal D'_U$ peut être décrit de la façon suivante~:
\begin{itemize}
\item ces objets sont les couples 
$T_U=(S_U,(\psi_h)_{h\in H})$, où $S_U$ est un objet de $\mcal D_U$ et $\psi_h:S_U\to
h^*S_U$ est un isomorphisme dans $\mcal D_U$ tel que pour tout $h,h'\in G$,
$(h^*\psi_h')\circ\psi_h=\psi_{h'h}$ (après avoir identifié $(h'h)^*$ et
$h^*{h'}^*$ par l'isomorphisme canonique pour simplifier les notations).
\item un morphisme $(S_U,(\psi_h))\to (S'_U,\psi'_h)$ est un morphisme $\phi:S_U\to
  S'_U$ dans $\mcal D_U$ tel que pour tout $g\in G$,
  $\psi'_h\phi=(h^*\phi)\psi_h$.\end{itemize}
Il existe un foncteur naturel $F_0: \mcal D'_U\to\mcal D_U$, qui envoie
$(S_U,(\psi_h))$ vers $S_U$.
Soit $F_U$ un foncteur fondamental $\mcal D_U\to \Set$, tel que $\Aut
F_U=\Pi_U$.\\
Soient $F:=F_UF_0$ et $\Pi'_U:=\Aut F$ (muni de la topologie discr\`ete).\\
\begin{prop}\label{extgp}
\begin{itemize}
\item Le foncteur naturel $\mcal D'_U\to\Pi'_U\tSet$ est une
  équivalence de catégories.
\item  Il existe une suite exacte naturelle \[1\to \Pi_U\to\Pi'_U\to
  H\to 1.\]
\end{itemize}
\end{prop}
\dem
En effet, on sait d\'ej\`a que $\mcal D'_U$ est une cat\'egorie
classifiante. Ainsi il suffit de prouver que son groupe fondamental est
discret, \ie qu'il existe un objet $T^\infty$ (appelé revêtement universel)
tel que pour tout
$t^\infty\in F(T^\infty)$ et pour  tout objet $T$ et tout $t\in F(T)$, il
existe un unique morphisme $f:T^\infty\to T$ tel que $F(f)(t^\infty)=t$. En
effet, si $\Pi'_U$ est un groupe discret, $\Pi'_U$ muni de l'action de $\Pi'_U$ par
translation à gauche est un revêtement universel de $\Pi'_U\tEns\simeq\mcal
D'_U$. Réciproquement, soient $(T^\infty,t^\infty)$ un revêtement universel
pointé et $\alpha\in \Stab^{\Pi'_U}_{t^\infty}$. Soit $(T,t)$ un objet
pointé et $f:T^\infty\to T$ qui envoie $t^\infty$ en $t$. Alors
$\alpha_{T}(t)=\alpha_{T}(F(f)(t^{\infty})=F(f)\alpha_{T^\infty}(t^\infty)=F(f)(t^\infty)=t$,
donc $\alpha$ est l'identité. L'élément neutre de $\Pi'_U$ est donc un
sous-groupe ouvert de $\Pi'_U$ et donc $\Pi'_U$ est discret.\\
Soit $S^\infty$ un objet universel de $\mcal D_U$ et soit  $S^0=\coprod_{g\in H}
g^*S^\infty$, et \[\psi_h:S^0=\coprod_{g\in H}
g^*S^\infty=\coprod_{gh\in H}
(gh)^*S^\infty\stackrel{\simeq}{\to}\coprod_{g\in H}
h^*g^*S^\infty=h^*(\coprod_{g\in H}
g^*S^\infty)=h^*S^0.\]
Cela définit un objet $T^\infty$ de $\mcal D'_U$. Soit $t^\infty$ un
élément de $F(t^\infty)=\coprod_g F_U(g^*S^\infty)$ et soit $g_0$ tel que
$F(t^\infty)$ soit dans $F_U(g_0^*S^\infty)$.\\

Soit $T$ un objet de $\mcal D'_U$ et $t\in F(T)$
$g_0^*S^\infty$ est également un objet universel de $\mcal D_U$, donc il
existe un unique morphisme $f_0:g_0^*S^\infty\to F_0(T)$ dans $\mcal D_U$ qui
envoie $t^\infty$ sur $t$. La restriction à $g_0^*S^\infty$ d'un morphisme
$f$ de $T^\infty$ vers $T$ doit donc être $f_0$.\\
Mais l'hypothèse d'équivariance des morphismes de $\mcal D'_U$ revient à
dire que pour tout $g$, la restriction à
$(g_0h)^*S^\infty=h^*g_0^*S^\infty$ est $\psi_{h,T}^{-1}\circ
h^*f_0:h^*g_0^*S^\infty\to F_0(h^*T)\to F_0(T)$, ce qui prouve l'unicité de $f$. Mais
cette formule définit bien un morphisme $\coprod h^*S^\infty\to F_0(T)$,
dont on vérifie facilement qu'elle vérifie la propriété d'équivariance
voulue, ce qui montre l'existence de $f$.\\

$F_0$ induit un morphisme $\Pi_U\to\Pi'_U$. Comme $F_0(T^\infty)$ est une
somme directe d'objets universels de $\mcal D_U$, $\Pi_U$ agit librement
sur $F(T^\infty)$ et donc $F_0$ est injectif.\\
Il existe aussi un foncteur exact naturel $F_1:H\tSet\to \mcal D'_U$ qui envoie
un $H$-ensemble $Y$ vers un couple $(Y=\coprod_{y\in Y}\{y\},(\psi_h))$ où $Y$ est
un objet constant de $\mcal D_U$ et $\psi_h$ envoie $y$ vers $h\cdot y$.\\
Ce foncteur est pleinement fidèle et $FF_1$ est canoniquement isomorphe au
foncteur oubli $H\tSet\to\Set$. $F_1$ induit donc un morphisme surjectif
$\Pi'_U\to H$. Comme $F_0F_1(Y)$ est un objet constant de $\mcal D_U$, le
composé $\Pi_U\to\Pi'_U\to H$ est trivial.\\
On a aussi un foncteur (non exact) $H_0:\mcal D'_U\to H\tEns$ qui envoie
$(S_U,(\psi_h))$ sur l'ensemble des composantes connexes de l'objet $S_U$ de
$\mcal D_U$, où l'action de $g\in G$ est induite par $\psi_h$. Si
$S=(S_U,(\psi_h))$ est un objet connexe de $\mcal D'_U$ tel que
$F_0(S)=S_U$ a une composante connexe triviale, $S\to F_1H_0(S)$ est un
isomorphisme et donc $S$ est dans l'image essentiel de $F_1$. Ainsi la
suite $\Pi_U\to\Pi'_U\to H$ est bien exacte en
$\Pi'_U$.
\findem

On peut faire une construction explicite de $\Pi'_U$ en terme d'un scindage
de la catégorie fibrée. Si $g$ est un
élément d'un groupe $G$, notons $\iota(g)\in\Aut(G)$ la conjugaison par $g$.\\ 
Un scindage de la
catégorie fibrée sur $H$ (vue comme catégorie avec un seul objet et
dont tous les morphismes sont inversibles) de fibre $\Pi\tSet$ et une famille d'isomorphismes de foncteurs $(h^*F_U\to F_U)_{h\in H}$ revient à se
donner~:
\begin{itemize}
\item $\alpha_h\in \Aut\Pi$ pour tout $h\in H$,
\item $g_{h,h'}\in \Pi$ pour tout $(h,h')\in H^2$ tels que
\begin{itemize}
\item
  $\alpha_h\alpha_{h'}=\iota(g_{h,h'})\alpha_{hh'}$
\item $g_{h,h'}g_{hh',h''}=\alpha_h(g_{h',h''})g_{h,h'h''}$ pour tout
  $(h,h',h'')\in H^3$
\end{itemize}
\end{itemize}
Une construction de Schreier permet d'associer \`a de telles donn\'ees une extension de $H$ par $\Pi$ (et toute extension est de cette forme-ci).\\
Définissons l'ensemble
\[\Pi\rtimes H=\{(\alpha,h,g)\in\Aut\Pi\times
H\times\Pi|\alpha=\iota(g)\alpha_h\}\]
muni de la multiplication
\[(\alpha,h,g)(\alpha',h',g')=(\alpha\alpha',hh',g\alpha_h(g')g_{h,h'}).\]
Elle est bien définie car
\[\iota(g\alpha_h(g')g_{h,h'})\alpha_{hh'}=\iota(g)\iota(\alpha_h(g'))\alpha_h\alpha_{h'}=\iota(g)\alpha_h\iota(g')\alpha_{h'}=\alpha\alpha'.\]
On vérifie que la multiplication est bien associative~:\\
$\big((\alpha,h,g)(\alpha',h',g')\big)(\alpha'',h'',g'')=(\alpha\alpha'\alpha'',hh'h'',g_1)$ avec 
\[g_1=g\alpha_h(g')g_{h,h'}\alpha_{hh'}(g'')g_{hh',h''}\]
et $(\alpha,h,g)\big((\alpha',h',g')(\alpha'',h'',g'')\big)=(\alpha\alpha'\alpha'',hh'h'',g_2)$ avec
\[\begin{array}{rcl} g_2 & = & g\alpha_h(g'\alpha_{h'}(g'')g_{h',h''})g_{h,h'h''}\\
& = & g\alpha_h(g')\alpha_h\circ\alpha_{h'}(g'')\alpha_h(g_{h',h''})g_{h,h'h''}\\
& = & g\alpha_h(g')i(g_{h,h'})\circ\alpha_{hh'}(g'')\alpha_h(g_{h',h''})g_{h,h'h''}\\
& = & g\alpha_h(g')g_{h,h'}\alpha_{hh'}(g'')g^{-1}_{h,h'}\alpha_h(g_{h',h''})g_{h,h'h''}\\
& = & g\alpha_h(g')g_{h,h'}\alpha_{hh'}(g'')g^{-1}_{h,h'}g_{h,h'}g_{hh',h''}\\
& = & g\alpha_h(g')g_{h,h'}\alpha_{hh'}(g'')g_{hh',h''}.\end{array}\]

$(1,1,g^{-1}_{1,1})$ est élément neutre de $\Pi\rtimes H$ et
\[(\alpha,h,g)^{-1}=(\alpha^{-1},h^{-1},g_{1,1}^{-1}g_{h^{-1},h}^{-1}\alpha_{h^{-1}}(g)^{-1}).\]
On obtient bien un groupe qui ne dépend essentiellement pas du choix des isomorphismes $h^*F_U\to F_U$.\\
Effectivement si l'on a une famille $(\gamma_h)_{h\in H}$ et que l'on
définit $(\beta_h),(\tilde g_{h,h'})_{(h,h')}$ par
$\beta_h=\iota(\gamma_h)\alpha_h$ et $\tilde
g_{h,h'}=\gamma_h\alpha_h(\gamma_{h'})g_{h,h'}\gamma_{hh'}^{-1}$,ce qui
fournit un nouveau scindage de la catégorie fibrée, et l'application
$(\alpha,h,g)\mapsto (\alpha,h,g\gamma_h)$ est un isomorphisme entre les
deux constructions de $\Pi\rtimes H$.\\

On a un morphisme injectif $\Pi\to\Pi\rtimes H$ qui envoie $g$ sur
$(\iota(g),1,gg^{-1}_{1,1})$ et un morphisme $\Pi\rtimes H\to H$ qui envoie
$(\alpha,h,g)$ sur $h$. D'où une suite exacte \[1\to\Pi\to\Pi\rtimes H\to
H\to 1.\]

\subsection*{}

Si $(U_i,u_i)_{i\in I}$ est un système projectif cofinal d'objet galoisiens pointés
(et soit $P$ l'objet correspondant de pro-$\mcal C$),
on définit alors $\Btemp(\mcal D/\mcal C,P)$ comme étant la catégorie $\injLim_i
\mcal D'_{U_i}$. Un isomorphisme de pro-objets $P\to P'$ induit une
équivalence $\Btemp(\mcal D/\mcal C,P')\to\Btemp(\mcal D/\mcal C,P)$ et
donc $\Btemp(\mcal D/\mcal C,P)$ ne dépend du choix de $(U_i)_i$ qu'à
isomorphisme près.
De plus si $h\in G_i=\Gal(U_i/e)$, l'endofoncteur $h^*:\mcal D'_{U_i}\to\mcal
D'_{U_i}$ envoie $T=(S_{U_i},\psi_g)$ vers
$h^*T=(h^*S_{U_i},\psi_{hg}\psi_h^{-1})$. Alors $\psi_h:S_{U_i}\to
h^*S_{U_i}$ définit un isomorphisme $T\to h^*T$ fonctoriellement en $T$. Donc
$h^*:\mcal D'_{U_i}\to \mcal D'_{U_i}$ est canoniquement isomorphe à
l'identité de $\mcal D'_{U_i}$. En particulier, tout isomorphisme du pro-objet $P$
induit un endofoncteur de $\Btemp(\mcal D/\mcal C,P)$ qui est canoniquement
isomorphe à l'identité (fonctoriellement sur $\Aut P$).\\

Soit $(F_i)_{i\in I}$ une famille de foncteurs fondamentaux $F_i:\mcal
D_{U_i}\to \Set$ et supposons qu'on ait une famille $(\alpha_f)_{f:U_i\to U_j}$, indexée par
l'ensemble des morphismes de $I$, d'isomorphismes de foncteurs $F_if^*\to F_j$ telle
que pour tout $U_i\stackrel{f}{\to} U_j\stackrel{g}{\to} U_k$,
$(\alpha_f\cdot g^*)\alpha_g=\alpha_{gf}$ (après avoir identifié $(gf)^*$ et
$f^*g^*$ pour alléger les notations) ; une telle famille existe par exemple
si $I$ est simplement $\mbf N$.\\
Alors, cela induit un système projectif  $(\Pi'_{U_i})_{i\in I}$ (unique à
isomorphisme près indépendemment de $(\alpha_f)$ si
$I$=$\mbf N$ et les foncteurs $\mcal D'_{U_i}\to\mcal D'_{U_j}$ sont
pleinement fidèles). On peut alors définir
\[\gtemp(\mcal D/\mcal C,(F_i))=\varprojlim\Pi'_{U_i}\]

Supposons que l'on ait un diagramme 2-commutatif à flèches verticales fibrées:
\[\begin{array}{ccc} \mcal D_1 & \to & \mcal D_2\\ \dar & & \dar\\ \mcal C_1 & \stackrel{f}{\to} &
  \mcal C_2\end{array}\]
tel que $f:\mcal C_1\to \mcal C_2$ soit exact, et $\mcal D_{1,U}\to\mcal
D_{2,f(U)}$ soit exact pour tout objet $U$ de $\mcal C_1$.
On obtient alors un foncteur $\Btemp(\mcal D_1/\mcal C_1)\to\Btemp(\mcal
D_2/\mcal C_2)$.\\

Par exemple, soient $X$ une $K$-variété analytique lisse, $\mcal C$ la catégorie des
revêtements étales finis de $X$ et $\mcal D\to \mcal C$ la catégorie fibrée
telle que $\mcal D_U$ soit la catégorie des revêtements topologiques de
$U$. Alors, puisque tout revêtement étale fini est un morphisme de descente
effectif pour les revêtements tempérés, $\mcal D'_U$ s'identifie
fonctoriellement \`a la sous-catégorie pleine de $\Covtemp(X)$ des
revêtements tempérés $S$ tels que $S_U$ soit un revêtement topologique de
$U$. Si $(U_i,u_i)$ est un système cofinal de revêtements galoisiens
pointés de $(X,x)$, alors $\Btemp(\mcal C/\mcal D)$ devient canoniquement
équivalente à $\Covtemp(X)$.\\

Appliquons la construction catégorique du groupe fondamental tempéré à
notre situation logarithmique.\\
Soit $\underline X:X\to X_{l-1}\to\cdots\to\Spec(k)$ une log fibration
polystable, et supposons $X$ connexe.\\
On a alors un foncteur $\Ctop:\KCov(X)\to \Ke$ obtenu en
composant le foncteur $\C$ de la proposition~\ref{ketpolysimpcplx} avec le
foncteur réalisation géométrique.\\
On peut définir une catégorie fibrée naturelle $\Dtop\to \KCov(X)$ telle que
la fibre en un revêtement két $Y$ de $X$ soit la catégorie des revêtements
topologiques de $\Ctop(Y)$ (qui est équivalente à
$\gtop(\Ctop(Y))\tEns$).\\
On définit une catégorie fibrée $\DDtemp\to\KCov(X)$ dont la fibre en un
revêtement két $f:Y\to X$ est la catégorie des données de descente de
$\Dtop\to \KCov(X)$ pour $Y\to X$ (cela correspond heuristiquement aux
revêtements ``tempérés'' de $X$ qui deviennent topologiques après
changement de base $Y\to X$).\\
Soit $x$ un point log géométrique de $X$ et soit $(Y,y)$ un revêtement két
galoisien connexe log géométriquement pointé de $(X,x)$.\\
Soit $\tilde y\to |\C(Y)|$ la cellule fermée  de $|\C(Y)|$ qui correspond à la strate de $Y$
contenant $y$. Elle
est contractile. Notons aussi $\mring y\to\tilde y\to |\C(Y)|$ la cellule
ouverte, qui est aussi contractile. Alors on a un foncteur fondamental $F_y:\Dtop_Y\to\Ens$
qui correspond au point base $\tilde y$ ($F_y(S)$ est l'ensemble des composantes
connexes de $S\times_{|\C(Y)|}\mring y$). De plus, pour tout morphisme $f:(Y',y')\to
(Y,y)$, les deux foncteurs $F_{y'}f^*$ et $F_y$ sont canoniquement isomorphes.\\
On peut alors considérer le foncteur
$F_{(Y,y)}:\DDtemp_Y\to\Ens$ qui associe à une donnée de descente $T$ l'ensemble
$F_y(T_Y)$.\\
Le foncteur induit \[\DDtemp_Y\to \Aut(F_{(Y,y)})\tEns\] est une équivalence
de catégorie et d'apr\`es la proposition~\ref{extgp} l'on a une suite exacte~:
\[1\to \gtop(|\C(Y)|,\tilde y)\to \Aut(F_{(Y,y)}) \to \Gal(Y/X)\to 1.\]
\begin{dfn}\label{defloggft}
Le groupe fondamental $\mbb L$-tempéré du log schéma pluristable $X$ en $x$ est
\[\gtemp(X,x)^{\mbb L}:=\varprojlim_{(Y,y)} \Aut(F_{(Y,y)}),\] où la limite
projective est prise suivant la catégorie filtrante $\LGalKCov(X,x)$ des
revêtements két galoisiens pointés $\mbb L$-finis de $(X,x)$.\end{dfn}

Si $x_1\to x_2$ est une spécialisation de points log géométriques de $X$,
elle induit une équivalence naturelle entre la catégorie des revêtements
pointés de $(X,x_1)$ et la catégorie des revêtements pointés de $(X,x_2)$
(on identifiera donc les deux catégories).\\
Si $Y$ est un revêtement pointé $(Y,y_1)$ de $(X,x_1)$, le revêtement
pointé correspondant de $(X,x_2)$ est $(Y,y_2)$ où $y_2$ est l'unique point
log géométrique au-dessus de $x_2$ tel qu'il existe une spécialisation $y_1\to
y_2$ au-dessus de $x_1\to x_2$ (et cette spécialisation est unique). Alors
il existe une application canonique
$\tilde y_2\to\tilde y_1$ telle que
\[\xymatrix{\tilde y_2 \ar[r] \ar[dr] & \tilde y_1 \ar[d]\\ &
  |\C(Y)|}\] commute.\\
Cela induit un isomorphisme canonique $F_{y_1}\simeq F_{y_2}$, fonctoriel en
$Y$, et l'on obtient donc un isomorphisme canonique $\gtemp(X,x_1)^{\mbb
  L}\to\gtemp(X,x_2)^{\mbb L}$. Si $X$ est connexe et $x_1,x_2$ sont deux
points log géométriques de $X$, il existe une suite de spécialisations et
de cospécialisations joignant $x_1$ à $x_2$. Ainsi $\gtemp(X,x_1)^{\mbb L}$
and $\gtemp(X,x_2)^{\mbb L}$ sont nécessairement isomorphes.\\

On a une équivalence de catégories entre
\[\BtempL_{(X,x)}=\injLim \DDtemp_Y
/\LGalKCov(X,x)\] et la catégorie
$\gtemp(X,x)^{\mbb L}\tEns$ des ensembles munis d'une action de
$\gtemp(X,x)^{\mbb L}$ qui se factorise à travers un quotient discret de $\gtemp(X,x)^{\mbb L}$.\\

Supposons maintenant que $X$ est log géométriquement connexe, \ie que
$X_{k'}$ est connexe pour toute extension két $k'$ de $k$.\\
Soit $\bar k$ un point log géométrique de $k$, soit $\bar x=(\bar x_{k'})$
un système compatible de points log géométriques de $X_{k'}$ où $k'$ parcourt
les extensions két de $(k,\bar k)$. Par exemple, pour construire un tel système,
on peut prendre un point géométrique de $X_{\tilde k}$ où $\tilde k$
est une clôture séparable stricte de $k$.  $\glog(\tilde k)$ est isomorphe
à $\Hom(\overline M_{\bar k}^{\gp},\widehat{\mbf Z}^{(p')}$ avec $\overline
M_{\bar k}^{\gp}$ de type fini, donc $\glog(\tilde k)$ est
topologiquement finiment engendré. On peut alors prendre un système
dénombrable de rev\^etements pointés $\tilde k_i$ de $\tilde k$, et choisir
des points géométriques sur $X_{\tilde k_i}$ par induction sur $i$.\\
On définit alors $\gtempgeom(X,\bar x)^{\mbb L}=\varprojlim_{k'}
\gtemp(X_{k'},\bar x_{k'})^{\mbb L}$, où
$k'$ parcourt les extensions kummerienne de $k$ dans un point log géométrique $\bar
k$.\\
Soit $\KCovgeom(X)=\injLim \KCov(X_{k'})$ où $k'$ parcourt les extensions
kummeriennes de $k$ dans $\bar k$. C'est la catégorie des revêtements log
géométriques de $X$.\\
On obtient alors une catégorie fibrée $\Dtopgeom \to \KCovgeom(X)$, dont la
fibre en $Y$ est la catégorie des revêtements topologiques de $|\Cgeom(Y)|$.\\
Si $Y\to X$ est un rev\^etement log géométrique, défini sur $k'$,
$\Cgeom(Y_{k'})$ ne dépend pas de $k'$, et l'on obtient donc un foncteur
$\KCovgeom(X)\to\Ke$ qui envoie $Y$ sur $|\Cgeom(Y)|$. Si $\bar x$ est un
système compatible de points, pour n'importe quel revêtement log
géométriquement pointé $(Y,\bar y)$ de $(X,\bar x)$, $\bar y$ définit un
foncteur fondamental $F_{\bar y}$ de
$\Dtopgeom_Y$. Pour tout morphisme $(Y',\bar
y')\to (Y,\bar y)$, les foncteurs correspondant sont canoniquement isomorphes.\\
Alors $\gtempgeom(X,\bar x)^{\mbb L}=\gtemp(\Dtopgeom/\KCovgeom(X), (F_{\bar y}))^{\mbb L}$.


\section{Comparaison pour le groupe fondamental tempéré pro-$(p')$}
Soit $K$ un corps complet \`a valuation discr\`ete. Soient $O_K$ son anneau d'entiers, $k$ sont corps r\'esiduel. $\Spec O_K$ est muni de la log structure fs induite par le mono\"ide $O_K\backslash\{0\}$.\\
Si $\underline X:X\to\cdots\to \Spec(O_K)$ est une log fibration polystable
propre, on veut comparer le groupe fondamental tempéré de la fibre générique
$X_\eta$ (qui est lisse sur $K$) avec le groupe fondamental tempéré de la fibre spéciale muni de sa
log structure naturelle. La théorie de la spécialisation du groupe
fondamental logarithmique nous donne déjà un foncteur de la catégorie des
revêtements két de la fibre spéciale vers la catégorie des revêtements
algébriques de la fibre générique. Pour étendre ceci aux groupes
fondamentaux tempérés, on doit comparer, pour tout revêtement két $T_s$
de la fibre spéciale, l'espace topologique $\C(T_s)$ avec l'espace de
Berkovich du revêtement correspondant $T_{\eta}$ de la fibre générique. Ainsi
on définit, comme dans~\cite{berk2}, une rétraction par déformation forte de
$T_{\eta}^{\an}$ sur un sous-espace canoniquement homéomorphe à
$|\C(T_s)|$. Nous construirons cette rétraction localement pour la
topologie étale, où $T$ a un revêtement galoisien $V'$ par une certaine log
fibration au-dessus d'une extension modérément ramifiée de $O_K$. Alors la
rétraction du tube de $T_s$ est obtenue
en descendant la rétraction du tube de $V'_s$ définie dans
\cite{berk2}. Nous vérifierons alors que la rétraction ne dépend pas du
choix de $V'$. Alors nous pourrons descendre la rétraction définie localement pour la topologie \'etale.

\subsection{Squelette d'un log schéma két au-dessus d'un log schéma pluristable}\label{retraction}

Soit $\underline X:X\to\cdots\to \Spec(O_K)$ une log fibration polystable
au-dessus de
$\Spec(O_K)$.
\begin{prop}\label{skelretract} Soit  $T\to X$ un morphisme két. Notons
  $\fk T_\eta$ la fibre générique, au sens de Berkovich, du complété formel
  de $T$ le long de sa fibre spéciale. Alors, il existe
  une application, fonctorielle en T,
  $|\C(T)_s)|\to \fk T_\eta$, qui identifie $|\C(T_s)|$ avec un sous-ensemble
  $S(T)$ de $\fk T_\eta$ sur lequel $\fk T_\eta$
  se rétracte par déformation forte.\end{prop}
\begin{rem}
$\fk T_\eta$ s'identifie naturellement à un domaine analytique de
$T^{\an}$, et si $T$ est propre (par exemple si $\underline X$ est propre
et $T$ est un revêtement két), alors $\fk T_\eta\to T^{\an}$ est un isomorphisme.
\end{rem}
\dem
Soit $f:T\to X$ un morphisme két.\\

Pour tout $x\in T_{s}$, 
soit $\underline U:U_l\to\cdots\to U_0$ une fibration polystable étale au
dessus de
$\underline X$ telle que $(U_l,x_l)$ soit un voisinage étale de $f(x)$,
et, pour tout $i$, $U_i$ ait une carte exacte  $P_i\to
A_i$ et des morphismes compatibles $P_i\to P_{i+1}$ tels que le morphisme
induit $U_{i+1}\to U_i\times_{\Spec \mbf Z[P_i]}\Spec \mbf Z[P_{i+1}]$ soit
étale strict.\\
On a un voisinage étale
$i:(V,x')\to (T,x)$ de $x$, un morphisme $(p')$-Kummer $P_l\to Q$ tel que
$V\to X$ se factorise à travers un morphisme strict étale $V\to U_l\times_{\Spec
  \mbf Z[P_l]}\Spec \mbf Z[Q]$.\\

 Soit $P_i \to \frac{1}{n} P_i$ l'injection canonique. Alors, par
 définition d'un morphisme $(p')$-kummerien, il existe
$n$ premier à $p$ tel que $P_l\to \frac{1}{n}P_l$ se factorise à travers $P_l\to Q$. Alors
$V$ a un revêtement galoisien két $V'$ qui provient d'une fibration polystable
$\underline U'=V'\to U'_{l-1}\to\cdots\to \Spec O_{K'}$, où $U'_i=U_i\times_{\Spec
  Z[P_i]}\Spec Z[\frac{1}{n}P_i]$ pour $i\leq l$ et $V'=V\times_{\mbf Z[Q]}\mbf
Z[\frac{1}{n}P]$ (d'où un morphisme étale strict $V'\to U'_l$) sur
$O_{K'}$ pour une certaine extension modérément ramifiée
$K'=K[\pi^{1/n}]$ de $K$ .
Appelons $G=(\frac{1}{n}P^{\gp}/Q^{\gp})^\vee$ le groupe de Galois de ce
revêtement két.\\

Notons $\fk U, \fk U_i, \fk V, \fk V'$ les complétés formels de
$U, U_i, V, V'$ le long de la fibre spéciale. Notons alors $\fk V_{\eta}$
la fibre générique de $\fk V$ au sens de Berkovich.\\
La rétraction de $\fk V'_{\eta}$ définie dans le théorème~\ref{berk81} est
$G$-équivariante, donc définit une rétraction de $\fk V_{\eta}$.\\
Soit $S(\ )$ l'image de la rétraction de $(\ )_{\eta}$.
Alors on a un isomorphisme naturel $S(\fk V_{\eta})=G\backslash S(\fk V'_{\eta})\simeq G\backslash
|C(V'_s)|=|G\backslash C(V'_s)|=|C(V_s)|$~(corollaire~\ref{berk85}).\\

Montrons que la rétraction de $\fk U_{\eta}$ que l'on vient de définir ne
dépend pas de $n$. Commençons par le cas d'un morphisme polystable.\\
Soit \[\psi:Z_1=\Spec A[P]/(p_i-\lambda_i)\to Z_2=\Spec A[P]/(p_i-\lambda_i^s) \]
où $P=\mbf N^{|\mbf r|}=\oplus_{(i,j)\in \mbf r}\mbf Ne_{ij}$ et
$p_i=\sum_j e_{ij}$ induit par la multiplication par
$s$ sur $P$, où $s$ est un entier premier à $p$ et où $\lambda\in A$.\\
Soit $\fk G=\prod_i\fGm^{(r_i)}$ où $\fGm^{(r_i)}$ est le noyau de la
multiplication $\fGm^{r_i+1}\to\fGm$. $\fk G$ agit sur $\fk Z_{1}$ et $\fk Z_{2}$. On a $\psi(g\cdot x)=g^s\cdot\psi(x)$.\\
Soit $T_{ij}$ les coordonnées de $G=\fk G_\eta$. Alors
$|T^s_{ij}-1|=|T_{ij}-1|$ si $|T_{ij}-1|<1$ car $s$ est premier à $p$. Ainsi pour $t<1$, $(\
)^s:G\to G$ induit un isomorphisme $(\ )^s:G_t\to G_t$, et $g^s_t=g_t$ (où
$G_t$ est le sous-groupe de $G$ défini par les inégalités $|T_{ij}-1|<t$ et
$g_t$ est son point maximal~; cf \S{}~\ref{berkspaces}).\\
Donc si $t<1$ (et aussi pour $t=1$ par
continuité), \[\psi(x_t)=\psi(g_t\ast
x)=g^s_t\ast\psi(x)=g_t\ast\psi(x)=\psi(x)_t,\]
o\`u $\ast$ est la multiplication d\'ej\`a utilis\'ee dans la discussion suivant le th\'eor\`eme~\ref{berk81}.\\
De plus le diagramme naturel d'homéomorphismes~:
\[\begin{array}{ccc} S(\fk Z_1) & \to & |\C(Z_{1,s})|\\ \dar & & \dar \\
  S(\fk Z_2) & \to & |\C(Z_{2,s})| 
\end{array}
\]
est commutatif.\\

Pour une fibration polystable standard, on obtient le même résultat par
r\'ecurrence sur la longueur, en utilisant que $\psi_n(r_i,t)^{1/s}=\psi_n(r_i^{1/s},t^{1/s})$
(revoir la discussion suivant~\ref{berk81} pour les notations).\\
Plus précisément, supposons qu'on ait le diagramme commutatif suivant~:
\[\xymatrix{B=B'[Y_{ij}]/(Y_{i0}\cdots Y_{in_i}-b_i) & B'\ar[l] \\ 
A=A'[X_{ij}]/(X_{i0}\cdots X_{in_i}-a_i) \ar[u]^\phi & A' \ar[u]^{\phi'} \ar[l]
}\]
où $\phi(X_{ij})=Y_{ij}^s$ et donc $\phi'(a_i)=b_i^s$, et $\tilde\phi':=\Spf\phi':\Spf
B'\to\Spf A'$ est un morphisme két de log fibrations polystables et
supposons par récurrence que l'on sache déjà que $\tilde\phi(x_t)=\tilde\phi(x)_t$.\\
Notons  $\Spf A$
(resp. $\Spf A'$, $\Spf B$, $\Spf B'$) $\fk X$ (resp. $\fk X'$, $\fk Y$, $\fk Y'$).\\
La première partie de la rétraction de $\fk X_\eta^{\an}$ et $\fk
Y^{\an}_{\eta}$ (consistant en la rétraction fibre par fibre) commute avec $\tilde\phi:=\Spf\phi$
d'après le cas précédent. Nous avons donc juste à étudier la deuxième
partie de la rétraction.\\
$\tilde\phi$ induit une application~:
\[\begin{array}{ccc} S_A & = & \{(x,r_{ij})\in(\fk X')_\eta^{\an}\times [0,1]^{[\mbf
    n]}|r_{i0}\cdots r_{in_i}=|a_i(x)|\}\subset\fk X_\eta^{\an} \\
\dar & & \\
S_B & = & \{(y,r_{ij})\in(\fk Y')_\eta^{\an}\times [0,1]^{[\mbf
    n]}|r_{i0}\cdots r_{in_i}=|b_i(y)|\}\subset\fk Y_\eta^{\an}
\end{array}\]
 qui envoie $(x,r_{ij})$ en $(\tilde\phi'(x),r_{ij}^{1/s})$ (remarquons que
 $|a_i(x)|=|b_i(\tilde\phi'(x))|^s$).\\
Alors, si $(x,r_{ij})\in S_A$ (notons $y:=\tilde\phi'(x)$~; par hypothèse
de récurrence, $\tilde\phi'(x_t)=y_t$)
\[\begin{array}{ccl}
\tilde\phi((x,r_{ij})_t) & = &
\tilde\phi((x_t,\psi_{n_i}(r_{ij},|a_i(x_t)|)_k))\\
& = & (y_t,\psi_{n_i}(r_{ij},|a_i(x_t)|)_k^{1/s})\\
& = & (y_t,\psi_{n_i}(r_{ij}^{1/s},|a_i(x_t)|^{1/s})_k)\\
& = & (y_t,\psi_{n_i}(r_{ij}^{1/s},|b_i(y_t)|)_k)\\
& = & (y,r_{ij}^{1/s})_t\\
& = & \tilde\phi(x,r_{ij})_t\end{array}\]
Si maintenant on a un diagramme de log schémas polystables~:
\[\xymatrix{X \ar[r]^g \ar[d]^f & X' \ar[d]^{f'} \\ Y 
\ar[r]^{g'} & Y'
}\]
où $g,g'$ sont étales, $f,f'$ sont két, et l'on sait déjà que $f'$ commute
avec les déformations fortes de $X'$ et $Y'$. Montrons-le pour $f$. Soit
$x\in\fk X_\eta$. Alors
$g'f(x_t)=f'g(x_t)=(f'g(x))_t=(g'f(x))_t=g'(f(x)_t)$. Soit $S_t:=\{y\in \fk
Y_\eta|g'(y)=g'f(x_{t'}), t'\in [0,t]\}={g'}^{-1}(\{(g'f(x))_{t'}\}_{t'\in
  [0,t]})$. Alors, la déformation de $Y$ induit par restriction une
déformation forte de $S_t$ sur ${g'}^{-1}((g'f(x))_t)$ (qui est discret)
grâce à~\ref{berk81}.(ii). En particulier toute composante
connexe de $S_t$ contient un unique élément de ${g'}^{-1}((g'f(x))_t)$. Or
$t'\mapsto f(x_{t'})$ et $t'\mapsto f(x)_{t'}$ sont deux fonctions
continues de $[0,t]$ dans $S_t$ et qui coincident en $0$. Comme $[0,t]$ est
connexe elles ont leurs images dans une même composante connexe et donc
$f(x_t)=f(x)_t$.\\
Donc $f$ commute bien aussi à la déformation forte. Il est clair que $f$
commute aussi aux homéomorphismes $S(\fk X)\to |\C(X_s)|$ et $S(\fk Y)\to
|\C(Y_s)|$. En effet si $x\in S(\fk X)$, les deux images de $x$ dans
$|\C(Y_s)|$ par les deux morphismes sont au-dessus du même point de
$|\C(Y'_s)|$, donc sont soit égaux, soit dans des strates différentes de
$|\C(Y_s)|$. Mais dans les deux cas, l'image de $x$ est dans la strate de
$|\C(Y_s)|$ correspondant à la strate de $Y_s$ dans laquelle se trouve
$\pi(g'f(x))$.\\
Par descente étale, tout morphisme két de fibration polystable commute à la
déformation forte et à l'homéomorphisme entre squelette et réalisation
géométrique du complexe polysimplicial.\\
 
On en déduit en particulier que la déformation forte de $\fk U_\eta$ et
l'isomorphisme $S(\fk U)\to |\C(U_s)|$ ne dépendent pas de $n$.\\

Soit $i:W\to T$ un autre voisinage de $x$ satisfaisant les mêmes propriétés
que $V$, et $W'$ défini de la même façon (on peut supposer grâce à la
remarque précédente que l'on a choisi le même $n$). Ainsi,
$W''=V'\times_{T}W'$ est étale sur $V'$ et $W'$ (les projections
canoniques sont notées $p$ et $p'$).\\
Soit $y\in \fk V'_{\eta}$ et $y'\in \fk W'_{\eta}$ ayant même image dans
$\fk T_{\eta}$. Soit $y''\in \fk W''_{\eta}$ au-dessus de $y$ et $y'$.
Alors, pour tout $t\in [0,1]$, $i(y_t)=i(p(\phi(y''))_t)=i(p(\phi(y'')_t))$
par définition de la rétraction de $\fk V_{\eta}$. En appliquant encore
le théorème~\ref{berk81}.(ii) à $\phi$, on obtient
$i(y_t)=ip\phi(y''_t)$. Par le même argument appliqué à $U'$ et puisque
$ip\phi=i'p'\phi'$, on obtient $i(y_t)=i'(y'_t)$. Ainsi, les rétractions
des différents $\fk V_{\eta}$ sont compatibles et on obtient une rétraction
bien définie de
 $\fk T_{\eta}$ (si $(V_i)$ est un recouvrement étale
 de $T$, la fonction obtenue en recollant les déformations des différents $\fk V_{i,\eta}$
 est continue puisque $\coprod \fk V_{i,\eta}\to \fk T_{\eta}$ est
 quasi-étale et surjectif et donc c'est un quotient topologique).\\

De plus, si $\phi:T_1\to T_2$
est un morphisme két, $\phi(x_t)=\phi(x)_t$. Comme dans le théorème~\ref{berk81}.(vi),
$\phi$ est aussi compatible aux extensions isométriques de $K$.\\

Si $T$ est recouvert par $\widetilde V=\bigcup_i V_i$ tel que tout $V_i$
satisfait aux mêmes propriétés que $V$,
on obtient un isomorphisme \[S(\fk T_{\eta})=\Coker(S(\fk{\widetilde V\times_T\widetilde
  V}_\eta)\rightrightarrows \fk{\widetilde
  V}))=\Coker(\bigcup_{i,j}|\C(V_{i,s}\times_{T_{s}}V_{j,s})|\rightrightarrows
\bigcup_i|\C(V_{i,s})|)=|\C(T_s)|\]
Cet isomorphisme ne dépend pas du choix de $\widetilde V$, et est
fonctoriel en $T$.
\findem
\begin{rem}
En fait, la même preuve permet de construire une déformation forte de $\fk
T_{O_{K'},\eta}$ pour n'importe quelle extension isométrique de $K$ (même à
valuation non discrète), fonctoriellement en $K'$. Si $K'$ contient une
clôture séparable de $K$ (et même une extension maximale modérément
ramifiée), alors $S(T_{O_{K'}})$ est fonctoriellement homéomorphe à $\Cgeom(T/k)$.
\end{rem}

\subsection{Théorème de comparaison des groupes fondamentaux temp\'er\'es}
Soit $K$ un corps complet à valuation discrète. Soit $p$ la caractéristique
résiduelle (qui peut éventuellement être 0).\\
Soit $\underline X:X\to\cdots\to \Spec(O_K)$ une log fibration polystable sur
$\Spec(O_K)$.\\
Comparons le groupe fondamental tempéré de la fibre générique, en tant que
$K$-variété lisse, et le groupe fondamental tempéré de sa fibre spéciale
logarithmique comme défini en~§\ref{tfgsp}.\\

\begin{thm}\label{isomfondtemp}
Soit $\bar x$ un point géométrique de $X^{\an}_\eta$, et soit $\tilde x$
sa log réduction.
On a un morphisme $\gtemp(X^{\an}_\eta,\bar x)^{\mbb L}\to\gtemp(X_s,\tilde
x)^{\mbb
  L}$ qui est un isomorphisme si $p\notin\mbb L$.\end{thm}
Cet isomorphisme est compatible avec les extensions finies de $K$.\\

\dem
On a les deux foncteurs exacts suivants~: \[\mbb L\text{-}\KCov(X)\to\mbb
L\text{-}\Covalg(X_{\eta}),\] qui est une équivalence de catégorie si $p\notin\mbb
L$, et
\[\mbb L\text{-}\KCov(X)\to\mbb L\text{-}\KCov(X_s)\] qui est une équivalence
de catégories (théorème~\ref{orgsp}).\\
Cela nous donne un foncteur exact de $\mbb L\text{-}\KCov(X_s)$ vers
$\mbb L\text{-}\Covalg(X_{\eta})$, qui est une équivalence de catégories si
$p\notin \mbb L$.\\
On a une catégorie fibrée $\Dtopan(X)$ sur $\mbb L\text{-}\KCov(X)$ dont la
fibre en un revêtement két $\mbb L$-fini $T$ est la catégorie des
revêtements topologiques de $T_\eta^{\an}$. On a aussi une catégorie fibrée
$\Dtops(X)$ sur $\mbb L\text{-}\KCov(X)$ obtenue par pullback de
la catégorie fibrée $\Dtop(X_s)\to\mbb L\text{-}\KCov(X_s)$ définie
en~\ref{tfgsp} le long de $\mbb L\text{-}\KCov(X)\to\mbb
L\text{-}\KCov(X_s)$.\\

La proposition~\ref{skelretract} induit une équivalence de catégories fibrées
$\Dtopan(X)\to \Dtops(X)$.\\

Mais $\mbb L\text{-}\KCov(X)\to\mbb
L\text{-}\Covalg(X_{\eta})$ induit un morphisme \[\gtemp(X_{\eta}^{\an})^{\mbb
  L}\to \gtemp(\Dtopan(X)/\mbb L\text{-}\KCov(X))\] qui est un isomorphisme
si $p\notin \mbb L$.\\
De même, $\mbb L\text{-}\KCov(X)\to\mbb
L\text{-}\KCov(X_s)$ induit un isomorphisme \[\gtemp(X_s)^{\mbb
L}\to\gtemp(\Dtops(X)/\mbb L\text{-}\KCov(X)).\]
\findem

\subsection{Théorème de comparaison des groupes fondamentaux temp\'er\'es géométriques}

\begin{thm}\label{isomtempgeom} Il existe un morphisme naturel
\[\gtemp(X_{\bar{\eta}})^{\mbb
  L} \to \gtempgeom(X_{s})^{\mbb L},\]
qui est un isomorphisme si $p\notin \mbb L$.\end{thm}
\dem
On a un foncteur $\KCov(s)\to\Covalg(K)$. En prenant la limite projective
des morphismes du théorème~\ref{isomfondtemp}, on obtient un morphisme
\[\gtemp(X_{\bar{\eta}})^{\mbb
  L}=\varprojlim_{K'}\gtemp(X_{K'})^{\mbb L} \to
\varprojlim_{s'}\gtemp(X_{s'})^{\mbb L}=\gtempgeom(X_{s})^{\mbb L}.\]
Pour montrer que c'est un isomorphisme si $p\notin \mbb L$, compte tenu des
isomorphismes de~\ref{isomfondtemp}, il suffit de montrer que $\gtemp(X_{\bar{\eta}})^{\mbb
  L}=\varprojlim_{K'}\gtemp(X_{K'})^{\mbb L}$ où $K'$ décrit les extensions
modérément ramifiées de $K$ (et non plus toutes les extensions séparable de
$K$). D'après~\cite[prop. 1.15]{kisin}, un revêtement $\mbb L$-fini $T_{\bar
  \eta}$ de $X_{\bar
\eta}$ est défini sur une extension modérément ramifiée $K'$ de $K$, et se
prolonge donc en un revêtement két $T$ de $X_{O_{K'}}$. On sait qu'il
existe une extension két $s''$ de $s'$ (\ie une extension modérément
ramifiée $K''$ de $K'$), telle que pour
toute extension $K'''$ de $K''$, $\C(T_{s'''})\to \C(T_{s''})$ est un
isomorphisme. On déduit de~\ref{skelretract} que $T^{\an}_{K'''}\to
T^{\an}_{K''}$ est une équivalence d'homotopie.\\
On en déduit que $\gtemp(X_{\bar\eta})^{\mbb L}\to
\varprojlim_{K_i}\gtemp(X_{K_i})^{\mbb L}$, où $K_i$ décrit les extensions
modérément ramifiées de $K$ dans $\overline K$ est un isomorphisme.\\
On obtient donc l'isomorphisme voulu en prenant la limite projective des
isomorphismes du théorème~\ref{isomfondtemp}.

\findem

Cet isomorphisme est $\Gal(\bar K,K)$-équivariant (puisque l'isomorphisme
pour chaque extension galoisienne $K_i$ de $K$ est $\Gal(K_i/K)$-équivariante).\\

\begin{rem} Nous pourrions aussi construire l'isomorphisme en prenant la
  limite projective sur toutes les extensions séparables de $K$ en
  remarquant que si $K_1\to K_2$ est totalement sauvagement ramifié, le
  morphisme de log points $s_2\to s_1$ est kuh et donc $\gtemp(X_{s_2})\to\gtemp(X_{s_1})$
  est un isomorphisme.\end{rem}

\section{Complexes classifiants associés et groupe fondamental algébrique
  de la fibre spéciale}

Dans cette section, nous associons à un log-schéma pluristable un complexe classifiant
dont le groupe fondamental tempéré est le groupe fondamental tempéré du
log-schéma de la définition~\ref{defloggft}.

\subsection{Complexe classifiant des strates}
Soit $S$ un log point fs, et soit $k$ son corps sous-jacent.\\ 
Soit $X\to S$ un morphisme strictement plurinodal de log schémas.\\
Soit $Y$ l'union disjointe des adhérences des strates de $X$ (chaque
adhérence de strate étant muni de la structure de sous-log-schéma réduit et
strict de $X$). Soit
$C(Y/X)$ le log-schéma simplicial dont la composante de degré $n$ est
$Y^n=Y\times_{X}\cdots\times_{X}Y$ (les morphismes de log-schémas étant
tous stricts ici, la formation des produits fibrés commute à la prise du
schéma sous-jacent) et les flèches sont évidentes. On obtient
alors un ensemble simplicial $\g0(C(Y/X))$ en prenant en chaque degré les
composantes connexes de $Y^n$. Remarquons que les composantes connexes de
$Y^n$ s'identifient à des adhérences de strates de $X$, ce qui donne un
foncteur $(\Delta/\g0(C(Y/X)))^{\op}\to \Str X$.\\
Considérons donc le complexe classifiant $\mcal G_{\Str}(X)$ sur $(\Str X)^{\op}$ dont la fibre en
une strate $x$ est la catégorie des revêtements két de l'adhérence de la
strate $\overline{\{ x\}}$. On a
alors un foncteur naturel $\KCov X\to \Ba(\mcal G_{\Str} X)$ qui, pour
tout objet $x$ de $(\Str X)^{\op}$, associe à un
revêtement de $X$ sa restriction à $\overline{\{ x\}}$.\\
On obtient donc un complexe classifiant $\mcal G_{\DD}(X)$ sur
$(\Delta/\g0(C(Y/X)))^{\op}$ par pullback. Si on tronque le complexe
classifiant à l'ordre 3, on a par définition une équivalence entre les
donnée de descentes de $Y/X$ et les revêtements de ce dernier complexe
classifiant (il n'est évidemment pas nécessaire de tronquer à l'ordre 2
pour obtenir la même catégorie).\\

Or, d'après~\cite[prop. 3.2.18]{stix}, le morphisme $Y \to X$ étant un
épimorphisme propre, c'est un morphisme de descente effectif (dans la
catégorie des fs log-schémas) pour les revêtements két. On a
donc une équivalence de catégories entre les revêtements két de
$X$ modérément ramifiés et les revêtements két de $Y$ munis
de données de descentes.\\
Le foncteur composé $\KCov(X)\to \Ba(\mcal G_{\Str}(X))\to\Ba(\mcal
G_{\DD}(X))$ est une équivalence de catégories.\\

Pour montrer que $\KCov(X)\to \Ba(\mcal G_{\Str}(X))$ est aussi une
équivalence de catégories, il suffit donc de montrer que $\Ba(\mcal
G_{\Str}(X))\to\Ba(\mcal G_{\DD}(X))$ est pleinement fidèle. Mais pour
cela il suffit d'appliquer la proposition~\ref{propcart}.(ii) à
$F:(\Delta/\g0(C(Y/X)))^{\op}\to \Str(X)$.\\
En effet, soit $x$ un objet de $\Str(X)$ (c'est-à-dire une strate de $X$), $x\!\!\dar\!\!F$ est la sous
catégorie pleine de $(\Delta/\g0(C(Y/X)))^{\op}$ constituée des objets $y$ tels
que $x\subset\overline{F(y)}$. Tout objet $y$ de cette catégorie admet une
flèche vers  un objet de $\g0(C(Y/X))_1$, correspondant à une composante
irréductible $Y_1$ de $X$ contenant $x$. Mais si $Y_1$ et $Y_2$ sont deux telles
composantes irréductibles de $X$, soit $Y_{12}$ la composante connexe de
$Y_1\cap Y_2$ qui contient $x$ et notons aussi $Y_{12}$ l'objet
correspondant de $(\Delta/\g0(C(Y/X)))^{\op}$, qui est aussi dans
$x\!\!\dar\!\!F$.\\
On a donc dans $x\!\!\dar\!\!F$ le diagramme suivant :
\[Y_1\gets Y_{12} \to Y_2,\]
ce qui montre que $x\!\!\dar\!\!F$ est bien connexe, et donc $\KCov(X)\to
\Ba(\mcal G_{\Str}(X))$ est bien une équivalence de catégories.\\

\subsection{Complexe classifiant polysimplicial d'un schéma strictement polystable}

Si $\underline X$ est une fibration strictement polystable, d'espace global $X$, on a
aussi un foncteur $F:(\mbf \Lambda/\mbf C(\underline X))^{\op}\to(\Str
X)^{\op}$ qui induit donc par pullback de $\mcal G_{\Str}(X)$ un
complexe classifiant $\mcal G_{\mbf C}(\underline X)$ sur $(\mbf
\Lambda/\mbf C(\underline X))^{\op}$.\\
Si $X$ est strictement polystable, $\mbf C(X)$ est un complexe
polysimplicial libre. Donc si $x$ est une strate de $X_s$, soit $y$ un
polysimplexe nondégénéré de $\mbf C(X)$ correspondant à $x$, alors
$x\!\!\dar\!\!F$ est la sous-catégorie pleine de $(\mbf
\Lambda/\mbf C(\underline X))^{\op}$ constituée des objets qui sont
l'image d'un morphisme de source $y$, qui est équivalente à $(\mbf
\Lambda/[\mbf n_y])^{\op}$ car $\mbf C(X)$ est libre. Or $\Lambda/[\mbf n_y]$
admet un objet final, donc on peut directement
appliquer~\ref{propcart}.(iii) à  $(\mbf \Lambda/\mbf C(\underline X))^{\op}\to(\Str
X_s)^{\op}$, et donc $\Ba(\mcal G_{\Str}(X))\to\Ba(\mcal G_{\mbf
  C}(\underline X))$ est une équivalence de catégories.\\

\subsection{Complexe classifiant $l$-polysimplicial d'une fibration
  strictement polystable}

Pour une fibration polystable $\underline X$ de longueur $l$, on peut
construire un objet $l$-polysimplicial $C(\underline X)$ de
$(\mbf \Lambda^l)^\circ\Ens$ tel que l'image de $C(\underline X)$ par le
foncteur $\sq^{l-1}_!(\mbf \Lambda^l)\circ\Ens\to\mbf \Lambda^\circ\Ens$,
qui prolonge en commutant aux limites directes $\sq^{l-1}:\mbf\Lambda^l\to\mbf\Lambda$, qui à $[\mbf
  n_1]\times\cdots\times[\mbf n_l]$ associe $[\mbf n_1]\sq\cdots\sq [\mbf
  n_l]$, est $\mbf C(\underline X)$ (en particulier, on a un foncteur $D_0:\mbf \Lambda^l/C(\underline
X_s)\to\mbf \Lambda/\mbf C(\underline X)$).\\
Supposons $C(\underline X)$ déjà construit pour les fibrations de
longueur $l$, et soit $\underline X=(X_{l+1}\to\underline X_l)$ une
fibration de longueur $l+1$. Considérons le foncteur :
\[D':(\mbf \Lambda^l/C(\underline X_{l}))^{\op}\stackrel{D_0^{\op}}{\to}(\mbf \Lambda/\mbf
C(\underline X_{l}))^{\op}\stackrel{D}{\to} \mbf \Lambda ^\circ\Ens.\]
On peut lui associer l'objet $C(\underline X_{l})\times D'$ dont
l'ensemble des objets en $([\mbf n_1],\cdots,[\mbf n_{l+1}])$ est
  \[\coprod_{x\in C(\underline X_{l})_{[\mbf n_1]\times\cdots\times[\mbf
        n_l]}}(D'_x)_{[\mbf n_{l+1}]}\]
L'image directe par
$\sq^{l-1}\times\id_{\mbf\Lambda}:\mbf\Lambda^{l+1}\to\mbf\Lambda\times\mbf\Lambda$
est $\mbf C(\underline X_{l})\times D$ par hypothèse de récurrence, et on
vérifie facilement que l'image directe par $\mbf \Lambda\times\mbf \Lambda\to\mbf \Lambda$
de $C\times D$ est $C\sq D$.\\
Par pullback de $\mcal G_{\mbf C}(\underline X)$ le long de $D_0^{\op}$, on en
déduit un complexe classifiant $\mcal G_{C}(\underline X)/(\mbf
\Lambda^l/C(\underline X_{l}))^{\op}$.\\
De m\^eme qu'on avait associ\'e \`a un morphisme k\'et $S\to X_l$ un ensemble polysimplicial $\C(S)$, on peut construire un ensemble $l$-polysimplicial $C(S)$ tel que $\sq^{l-1} C(S)=\C(S)$.\\

Montrons par récurrence sur la longueur de $\underline X$ que la conclusion
de \ref{propcart}.{iii} est vraie pour $(\mbf \Lambda^l/C(\underline X))^{\op}\to(\Str
X)^{\op}$. Supposons-la vraie pour $\underline X_l$ et montrons-la pour
$\underline X=(X_{l+1}\to\underline X_l)$.\\
On a le diagramme commutatif suivant :
\[\xymatrix{(\mbf \Lambda^{l+1}/C(\underline X))^{\op} \ar[r] \ar[d]
  & \Str(X_{l+1})^{\op} \ar[d]\\(\mbf \Lambda^l/\mbf C(\underline X_{l}))^{\op}
  \ar[r] & \Str(X_{l})^{\op}},\]
où les flèches verticales sont cofibrantes. La flèche du bas vérifie la
  conclusion de~\ref{propcart}.(iii), et les fibres aussi d'après le cas
  où $X$ est strictement polystable (car si $x$ est un objet de $\mbf
  C(\underline X_{l})$,  $(\mbf \Lambda^{l+1}/C(\underline X))^{\op}_x\to(\Str
X)^{\op}_x$ s'identifie à $(\mbf \Lambda^l/\mbf C((X_{l+1})_x)^{\op}\to
  (\Str (X_{l+1})_x)^{\op}$ par définition de $C(\underline
  X)$).\\
Donc la conclusion de~\ref{propcart}.(iii) est aussi vraie pour $(\mbf
  \Lambda^{l+1}/C(\underline X))^{\op}\to(\Str X)^{\op}$.\\

Donc \[\KCov(X)\to \Ba(\mcal G_{C}(\underline X))\]
est une équivalence de catégories.\\

\subsection{Descente et complexe classifiant $l$-polysimplicial d'une
  fibration polystable}\label{cpgpsch}
Soit $\underline X$ une $l$-fibration polystable sur $S$. Soit $\underline
X'\to\underline X$ un morphisme étale surjectif de $l$-fibrations
polystables avec $\underline X'$ strictement polystable. Notons
$\underline X'^{(n)}=\underline X'\times_{\underline
  X}\cdots\times_{\underline X} \underline X'$\\
Alors $C(\underline X)=\Coker(C(\underline X'^{(2)})\rightrightarrows
C(\underline X'))$ ne dépend pas
du choix de $\underline X'$.\\
Le foncteur $(\mbf\Lambda^l/C(\underline
X'^{(n)}))^{\op}\to(\mbf\Lambda^l/C(\underline X))^{\op}$ est en fait une
catégorie cofibrée en catégories discrètes, et en recollant les $\underline
X'^{(n)}$ pour $n\geqslant 1$, on obtient ainsi une catégorie cofibrée
$\mcal C_{\DD}(\underline X'/\underline X)$ en ensembles
simpliciaux sur $(\mbf\Lambda^l/C(\underline X))^{op}$ (correspondant à la
catégorie des objets d'un foncteur $\Delta\times\mbf\Lambda^l\to\Ens$).\\
Sur $\mcal C_{\DD}(\underline X'/\underline X)$, on peut mettre une
structure de complexe classifiant $\mcal G_{DD}(\underline X'/\underline
X)$.\\

En descendant le long de $F_{\underline X'}:\mcal C_{\DD}(\underline X'/\underline
X)\to\mbf\Lambda^l/C(\underline X)$, on obtient un complexe classifiant
$\mcal G(\underline X'/\underline X)=F_{\underline X'*}\mcal G_{DD}(\underline X'/\underline
X)$ sur $\mbf\Lambda^l/C(\underline X)$.\\
On a un isomorphisme
\[\Ba(\mcal G(\underline X'/\underline X))\simeq\Ba(\mcal G_{\DD}(\underline
X'/ \underline X))\] d'après la proposition~\ref{propdesc}.\\

Mais on peut aussi ``calculer'' $\Ba(\mcal G_{\DD}(\underline X'/\underline
X))$ en le descendant le long du foncteur cofibrant $F':\mcal C_{\DD}(\underline X'/\underline
X)\to\Delta^{\op}$.\\
Or la fibre en $[n]$ du complexe classifiant $\mcal G_{\DD}(\underline X'/\underline
X)$ est $\mcal G_C(\underline X'^{(n)})$. Donc $F'_*(\mcal G_{\DD}(\underline X'/\underline
X))$ est l'objet simplicial des données de descentes de $\underline
X'/\underline X$, et donc par descente étale des revêtements két,
$\Ba(\mcal G_{\DD}(\underline X'/\underline X))\simeq\KCov(X)$. Donc
:
\[\Ba(\mcal G(\underline X'/\underline X))\simeq\KCov(X).\]

Si $\underline X'$ et $\underline X''$ sont maintenant deux $l$-fibrations
étales surjectives sur $X$, on peut construire comme
précédemment un objet bicosimplicial en complexes $l$-polysimpliciaux
classifiant dont la fibre en $([m],[n])$ est le complexe $l$-polysimplicial
classifiant $\mcal G_C(\underline X'^{(m)}\times_{\underline X}\underline
X''^{(n)})$. Notons-le $\mcal G_{\DD}(\underline X',\underline X''/\underline X)/\mcal
C_{\DD}(\underline X',\underline X''/\underline X)$.\\
On a de plus un diagramme commutatif de foncteurs :
\[\xymatrix{ & \mcal C_{\DD}(\underline X',\underline X''/\underline X)
  \ar^{G_{\underline X'}}[dl] \ar_{G_{\underline X''}}[dr] & \\ \mcal
  C_{\DD}(\underline X'/\underline X) \ar^{F_{\underline X'}}[dr] & & \mcal
  C_{\DD}(\underline X''/\underline X) \ar_{F_{\underline X''}}[dl] \\
   & (\mbf\Lambda^l/C(\underline X))^{\op} & }\]
Or en descendant $\mcal G_{\DD}(\underline X',\underline X''/\underline X)$
  le long de $G_{\underline X'}$, on obtient $\mcal G_{\DD}(\underline
  X'/\underline X)$.\\
En effet la fibre en un objet $x$ de $\mcal G_{\DD}(\underline
  X'/\underline X)$ correspondant à une strate $x$ de $X'^n_l$ est le
  complexe classifiant des revêtements két du log-schéma simplicial ($\bar
  x\times_{X_l}X''^{(n)}_l$) ; or comme $\bar x\times_{X_l}X''_l\to \bar x$ est
  un morphisme de descente pour les revêtements két
  d'après~\cite[prop. 3.2.18]{stix}, on en déduit le résultat voulu.\\
On a évidemment le résultat analogue pour $X''$. En appliquant le résultat
  de~\ref{propdesc}, on en déduit que $\mcal G_{\DD}(\underline
  X'/\underline X)$ ne dépend pas  essentiellement pas de $\underline
  X'$.\\
On le notera donc $\mcal G_{\underline X}$.

\subsection{Comparaison des groupes fondamentaux temp\'er\'es}\label{cclogschtemp}
On a d\'ej\`a une \'equivalence de cat\'egories~:
\[\Ba(\mcal G_{\underline X})\simeq\KCov(X).\]
Soit $S$ un rev\^etement k\'et de $X$. Soit $T$ l'objet correspondant de $\Ba(\mcal G(\underline X'/\underline X))$.
Alors le complexe $l$-polysimplicial $C(S)$ associ\'e \`a $S$ est canoniquement isomorphe \`a $C(X)\sq D_T$.\\
En effet, il suffit de construire l'isomorphisme localement pour la
topologie \'etale sur $X$, car si $X'\to X$ est \'etale surjectif,
$C(S)=\Coker(C(S_{X''})\rightrightarrows C(S_{X'}))$ et $C(X)\sq
D_T=\Coker(C(X'')\sq D_{T_{X''}}\rightrightarrows C(X')\sq D_{T_{X'}})$. On
peut donc supposer $X$ strictement polystable, o\`u l'isomorphisme est \'evident par d\'efinition de $C(S)$.\\
L'isomorphisme $C(S)\simeq C(X)\sq D_T$ induit une \'equivalence de cat\'egories \[\Covtop(C(X)\sq D_T)\to\Covtop(C(S))\] 2-fonctorielle en $S$, d'o\`u une \'equivalence des cat\'egories fibr\'ees correspondantes sur $\Ba(\mcal G_{\underline X})$ et $\KCov(X)$. En consid\'erant les champs associ\'es et en prenant les sections globales, on en d\'eduit une \'equivalence de cat\'egories \[\Btemp(\mcal G_{\underline{X}})\simeq\Covtemp(X),\]
qui induit un isomorphisme
\[\gtemp(\mcal G_{\underline X})\to\gtemp(X).\]

\chapter{Cospécialisation et groupe fondamental tempéré}

Pour un morphisme propre de schémas $f: X\to Y$ de fibres géométriquement
connexes et une spécialisation  $ 
\bar y_1 \to \bar y_2$ de points géométriques de $Y$, A. Grothendieck  
a construit un morphisme de spécialisation  $\ga(X_{\bar y_1})\to  
\ga(X_{\bar y_2})$. Le théorème de spécialisation de Grothendieck 
dit que ce morphisme est surjectif si $f$ est séparable et induit un isomorphisme  
entre les quotients premiers à $p$ si $f$ est lisse  
(o\`u $p$ d\'esigne la caractéristique de $\bar y_2$, qui est éventuellement $0$), cf. \cite[cor. X.2.4,
cor. X.3.9]{sga}.\\
En géométrie analytique complexe, un morphisme propre et lisse est,
localement sur $Y$, une fibration triviale de variétés différentielles
réelles~; en particulier toutes les fibres sont homéomorphes et
ont donc des groupes fondamentaux (topologiques) isomorphes.\\

Peut-il y avoir un résultat analogue pour le groupe fondamental tempéré (ou
plus raisonnablement $(p')$-tempéré)~?\\
Considérons le cas d'une famille de courbes elliptiques. Le groupe
fondamental tempéré géométrique d'une courbe elliptique est $\widehat{\mbf
  Z}^2$ si elle a bonne réduction et $\widehat{\mbf Z}\times \mbf Z$ si
c'est une courbe de Tate. En particulier, dans une famille modulaire
de courbes elliptiques avec structure de niveau\footnote{pour éviter les
  champs. Cependant, les morphismes de cospécialisation que nous
  construirons seront locaux pour la topologie étale de la fibre spéciale
  de la base. Ainsi le fait de travailler avec un champ de Deligne-Mumford n'est
  pas vraiment un problème.}, le groupe fondamental tempéré ne peut être
constant (ni même le groupe fondamental $(p')$-tempéré).\\
De plus, dans le contexte $p$-adique, si $E_1$ est une courbe elliptique \`a mauvaise réduction et $E_0$ est une courbe elliptique \`a r\'eduction g\'en\'erique (et donc bonne réduction), il ne peut pas y avoir de morphisme
$\gtemp(E_0)^{(p')}\to\gtemp(E_1)^{(p')}$ qui induise le morphisme de spécialisation de
Grothendieck sur les complétés profinis, bien que  la réduction du point
correspondant à
$E_1$ se spécialise sur la réduction du point correspondant à $E_0$. Il ne
peut donc pas y avoir de spécialisation raisonnable.\\
Si n\'eanmoins, $\eta_1$ et $\eta_2$ sont des points g\'eom\'etriques de
l'espace analytique de modules tels que la réduction de $\eta_1$ se spécialise sur la
réduction de $\eta_2$,
alors $E_{\eta_1}$ a "meilleure réduction" que $E_{\eta_2}$ et
il existe des morphismes
$\gtemp(E_{\eta_2})\to\gtemp(E_{\eta_1})$ qui induisent un isomorphisme sur
les complétés profinis. C'est pourquoi nous chercherons à construire des morphismes de
\emph{cosp\'ecialisation} (plut\^ot que de sp\'ecialisation) du groupe fondamental tempéré.\\

Un problème qui apparaît tout de suite est qu'il n'y a pas en général de spécialisation non triviale
entre les points d'un espace de Berkovich~: par exemple un espace de
Berkovich séparé a un espace topologique sous-jacent séparé, et donc s'il y
a une spécialisation (pour la topologie de Berkovich, la topologie étale\dots)
entre deux points géométriques d'un espace de Berkovich séparé, les deux
points géométriques doivent avoir même point sous-jacent. Ainsi nous fixerons un mod\`ele de $Y$ sur $O_K$ ayant d'assez bonnes propriétés, et les "morphismes de sp\'ecialisation" entre points g\'eom\'etriques de $Y^{\an}$ que nous consid\'ererons seront en fait des morphismes de sp\'ecialisation entre leurs r\'eductions.\\

Nous allons étudier la situation suivante. Soit $K$ un corps à valuation
discrète, $O_K$ son anneau des entiers, $k$ son corps résiduel et $p$
la caractéristique de $k$ (qui peut être $0$). Soit $\mbb L$ un ensemble de
nombres premiers ne contenant pas $p$. Soit $X\to Y$ un morphisme pluristable de log schémas sur $O_K$.\\
Si $\eta_1$ est un point de $Y_{\tr}^{\an}$ tel que $\mcal H(\eta_1)$ soit
à valuation discrète, on sait décrire le
groupe fondamental $(p')$-tempéré de $X_{\eta_1}$  en terme de $X_{s_1}$
où $s_1$ est la log réduction de $\eta_1$. Pour que cette réduction existe,
on supposera $Y$ propre sur $O_K$~; sinon il faudrait se limiter à considérer
les points  $\eta_1$ dans le tube de la fibre spéciale $Y$ et remplacer
$Y_{tr}^{\an}$ par son intersection avec ce tube (les résultats qui
suivront sont encore vrai dans ce contexte-ci plus général).\\

On peut alors reformuler notre problème de cospécialisation en terme de la
fibre spéciale de $Y$. Nous supposerons alors $Y$ log lisse sur $O_K$ (ceci munit $Y$ d'une
stratification canonique). Nous ferons également une hypothèse sur le comportement combinatoire des
fibres géométriques $X\to Y$. Plus précisément, on supposera que l'ensemble
polysimplicial associé aux fibres géométriques de $X\to Y$ est
intérieurement libre~; c'est par exemple le cas si $X\to Y$ est strictement
polystable ou si $X\to Y$ est de dimension relative 1.\\
Nous prouverons
alors le résultat suivant~:
\begin{thm}\label{sp} Soit $\eta_1$ et $\eta_2$ deux points
  de $Y_{\tr}^{\an}$ à valuation discrète, et soit $\bar \eta_1,\bar
  \eta_2$ des points géométriques au-dessus de ces points. 
Soit $\bar s_2\to\bar s_1$ une spécialisation entre leurs log réductions. Alors
il existe un morphisme naturel dit de \emph{cospécialisation}\index{Morphisme!de cosp\'ecialisation} \[\gtemp(X_{\bar\eta_1})^{\mbb
    L}\to\gtemp(X_{\bar\eta_2})^{\mbb L},\] qui est un isomorphisme si
  $\bar s_1$ et $\bar s_2$ sont dans la même strate de
  $Y$.\end{thm}
Si l'on ne suppose plus que les fibres géométriques de $X\to Y$ ont des
ensembles polysimpliciaux intérieurement libres, il y a encore un morphisme
de cospécialisation si $\bar s_2$ est le point générique d'une strate.\\

La première chose dont nous aurons besoin pour construire le morphisme de
cospécialisation pour le groupe fondamental $(p')$-tempéré est un morphisme
de spécialisation pour les groupes fondamentaux logarithmiques pro-$(p')$ de
$X_{\bar s_1}$ et $X_{\bar s_2}$ qui étend tout revêtement két $(p')$ de
$X_{s_1}$ sur un voisinage két de $s_1$ (si l'on a 
un tel morphisme de spécialisation, en les comparant aux groupes
fondamentaux algébriques de $X_{\bar\eta_1}$ et $X_{\bar\eta_2}$ et en
utilisant le théorème de spécialisation de Grothendieck, nous obtiendrons
facilement que c'est un isomorphisme). Ce morphisme de spécialisation se
déduit aisément de~\cite{org2} si $s_1$ est un point strict de $Y$. Ainsi
nous étudierons l'invariance du groupe fondamental log géométrique par
changement de point base.\\
Nous nous int\'eresserons alors au cas, plus simple, où $X\to Y$ est une
courbe relative. Dans ce contexte nous n'aurons pas à utiliser les
ensembles polysimpliciaux de Berkovich. Nous obtiendrons également un
résultat meilleur que dans le cas de dimension supérieure~: nous n'aurons
pas à supposer $Y$ log lisse et nous n'aurons pas non plus à supposer $X\to
Y$ vertical.\\
Nous reviendrons alors au cas général. Nous étudierons le comportement combinatoire
d'un revêtement két vis-à-vis de la cospécialisation. Des morphismes de
cospécialisations d'ensembles polysimpliciaux des fibres d'une fibration polystable
sont déjà donnés par Berkovich dans~\cite[cor. 6.2]{berk2} quand
$Y$ est plurinodal et $s_1,s_2$ sont les points génériques topologiques
d'une strate de $Y$ sans hypothèse de propreté. La construction s'étend
facilement, localement pour la topologie étale, à notre situation si l'on
suppose encore que $s_2$ est le point générique d'une strate. Pour obtenir
un morphisme de cospécialisation d'ensembles polysimpliciaux des fibres géométriques, 
nous devrons prouver qu'après une nouvelle localisation két en $\bar s_1$,
les strates de $X_{s_2}$ dont l'adhérence rencontre $X_k$ sont
géométriquement connexes. Cela proviendra du fait qu'après localisation supplémentaire,
 l'adhérence de la clôture de ces strates sont plates sur leur image dans
 $Y$ et ont des fibres géométriques réduites. On descend alors ces
morphismes de cospécialisation que l'on avait localement pour la topologie
étale. Dans le cas propre initial, le morphisme de cospécialisation est un
isomorphisme si $s_1$ et $s_2$
sont dans la même strate de $Y$ et l'ensemble polysimplicial de $X_{\bar s_2}$
est intérieurement libre. Ces morphismes de cospécialisation commutent
avec les revêtements két, et donc nous donneront les morphismes de
cospécialisation voulus.\\

\section{Spécialisation des groupes fondamentaux log géométriques}

Etudions la spécialisation des groupes fondamentaux log géométriques (\ie
la limite projective des groupes fondamentaux logarithmiques par extension
két du log point de base). Notre principal résultat
sera l'invariance du groupe fondamental log géométrique d'un log schéma fs
 $X$, saturé et de type fini sur un log point fs $s$ de corps séparablement clos,
par changement de base fs $s'\to s$ qui est un isomorphisme sur le schéma sous-jacent.
L'hypothèse implique que le changement de base $X_{s'}\to X$ dans la cat\'egorie des log sch\'emas fs induit un isomorphisme $\mring{X_{s'}}\to\mring X$ entre les sch\'emas sous-jacents. En travaillant localement sur la topologie étale
de $\mring X$, nous nous réduirons au cas où $\mring X$ est strictement local
et hensélien, cas dans lequel le groupe fondamental log géométrique peut
être explicitement décrit en terme du morphisme de monoïdes
$\overline M_{X}\to\overline M_s$.\\
Combinant ce résultat d'invariance par changement de base
du groupe fondamental log géométrique pro-$(p')$ et la spécialisation
du groupe fondamental log géométrique pro-$(p')$ dans le cas strict (\cite{org2}),
nous obtiendrons que si $X\to S$ est un morphisme propre, log lisse et saturé, et
si $s_2,s_1$ sont des points fs de $S$ et $\bar s_2\to\bar s_1$ est une
spécialisation de points log géométriques de $S$ au-dessus de $s_2$ et $s_1$, alors
il y a un morphisme de spécialisation
\[\ggeom(X_{s_2})^{(p')}\to\ggeom(X_{s_1})^{(p')}\] (ce sera le seul résultat
dont nous aurons besoin par la suite).\\
 
\begin{lem}\label{lemchgmtbase} Soit $s'\to s$ un morphisme strict de log
  points tels que $\mring s'$ et $\mring s$ soient des points géométriques
  de caractéristique $p$. Soit $X\to s$ un morphisme
 de log schémas fs tel que $\mring X\to\mring s$ soit de type fini.\\
Alors $F:\KCov(X)^{(p')}\to \KCov(X_{s'})^{(p')}$ est une équivalence de catégories.\end{lem}
\dem
Si $T$ est un revêtement két connexe de  $X$, $\mring{T\times_ss'}\to\mring
T\times_{\mring s}\mring s'$ est un isomorphisme puisque $s'\to s$ est
strict. $\mring T\times_{\mring s}\mring s'$ est aussi connexe, donc $F$
est bien pleinement fidèle.\\
Comme on sait déjà que $F$ est pleinement fidèle pour tout $X$, et comme
les recouvrements étales stricts sont de descente effective pour les
revêtements két, on peut prouver le résultat localement pour la topologie
étale, et ainsi supposer que $X$ a une carte globale
$X\to\Spec \mbf Z[P]$.\\
Soit $S'$ un revêtement két de $X_{s'}$. Alors il existe un morphisme
$(p')$-Kummer de monoïdes $P\to Q$ tel que \[S'_Q:=S'\times_{\Spec \mbf
  Z[P]}\Spec\mbf Z[Q]\to X_{s',Q}:=X_{s'}\times_{\Spec \mbf
  Z[P]}\Spec\mbf Z[Q]\] soit strict étale (et surjectif).\\
Mais, puisque $\mring X_{s',Q}\to\mring X_Q\times_{\mring s}\mring s'$ est un
isomorphisme de schémas, $\Covalg(\mring X_{s',Q})\to\Covalg(\mring X_Q)$
est une équivalence de
catégories~(\cite[cor 4.5]{org}). Ainsi, il existe un revêtement étale strict
$S_Q$ de $X_Q$ (et donc $S_Q\to X$ est un revêtement két) tel que $S'_Q$ soit
$X_{s',Q}$-isomorphe à $S_Q\times_ss'$. Donc $F$ est une équivalence de catégories.
\findem

Soit $X\to s$ un morphisme de log schémas fs, où $s$ est un log point fs. Soit
$\bar x$ un point log géométrique de $X$ et soit $\bar s$ son image dans
$s$.\\
Définissons le \emph{groupe fondamental log géométrique} de $X$ en $\bar x$
comme étant
\[\ggeom(X/s,\bar x)=\Ker(\glog(X,\bar x)\to\glog(s,\bar s)).\]
Si $(t,\bar t)\to (s,\bar s)$ est un revêtement két connexe galoisien pointé de $(s,\bar
s)$, notons $X_t=X\times_st$ et $\bar x_t=(\bar x,\bar t)$. On a alors la suite exacte:
\[1\to \glog(X_t,\bar x_t)\to\glog(X,\bar x)\to\Gal(t/s),\]
et le morphisme de droite est surjectif si $X_t$ est connexe (nous dirons que $X$ est
\emph{log g\'eom\'etriquement connexe} si $X_t$ est connexe pour tout
revêtement két connexe $t$ de $s$).\\
En prenant la limite projective de la suite exacte précédente quand $(t,\bar
t)$ parcourt la catégorie dirigée des revêtements galoisiens pointés de
$(s,\bar s)$, on obtient un isomorphisme canonique $\ggeom(X/s,\bar
x)=\varprojlim_{(t,\bar t)} \glog(X_t,\bar x_t)$.\\
En fait, si $\tilde t\to t$ est le sous-schéma réduit de $t$ muni de
la log structure image inverse, le morphisme $\glog(X_{\tilde t},\bar
x_t)\to\glog(X_t,\bar x_t)$ est un isomorphisme. On peut donc remplacer $X_t$
par $X_{\tilde t}$ dans la limite projective précédente.\\

Si l'on a un carré commutatif de log schémas fs pointés:
\[\xymatrix{(X',\bar x') \ar[r]^{\phi} \ar[d] & (X,\bar x) \ar[d]\\
(s',\bar s')\ar[r]^{\psi} & (s,\bar s)}\]
où $s'$ et $s$ sont des log points, on obtient un diagramme commutatif de
groupes profinis~:
\[\xymatrix{\glog(X',\bar x') \ar[r] \ar[d] & \glog(X,\bar x) \ar[d]\\
\glog(s',\bar s')\ar[r] & \glog(s,\bar s)}\]
En prenant le noyau des flèches verticales, on obtient un morphisme
$\ggeom(X'/s',\bar x')\to\ggeom(X/s,\bar x)$, qui est fonctoriel en $(\phi,\psi)$.\\
On a aussi, par définition du pro-revêtement pointé universel de
$(s',\bar s')$, un morphisme canonique de pro-revêtements két pointés de
$(s',\bar s')$
\[\varprojlim_{(t',\bar t')} (t',\bar t')\to \varprojlim_{(t,\bar t)}
\psi^*(t,\bar t),\]où $(t',\bar t')$ parcourt les revêtements két
galoisiens connexes de $(s',\bar s')$ et $(t,\bar t)$ parcourt les
revêtements két galoisiens
connexes de $(s,\bar s)$, d'où un morphisme de pro-log schémas fs pointés \[\varprojlim_{(t',\bar t')} (t',\bar
t')\to \varprojlim_{(t,\bar t)} (t,\bar t).\]
Cela induit un morphisme de pro-log schémas fs pointés
\[\varprojlim (X'_{t'},\bar x'_{t'})\to \varprojlim (X_t,\bar x_t),\]
d'où un morphisme de groupes profinis \[\varprojlim \glog(X'_{t'},\bar
x'_{t'})\to\varprojlim \glog(X_t,\bar x_t),\]
tel que le carré commutatif de groupes profinis soit commutatif:
\[\xymatrix{\varprojlim\glog(X'_{t'},\bar x'_{t'})\ar[r]\ar@{=}[d] & \varprojlim\glog(X_t,\bar
x_t)\ar@{=}[d]\\ \ggeom(X'/s',\bar x')\ar[r] & \ggeom(X/s,\bar x)}\]

Consid\'erons maintenant $s'\to s$ un morphisme de log points fs, tel que le
morphisme sous-jacent de schémas $\mathring s'\to
\mathring s$ soit un isomorphisme de points géométriques, et soit $X\to s$
un morphisme saturé de log schémas fs avec $\mathring X$
nœth\'erien et $\mathring{X}\to\mring{s}$ (géométriquement) connexe. Puisque
$X\to s$ est saturé, $X$ est log géométriquement connexe.\\
Posons $X':=X\times_ss'$. Soit $\bar x'$ un point log géométrique de $X'$ et soient $\bar
x$, $\bar s'$ et $\bar s$ l'image de $\bar x'$ dans $X$, $s'$ et $s$ respectivement.\\
On a un diagramme commutatif
\[\xymatrix{\glog(X',\bar x') \ar[r] \ar[d] & \glog(X,\bar x) \ar[d]\\
\glog(s',\bar s)\ar[r] & \glog(s,\bar s)}\]

\begin{thm}\label{chgmtbase1} Le morphisme $\ggeom(X'/s',\bar x')\to\ggeom(X/s,\bar x)$ est un isomorphisme.\end{thm}
\dem
Soit $(s_i,\bar s_i)_{i\in I}$ un système cofinal de revêtements két
connexes galoisiens pointés de $(s,\bar s)$ et soit $\tilde s_i$ le sous
schéma réduit de $s_i$ muni de la log structure image inverse. Posons $(X_i,\bar x_i):=(X\times_s\tilde s_i,(\bar x,\bar s_i))$.\\
Soit $(s'_j,\bar s'_j)_{j\in J}$ un système cofinal de revêtements két
galoisiens connexes pointés de $(s',\bar s')$ (et soit $\tilde s'_j$ le
sous-schéma réduit de $s'_j$ muni de la log structure image inverse. Posons $(X'_j,\bar x'_j)=(X\times_{s'}s'_j,\bar
x'\times_{\bar s'}\bar s'_j)$.\\
Il faut prouver que
 \[\varprojlim_j\glog(X'_j,\bar
 x'_j)\to\varprojlim_i\glog(X_i,\bar
x_i)\]
est un isomorphisme.\\

Si $Y$ est un log schéma fs, notons $\epsilon_Y:Y_{\ket}\to \mathring
Y_{\etale}$  le morphisme usuel de topoi du topos két de $Y$
au topos étale du sous-schéma sous-jacent à $Y$.\\
$\tilde s_i$ a même schéma sous-jacent que $s$, donc $\mring X_i\to\mring
X$ est un isomorphisme de schémas puisque $X\to s$ est saturé.
De même 
$\mring X'_j\to \mring X'$ est un isomorphisme, et sont isomorphes à $\mring X$.\\
Plus précisément, pour tout $i$, il existe $j_0$ tel que pour $i\geq i'$ et $j\geq j'\geq j_0$, on
ait un diagramme 2-commutatif~:
\[\xymatrix {
    X'_{j,\ket} \ar[rr] \ar[dd] \ar[dr] && X'_{j',\ket} \ar[dr] \ar[dd] |\hole \\
    & X_{i,\ket} \ar[rr] \ar[dd] && X_{i',\ket} \ar[dd] \\
    \mathring X'_{j,\etale} \ar[rr] |\hole \ar[dr] && \mathring X'_{j',\etale} \ar[rd] \\
    & \mathring X_{i,\etale} \ar[rr] && \mathring X_{i',\etale} \\
  }\]
où tous les morphismes de schémas du carré du bas sont des isomorphismes.\\
Si $G$ est un groupe fini, on a un diagramme 2-commutatif de champs sur $\mathring
X_{\etale}$:
\[\xymatrix{\epsilon_{X'_{j'}*}\Tors_{X'_{j',\ket}}(G) \ar[r]\ar[d] & \epsilon_{
    X'_{j}*}\Tors_{X'_{j,\ket}}(G) \ar[d]\\
  \epsilon_{X_{i'}*}\Tors_{X_{i',\ket}}(G) \ar[r] &
  \epsilon_{X_{i}*}\Tors_{X_{i,\ket}}(G)} \]
Plus précisément, notons $IJ=I\coprod J$ muni de la relation d'ordre suivante~:
\begin{itemize}
\item la restriction de l'ordre à $I$ (resp. $J$) est l'ordre usuel,
\item si $i\in I$ et $j\in J$, elors $i\leq j$ si et seulement si il existe
  un (nécessairement unique)
morphisme de log schémas fs pointés $(s'_j,\bar s'_j)\to (s_i,\bar s_i)$
qui rende le carré \[\begin{array}{ccc}(s'_j,\bar s'_j) & \to &
  (s_i,\bar s_i)\\ \dar & & \dar \\ (s',\bar s') & \to & (s,\bar
  s)\end{array}\] commutatif,
\item si $j\in J$ et $i\in I$, alors $j\nleq i$.\end{itemize}

Notons aussi $IJ$ la catégorie correspondante.
On a un catégorie fibrée sur $IJ^{\op}\times \mring X_{\etale}$, dont la
fibre en $(i,U)$ est $\epsilon_{X_{i}*}\Tors_{X_{i,\ket}}(G)(U)$ et la
fibre en $(j,U)$ est $\epsilon_{X'_{j}*}\Tors_{X'_{j,\ket}}(G)(U)$.\\

En prenant la limite inductive pour $i\in I$, qui est filtrante, on
obtient, grâce à~\cite[I.1.10]{giraud}, une catégorie fibrée sur
$\mathring X_{\etale}$ dont la fibre en $U$ est $\injLim_i
\epsilon_{\widetilde X_{i}*}\Tors_{\widetilde X_{i,\ket}}(G)(U)$. Notons $\injLim_i \epsilon_{\widetilde X_{i}*}\Tors_{\widetilde
  X_{i,\ket}}(G)$ le champ associé à cette catégorie fibrée (et faisons la
même chose pour $X'$ et $J$).\\
Puisque $\mring X$ est supposé nœtherien (et donc $\mring X_{\et}$ est un
topos cohérent), le foncteur de  $\injLim_j \Tors(G,
X'_{j,\ket})$ vers la catégorie des sections globales de $\injLim_j\epsilon_{X'_j*}\Tors_{
  X'_{j,\ket}}(G)$ est une équivalence de catégories. De plus \[\injLim_j \Tors(G,
X'_{j,\ket})\simeq\injLim_j \Tors(G,\glog(X'_j,\bar x'_j)-\Ens)\simeq\Tors(G,\varprojlim \glog(X'_j,\bar x'_j)-\Ens).\]
Les diagrammes 2-commutatifs induisent un morphisme de champs:
\begin{equation}\label{morphisme_champ}\injLim_j\epsilon_{X'_j*}\Tors_{X'_{j,\ket}}(G)\to\injLim_i\epsilon_{X_i*}\Tors_{X_{i,\ket}}(G)\end{equation}

Nous n'avons donc qu'à prouver que~(\ref{morphisme_champ}) est une
équivalence de champs, ce qui se prouve fibre à fibre (puisque les points
de $\mring X_{\et}$ forment un système conservatif). Ceci ne dépend plus
des points de base initiaux que l'on avait dans l'énoncé du théorème~; nous
pouvons donc les oublier, ainsi nous utiliserons les mêmes notations pour
d'autres points.\\
Soit $x$ un point de $\mring X_{\et}$, $\bar x'$ un point de $X'_{\ket}$ au
dessus de $x$
d'image $\bar x$ dans $X_{\ket}$, et soit $V_x$ la catégorie des voisinages
étales de $x$ dans $\mring X$. On a alors (par cohérence du morphisme de
topoi $\epsilon_{X_i}$, comme dans~\cite[dem of 2.4]{org2}):
\[\begin{array}{rl}\injLim_{U\in
    V_x}\injLim_i\epsilon_{X_i*}\Tors_{X_{i,\ket}}(G)(U) &
  =\injLim_i\injLim_{U\in V_x}\epsilon_{X_i*}\Tors_{X_{i,\ket}}(G)(U)\\
 & =\injLim_i\Tors(G,X(x)_{i,\ket})\\
 & =\Tors(G,\varprojlim_i\glog(X(x)_{i,\ket},\bar x\times s_i)-\Ens),\end{array}\]
et l'on a un résultat similaire pour $X'$.\\
On a donc juste à prouver que \[\varprojlim_j\glog(X'(x)_{j,\ket},\bar
x'\times_{\bar s'} \bar s'_j)\to\varprojlim_i\glog(X(x)_{i,\ket},\bar
x\times_{\bar s}\bar s_i)\] est un isomorphisme.\\
Nous nous sommes donc réduit au cas où $\mathring X$ est un schéma
nœthérien, strictement local et hensélien.
Mais alors, pour un schéma nœthérien,
strictement local et hensélien $\mathring X$,  on a d'apr\`es~\cite[prop. 3.1.11]{stix} 
\[\varprojlim_i\glog(X_{i,\ket},\bar x\times_{\bar s}\bar s_i)=\varprojlim_i \overline M_{X_i}^{\gp \vee}\otimes \widehat{\mbf Z}^{(p')}=\Coker(\overline M_s^{\gp}\to \overline M_X^{\gp})^{\vee}\otimes\widehat{\mbf Z}^{(p')},\]
et l'on a un résultat similaire pour $X'$.
$\Coker(\overline M_s^{\gp}\to \overline M_X^{\gp})\to\Coker(\overline
M_{s'}^{\gp}\to \overline M_{X'}^{\gp})$ est un isomorphisme.\\
On obtient donc le résultat voulu.
\findem

Supposons maintenant seulement que $(s',\bar s')\to (s,\bar s)$ soit un morphisme de
log points, que $Y\to s$ soit un morphisme saturé et que $X\to Y$ soit un morphisme két
avec $\mring X$ de type fini sur $s$.
\begin{cor}\label{loginvar} Le morphisme de groupes
  profinis \[\ggeom(X/s,\bar x)^{(p')}\to\ggeom(X'/s',\bar x')^{(p')}\] est
  un isomorphisme.\end{cor}
\dem
Quitte à remplacer $s$ (resp. $s'$) par le sous-schéma réduit d'un
revêtement két connexe de $s$ (resp. $s'$), on peut supposer que $X\to s$ est saturé
($\mring X\to\mring s$ sera encore de type fini).\\
Si $(t,\bar t)\to (s,\bar s)$ est un revetement étale strict, alors $\ggeom(X_t/t,\bar
x_t)\to\ggeom(X/s,\bar x)$ est un isomorphisme. Ainsi, en notant $s_0$ la
clôture séparable de $s$ et en prenant la limite projective sur les
revêtements étales stricts pointés de $s$ (puisque
$\glog(X_{s_0})=\varprojlim\glog(X_t)$, où $t$ parcours les revêtements
étales stricts pointés de $s$), on obtient $\ggeom(X_{s_0}/s_0,\bar
x_0)\to\ggeom(X/s,\bar x)$ est un isomorphisme. On peut alors supposer que $\mring s$ et $\mring s'$ sont des points géométriques. Considérons le log schéma fs $s''$ de schéma sous-jacent $\mring
s'$ et dont la log structure est l'image inverse de la log structure de
$s$.\\
Ainsi, on a des morphismes $s'\to s''\to s$, où $s'\to s''$ est un
isomorphisme de schémas sous-jacents et $s''\to s$ est strict.\\
Grâce au lemme~\ref{lemchgmtbase}, $\glog(X_{s''})^{(p')}\to
\glog(X)^{(p')}$ et $\glog(s'')^{(p')}\to\glog(s)^{(p')}$ sont des isomorphismes. Ainsi,
\[\ggeom(X_{s''}/s'')\to\ggeom(X/s)\] est un isomorphisme.\\
Par~\ref{chgmtbase1}, $\ggeom(X_{s'}/s')\to\ggeom(X_{s''}/s'')$ est aussi
un isomorphisme.
\findem

Rappelons que si $S$ est un schéma strictement local de point fermé $s$ et si
$X$ est un log schéma fs connexe tel que
  $\mring X$ est propre sur $S$, alors
\[\KCov(X)\to\KCov(X_s)\]est une équivalence de cat\'egories (théorème~\ref{orgsp}).\\
Soit $X\to S$ un morphisme propre et saturé de log schemes fs,
  et soit $Y\to X$ un revêtement két. Soit $s$ et $s'$ deux points de $S$
  et supposons que l'on ait une spécialisation $\bar
    s'\to\bar{s}$ (où $\bar s$ et
  $\bar s'$ sont des points log géométriques sur $s$ et
  $s'$).\\
Soit $Z$ le localisé strict de $S$ à $s$ muni de la log structure image
inverse, et soit $z$ son point fermé, muni de la log structure image inverse.\\
On a les morphismes suivants (définis à automorphismes intérieurs près):
\[\ggeom(Y_{s}/s)^{(p')}\stackrel{\simeq}{\to}\ggeom(Y_z/z)^{(p')}\stackrel{\simeq}{\to}\ggeom(Y_Z/Z)^{(p')}\leftarrow\ggeom(Y_{s'}/s')^{(p')}\]
où les deux premiers morphismes sont des isomorphismes d'après le
corollaire~\ref{loginvar} et le théorème~\ref{orgsp}.\\
\begin{cor}\label{logsp}On a un morphisme de spécialisation
\[\ggeom(Y_{s'}/s')^{(p')}\to\ggeom(Y_{s}/s)^{(p')}\] qui se factorise à travers $\ggeom(Y_Z/Z)^{(p')}$.\end{cor}

\section{Cas des courbes}
\label{cospcourbes}
\subsection{D\'efinitions}
\subsubsection{Graphes}
Rappelons qu'un \emph{graphe} $\mbb G$ consiste en la donn\'ee d'un ensemble $\mcal E$ d'"ar\^etes", d'un ensemble 
$\mcal V$ de "sommets", et pour tout $e\in \mcal E$, d'un ensemble de branches $\mcal B_e$ de
cardinal 2 et d'une fonction $\psi_e:\mcal B_e\to\mcal V$.\\
On peut, de fa\c con \'equivalente remplacer la donn\'ee des ar\^etes et des branches de chaque ar\^ete par la donn\'ee de l'ensemble de toutes les branches $\mcal B=\coprod_e\mcal B_e$ muni d'une involution
$\iota$ sans point fixe et d'une fonction $\psi:\mcal B\to \mcal V$. $E$ est
alors l'ensemble des orbites de branches pour $\iota$.\\
Un \emph{vrai morphisme de graphe} $\phi:\mbb G\to \mbb G'$ est donn\'e par une fonction
$\phi_{\mcal E}:\mcal
E\to \mcal E'$, une fonction $\phi_{\mcal V}:\mcal V\to\mcal V'$ et pour tout
$e\in \mcal E$ une bijection $\phi_e:\mcal B_e\to\mcal B'_{\phi_{\mcal
    E}(e)}$ telles que le diagramme suivant commute~:
\[\xymatrix{\mcal B_e \ar[d] \ar[r] & \mcal B'_{\phi_{\mcal E}(e)} \ar[d]\\
  \mcal V \ar[r] & \mcal V'.}\]
Cependant la cosp\'ecialisation topologique des courbes semistables sera donn\'ee par des applications entre graphes qui ne sont pas en g\'en\'eral des vrais morphismes.\\
Un \emph{morphisme g\'en\'eralis\'e de graphe} $\phi:\mbb G\to\mbb G'$ est donn\'e par:
\begin{itemize}
\item une fonction $\phi_{\mcal V}:\mcal V\to\mcal V'$,
\item une fonction $\phi_{\mcal E}:\mcal E\to\mcal E'\coprod\mcal V'$,
\item pour tout $e\in \mcal E$ tel que $\phi_{\mcal E}(e)\in\mcal
  E'$, une bijection $\phi_e:\mcal B_e\to\mcal B'_{\phi_{\mcal
    E}(e)}$ qui fait commuter le diagramme \'evident (le m\^eme que pour
les vrais morphismes),
\item pour tout $e\in\mcal E$ tel que $\phi_{\mcal E}(e)\in\mcal V'$ et
  toute branche $b\in \mcal B_e$ de $e$, $\phi_{\mcal
    V}(\psi_e(b))=\phi_{\mcal E}(e)$.
\end{itemize}
On peut remplacer les deux derni\`eres donn\'ees par une fonction $\phi_{\mcal B}:\mcal
B\to\mcal B'\coprod\mcal V'$ telle que si $\phi_{\mcal B}(b)\in\mcal B'$
then $\phi_{\mcal B}(\iota(b))=\iota'(\phi_{\mcal B})$ et si $\phi_{\mcal
B}(b)\in\mcal V'$, alors $\phi_{\mcal B}(\iota(b))=\phi_{\mcal B}(b)$.\\
En particulier, un vrai morphisme est un morphisme g\'en\'eralis\'e. Il y a une notion \'evidente de composition des morphismes g\'en\'eralis\'es de graphes.\\
On obtient ainsi une cat\'egorie $\Graph$ de graphes pour laquelle les morphismes sont des vrais morphismes de graphes et une
categorie $\GenGraph$ de graphes pour laquelle les morphismes sont les morphismes g\'en\'eralis\'es de graphes.\\

Il y a un foncteur r\'ealisation g\'eom\'etrique $|\ |:\GenGraph\to\Top$
qui envoie un graphe $\mbb G$ vers \[|\mbb G|:=\Coker(\coprod_{b\in\mcal B} \pt_{1,b}\amalg\pt_{2,b}
\rightrightarrows\coprod_{v\in\mcal V} \pt_v\amalg\coprod_{b\in\mcal B}[1/2,1]_b),\]
o\`u \begin{itemize}\item la fl\`eche du haut envoie~: 
\begin{itemize}\item $\pt_{1,b}$ sur $1/2\in [1/2,1]_b$ \item
$\pt_{2,b}$ sur $1\in [1/2,1]_b$,\end{itemize} \item la fl\`eche du bas envoie~: \begin{itemize}\item$\pt_{1,b}$ sur
$1/2\in [1/2,1]_{\iota(b)}$\item $\pt_{2,b}$ sur $\pt_{\psi(b)}$.\end{itemize}\end{itemize}
Si $\phi:\mbb G\to\mbb G'$ est un morphisme g\'en\'eralis\'e de graphes, $|\phi|$ est obtenu
en envoyant \begin{itemize}\item$\pt_v$ sur $\pt_{\phi_{\mcal V}(v)}$,\item $[1/2,1]_b$ sur
$[1/2,1]_{\phi_{\mcal B}(b)}$ si $\phi_{\mcal B}(b)\in \mcal B'$ (par l'identit\'e
de $[1/2,1]$),\item
$[1/2,1]_b$ sur
$\pt_{\phi_{\mcal B}(b)}$ si $\phi_{\mcal B}(b)\in \mcal V'$.\end{itemize}

\subsubsection{Courbes semistables}
\begin{dfn}\label{courbessemistables}
Un morphisme de sch\'ema $X\to S$ est une \emph{courbe strictement semistable} (resp. \emph{semistable}) si pour tout point g\'eom\'etrique $x$ de $X$, d'image $s$ dans $S$, il existe un voisinage \'etale $U$ de $s$ et un voisinage de Zariski (resp. \'etale) $V$ de $x$ dans $X\times_SV$ tels que $V\to U$ se factorise \`a travers un morphisme \'etale $V\to\Spec A[X,Y]/(XY-a)$ ou $V\to \Spec A[X]$.
\end{dfn}
De m\^eme, il y a un pendant logarithmique \`a cette d\'efinition.
\begin{dfn}
Un morphisme $X\to S$ de log sch\'emas fs est une \emph{courbe log strictement semistable} (resp. \emph{log semistable} si pour tout point g\'eom\'etrique $x$ de $X$, d'image $s$ dans $S$, il existe un voisinage \'etale $U$ de $s$, un voisinage de Zariski (resp. \'etale) $V$ de $x$ dans $X\times_SV$ et une carte $U\to \Spec P$ de $U$ telle que l'une des propri\'et\'es suivantes soit v\'erifi\'ee:
\begin{itemize}
\item $V\to U$ est un courbe lisse stricte,
\item $V\to U$ se factorise \`a travers un morphisme \'etale  $V\to U\times_{\Spec\mbf Z[P]}\Spec\mbf
    Z[Q]$ avec $Q=(P\oplus<u,v>)/(u\cdot v=p)$ avec $p\in P$ et où la log
    structure de $V$ est celle associée à $Q$,
\item $V\to U$ se factorise \`a travers un morphisme \'etale $V\to U\times_{\Spec\mbf Z[P]}\Spec\mbf
    Z[P\oplus\mbf N]$ et où la log structure de $V$ est celle associée à
    $P\oplus\mbf N$.
\end{itemize}
\end{dfn}

Remarquant que, pour cette d\'efinition, une courbe log semistable n'est pas n\'ecessairement un morphisme semistable de log sch\'emas, car la log structure n'est pas n\'ecessairement suppos\'ee verticale ici.\\
Une courbe log semistable est log lisse et satur\'ee (\^etre une courbe log semistable \'equivaut en fait \`a \^etre log lisse, satur\'e et purement de dimension relative 1).\\

Le morphisme de sch\'ema sous-jacent $\mring X\to \mring S$ \`a une courbe log semistable (resp. strictement log semistable) est une courbe semistable (resp. strictement semistable). En particulier, si $\mring S$ est un point g\'eom\'etrique, on peut associer \`a $X$ un graphe $\mbb G(X)$ de la fa\c con
suivante~: les sommets sont les composantes irr\'eductibles de $\mcal X_s$,
les ar\^etes sont les points doubles. Les deux branches d'une ar\^ete $e$ aboutissent aux sommets correspondant
aux composantes irr\'eductibles contenant le point double correspondant \`a $e$.\\

Soit $X\to S$ une courbe propre log semistable et soit $X'\to X$ un rev\^etement két. Alors pour tout point log g\'eom\'etrique $\bar s$ de $S$, il existe un voisinage k\'et $U$ de $\bar s$ tel que $X'_U\to U$ soit satur\'e. Alors
$X'_U\to U$ est aussi une courbe log semistable.\\
Le morphisme $\mring{\mcal X}'_{\bar s}\to\mring{\mcal X}_{\bar s}$ induit un vrai
morphisme $\mbb G(X_{\bar s})\to\mbb G(X_{\bar s})$ de graphes.

Soit $K$ un corps complet non archim\'edien de corps r\'esiduel s\'eparablement clos. Soit $O_K$ l'anneau des entiers de $K$. Soit $X\to
O_K$ une courbe semistable propre avec fibre g\'en\'erique lisse, il existe un plongement canonique $|\mbb G(X)|\to
X^{\an}_\eta$ qui est une \'equivalence d'homotopie. Il est compatible aux extensions isom\'etriques de $K$.\\
De plus, si $U$ est un ouvert de Zariski dense de $X_{\eta}$, $|\mbb
G(X)|$ s'envoie dans $U^{\an}$ et $|\mbb G(X)|\to U^{\an}$ est encore une \'equivalence d'homotopie.\\
Si $X\to O_K$ est une courbe log semistable et $X'\to X$ est un morphisme két 
tel que $X'$ soit encore une courbe log semistable, le diagramme suivant est commutatif~:
\[\begin{array}{ccc} |\mbb G(X')| & \to & {X'}_{\eta}^{\an}\\ \dar & & \dar\\
|\mbb G(X)| & \to & X_{\eta}^{\an}\end{array}\]

\subsection{Cospécialisation topologique des courbes semistables}

Soit $f:X\to Y$ une courbe semistable, et soit $\bar y_2\to \bar y_1$ une
spécialisation de points géométriques de $Y$. Dans cette section, nous
définirons un morphisme de cospécialisation  entre graphes associ\'es $\mbb G(X_{\bar y_1})\to \mbb G(X_{\bar y_2})$.\\

Si $x$ est un point de $X$, nous noterons $X(x)$ le localisé de $X$ en $x$.
\begin{lem}Supposons $Y$ strictement local de point fermé $y_1$, et $X\to
  Y$ strictement semistable. Soit $x$
  un nœud ou un point générique de $X_{y_1}$. Alors $X(x)_{y_2}$ est soit
  contenu dans le lieu lisse d'une composante géométriquement irréductible
  (qu'on notera $F(x)$) de
  $X_{y_2}$ ou contient un unique point double (qu'on notera $F(x)$) de $X_{y_2}$, qui est
  un point rationnel.\\
\end{lem}
\dem
Quitte à remplacer $Y$ par un sous-schéma fermé , on peut supposer que $Y$ est intègre et que
$y_2$ est le point générique de $Y$.\\
\begin{enumerate}[(i)]
\item
Si $x$ est dans le lieu lisse de $X_{y_1}$, $X\to Y$ est lisse en $x$, et
$X(x)_{y_2}$ est géométriquement connexe par 0-acyclicité locale des
morphismes lisses.
\item
Si $x$ est un nœud, on peut supposer que $Y=\Spec A$ et $f$ se factorise
à travers un morphisme étale $X\to\Spec B$ avec $B=A[u,v]/(uv-a)$ et
$a(y_1)=0$.\\
Si $a=0$, posons $Z:=X\times_{\Spec B}\Spec A$ où $g:B\to A$ est défini par
$g(u)=g(v)=0$ ($Z$ est le sous-schéma fermé réduit de $X$ défini par les
nœuds; en particulier $Z_{y_2}$ est l'union de tout les nœuds de
$X_{y_2}$). $Z\to Y$ est étale et donc $Z(x)\to Y$ est un isomorphisme
puisque $Y$ est strictement local. Ainsi
$Z(x)_{y_2}$ est juste un point rationnel $F(x)$.\\
Si $a\neq 0$, $X_{y_2}$ est lisse. Le morphisme $X\to Y$ est séparable (\ie
plat à fibre géométrique réduite) et donc localement
0-acyclique. $X(x)_{y_2}$ est donc géométriquement connexe.
\end{enumerate}
\findem
Puisque $F(x)$ est géométriquement irréductible, nous écrirons aussi $F(x)$ pour
la composante irréductible ou le nœud correspondant de $X_{\bar y_2}$.\\

Si $\phi:X'\to X$ est un morphisme fini et ouvert de courbes strictement semistables sur
$Y$ qui, fibre à fibre, envoie nœud sur nœud, alors $\phi F'=F\phi$. En effet $\phi(x)$ est dans l'adh\'erence de $\phi F'(x)$, donc $F\phi(x)$ est dans l'adh\'erence de $\phi F'(x)$. Il suffit donc de v\'erifier que si $F\phi(x)$ est un point double, $\phi F'(x)$ aussi. Supposons que $F \phi(x)$ soit un point double de $X_{y_2}$. Soient $z_1$ et $z_2$ les deux points g\'en\'eriques des composantes irr\'eductibles de $X_{y_2}$ dont l'adh\'erence contient $F \phi(x)$ (et donc aussi $\phi(x)$). Comme $\phi$ est ouvert, il existe $z'_1$ et $z'_2$ dans $X(x)_{y_2}$ tels que $\phi(z'_1)=z_1$ et $\phi(z'_2)=z_2$. Donc $X(x)_{y_2}$ ne peut pas \^etre contenu dans une seule composante irr\'eductible de $X_{y_2}$, et donc $F'(x)$ est un point double de $X'_{y_2}$. Par hypoth\`ese, $\phi F'(x)$ est bien un point double de $X_{y_2}$.

\begin{cor}\label{cosptopcourbes}
Il existe un unique morphisme généralisé de graphes \[\psi:\mbb G(X_{\bar
  y_1})\to\mbb G(X_{\bar y_2})\] qui est
\begin{itemize}
\item fonctoriel pour les morphismes étales $X'\to X$,
\item compatible aux changements de bases $Y'\to Y$,
\item tel que si $f:X\to Y$ est strictement semistable et $Y$ est
strictement local de point spécial $\bar y_1$, $\psi(x)=F(x)$ pour tout
nœud ou point générique $x$.
\end{itemize}
\end{cor}
\dem
Soit $f:X\to Y$ une courbe strictement semistable, et soit $\bar y_2\to
\bar y_1$ une spécialisation de points géométriques de $Y$.\\
On peut construire un morphisme généralisé $\psi_{X/Y}:\mbb
  G(X_{y_1})\to\mbb G(X_{y_2})$ en remplaçant $Y$ par son localisé strict
  à $\bar y_1$ et en posant $\psi(x):=F(x)$. Evidemment, si $b$ est une branche
  de $e$ dans $\mbb G(X_{y_1})$ qui aboutit en $v$, $\psi(e)\subset\psi(v)$,
  d'où $\psi(b)$ (qui est bien défini puisque $X_{\bar y_2}$ est
  strictement semistable).\\
Ce morphisme généralisé est clairement compatible aux morphismes étales de courbes
strictement semistables.\\

Si $X\to Y$ est une courbe semistable, on peut choisir un recouvrement étale 
$X'\to X$ tel que $X'\to Y$ soit strictement semistable. Posons $X'':=X'\times_XX'$.\\
On a un diagramme commutatif~:
\[\begin{array}{ccccc} \mbb G(X''_{\bar y_1}) & \rightrightarrows & \mbb
  G(X'_{\bar y_1}) & \to & \mbb G(X_{\bar y_1})\\ \dar & & \dar & & \\
\mbb G(X''_{\bar y_2}) & \rightrightarrows & \mbb
  G(X'_{\bar y_2}) & \to & \mbb G(X_{\bar y_2})\end{array}\]
Il existe un unique morphisme g\'en\'eralis\'e de graphes $\psi:\mbb G(X_{\bar
  y_1})\to\mbb G(X_{\bar y_2})$ rendant le diagramme commutatif.
\findem

Ce morphisme généralisé de graphes n'est en général pas un vrai
morphisme. En effet, si $X$ est localement isomorphe à $\Spec
A[u,v]/(uv-a)$ avec $a(\bar y_1)=0$ et $a(\bar y_2)\neq 0$, l'arête de
$\mbb G(X_{\bar y_1})$ correspondant au point double $u=v=0$ s'envoie sur
un sommet de $\mbb G(X_{\bar y_2})$.\\
 
Les morphismes de cospécialisation de graphes sont aussi compatibles
avec les morphismes ouverts et finis de courbes semistables qui envoient
fibre à fibre les nœuds sur les nœuds (cela d\'ecoule de la compatibilit\'e avec $F$).\\

De plus si $\bar y_3\to\bar y_2\to\bar y_1$ sont des spécialisations,
$\psi_{\bar y_3\to\bar y_1}=\psi_{\bar y_3\to\bar y_2}\psi_{\bar
y_2}\psi_{\bar y_1}$.

On veut savoir quand ce morphisme généralisé de graphes est un isomorphisme.

\begin{prop}
Si $\psi:\mbb G(X_{\bar y_1})\to\mbb G(X_{\bar y_2})$ est un vrai morphisme
de graphes et $f$ est propre, alors $\psi$ est un isomorphisme.
\end{prop}
\dem
On peut supposer que $Y=\Spec A$ est strictement local et intègre de point spécial
$y_1$ et de point générique $y_2$.\\
L'hypothèse signifie que localement pour la topologie étale de la fibre
spéciale (et donc sur
$X$ par propreté), $X$ est isomorphe à $\Spec A[u,v]/uv$ ou est lisse. En effet si $X$ est localement isomorphe \`a $\Spec A[u,v]/uv-a$ avec $a\notin A^*$, la fibre sp\'eciale a un point double, correspondant \`a une ar\^ete $e$. Alors $\psi(e)$ est une ar\^ete, correspondant a un point double et donc $\Spec \Frac(A)[u,v]/uv-a$ est non lisse, et donc $a=0$.\\

Soit $Z\subset X$ le lieu non lisse de $X\to Y$, muni de la structure de
sous-schéma réduit. $Z\to Y$ est
étale (comme on peut le voir localement sur la topologie étale de $X$), et
propre. On obtient ainsi
que $F$ induit une bijection entre les nœuds de $X_{\bar y_1}$ et ceux de
$X_{\bar y_2}$.\\
Soit $\widetilde X$ l'éclaté de $X$ le long de $Z$. 
Quand $X=\Spec A[u,v]/(uv)$, $Z$ est défini par l'idéal engendré par $u$
et $v$, et $\widetilde X=\Spec A[u]\coprod \Spec A[v]$.\\
Ainsi en regardant localement sur $X$ pour la topologie étale, on voit que
$\widetilde X$ est lisse sur $Y$, et que $\widetilde X_y$ est simplement le
normalisé de $X_y$.\\
Puisqu'on suppose $X\to Y$ propre, $\widetilde X\to Y$ est lisse et
propre, donc sa factorisation de Stein est étale sur $Y$ et donc induit une
bijection entre les composantes connexes de $\widetilde X_{\bar y_1}$ et
celles de $\widetilde X_{\bar
  y_2}$, et donc la fonction entre les composantes irréductibles de $\widetilde
X_{\bar y_1}$ et celles de $\widetilde X_{\bar y_2}$ est aussi une bijection.
\findem

\begin{prop}\label{invstratescourbes}
Soit $f:X\to Y$ une courbe log semistable, et $\bar y_2\to
\bar y_1$ une spécialisation de points log géométriques.
Supposons que $\overline M_{\bar y_1}\to\overline M_{\bar y_2}$ est un isomorphisme.
Alors $\psi:\mbb G(X_{y_1})\to\mbb G(X_{y_2})$ est un vrai morphisme de graphes. \end{prop}
\dem
On peut supposer $Y$ strictement local, intègre de point générique $y_2$~:
$Y=\Spec A$, et soit $P\to A$ une carte.\\

Pour montrer que c'est un vrai morphisme, il suffit de montrer que $\psi(e)$
est une arête si $e$ est une arête de $\mbb G(X_{y_1})$. Cette propriété
n'est pas modifiée par un morphisme étale, et donc on peut supposer $X=\Spec A\otimes_{\mbf
  Z[P]}\mbf Z[Q]$ avec $Q=(P\oplus<u,v>)/(u\cdot v=p)$ et $p\in P$, tel que
l'image de $p$ dans $M_{\bar y_1}$ soit non inversible. Alors l'image de $p$
dans $M_{\bar y_2}$ est non inversible et donc $X=\Spec A[u,v]/(uv=0)$, ce
qui donne le résultat voulu.\\
\findem

\subsection{Cospécialisation topologique et revêtements két}
Soit $f:X\to Y$ une courbe log semistable et propre. Soit $\bar y_2\to\bar y_1$
une spécialisation de points log géométriques.\\
Soit $S_0\to X_{\bar y_1}$ un revêtement két log géométrique.\\
Alors il s'étend au-dessus d'un voisinage két $U$ de $\bar y_1$ en un revêtement két
$S\to X_U$. Quitte à localiser encore pour la topologie két, on peut même
supposer que $S\to X_U$
est une courbe log semistable.\\
On peut donc utiliser le corollaire~\ref{cosptopcourbes} et obtenir un morphisme
généralisé de graphes $\mbb G(S_{\bar y_1})\to\mbb G(S_{\bar y_2})$. Cela
ne dépend pas du choix de $U$ ni de $S$ puisque deux extensions de $S_0$
sur un voisinage két de $\bar y_1$ sont isomorphes sur un voisinage két
plus petit et cet isomorphisme est unique à localisation két supplémentaire
près.\\
Si $S'_0\to S_0$ est un revêtement két, le diagramme suivant est commutatif~:
\[\begin{array}{ccc} \mbb G(S'_{\bar y_1}) & \to  & \mbb G(S'_{\bar y_2})\\
\dar & & \dar \\ \mbb G(S_{\bar y_1}) & \to & \mbb G(S_{\bar
  y_2})\end{array}\]
Si $\overline M_{Y,y_1}\to\overline M_{Y,y_2}$ est un isomorphisme,
$\overline M_{U,y_1}\to\overline M_{U,y_2}$ est encore un isomorphisme,
donc on peut encore appliquer la proposition~\ref{invstratescourbes} à $S$.\\
 
\subsection{Cospécialisation du groupe fondamental tempéré}
Supposons $p\notin\mbb L$.
Soit $Y\to O_K$ un morphisme de fs log schémas et $X\to Y$ une courbe
propre log semistable.\\
Soit $Y_{\tr}$ (resp. $U:=X_{\tr}$) le lieu ouvert de $Y$ (resp. $X$) où la log
structure est triviale
($Y_{\tr}\subset Y_\eta$). Soit $\fk Y$ le complété de $Y$ le long de sa
fibre spéciale.\\
Soit $\mcal Y=\fk Y_\eta\cap Y_{\tr}^{\an}$, c'est un sous-domaine
analytique de $Y^{\an}$.\\

Rappelons qu'une spécialisation (két) $\bar x_1\to\bar x_2$ de points log
géométrique d'un log schéma fs $X$ est un morphisme entre les log strictes
localisations $X(\bar x_1)\to X(\bar x_2)$ (cf. \S{}~\ref{revket}).
\begin{dfn} On note $\Pt^{\an}(Y)$ la catégorie dont les objets sont les points
géométriques $\bar y$ de $Y_{\tr}^{\an}$ tels que $\mcal H(y)$  soit à valuation
  discrète (où $y$
  est le point sous-jacent à $\bar y$) et o\`u $\Hom_{\Pt^{\an}(Y)}(\bar y,\bar
  y')$ est l'ensemble des spécialisations két de $Y_k$ de la log réduction
  $\bar y_s$ de $\bar y$ à la log réduction $\bar y_s'$ de $\bar y'$.\end{dfn}

\begin{thm}\label{cospcourbes}
Pour tout morphisme $\bar y_2\to\bar y_1$ dans $\Pt^{\an}(Y)$, il existe un
morphisme de cospécialisation
\[\gtemp(U_{\bar y_1})^{\mbb L}\to\gtemp(U_{\bar y_2})^{\mbb L}\]
qui est un isomorphisme si $\overline M_{\bar y_{1,s}}\to\overline M_{\bar y_{2,s}}$
est un isomorphisme.
\end{thm}
\dem
Soit $\bar
y_2\to\bar y_1$ un morphisme de $\Pt^{\an}(Y)$.\\
Il y a un diagramme 2-commutatif de foncteurs de cospécialisation~:
\[\begin{array}{ccc}\KCovgeom(X_{\bar y_{1,s}})^{\mbb L} & \to & \KCovgeom(X_{\bar
    y_{2,s}})^{\mbb L}\\
\dar & & \dar \\  \Covalg(U_{\bar y_{1}})^{\mbb L} & \to & \Covalg(U_{\bar
    y_{2}})^{\mbb L}\end{array}\]
Les flèches verticales sont des équivalences~(\cite[prop. 1.15]{kisin}) et
la flèche du bas aussi~(\cite[XIII.2.10]{sga}).\\
Soit $S_{1,s}$ un revêtement két log géométrique de $X_{\bar y_{1,s}}$ et
    soit $S_{2,s}$ (resp. $S_1$, $S_2$)
le revêtement correspondant de $X_{\bar y_{2,s}}$ (resp. $U_{\bar
  y_1}$, $U_{\bar y_2}$).\\

Il y a des applications fonctorielles en $S$~:
\[|S^{\an}_{1}|\leftarrow |\mbb G(S_{1,s})|\to |\mbb G(S_{2,s})|\to
|S^{\an}_{2}|\]
la première et la troisième flèche \'etant les plongements du squelette
d'une courbe semistable. Ce sont des équivalences d'homotopie.\\
On obtient donc un morphisme de types d'homotopie
$|S^{\an}_{1}|\to|S^{\an}_{2}|$, fonctoriel en $S$.\\

Rappelons que, si $Z$ est une variété analytique,
$\Dtop(Z)/\Covalg(Z)^{\mbb L}$ est la catégorie fibrée dont la fibre
en un revêtement $\mbb L$-fini $T\to Z$ est la catégorie des revêtements
topologiques de $T$ et $\Dtemp(Z)^{\mbb L}/\Covalg(Z)^{\mbb L}$ est la
catégorie fibrée dont la fibre en $T$ est la catégorie des revêtements
tempérés de $T$. Rappelons également que $\Dtemp(Z)^{\mbb
  L}/\Covalg(Z)^{\mbb L}$ est le champ associé à $\Dtop(Z)/\Covalg(Z)^{\mbb
  L}$.\\

Le morphisme fonctoriel $|S^{\an}_{1}|\to|S^{\an}_{2}|$ de types
d'homotopie fournit un foncteur de catégories fibrées~:
\[\begin{array}{ccc} \Dtop(U_{\bar y_2}) & \to & \Dtop(U_{\bar y_1})\\
\dar & & \dar\\
 \Covalg(U_{\bar y_{2}})^{\mbb L} & \simeq & \Covalg(U_{\bar
    y_{1}})^{\mbb L}\end{array}\]
qui induit un foncteur de champs~:
\[\begin{array}{ccc} \Dtemp(U_{\bar y_2})^{\mbb L} & \to & \Dtemp(U_{\bar
    y_1})^{\mbb L}\\
\dar & & \dar\\
 \Covalg(U_{\bar y_{2}})^{\mbb L} & \simeq & \Covalg(U_{\bar
    y_{1}})^{\mbb L}\end{array}\]
En prenant les sections globales, on obtient un foncteur~:
\[\Covtemp(U_{\bar y_2})^{\mbb L}\to\Covtemp(U_{\bar y_1}))^{\mbb L}.\]
Cela induit un morphisme de cospécialisation de groupes fondamentaux tempérés~:
\[\gtemp(U_{\bar y_1})^{\mbb L}\to\gtemp(U_{\bar y_2})^{\mbb L}.\]
\findem

\section{Cospécialisation du groupe fondamental $(p')$-tempéré}

Revenons au cas d'une log fibration polystable propre $X\to Y$, telle que $Y$ soit log
lisse et propre sur $O_K$ (la propreté de $Y\to O_K$ est seulement supposée
ici afin que tout point de $Y_\eta$ ait une réduction dans $Y_s$, mais la
construction que nous ferons sera locale sur $Y$).
Dans ce paragraphe, nous construisons des morphismes de cospécialisation
pour le groupe fondamental $(p')$-tempéré des fibres géométriques de
$X_\eta\to Y_\eta$. Grâce au théorème~\ref{isomtempgeom}, nous nous sommes
réduits à construire des morphismes de cospécialisation pour les groupes
fondamentaux $(p')$-tempérés des fibres log géométriques de $X_s\to
Y_s$. Soit donc $\bar s_2\to\bar s_1$ une spécialisation de points log
géométriques de $Y$, où $\bar s_1$ et $\bar s_2$ sont
les réductions de points géométriques $\bar \eta_1,\bar \eta_2 $ de
$Y_\eta$.\\
Pour comparer les groupes fondamentaux tempérés, nous aurons d'abord besoin
de comparer les groupes fondamentaux logarithmiques. Ainsi nous montrerons
que tout revêtement két géométrique de $X_{s_1}$ s'étend à $X_U$ pour un
certain voisinage két $U$ de $s_1$ dans le localisé $Z$ en $s_1$ de $Y$.\\
Ceci
découle de~\cite{org2} si $s_1$ est strict sur $Y$. Ainsi nous n'avons qu'à montrer
l'invariance du groupe fondamental log géométrique pro-$(p')$ par changement
de point de base. Ceci nous donnera un morphisme de spécialisation du
groupe fondamental log géométrique de $X_{s_1}$ vers celui de $X_{s_2}$, qui sera
un isomorphisme (par comparaison aux groupes fondamentaux de
$X_{\bar \eta_1}$ et de $X_{\bar \eta_2}$).\\
Nous avons maintenant une équivalence entre les revêtements két
géométriques pro-$(p')$ de $X_{\eta_1}$ et de
$X_{\eta_2}$. Nous devons comparer, pour tout tel revêtement két  $Z_{s_1}$
correspondant à $Z_{s_2}$ (qui s'étend à $Z_U$), leurs ensembles
polysimpliciaux, définis dans~\ref{ketpolysimpcplx}. Nous supposerons que
$s_2$ est le point générique de sa strate (si $s_1$ et $s_2$ sont dans la
même strate et $\Cgeom(X_{s_2})$ est intérieurement libre, il s'avérera que
$\Cgeom(X_{s_1})\to\Cgeom(X_{s_2})$ est un isomorphisme. Ainsi on pourra
remplacer $s_2$ par le point générique de sa strate). Nous construisons le
morphisme de cospécialisation d'ensembles polysimpliciaux localement pour
la topologie étale. Ainsi nous pourrons supposer $X$ strictement polystable
(la propreté n'est pas nécessaire ici). Ce morphisme de cospécialisation
d'ensembles polysimpliciaux sera construit en associant, après localisation
két de la base pour que
$Z_U$ devienne saturé, à une strate $z$ de $Z_{s_1}$ la strate minimale
$z'$ de $Z_{s_2}$ telle que $z$ soit dans la clôture de $z'$ (comme dans le
lemme~\ref{berk62}). Alors la clôture de $z'$ dans le localisé strict du
point générique de $z$ est séparable sur son image. Ainsi $z'$ est géométriquement
connexe et définit donc une strate géométrique de $X_{s_2}$.\\
Nous terminerons en recollant notre isomorphisme de spécialisation du
groupe fondamental log géométrique pro-$(p')$ avec notre morphisme de
cospécialisation d'ensembles polysimpliciaux en un morphisme de
cospécialisation de groupes fondamentaux tempérés.

\subsection{Cospécialisation d'ensembles polysimpliciaux}

Dans cette section, nous construisons un morphisme de cospécialisation
d'ensembles polysimpliciaux pour la composée d'un morphisme két et d'une
log fibration polystable.\\

On a un résultat analogue à~\cite[prop 2.9]{berk2}:
\begin{prop} Soit $Z' \to Z$ un morphisme strictement plurinodal de log
  schémas fs, et $Z''\to Z'$ un morphisme két de log schémas. Soit
  $z_1$ et $z_2$ deux strates de $Z$ (muni de la log structure image
  inverse $Z$), tel que $z_2\leq z_1$ (\ie $z_1 \in \overline{\{ z_2 \} }$). 
Alors on a un morphisme de cospécialisation $\Str(Z''_{z_1})\to \Str(Z''_{z_2})$
qui envoie une strate $x_1$ de $\Str(Z''_{z_1})$ vers l'unique élément
maximal de $\{ x_2\in \Str(Z''_{z_2}) | x_1\in \bar x_2 \}$.\\
Le morphisme de cospécialisation envoie les points minimaux vers les points
minimaux.\\
Si $z_3\leq z_2\leq z_1$, le triangle évident de morphismes de
cospécialisation est commutatif.\\
\end{prop}
\dem
Comme pour~\cite[prop 2.9]{berk2}, si le résultat est vrai pour deux
morphismes 
$\phi:Z'' \to Z'$ et $\psi: Z' \to Z$, alors il est aussi vrai pour $\psi
\circ \phi$ car $\Str(Z''_z)=\coprod_{z' \in \Str(Z'_z)}
\Str(Z''_{z'})$. Mais il est vrai si $\psi$ est strictement
plurinodal~(\cite[prop 2.9]{berk2}), et c'est aussi vrai pour $\phi$ két d'après~(\ref{lemkum}).\findem

\begin{dfn}
Un couple de points $(z_2 \leq z_1)$ d'un Zariski
log schéma fs $Z$ est un \emph{bon couple} si il existe un voisinage $U$ de $z_1$ et
une carte fs $U\to \Spec \mathbf Z[P]$ telle que si $\mathfrak p$ est
l'image de $z_2$ dans $\Spec P$ par $U \to \Spec P$, et si $F=P\backslash
\mathfrak p$, le schéma réduit $\overline{\{z_2\} }$ muni de la log structure
associée à $F$ (ce log schéma sera noté $\overline{\{z_2\} }_F$) par le morphisme
\[\overline{\{z_2\}  }\to \Spec k[P]/k[\mathfrak p] \simeq \Spec k[F]\]
est log régulier.\\
Une spécialisation de points géométriques (resp. de points log géométriques) $(\bar z_2\to \bar
z_1)$ d'un log schéma fs $Z$  est un \emph{bon couple} si il existe un voisinage
étale (resp. két) $U$ de $\bar z_1$ tel que $U$ ait une carte globale
(et donc soit Zariski) et le couple $(z_2\leq z_1)$ de points
sous-jacents de $U$ est un bon couple.\end{dfn} 
Un couple de points $(z_2\leq z_1)$ est un bon couple si $z_2$
est le point générique d'une strate d'un log schéma fs Zariski log régulier
(\cite[prop. 7.2]{kato2}).\\

\begin{lem} Soit $Z' \to Z$ un morphisme strictement plurinodal de log schémas fs, soit $(z_2,z_1)$ un
  bon couple de points de $Z$. Soit $z'_2$ (respectivement $z'_1$) le
  point générique d'une strate de $Z'_{z_2}$ (respectivement $Z'_{z_1}$), tel
  que $z'_2\leq z'_1$. Alors
  $(z'_2,z'_1)$ est un bon couple de points de $Z'$. \label{propboncouples}\end{lem}
\dem Comme l'assertion est locale sur $Z$ et comme la log régularité est
stable par changement de base étale, on
peut supposer que l'on a une carte fs ($\phi: P \to Q$) de $Z\to Z'$ telle
que le carré
\[\begin{array}{ccc}Z' & \to & \Spec \mathbf Z[Q]\\ \dar & \square & \dar \\ Z & \to & \Spec \mathbf Z[P]\end{array}\]
soit cartésien.\\
Soit $\mathfrak p$ (respectivement $\mathfrak q$) l'idéal de $P$
(respectivement $Q$) qui est l'image de $z_1$ (respectivement $z_2$) et soit
$F=P\backslash \mathfrak p$ (respectivement $F'=Q\backslash \mathfrak q$) la
face associée. On veut montrer que $\overline{ \{z'_1\} }$ muni
de la log structure associée à $F'$ est log régulier.\\
On a $F=\phi^{-1}(F')$, $\Ker \phi_{|F}^{\gp} \subset \Ker \phi^{\gp}$
et $\Coker \phi_{|F}^{\gp} \subset \Coker \phi^{\gp}$, donc $\phi_{|F}: F
\to F'$ est aussi un morphisme log lisse de monoïdes.\\
Notons $\mbf Z[\fk p]:=\{\sum_{p\in\fk p}a_pp,a_p\in\mbf Z\}\subset\mbf
Z[P]$ (c'est un idéal premier de $\mbf Z[P]$).\\
Le diagramme suivant de schémas est commutatif~:
\[\begin{array}{ccccc}\Spec \mathbf Z[F'] & \simeq & \Spec \mathbf Z[Q]/\mathbf Z[\mathfrak q] & \to & \Spec \mathbf Z[Q] \\ \dar & & \dar & & \dar  \\ \Spec \mathbf Z[F] & \simeq & \Spec \mathbf Z[P]/\mathbf Z[\mathfrak p] & \to & \Spec \mathbf Z[P] \end{array}\]
Ainsi
\[\begin{array}{ccc} \overline{\{ z'_1\} } & \to & \Spec \mathbf Z[F'] \\ \dar & & \dar \\ \overline{\{ z_1\} } & \to & \Spec \mathbf Z[F] \end{array}\]
est commutatif. Soit $Z''=\overline{\{ z_1\} } \times_{\Spec \mathbf Z[F]}
\Spec \mathbf Z[F']$ et munissons-le de la log structure associée à
$F'$. Puisque $Z'' \to \overline{\{ z_1\} }_F$ est log lisse et
$\overline{\{ z_1\} }_F$ est log régulier, $Z''$ est log régulier
(\cite[th. 8.2]{kato2}).\\
On a un morphisme de log schémas $\overline{\{ z'_1\} } \to Z''$ qui est
l'immersion fermée d'une composante irréductible (car $Z''$ est la préimage
de $\overline{\{ \mathfrak q\} }$ dans $Z'_{\overline{\{ z_1\} }}$ et $z'_1$
est un point générique de $\overline{\{ \mathfrak q\} }$ dans
$Z'_{z_1}$ par définition d'une strate), qui induit un morphisme strict de
log schémas $\overline{\{ z'_1\} }_{F'} \to Z''$. Ainsi, comme $Z''$ est log
régulier (et donc normal), $\overline{\{ z'_1\} }_{F'}$ est une
composante connexe de $Z''$ et donc est aussi log régulier.\findem

\begin{lem}\label{cosp1} Soit $\phi: Z' \to Z$ un morphisme strictement plurinodal
 de log schémas fs. Soit $(z_2,z_1)$ un bon couple de points de $Z$.\\ 
On a un morphisme de cospécialisation $\Str(Z'_{z_1})\to \Str(Z'_{z_2})$
qui envoie une strate $x_1$ de $\Str(Z'_{z_1})$ vers l'unique élément
maximal de $\{ x_2\in \Str(Z'_{z_2}) | x_1\in \bar x_2 \}$.\\
Le morphisme de cospécialisation envoie les éléments minimaux vers les
éléments minimaux.\\
Si $z_3\leq z_2\leq z_1$, le diagramme évident de morphismes de
cospécialisation est commutatif.\\
\end{lem}
\dem C'est clairement vrai si $\phi$ est un morphisme plurinodal standard $\Spec B \to \Spec
A$ avec $f:  P\to A$ une carte $\Spec A$ et $B=A[u,v]/(uv-f(a))$ où $a
\in M$ (on peut par exemple utiliser~\cite[lem 2.3]{berk2}).\\
D'après le lemme~(\ref{propboncouples}), on voit, comme dans la preuve de~\cite[prop. 2.9]{berk2},
que si la proposition est vraie pour $\phi$ et $\phi'$, elle est vraie pour
$\phi \circ \phi'$. De plus le résultat est local pour la topologie de
Zariski de $Z'$, donc il n'y a à montrer le résultat que pour $\phi$ étale, mais
cela découle alors de~(\ref{lemetale}) et du fait que $\overline{\{z_2\} }$ est
normal dans un voisinage de $z_1$ et les adh\'erences dans
$Z'\times_Z\overline{\{z_2\}}$ de deux points différents de $Z_{z_2}$ (qui
sont deux composantes irréductibles de $Z'\times_Z\overline{\{z_2\}}$, qui
est normal)
ont une intersection vide. \findem

\begin{rem} Supposons toujours que $(z_2,z_1)$ est un bon couple. Si l'on a un morphisme két $Z''\to Z'$ tel que $Z''\to Z$ soit
  aussi strictement plurinodal, alors, comme pour~\cite[cor. 2.11]{berk2},
  le diagramme suivant
\[\begin{array}{ccc}\Str(Z''_{z_1}) & \to & \Str(Z''_{z_2})\\ \dar & 
  & \dar \\ \Str(Z'_{z_1}) & \to & \Str(Z'_{z_2})\end{array}\]
est commutatif.\end{rem}
Si $Z'\to Z$ est un morphisme strictement polystable, il induit (comme 
dans~\cite[lem. 6.1]{berk2}) un morphisme de complexes polysimpliciaux
\[\C(Z'_{z_1})\to \C(Z'_{z_2}).\]
Si maintenant $Z'\to \cdots\to Z$ est une fibration strictement polystable,
en utilisant le lemme~\ref{propboncouples}, on construit par récurrence sur
la longueur de la fibration un morphisme de complexes polysimpliciaux~:
\[\C(Z'_{z_1})\to \C(Z'_{z_2}).\]

\begin{rem}Le morphisme $Z'\to Z$ est saturé, donc pour toute extension
  $z'_1$ két de
  $\bar z_1$, $Z'_{\tilde z'_1}\to Z'_{\bar z_1}$ est un isomorphisme sur
  les schémas sous-jacents et donc on même strates et donc même complexe
  polysimplicial. En particulier, si $s_1$ est un log
  point géométrique de $Z$ au-dessus de $\bar z_1$, le morphisme
  $\Cgeom(Z'_{\bar z_1}/s_1)\to \C(Z'_{\bar z_1})$ est un isomorphisme (et
  le même résultat est vrai pour $z_2$).\end{rem}

Soit $(\bar z_2\to\bar z_1)$ un bon couple de points log géométriques de
$Z$.\\
Quitte à remplacer $Z$ par sa stricte localisation $Z_1$ en $\bar z_1$ (et choisissons
une bonne carte modelée sur $P$ de $Z_1$ en $\bar z_1$), on obtient un morphisme
$\psi:\Str(Z'_{\bar z_1})\to\Str(Z'_{z'_2})$,
où $z'_2$ est l'image de $\bar z_2$ dans le localisé strict de $Z$
en $\bar z_1$.\\
Soit $x_1$ le point générique d'une strate $\tilde x_1$ de $Z'_{\bar z_1}$, et
soit $Z'_1$ le localisé de $Z'_{Z_1}$ en $x_1$ et soit $Z''_1$ la
clôture de $\psi(x_1)$ dans $Z'_1$ ($x_1$ est encore dans $Z''_1$).\\
Etale localement dans un voisinage de $x_1$, $Z'\to Z$ est isomorphe au
pullback à $Z$ de $\Spec \mbf Z[Q]\to
\Spec \mbf Z[P]$ où $P\to Q$ est un morphisme saturé. Donc le morphisme de
la clôture d'une strate $\overline{\{x_2\}}$ de $Z'_{z'_2}$ vers son image
$Z_0$ est étale localement isomorphe au pullback à $Z_0$ de $\Spec \mbf
Z[F']\to\Spec \mbf Z[F]$ où $F'$ est la face de $Q$ correspondant à 
$x_2$ et $F$ est la face préimage de $F'$ dans $P$. Alors $F\to F'$ est aussi
un morphisme saturé de monoïdes grâce à~\ref{satface}. En particulier $\Spec \mbf
Z[F']\to\Spec \mbf Z[F]$ est un morphisme séparable de schémas.\\
D'après~\cite[cor. 18.9.8]{ega4}, les fibres de $Z''_{1}\to
Z_{0}$ sont géométriquement connexes.
En particulier la strate
$x_{2}$, image de $x_1$ par $\psi:\Str(Z'_{\bar
  z_1})\to\Str(Z'_{z'_{2}})$ est géométriquement connexe, définissant ainsi
  une strate de $\Str(Z'_{\bar z_2})$.\\
On obtient ainsi un morphisme canonique $\Str(Z'_{\bar
  z_1})\to\Str(Z'_{\bar z_2})$ qui rend le diagramme commutatif:
\[\xymatrix{ & \Str(Z'_{\bar z_2})\ar[d] \\ \Str(Z'_{\bar z_1})\ar[ur]\ar[r]
  & \Str(Z'_{z'_2})}\]
On obtient donc un morphisme de cospécialisation
\[\Cgeom(Z'_{z_1}/z_1)\to\Cgeom(Z'_{ z_2}/z_2).\]

Si $\underline Z':Z'\to \cdots\to Z$ est maintenant une log fibration polystable, quitte
à changer $Z$ par un voisinage étale, il existe un morphisme étale
surjectif $\underline Z''\to
\underline Z'$ de log fibrations polystables au-dessus de $Z$ tel que
$\underline Z''$
soit strictement polystable. Notons $\underline Z'''=\underline
Z''\times_{\underline Z'}\underline Z$. Comme $\C(Z'''_z)=\Coker(\C(Z''_z)\rightrightarrows\C(Z'_z))$, en prenant le conoyau
des flèches horizontales du diagramme commutatif~:
\[\begin{array}{ccc}C(Z'''_{\bar z_1}) & \rightrightarrows &
  \C(Z''_{\bar z_1})\\ \dar & & \dar \\ \C(Z'''_{\bar z_2}) &
  \rightrightarrows & \C(Z''_{\bar z_2})\end{array}\] 
on obtient un morphisme de cospécialisation $\C(Z'_{\bar z_1})\to \C(Z'_{\bar
  z_2})$ (qui est fonctoriel vis-à-vis des morphismes étales).\\
Ces morphismes de cospécialisation commutent aux morphismes két de log
fibrations polystables.\\

Si $Z''\to Z'$ est un morphisme két et $Z'\to Z$ est une log fibration
polystable et soit $(\bar z_2\to \bar z_1)$ un bon couple de points log
géométriques.
\begin{prop}\label{cospcc} Il existe un morphisme de cospécialisation canonique
  \[\Cgeom(Z''_{z_1}/z_1)\to\Cgeom(Z''_{z_2}/z_2),\] fonctoriel en $Z''$.\\

Si l'on dispose d'un morphisme $Z'''\to Z''$
 \end{prop}
\dem 
Supposons $Z'\to Z$ strictement polystable.
Quitte à remplacer $Z$ par un voisinage két de $\bar z_1$ et
à remplacer $z_1$ et $z_2$ par le sous-schéma réduit d'un revêtement
galoisien két connexe dans $\bar z_1$ et $\bar z_2$, on peut supposer que
$Z$ a une carte globale modelée sur $M$, qu'il existe un recouvrement étale
$(U'_i)_{i\in I}$ de $Z'$ avec $I$ fini, que $U'_i$ est étale sur
$Z\times_{\Spec \mbf Z[M]}\Spec \mbf Z[P_i]$, $Z''$ a un recouvrement étale
$(U''_i)_{i\in I}$ où $U''_i$ est étale sur $U'_i\times_{\Spec \mbf
  Z[P_i]} \Spec\mbf Z[Q_i]$ et où $M\to Q_i$ est saturé. Ainsi il existe un
morphisme $(p')$-Kummer de monoïdes $Q_i\to P'_i$ tel que
$V_i=U''_i\times_{\Spec\mbf Z[Q_i]}\Spec\mbf Z[P'_i]$ s'insère dans une
fibration polystable $\underline V_i:V_i\to\cdots \to Z$ and $\underline
V_i\to\underline U'_i$ est un morphisme két de log fibrations
polystables (en particulier $V_i\to U''_i$ est un revêtement galoisien
de groupe $G=((P'_i)^{\gp}/Q_i^{\gp})^{\vee}$, et $Z''\to
Z$ est saturé).\\ 

$\C(V_{i,z_1})\to\C(V_{i,z_2})$ est $G$-équivariant, et donc
induit un morphisme
$\C(U''_{i,z_2})=\C(V_{i,z_2})/G\to\C(V_{i,z_1})/G=\C(U''_{i,z_1})$.
On déduit du fait que les morphismes de cospécialisation commutent avec les
morphismes két
de log fibrations polystables qu'il ne dépend pas du choix de $V_i$ et
qu'il se descend en un morphisme $\C(Z''_{z_1})\to\C(Z''_{z_2})$.\\
En prenant la limite projective sur les voisinages étales stricts de
$\mring{\bar{z_1}}$, on obtient un morphisme $\Cgeom(Z''_{z_1}/z_1)\to
\C(Z''_{z'_2})$, où $z'_2$ est l'image de $\bar z_2$ dans le localisé
strict de $Z$ en $\mring{\bar{z_1}}$.\\
Si l'on a une strate log géométrique de $Z''_{z_1}$, en utilisant
\cite[cor. 18.9.8]{ega4} comme précédemment, on obtient que la strate image
de $Z''_{z_2}$ est géométriquement connexe.\\
On obtient alors le morphisme voulu
$\Cgeom(Z''_{z_1}/z_1)\to\Cgeom(Z''_{z_2}/z_2)$.\\

Si l'on a un morphisme $Z'''\to Z''$, on peut recouvrir $Z'''$ par des
ouverts étales
$U'''_i$ et $Z''$ par des ouverts étales  $U''_i$ tels que $U'''_i\to U''_i$ soient des
revêtments két. Quitte à raffiner les familles $(U''_i)$ et $(U'''_i)$, on
peut supposer qu'il existe une fibration strictement polystable $\underline
V_i$ et un revêtement két $V_i\to U'''_i$. De plus, quitte à raffiner
encore les recouvrements, on peut supposer les revêtements galoisiens
(puisque le groupe fondamental logarithmique d'un log schéma de schéma
sous-jacent strictement local est commutatif). Soit $G=\Gal(V_i/U''_i)$ et
$H=\Gal(V_i/U'''_i)$. Alors le diagramme 
\[\begin{array}{ccc}\C(V_{i,z_2})/H & \to  & \C(V_{i,z_1})/H\\ \dar & & \dar\\
  \C(V_{i,z_2})/G & \to & \C(V_{i,z_1})/G\end{array}\]
est commutatif. On en déduit la fonctorialité par descente le long des
recouvrements étales.\\

Si $Z'\to Z$ n'est plus supposé strictement polystable, quitte à changer
$Z$ par un voisinage étale, il existe un morphisme étale $Z'_0\to Z'$ de
log fibrations polystable sur $Z$ telles que $Z'_0\to\cdots \to Z$ soit une
log fibration strictement polystable. Soit
$Z'_1=Z'_0\times_{Z'}Z'_0$. Alors, le morphisme voulu est obtenu en prenant
le conoyau des flèches horizontales du carré commutatif suivant~:
\[\begin{array}{ccc}\Cgeom(Z''_{1,z_1}/z_1) & \rightrightarrows &
  \C(Z''_{0,z_1}/z_1)\\ \dar & & \dar \\ \Cgeom(Z''_{1,z_2}/z_2) &
  \rightrightarrows & \Cgeom(Z''_{1,z_2}/z_2)\end{array}\]
\findem

Supposons maintenant que $Z''\to Z$ soit propre, que $Z$ soit log régulier et
que $\bar z_1$ et $\bar z_2$ soient dans la même strate de $Z$ (\ie que le morphisme
de cospécialisation $M_{Z,z_1}\to M_{Z,z_2}$ est un isomorphisme). On peut remplacer
$Z$ par son localisé strict en $z_1$ (ceci n'affecte pas notre morphisme de
cospécialisation). En particulier $Z$ est log Zariski. Quitte à prendre
encore un voisinage két, on peut supposer que $Z''\to Z$ est saturé.\\
Le fait que $z_1$ et $z_2$ soient dans la même strate de $Z$ implique
que $\Cgeom(Z''_{z_1})\to\Cgeom(Z''_{z_2})$ envoie les polysimplexes non dégénérés 
vers les polysimplexes non dégénérés (il suffit de regarder localement pour la topologie \'etale).\\
Soit ${Z''}^{(i)}$ l'adhérence de $(Z''_{\bar z_2})^{(i)}$ dans $Z''$
muni de la structure de sous-schéma réduit, et soit
$(\widetilde{Z}'')^{(i)}$ sa normalisation. En regardant localement pour la topologie \'etale sur
$Z''$, et grâce au fait que $z_2$ et $z_1$ sont dans la même strate,
on voit que $({Z''}^{(i)})_{z_1}$ est juste $(Z''_{z_1})^{(i)}$
et que $((\widetilde{Z}'')^{(i)})_{z_1}$ est juste la normalisation de 
$(Z''_{z_1})^{(i)}$. Ainsi les composantes connexes de $((\widetilde
  Z'')^{(i)})_{z_1}$ et de $((\widetilde Z'')^{(i)})_{z_2}$ sont en bijection
avec les strates de $Z''_{z_1}$ et de $Z''_{z_2}$ de rang $i$. Comme
expliqué précédemment, puisque $Z''\to Z$ est saturé, $(\widetilde Z'')^{(i)}$ est
séparable sur la clôture de $z_2$ et donc la factorisation de Stein de
$(\widetilde Z'')^{(i)}$ pour tout $i$ nous dit que $\Str(Z_{\bar
  z_1})\to\Str(Z_{\bar z_2})$ est bijectif.\\

Si l'on suppose de plus que $\Cgeom(Z''_{z_2}/z_2)$ est intérieurement libre
(c'est le cas si $\Cgeom(Z'_{z_2}/z_2)$ est intérieurement libre),
\[\Cgeom(Z''_{{z}_1}/z_1)\to\Cgeom(Z''_{z_2}/z_2)\] est aussi un isomorphisme.\\

Supposons $Z''\to Z$ propre, $Z$ log régulier, et supposons que pour toute strate
de $Z$ de point générique log géométrique $\bar z$, $\Cgeom(Z'_z)$ est
intérieurement libre.\\
Soit $z_2\to z_1$ une spécialisation de points log géométriques $Z$. Soit
$x\to z_2$ une spécialisation où $x$ est un point log géométrique au-dessus
du point générique de la strate de $z_2$ (il existe une telle
spécialisation). Rappelons que $(x\to z_2)$ et
$(x\to z_1)$ sont alors de bons couples (cf. \cite[prop. 7.2]{kato2}) . On a alors les morphismes~:
\[\Cgeom(Z_{z_2})\stackrel{\simeq}{\leftarrow}\Cgeom(Z_x)\to\Cgeom(Z_{z_1}).\]
On obtient donc un morphisme $\Cgeom(Z_{z_2})\to\Cgeom(Z_{z_1})$ qui ne
dépend pas de $x\to z_2$ (puisque tout autre morphisme $x\to z_2$ se
factorise à travers notre spécialisation précédente).\\ 
De plus, si $z_1$ et $z_2$ sont dans la même strate, le morphisme de
cospécialisation est un isomorphisme.

\subsection{Morphisme de cospécialisation du groupe fondamental $(p')$-tempéré}
Sous les mêmes hypothèses que précédemment, soit $\bar y_2\to \bar y_1$ une
  spécialisation de points log géométriques sur des log points fs $y_2\to y_1$ de $Z'$ dont
les images $z_2\to z_1$ dans $Z$ se trouvent dans la même strate de $Z$.
 On a alors un foncteur de cospécialisation
$F:\KCovgeom(Z'_{z_1}/z_1)^{\mbb L}\to\KCovgeom(Z'_{z_2}/z_2)^{\mbb L}$ si
$\mbb L$ ne contient pas la caractéristique résiduelle $p$ en $z_1$. Si
$Z''_{z_1}$ est un revêtement két géométrique de $Z'_{z_1}$, il s'étend
grâce au corollaire~\ref{logsp} à un voisinage két $U$ de $\bar z_1$
dans $Z$. Soit $Z''_U\to U$ cette extension (unique quitte à remplacer $U$ par
un voisinage plus petit de $\bar z_1$), et $F(Z''_{z_1})$ est juste la
fibre de $Z'_0\to Z'$  dans $\bar z_2$ $Z''_U$. On a alors un isomorphisme
$\Cgeom(Z''_1)\simeq\Cgeom(Z''_2)$, qui induit un foncteur de catégories fibrées~:
\[\begin{array}{ccc}\Dtopgeom(Z'_1) & \to & \Dtopgeom(Z'_2)\\ \dar & &
  \dar\\ \KCovgeom(Z'_{z_1}/z_1)^{\mbb L} & \to & \KCovgeom(Z'_{z_2}/z_2)^{\mbb
    L} \end{array}\]
et donc un morphisme de spécialisation $\gtempgeom(Z'_{z_2},y_2)^{\mbb L}
\to\gtempgeom(Z'_{z_1},y_1)^{\mbb L}$.\\

Soit maintenant $K$ un corps complet pour une valuation discrète. $\Spec O_K$ est muni de sa
log-structure usuelle. Soit $\mbb L$ un ensemble de nombres premiers ne contenant pas la caractéristique résiduelle $p$ de
$K$ ($p$ peut \^etre nul).\\ 
Soit $X\to Y$ une log fibration polystable propre telle que $Y\to \Spec O_K$
soit log lisse, et supposons que pour tout point géométrique $\bar
y$ de $Y_s$,
$\C(X_{\bar y})$ est \emph{intérieurement libre} (c'est par exemple le cas si $X\to
Y$ est strictement polystable ou si les fibres de $X\to Y$ sont des courbes
semistables, au sens de la définition~\ref{courbessemistables}).\\
Rappelons que la fibre générique $\fk Y_{\eta}$ du complété formel de $Y$
le long de sa fibre spéciale s'identifie à un sous-domaine analytique de
$Y^{\an}$ (de plus $\fk Y_{\eta}\to Y^{\an}$ est un isomorphisme si $Y$ est
supposé propre).\\
Soient $y_1$ et $y_2$ deux points de
$Y_{\triv}^{\an}\cap \fk Y_{\eta}$ à valeur dans un corps à valuation discrète (quitte à
remplacer $\mcal H(y_1)$ par une extension isométrique, ce qui ne change
pas le groupe fondamental tempéré géométrique, on supposera que
$\mcal H(y_1)$ a un corps résiduel algébriquement clos). On a des
morphismes canoniques de log schémas fs $\Spec O_{\mcal
  H(y_i)}\to Y$ pour $i=1,2$. Soit $s_i$ le log point fs de $Y$
correspondant au point spécial de $\Spec O_{\mcal
  H(y_i)}$ muni de la log structure image inverse. Soit $s'_i$ le log
point fs de $Y$ qui a même schéma sous-jacent que $s_i$ mais muni de la log
structure image inverse.\\
Un point géométrique $\bar y_i$ (de l'espace de Berkovich $\fk Y_{\eta}$) au
dessus de $y_i$ induit un point log géométrique $\bar s_i$ au-dessus de
$s_i$ ($\bar y_i$ peut aussi être vu comme un point log géométrique de $Y$
puisque la log structure de $Y$ est triviale en $y_i$. Soit $\bar s'_i$ un
point log géométrique de
$s'_i\times_{s_i}\bar s_i$). \\
Considérons une spécialisation két $\bar s'_2\to\bar
s'_1$ (il en existe dès qu'il existe une spécialisation entre les points
géométriques sous-jacents de $\mring Y$).\\

Plus précisément, définissons la catégorie $\Pt^{\an}(Y)$ dont les objets
sont les points géométriques $\bar y$ de $Y_{\tr}^{\an}\cap \fk Y_{\eta}$ tels que
$\mcal H(y)$ est à valuation discrète (où $y$ est le point sous-jacent à
$\bar y$) et $\Hom(\bar y,\bar y')$ est l'ensemble des spécialisations két
de $\bar s$ vers $\bar s'$ où $\bar s$
et $\bar s'$ sont les log-réductions de $\bar y$ et
$\bar y'$.\\
Définissons aussi la catégorie $\Pt^{\an}_0(Y)$ obtenu à partir de
$\Pt^{\an}(Y)$ en inversant les morphismes $\bar
y\to\bar y'$ tels que $\bar s$ et $\bar s'$ soient dans la même strate $Y$.\\  

\begin{thm}\label{cospthm} Soit $X\to Y$ une log fibration polystable
  propre et $Y\to\Spec O_K$ un morphisme log lisse. Supposons que, pour tout
  point géométrique $\bar y$ de $Y_s$, $\C(X_{\bar y})$ est intérieurement
  libre. Alors, pour tout morphisme $\bar y_2\to\bar y_1$ dans
  $\Pt^{\an}(Y)$, il existe un morphisme extérieur
\[\gtemp(X_{\bar y_1}^{\an})^{\mbb L}\to\gtemp(X_{\bar
  y_2}^{\an})^{\mbb L},\]
qui est un isomorphisme si $s_1$ et $s_2$ sont dans la même strate de $Y$.\end{thm}
\dem
On a un foncteur de cospécialisation
\[F:\KCovgeom(X_{s_1}/s_1)^{\mbb L}\to\KCovgeom(X_{s_2}/s_2)^{\mbb L}\] qui
se factorise à travers
$\KCovgeom(X_{Z_0}/Z_0)^{\mbb L}$ où $Z_0$ est le localisé strict en $s_1$.
Soit $\eta$ un point générique au-dessus du point générique de $Y$.\\
Comme les foncteurs de cospécialisation
$\KCovgeom(X_{s_i}/s_i)^{\mbb L}\to\KCovgeom(X_{y_i}/y_i)^{\mbb L}$ et
$\KCovgeom(X_{y_i}/y_i)^{\mbb L}\to\KCovgeom(X_{\eta}/\eta)^{\mbb L}$ sont
des équivalences (\cite[prop. 1.15]{kisin}), on obtient que
$\KCovgeom(X_{s_1}/s_1)^{\mbb L}\to\KCovgeom(X_{s_2}/s_2)^{\mbb L}$ est une
équivalence.\\

Si $Z_{s_1}$ est un revêtement két de $X_{s_1}$, il s'\'etend grâce au
corollaire~\ref{logsp} à un voisinage két $U$ de $\bar s_1$
dans $Z$. Soit $Z_U\to U$ cette extension (unique quitte à remplacer $U$ par
un voisinage két plus petit de $\bar s_1$), et $F(Z_{s_1})$ est juste
la fibre en $\bar z_2$ de $Z_U$. On a alors un morphisme de cospécialisation
$\Cgeom(Z_{s_1})\to\Cgeom(Z_{s_2})$, qui induit un foncteur de spécialisation
\[(\Dtopgeom)_{X_{s_2}}(Z_{s_2})\to (\Dtopgeom)_{X_{s_1}}(Z_{s_1}),\]
qui est une équivalence de catégories si $\bar s_1$ et $\bar
s_2$ sont dans la même strate de $Y$ ($(\Dtopgeom)_{X_{s_i}}$ est
la catégorie fibrée sur $\KCovgeom(X_{s_i})^{\mbb L}$ des revêtements topologiques).\\
On a donc un diagramme 2-commutatif~:
\[\begin{array}{ccc}(\Dtopgeom)_{X_{s_2}} & \to & (\Dtopgeom)_{X_{s_1}}\\ \dar & &
  \dar\\ \KCovgeom(X_{s_2}/s_2)^{\mbb L} & \stackrel{F^{-1}}{\to} &
  \KCovgeom(X_{s_1}/s_1)^{\mbb L} \end{array}\]
où $F^{-1}$ est un quasi-inverse de $F$.
Cela induit un morphisme de cospécialisation~:
\[\gtempgeom(X_{s_1}/s_1)^{\mbb L}\to\gtempgeom(X_{s_2}/s_2)^{\mbb L}\]
Les morphismes de comparaison du théorème~\ref{isomtempgeom} nous donnent
le morphisme voulu.
\findem

On obtient donc un foncteur $\Pt^{\an}_0(Y)^{\op}\to \OutGptop$ (où
$\OutGptop$ est la catégorie des groupes topologiques à automorphisme
intérieur près)
qui envoie $\bar y$ vers $\gtemp(X_{\bar y}^{\an})^{\mbb L}$.\\

\begin{rem}
Un tel foncteur ne peut pas exister si $p\neq 0$ et $\mbb L$ est l'ensemble
de tous les nombres premiers. En effet, si $X_1$ et $X_2$ sont deux courbes
de Mumford avec réduction isomorphe mais une métrique différente sur le
graphe de leur modèle stable, alors leurs groupes fondamentaux tempérés
sont non isomorphes. Considérons un espace de module de courbes stables
avec structure de niveau muni de sa log structure canonique sur $\Spec \mbf Z_p$
(lui même muni de sa log structure induite par $\mbf Z_p^*\to\mbf Z_p$) et un point
géométrique $\bar s$ dans la fibre spéciale de l'espace de modules
correspondant à une courbe stable totalement dégénérée. En particulier, elle
a au moins un point double, et donc le rang de $\overline M_{\bar
  s}^{\gp}$ est au moins deux (rappelons que $\overline M_{\bar
  s}=\mbf N |p|\oplus\bigoplus_{e\in \mcal E}\mbf Nm_e$ où $\mcal E$ est
l'ensemble des arêtes du graphe
dual de la courbe correspondant à $\bar s$). Prenons deux log points (de
corps sous jacent séparablement clos) fs $s_1$ et $s_2$ valuatifs (\ie $\overline
M_{s_i}\simeq\mbf N$) tels que les morphismes correspondants $\overline M_{\bar
  s}\to\mbf N$ soient
linéairement indépendants. Il existe une unique normalisation $M_{s_i}\to\mbf
Q$ tel que l'image de $|p|$ par la composée $\overline M_{\mbf Z_p,\bar{
  \mbf F}_p}\to\overline M_{s_i}\to\mbf Q$ soit $1$. Soient $\eta_1$ et $\eta_2$ des
points discr\`etement valués de la fibre géométrique analytique dont les
log-réductions sont $s_1$ et $s_2$. Alors les
deux courbes de Mumford géométriques correspondantes $C_i$ ont des métriques
différentes sur le graphe de leur modèle stable (car la longueur d'une
arête $e$ pour la distance de $C_i$ est égale à l'image de $m_e\in
\overline M_{\bar s}$ par la 
composée $\overline M_{\bar s}\to\overline M_{s_i}\to\mbf Q$), et donc ont
des groupes 
fondamentaux tempérés différents d'apr\`es le th\'eor\`eme~\ref{mumford}. Mais deux log points géométriques au
dessus de $s_1$ et $s_2$ sont isomorphes vis-à-vis de la spécialisation
pour la topologie két.\end{rem}

Si l'on enlève l'hypothèse de liberté intérieure, on obtient quand même le
résultat suivant, avec la même preuve~:
\begin{thm} Pour tout couple de points géométriques $\bar y_1$ et
  $\bar y_2$ de $Y^{\an}_{\tr}\cap \fk Y_\eta$ au-dessus de points à valuation discrète $y_1, y_2$ de
  $Y$. Soit $\bar s_1,\bar s_2$ leurs log réductions et supposons que
  $\bar s_1$ est le point générique d'une strate de $Y$. Alors il existe un
  homomorphisme extérieur de cospécialisation, fonctoriel~:
\[\gtemp(X_{\bar y_1}^{\an})^{\mbb L}\to\gtemp(X_{\bar
  y_2}^{\an})^{\mbb L}.\]\end{thm}
La preuve est similaire à celle du  théorème \ref{cospthm}, si ce n'est que
le foncteur \[(\Dtopgeom)_{X_{s_2}}(Z_{s_2})\to
(\Dtopgeom)_{X_{s_1}}(Z_{s_1})\] n'a était défini que pour $\bar s_1$ point
générique d'une strate (cf. proposition \ref{cospcc})

\printindex

\providecommand{\bysame}{\leavevmode\hbox to3em{\hrulefill}\thinspace}
\providecommand{\MR}{\relax\ifhmode\unskip\space\fi MR }
\providecommand{\MRhref}[2]{%
  \href{http://www.ams.org/mathscinet-getitem?mr=#1}{#2}
}
\providecommand{\href}[2]{#2}

\end{document}